\def\no{\if01}
\def\iftwelvept{\no}

\def\ifusepdf{\no}
\def\ifpsfont{\no}

\iftwelvept
\documentclass[leqno,12pt]{amsart}
\else
\documentclass[leqno,11pt]{amsart}
\fi
\usepackage{amssymb}
\usepackage{amscd}
\usepackage{latexsym}
\usepackage{verbatim}
\usepackage[all]{xy}

\ifusepdf
\usepackage{hyperref}
\else\fi
\ifpsfont
\usepackage[T1]{fontenc}
\usepackage{times}
\else\fi

\iftwelvept
\else\fi

\theoremstyle{plain}
\newtheorem{Theorem}{Theorem}[section]

\newtheorem{Proposition}[Theorem]{Proposition}
\newtheorem{Lemma}[Theorem]{Lemma}
\newtheorem{Corollary}[Theorem]{Corollary}

\theoremstyle{definition}

\newtheorem{Definition}[Theorem]{Definition}
\newtheorem{Remark}[Theorem]{Remark}
\newtheorem{Example}[Theorem]{Example}

\renewcommand{\theTheorem}{\arabic{section}.\arabic{Theorem}}

\newcommand{\ZZ}{{\mathbb{Z}}}
\newcommand{\QQ}{{\mathbb{Q}}}
\newcommand{\RR}{{\mathbb{R}}}
\newcommand{\RRR}{\mathsf{R}}
\newcommand{\CC}{{\mathbb{C}}}
\newcommand{\PP}{{\mathbb{P}}}

\newcommand{\NNNN}{\operatorname{N}}
\newcommand{\GG}{{\overline{G}}}
\newcommand{\KK}{K}

\newcommand{\HOM}{\mathsf{Hom}}

\newcommand{\MG}{\mathsf{MG}}

\newcommand{\DIS}{\textup{dis}}
\newcommand{\Sm}{\textup{Sm}}

\newcommand{\uni}{\mathbf{1}}
\newcommand{\Sh}{\operatorname{Sh}}

\newcommand{\Chow}{\textup{Chow}}
\newcommand{\DDD}{{\mathcal{D}}}
\newcommand{\wSSS}{\widehat{\mathcal{S}}}

\newcommand{\CCC}{{\mathcal{C}}}

\newcommand{\PR}{\operatorname{Pr}^{\textup{L}}}
\newcommand{\smCat}{\operatorname{Cat}_{\infty}^{\textup{sMon}}}
\newcommand{\wsmCat}{\widehat{\operatorname{Cat}}_{\infty}^{\textup{sMon}}}

\newcommand{\Rep}{\operatorname{Rep}}

\newcommand{\OO}{{\mathcal{O}}}

\newcommand{\oDTM}{\overline{\mathsf{DTM}}}

\newcommand{\Vect}{\operatorname{Vect}}

\newcommand{\DM}{\mathsf{DM}}

\newcommand{\DTM}{\mathsf{DTM}}
\newcommand{\CH}{\operatorname{CH}}
\newcommand{\Ch}{\mathsf{ch}}

\newcommand{\Alb}{\textup{Alb}}
\newcommand{\TM}{\mathsf{TM}}

\newcommand{\Hom}{\operatorname{Hom}}

\newcommand{\Ker}{\operatorname{Ker}}
\newcommand{\Coker}{\operatorname{Coker}}
\newcommand{\Comp}{\operatorname{Comp}}
\newcommand{\Spec}{\operatorname{Spec}}

\newcommand{\Gal}{\operatorname{Gal}}

\newcommand{\SP}{\operatorname{Sp}}

\newcommand{\Mod}{\operatorname{Mod}}

\newcommand{\SSS}{\mathcal{S}}
\newcommand{\Free}{\mathbb{F}}

\newcommand{\Cat}{\textup{Cat}_{\infty}}

\newcommand{\Map}{\operatorname{Map}}
\newcommand{\Fun}{\operatorname{Fun}}

\newcommand{\End}{\operatorname{End}}

\newcommand{\sSet}{\operatorname{Set}_{\Delta}}

\newcommand{\Aff}{\operatorname{Aff}}
\newcommand{\GL}{\textup{GL}}
\newcommand{\Grp}{\operatorname{Grp}}
\newcommand{\Fin}{\operatorname{Fin}}

\newcommand{\wCat}{\widehat{\textup{Cat}}_{\infty}}

\newcommand{\CAlg}{\operatorname{CAlg}}

\newcommand{\QC}{\operatorname{QC}}

\newcommand{\hhh}{\operatorname{h}}
\newcommand{\Ind}{\operatorname{Ind}}

\newcommand{\Pic}{\operatorname{Pic}}

\newcommand{\Aut}{\operatorname{Aut}}

\newcommand{\Sym}{\operatorname{Sym}}
\newcommand{\Proof}{{\sl Proof.}\quad}
\newcommand{\QED}{{\unskip\nobreak\hfil\penalty50\quad\null\nobreak\hfil
{$\Box$}\parfillskip0pt\finalhyphendemerits0\par\medskip}}

\begin{document}

\title{Motivic rational homotopy type}

\author{Isamu Iwanari}


\address{Mathematical Institute, Tohoku University, Sendai, Miyagi, 980-8578 Japan}

\thanks{The research is partially supported by the grant of JSPS}

\email{iwanari@math.tohoku.ac.jp}

\maketitle

\setcounter{tocdepth}{1}
\tableofcontents

\section{Introduction}

In this paper, our interest lies in motives for rational homotopy
types of algebraic varieties.
Rational homotopy theory originated from Quillen
\cite{Q} and Sullivan \cite{S}.
In both approaches,
the main object of interest is
an algebraic invariant associated to
a topological space, that encodes a rational homotopy type of the space
under suitable conditions.
In Quillen's theory, the algebraic invariant
is a differential graded Lie algebra
obtained from a simply connected
topological space.
On the other hand, to a topological space $S$,
Sullivan associated a commutative differential graded (dg) algebra
$A_{PL}(S)$ of polynomial differential forms on $S$ with rational coefficients.
The cohomology ring of $A_{PL}(S)$ is isomorphic to the
graded-commutative ring $H^*(S,\QQ)$ of the singular cohomology.
In his approach, the main algebraic invariants of $S$ are
$A_{PL}(S)$ and its (so-called) Sullivan model.

We now turn to our attention to algebraic varieties.
One of motivating sources of motives is Hodge theory.
When $S$ is a complex algebraic variety,
thanks to the works of Morgan \cite{M} and Hain \cite{Ha},
a suitable model of $A_{PL}(S)$ admits a mixed Hodge structure
in an appropriate setting.
Their work generalized the classical Hodge theory to
Hodge theory for
higher rational
homotopy groups and unipotent fundamental groups, i.e.,
the pro-unipotent completion of fundamental group.
Meanwhile, in 80's, a notion of motivic homotopy type was envisaged by Grothendieck
\cite{G}.
Deligne and Gonchalov developed
a motivic theory for the pro-unipotent completions of fundamental groups
in the setting of mixed Tate (and Artin-Tate) motives over a number field
and its ring of integers \cite{DG}.

Our investigation is an attempt to define and study
a motivic generalization of $A_{PL}(S)$.
In order to get a feeling for invariants we will study,
let us compare the homotopy (triangulated) category arising from topological
spaces and the category of motives.
Let $DM^\otimes(k)$ be the symmetric monoidal
triangulated category of Voevodsky motives over a perfect field $k$,
\cite{MVW}, \cite{Voe} (here $DM^\otimes(k)$ is allowed to admit infinite coproducts).
One of pleasant features of $DM^\otimes(k)$ is
that the construction is given by ``doing homotopy theory'' of schemes,
so that the analogy is quite transparent,
while motivic cohomology groups appear as the hom sets in $DM(k)$.
By analogy with homotopy theory,
$DM^\otimes(k)$ should be thought of as an analogue
of the homotopy category
of module spectra over the Eilenberg-MacLane ring spectrum $H\ZZ$.
The motive $M(X)\in DM(k)$ associated to $X$ \cite{MVW} plays the role of
the singular chain complex of a topological space.
We now work with rational coefficients instead of $\ZZ$,
and take a point of view that
a topological counterpart of $DM^\otimes(k)$ is the derived category of $\QQ$-vector spaces.
Remember that for a topological space $S$,
$A_{PL}(S)$ is a commutative dg algebra with rational coefficients
whereas
the singular cochain complex $C^*(S,\QQ)$ is only a dg algebra that is
not necessarily commutative.
We can think that the commutative dg algebra
$A_{PL}(S)$ amounts to the (underlying) complex $C^*(S,\QQ)$ endowed with
an $E_{\infty}$-algebra structure, that is, a commutative algebra structure
in the operadic or $(\infty,1)$-categorical sense.
This structure is crucial for rational homotopy theory.
(Also, the integral singular
cochain complex $C^*(S,\ZZ)$ admits an $E_{\infty}$-algebra structure \cite{FrC},
\cite{MS}, and it is important to generalizations of rational homotopy theory
such as integral homotopy theory \cite{Man}.)
To incorporate such structures
and to pursue the comparison, we need to replace the derived category
of $\QQ$-vector spaces with its $(\infty,1)$-categorical enhancement,
i.e., the derived $(\infty,1)$-category $\mathsf{D}(\QQ)$
of $\QQ$-vector spaces, that inherits a symmetric monoidal structure
given by the tensor product of complexes.
For the introductions to the $(\infty,1)$-categorical language,
we refer to \cite[Chapter 1]{HTT}, \cite{Ber}, \cite{GrI} for instance.
Then $A_{PL}(S)$ may be viewed as a commutative algebra object
of the symmetric monoidal $(\infty,1)$-category $\mathsf{D}^\otimes(\QQ)$
in the $(\infty,1)$-categorical sense.
Let $\DM^\otimes(k)$ be a symmetric monoidal
$(\infty,1)$-category of motives, that
is an $(\infty,1)$-categorical enhancement of $DM^\otimes(k)$.
Let $\CAlg(\DM^\otimes(k))$ be the $(\infty,1)$-category of
commutative algebra objects of $\DM^\otimes(k)$.
The analogy suggests that it is natural to think that
a motivic generalization of $A_{PL}(-)$ should be defined as an object
of $\CAlg(\DM^\otimes(k))$ whose underlying object in $\DM(k)$
is equivalent to the (weak) dual of $M(X)$.
There are (at least) two approaches to constructing this:
\begin{enumerate}
\renewcommand{\labelenumi}{(\roman{enumi})}

\item If $\Sm_k$ denotes the category of smooth schemes over $k$,
equipped with the symmetric monoidal structure given by
the product $X\times_kY$, then an object $X$ of $\Sm_k$ can be viewed as a
cocommutative coalgebra
object such that the comultiplication is the diagonal $X\to X\times_kX$,
and the counit is the structure morphism $X\to \Spec k$.
If we regard the assignment $X\mapsto M(X)$ as a symmetric monoidal
functor $\Sm_k\to \DM^\otimes(k)$, then $M(X)$ is a cocommutative
coalgebra object in $\DM(k)$. Let $\uni_k$ be a unit object in $\DM(k)$.
Then the internal hom object $\mathsf{Hom}_{\DM(k)}(M(X),\uni_k)$
inherits a commutative algebra structure in the $(\infty,1)$-categorical sense
(i.e., an $E_{\infty}$-algebra structure) from $M(X)$.

\item Let $X$ be an object of $\Sm_k$ and let $f:X\to \Spec k$ be
the structure morphism.
Suppose that a symmetric monoidal $(\infty,1)$-category $\DM^\otimes(X)$
of motives over $X$ is available and there is an adjoint pair
$f^*:\DM(k)\rightleftarrows \DM(X):f_*$.
If $f^*$ is symmetric monoidal, then the right adjoint $f_*$ is a
lax symmetric
monoidal functor, so that $f_*$ sends a commutative algebra object in $\DM(X)$
to a commutative algebra object in $\DM(k)$. We denote by $\uni_X$
a unit object of $\DM(X)$ and think of it as a commutative algebra object.
We then have a commutative algebra object $f_*(\uni_X)$, that is a natural
candidate.
\end{enumerate}

The approach (i) is reminiscent of the setup in topology: singular
chain complexes and singular cochain complexs 
(but, the assignment $S\mapsto C_*(S,\ZZ)$ is only oplax monoidal).
We will adopt the approach (ii) since it gives
a clear relationship with the relative situation. We will
use the formalism of motives over $X$, extensively developed by Cisinski and D\'eglise. For a smooth scheme $X$, we define an object $M_X$ of $\CAlg(\DM^\otimes(k))$, which we shall
refer to as the cohomological motivic algebra of $X$.
The definition will be given in Section~\ref{CMAsection}. Actually,
in Section~\ref{CMAsection},
we work with not only rational coefficients but an arbitrary coefficient ring.

The first important property of $M_X$ is that
a (topological) realization of $M_X$ 
gets identified with the commutative dg algebra $A_{PL}(X^t)$ of polynomial
differential forms on the underlying topological space $X^t$ of $X\times_{k}\Spec\CC$
when $k\subset \CC$.
To Weil cohomology theory such as singular cohomology, analytic/algebraic de Rham cohomology, $l$-adic \'etale cohomology, one can associate a symmetric monoidal functor called a realization functor:
\[
\mathsf{R}:\DM^\otimes(k)\to \mathsf{D}^\otimes(K)
\]
where $K$ is a coefficient field of cohomology theory,
and $\mathsf{D}^\otimes(K)$ is the symmetric monoidal derived $(\infty,1)$-category of $K$-vector spaces. The field $K$ is assumed to be
of characteristic zero.
For example, when $k$ is embedded in $\CC$, the realization functor 
$\mathsf{R}:\DM^\otimes(k)\to \mathsf{D}^\otimes(\QQ)$ associated to singular cohomology theory (with rational coefficients) carries $M(X)$ to a complex quasi-isomorphic to the singular chain complex $C_*(X^t,\QQ)$ of the underlying topological space $X^t$. Notice that the realization functor is symmetric monoidal.
It gives rise to a functor
\[
\CAlg(\DM^\otimes(k))\to \CAlg(\mathsf{D}^\otimes(K)),
\]
which we call the multiplicative realization functor,
where $\CAlg(\mathsf{D}^\otimes(K))$ is the $(\infty,1)$-category of commutative algebra objects in $\mathsf{D}^\otimes(K)$. One can naturally identify
$\CAlg(\mathsf{D}^\otimes(K))$ with the $(\infty,1)$-category obtained from the category of commutative dg algebras over $K$ by inverting quasi-isomorphisms (cf. Section~\ref{convention}).
In the case of singular cohomology, we have $\CAlg(\DM^\otimes(k))\to \CAlg(\mathsf{D}^\otimes(\QQ))$.
The commutative dg algebra $A_{PL}(X^t)$ appears as the image of $M_X$
under the multiplicative realization functor. This property is proved in Section~\ref{realization}. Thus, along with merely an analogy,
the multiplicative realization functor
relates $M_X$ with $A_{PL}(X^t)$.
It is worth emphasizing that
it allows one to promote many operations on $A_{PL}(X^t)$ to a motivic level
(even as the elementary nature tends to obscure the significance).
For example,
the multiplicative realization functor preserves (small) colimits.
Suppose that $x$ is a 
$k$-rational point on $X$.
Let $\epsilon:A_{PL}(X^t)\to \QQ$ be the augmentation induced by
the point $x$ on $X^t$. The bar construction of the augmented commutative
dg algebra
can be described in terms of (a cosimplicial diagram of)
colimits. Thus, it is possible to promote the bar construction
of $A_{PL}(X^t)\to \QQ$ to a bar construction of $M_X\to \uni_k$
in $\CAlg(\DM^\otimes(k))$.

\vspace{2mm}

{\it Tannakian aspect.}
One of descriptions of ``motivic stuctures'' is
a tannakian formalism.
We discuss a tannakian aspect in Section~\ref{MGA}.
Recall that various ``topological invariants'' of algebraic varieties
are equipped with actions of groups. For instance,
an $l$-adic \'etale cohomology group has an action of the absolute
Galois group, and
a Hodge structure can be described by an action of a Mumford-Tate group.
In our context, the groups will be
the derived motivic Galois group $\MG$, introduced in \cite{Tan},
and
the associated pro-algebraic group $MG$ we call the motivic Galois group
(see the beginning of Section~\ref{MGA}, Section~\ref{underlyinggroup},
and \cite{Tan}).
By using $M_X$ we construct a canonical action of $\MG$ on $A_{PL}(X^t)$
(when $k\subset \CC$).
It is a tannakian representation of the motivic structure
on the rational homotopy type.
When $X$ has a base point,
it is possible to deduce the pro-unipotent completions 
$\pi_i(X^t,x)_{uni}$ of homotopy
groups $\pi_i(X^t,x)$ ($i\ge1$) from $A_{PL}(X^t)$
with the augmentation. 
We obtain canonical actions of the motivic Galois group $MG$ on pro-unipotent
groups $\pi_i(X^t,x)_{uni}$ from the action of $\MG$ on $A_{PL}(X^t)$
(cf. Theorem~\ref{actiononhomotopy}, Corollary~\ref{actiononhomotopy2}).
Thus, from the tannakian viewpoint, our study may be regarded as
a generalization of motivic structures on (co)homology groups to
motivic structures on
the unipotent non-abelian fundamental groups and
higher rationalized homotopy groups.

\vspace{2mm}

{\it Structure of cohomological motivic algebras.}
In order to understand things more explicitly,
it is natural to attempt to understand the structure of $M_X$,
that is, what a cohomological motivic algebra looks like.
In Section~\ref{Sullivanmodel}, as a first step towards the understanding,
we describe an explicit structure of the cohomological motivic algebra
in several cases such as a projective space over a field.
To do this, we recall an approach that traces back to Sullivan's work.
A (minimal) Sullivan model of $A_{PL}(S)$ is given
by an iterated homotopy pushout of free commutative dg algebras
(see e.g. \cite{He}, \cite{Hin}, \cite{FHT} or the beginning of Section~\ref{Sullivanmodel}).
Based on this idea, we describe $M_X$ as
a colimit of an analogous diagram of
free commutative algebra objects
in $\CAlg(\DM^\otimes(k))$ in an explicit way.
Unlike the classical rational homotopy theory,
the study of $M_X$ is not so simple even in relatively elementary
cases: we need some devices and deep results.
This difference may be regarded as a reflection of the fact that
$M_X$ has rich and interesting structures.
For instance,
suppose that $C$ is a proper smooth curve of genus $g>1$
with a base $k$-rational point $c$.
Let $J_C$ be the Jacobian variety and let $u:C\to J_C$ be
the Abel-Jacobi morphism.
We here take a viewpoint that
the Abel-Jacobi morphism is an ``algebraic abelianization'' of $C$:
when $k=\CC$, the map
$C^t\to J_{C}^t$ of the underlying toplogical spaces
induces an abelianization $\pi_1(C^t,c)\to \pi_1(J_{C}^t,u(c))\simeq \pi_1(C^t,u)^{ab}$.
The Abel-Jacobi morphism $u$ induces
a morphism $u^*:M_{J_C}\to M_C$ of cohomological
motivic algebras.
Then it gives rise to an inductive sequence in $\CAlg(\DM^\otimes(k))$:
\[
M_{J_C}=M_1\to M_2\to \cdots \to M_n\to M_{n+1}\to \cdots \to M_C
\]
that decomposes $u^*:M_{J_C}\to M_C$ such that
$M_C$ is a filtered colimit $\varinjlim_{n\ge 1}M_n$
(cf. Section~\ref{Sull1}, Section~\ref{curveexp}).
One can think of this sequence (or co-tower) starting with $M_{J_C}$
as a structure of $M_C$ or a refined Abel-Jacobi morphism.
It is notable that it does not exist in the category of schemes
and does not arise from $\DM(k)$ .
Roughly speaking, this sequence gives a step-by-step description of
the non-abelian nature of $C$ that starts with
its ``abelian part'' $M_{J_C}$.
From a perspective of the formality,
it is not reasonable to expect a formality of $M_X$
of a smooth projective variety $X$ in general
(even if
one can define a formality by using a motivic $t$-structure).
Actually, there is a counterexample to
the formality at the Hodge level (see \cite{CCM}).

The large class is yet to be explored and remains mysterious, so that
one may expect more to understand
structures of cohomological motivic algebras.

\vspace{2mm}

{\it Cotangent motives.}
In Section~\ref{cotangenthomotopy}, we introduce
a new invariant of a pointed smooth scheme $(X,x)$ over a perfect field,
that lies in $\DM(k)$.
The invariant $LM_{(X,x)}$ in $\DM(k)$ is defined by means of cotangent
complex of $M_X$ endowed with the augmentation induced by $x$.
We shall
call $LM_{(X,x)}$
the cotangent motive of $X$ at $x$ (cf. Definition~\ref{cotangentmotivedef}). For the definition,
we apply the theory of cotangent complexes
in a very general setting, developed by Lurie.
We prove that the rationalized homotopy group appears as the
realiziation
of $LM_{(X,x)}$ (cf. Theorem~\ref{motivichomotopy}, Theorem~\ref{motivichomotopy2}). Namely, when $k$ is embedded in $\CC$ and the underlying
topological space $X^t$ is simply connected,
$H^i(\mathsf{R}(LM_{(X,x)}))$ is the dual of the $i$-th
rationalized homotopy
group of $X^t$. In addition, $H^1(\mathsf{R}(LM_{(X,x)}))$
can be identified with the cotangent space of the origin of
the pro-unipotent completion
of the fundamental group, that is, the 
``linear data'' of the fundamental group.
By using Hodge realization of $LM_{(X,x)}$
one can obtain a mixed  Hodge structure on the rational homotopy group
(in the simply conneced case).
Intuitively,
we may consider $LM_{(X,x)}$ to be a motive
for (the dual of) rational higher homotopy groups and the linear data of
 the fundamental group.
Though $LM_{(X,x)}$
has less information than $M_X$,
the motive $LM_{(X,x)}$ has the relation with
homotopy groups in a more direct way than $M_X$,
and furthermore one can consider motivic cohomology of $LM_{(X,x)}$
since it belongs to $\DM(k)$.
We apply the (explicit) study of $M_X$
in Section~\ref{Sullivanmodel} to compute $LM_{(X,x)}$.
Indeed, one of motivations for it is computation of the cotangent motives.
For instance, if $\mathbb{P}^n$ is the $n$-dimensional projective space (over a perfect field) endowed with a base point $x$, then
\[
LM_{(\mathbb{P}^n,x)}\simeq \uni_k(-1)[-2]\oplus \uni_k(-n-1)[-2n-1],
\]
where ``$(s)$'' and ``$[t]$'' indicate the Tate twist and the shift, respectively. This means that $\uni_k(1)$ is a ``motive for the second rational homotopy group'', and 
$\uni_k(n+1)$ is a ``motive for the $(2n+1)$-th rational homotopy group'' (cf. Remark~\ref{motcpsR}).

\vspace{2mm}

{\it Homotopy exact sequence.}
Remember the homotopy exact sequence for \'etale
fundamental groups
\[
1 \to \pi_1^{\textup{\'et}}(X\times_{k}\Spec \bar{k},\bar{x})\to \pi_1^{\textup{\'et}}(X,\bar{x}) \to \Gal(\bar{k}/k)\to 1
\]
where $\bar{k}$ is a separable closure of $k$.
This plays a central role in the theory of \'etale fundamental groups.
In Section~\ref{mhes}, by means of a tannakian theory developed in \cite{DTD},
when $X$ is an algebraic curve we formulate and prove a version of the homotopy exact sequence
in which the derived motivic Galois group (or stack) instead of $\Gal(\bar{k}/k)$ (cf. Proposition~\ref{MHES}).

\vspace{2mm}

From a conceptual point of view,
subjects in Section~\ref{MGA}, \ref{Sullivanmodel}, \ref{cotangenthomotopy},
\ref{mhes} are interconnected with each other.
Nevertheless, it is possible to read these Sections in any order with some exceptions.

It is desirable and important to study various
(multiplicative) realizations of cohomological motivic algebras $M_X$
and $LM_{(X,x)}$ in detail: \'etale, de Rham, Hodge, crystalline realizations,
a relation with Chen's theory of iterated integrals, etc.
These issues remain untouched and are beyond the scope of this paper. We hope to return to subjects in the future.

In Appendix, in the case of mixed Tate motives over a number field
we compare our approach and an approach to motivic fundamental groups
due to Deligne and Gonchalov. We hope that the comparison is helpful
for understanding circle of ideas from the viewpoint of their work.

\section{Notation and Convention}
\label{convention}

\subsection{}
We shall use the theory of {\it quasi-categories} extensively developed by
Joyal and Lurie from the viewpoint of $(\infty,1)$-categories.
This theory provides us with powerful tools and adequate language for
our purpose, though a part of contents might be reformulated
in term of other languages such as model categories or the like.
Following \cite{HTT}, we shall refer to quasi-categories
as {\it $\infty$-categories}.
Our main references are \cite{HTT} and \cite{HA}.
To an ordinary category $\mathcal{C}$,
one can assign an $\infty$-category by taking
its nerve $\NNNN(\mathcal{C})$.
Such simplicial sets $\NNNN(\mathcal{C})$ arising from ordinary categories
naturally constitute a full subcategory of
the simplicial category of $\infty$-categories.
Therefore,
when we treat ordinary categories we often omit the nerve $\NNNN(-)$
and think of them directly as $\infty$-categories.
We often refer to a map $S\to T$ of $\infty$-categories
as a functor. We call a vertex in an $\infty$-category $S$
(resp. an edge) an object (resp. a morphism).
We use Grothendieck universes $\mathbb{U}\in \mathbb{V}\in \mathbb{W}\in\ldots$
and usual mathematical objects such as groups, rings, vector spaces
are assumed to belong to $\mathbb{U}$.
Here is a list of (some) of the convention and notation that we will use:

\begin{itemize}

\item $\Delta$: the category of linearly ordered finite sets (consisting of $[0], [1], \ldots, [n]=\{0,\ldots,n\}, \ldots$)

\item $\Delta^n$: the standard $n$-simplex as the simplicial set represented by
$[n]$,

\item $\sSet$: the category of simplicial sets,

\item $\textup{N}$: the simplicial nerve functor (cf. \cite[1.1.5]{HTT})

\item $\Gamma$: the nerve of the category of pointed finite sets, $\langle 0\rangle=\{*\}, \langle 1\rangle=\{*,1\},\ldots,\langle n\rangle=\{*,1,\ldots,n\},\ldots$

\item $\mathcal{C}^{op}$: the opposite $\infty$-category of an $\infty$-category $\mathcal{C}$. For a functor $F:\mathcal{C}\to \mathcal{D}$, we denote by
$F^{op}:\mathcal{C}^{op}\to \mathcal{D}^{op}$ the induced functor

\item Let $\mathcal{C}$ be an $\infty$-category and suppose that
we are given an object $c$. Then $\mathcal{C}_{c/}$ and $\mathcal{C}_{/c}$
denote the undercategory and overcategory, respectively (cf. \cite[1.2.9]{HTT}).

\item $\CCC^\simeq$: the largest Kan subcomplex (contained) in an $\infty$-category
$\CCC$,
that is, the Kan complex obtained from $\CCC$
by restricting
morphisms (edges) to equivalences.

\item $\operatorname{Cat}_\infty$: the $\infty$-category of small $\infty$-categories, Similarly,
$\wCat$ denotes $\infty$-category of large $\infty$-categories (i.e.,
$\infty$-categories that belong to $\mathbb{V}$),

\item $\SSS$: $\infty$-category of small spaces. We denote by $\widehat{\SSS}$
the $\infty$-category of large $\infty$-spaces (cf. \cite[1.2.16]{HTT})

\item $\textup{h}(\mathcal{C})$: homotopy category of an $\infty$-category (cf. \cite[1.2.3.1]{HTT})

\item $\Fun(A,B)$: the function complex for simplicial sets $A$ and $B$

\item $\Fun_C(A,B)$: the simplicial subset of $\Fun(A,B)$ classifying
maps which are compatible with
given projections $A\to C$ and $B\to C$.

\item $\Map(A,B)$: the largest Kan subcomplex of $\Fun(A,B)$ when $B$ is an $\infty$-category.

\item $\Map_{\mathcal{C}}(C,C')$: the mapping space from an object $C\in\mathcal{C}$ to $C'\in \mathcal{C}$ where $\mathcal{C}$ is an $\infty$-category.
We usually view it as an object in $\mathcal{S}$ (cf. \cite[1.2.2]{HTT}).
If $\mathcal{C}$ is an ordinary category, we write
$\Hom_{\mathcal{C}}(C,C')$ for the hom set.

\item $C^\vee$:
For an object $C$ of a symmetric monoidal $\infty$-category $\mathcal{C}$,
we write $C^\vee$ for a dual of $C$ when $C$ is a dualizable object.
If there are internal objects,
we write $C^\vee$ also for the weak dual, that is,
the internal hom $\mathsf{Hom}_{\mathcal{C}}(C,\uni_{\mathcal{C}})$ with $\uni_{\mathcal{C}}$ a unit object.

\item $\Ind(\mathcal{C})$; $\infty$-category of Ind-objects
in an $\infty$-category $\mathcal{C}$
(see \cite[5.3.5.1]{HTT}, \cite[4.8.1.13]{HA} for the symmetric monoidal setting).

\item $\PR$: the $\infty$-category of presentable $\infty$-categories
whose morphisms are left adjoint functors.

\end{itemize}

\subsection{}
{\it From model categories to $\infty$-categories.}
We recall Lurie's construction
by which one can obtain
$\infty$-categories
from a category (more generally
$\infty$-category) endowed with a prescribed collection of
morphisms (see \cite[1.3.4, 4.1.3, 4.1.4]{HA} for details).
It can be viewed as an alternative of the Dwyer-Kan
hammock localization.
Let $\mathcal{D}$ be a category and let $W$ be a collection
of morphisms in $\mathcal{D}$ which is closed under composition
and contains all isomorphisms.
A typical example of $(\mathcal{D},W)$
which we have in mind is
$(\mathbb{M},W_{\mathbb{M}})$ such that
$\mathbb{M}$ is a model category (see e.g. \cite[Appendix]{HTT}, \cite{Hov})
and $W_{\mathbb{M}}$ is the collection of all weak equivalences.
For $(\mathcal{D},W)$, there is an $\infty$-category $\NNNN(\mathcal{D})[W^{-1}]$
and a functor
$\xi:\textup{N}(\mathcal{D})\to \NNNN(\mathcal{D})[W^{-1}]$
such that for any $\infty$-category $\CCC$ the composition
induces a fully faithful functor
\[
\Map(\NNNN(\mathcal{D})[W^{-1}],\CCC)\to \Map(\textup{N}(\mathcal{D}),\CCC)
\]
whose essential image consists of those functors $F:\textup{N}(\mathcal{D})\to\CCC$ such that $F$ carry morphisms lying in $W$
to equivalences in $\CCC$.
We shall refer to $\NNNN(\mathcal{D})[W^{-1}]$
as the $\infty$-category obtained from $\mathcal{D}$
by inverting morphisms in $W$.
Consider $(\mathbb{M},W_{\mathbb{M}})$ such that
$\mathbb{M}$ is a combinatorial model category
and $W_{\mathbb{M}}$ is the collection of weak equivalences.
The $\infty$-category
$\mathbb{M}^c[W^{-1}]:=\NNNN(\mathbb{M}^c)[(\mathbb{M}^c\cap W_{\mathbb{M}})^{-1}]$ is presentable where $\mathbb{M}^c$ is the full subcategory of cofibrant objects.
(When $\mathbb{M}$ is a monoidal model category,
it is convenient to work with
 the full subcategory of cofibrant objects $\mathbb{M}^c\subset \mathbb{M}$
instead of $\mathbb{M}$.)
If $\mathbb{M}$ is a stable model category, then $\mathbb{M}^c[W^{-1}]$ is 
a stable $\infty$-category (cf. \cite{Tan}).
The homotopy category of $\mathbb{M}^c[W^{-1}]$ coincides
with the homotopy category of the model category $\mathbb{M}$.
If $\mathbb{M}$ is a symmetric monoidal model
category (whose unit object is cofibrant),
$\mathbb{M}^c[W^{-1}]$ is promoted to a
symmetric monoidal $\infty$-category $\mathbb{M}^c[W^{-1}]^\otimes:=\NNNN(\mathbb{M}^c)[(\mathbb{M}^c\cap W_{\mathbb{M}})^{-1}]^{\otimes}$
(see below for symmetric monoidal
$\infty$-categories). 
In addition, there is a symmetric monoidal
functor $\tilde{\xi}:\NNNN(\mathbb{M}^c)^\otimes \to \mathbb{M}^c[W^{-1}]^\otimes$ which has $\xi$ as the underlying functor and
satisfies a similar universal property.
If $\mathbb{M}$ is combinatorial, then
the tensor product $\otimes:\mathbb{M}^c[W^{-1}]\times \mathbb{M}^c[W^{-1}]\to \mathbb{M}^c[W^{-1}]$ preserves
small colimits separately in each variable.
Let $\mathbb{L}$ be another symmetric monoidal model category
and let $\phi:\mathbb{M}\to \mathbb{L}$ be a symmetric monoidal functor.
If $\phi$ carries cofibrant objects to cofibrant objects and preserves
weak equivalences between them (e.g. symmetric monoidal left
Quillen functors), it induces a symmetric monoidal functor
$\mathbb{M}^c[W^{-1}]^\otimes\to \mathbb{L}^c[W^{-1}]^\otimes$
of symmetric monoidal $\infty$-categories.

\subsection{}
{\it Symmetric monoidal $\infty$-categories, modules and algebras.}
We use the theory of (symmetric) monoidal $\infty$-categories
developed in \cite{HA}.
A symmetric monoidal $\infty$-category is a coCartesian fibration
$\mathcal{C}^\otimes\to \Gamma$ that satisfies a ``symmetric monoidal
condition'', see \cite[2.1.2]{HA}.
For a symmetric monoidal $\infty$-category $\mathcal{C}^\otimes\to \Gamma$,
we often write $\mathcal{C}$ for the underlying $\infty$-category.
Also, by abuse of notation, we usually use the superscript in $\mathcal{C}^\otimes$ to indicate a symmetric monoidal structure on an $\infty$-category.
For a symmetric monoidal $\infty$-category $\mathcal{C}^\otimes$,
we write $\CAlg(\mathcal{C}^\otimes)$ (or simply
$\CAlg(\mathcal{C})$) for the $\infty$-category
of commutative algebra objects in $\mathcal{C}^\otimes$.
Let $A$ be a commutative ring spectrum, that is, a commutative algebra
object in the category $\SP$ of spectra.
We write $\Mod^\otimes_A$ for the symmetric monoidal $\infty$-category
of $A$-module spectra, (see e.g. \cite{HA}).
We put $\CAlg_A=\CAlg(\Mod_A^\otimes)$.
For an ordinary commutative ring $K$,
we put $\Mod^\otimes_{K}:= \Mod^\otimes_{HK}$ and
$\CAlg_{K}:=\CAlg_{HK}$
where $HK$ is the Eilenberg-MacLane ring spectrum.

Let $K$ be a field of characteristic zero.
Let $\Comp^\otimes(K)$ be the symmetric monoidal
category of cochain complexes of $K$-vector spaces
(the symmetric monoidal structure is given by the tensor product
of cochain complexes).
This category admits a projective combinatorial symmetric monoidal
model structure, whose weak equivalences are quasi-isomorpisms,
and whose cofibrations (resp. fibrations) are monomorphisms
(resp. epimorphisms), see e.g. \cite[Section 2.3]{Hov} or \cite[7.1.2.8]{HA}.
We shall write $\mathsf{D}^\otimes(K)$ for the symmetric monoidal
stable presentable $\infty$-category obtained from $\Comp^\otimes(K)$
by inverting weak equivalences.
According to \cite[7.1.2.12, 7.1.2.13]{HA}, there is a canonical equivalence
$\mathsf{D}^\otimes(K)\simeq \Mod_K^\otimes$.
We refer to $\mathsf{D}^\otimes(K)$ and $\Mod_K^\otimes$
as the (symmetric monoidal) derived $\infty$-category of $K$-vector spaces.
The equivalence $\mathsf{D}^\otimes(K)\simeq \Mod_K^\otimes$
induces $\CAlg(\mathsf{D}^\otimes(K))\simeq \CAlg_{K}=\CAlg(\Mod_K^\otimes)$. Let $\CAlg_K^{dg}$ be the category of commutative differential graded
$K$-algebras. A commutative differential graded
$K$-algebras is a commutative algebra object in $\Comp^\otimes(K)$.
There is a natural forgetful functor $U:\CAlg_K^{dg}\to \Comp(K)$.
The category $\CAlg_K^{dg}$ admits a combinatorial model structure
such that a morphism $f$ is
a weak equivalences (resp. a fibration)
if and only if $U(f)$ is a quasi-isomorphism
(resp. a epimorphism) (here, we use the assumption of characteristic zero) .
If we write $\NNNN(\CAlg_{K}^{dg})[W^{-1}]$ for the $\infty$-category
obtained from $\CAlg_{K}^{dg}$ by inverting weak equivalences,
then there is a canonical equivalences
$\NNNN(\CAlg_{K}^{dg})[W^{-1}]\simeq \CAlg_K$
(see \cite[7.1.4.10, 7.1.4.11]{HA}, \cite[4.5.4.6]{HA}).
We often use these equivalences
\[
\NNNN(\CAlg_{K}^{dg})[W^{-1}]\simeq \CAlg_K \simeq \CAlg(\mathsf{D}^\otimes(K)).
\]

A variety is a geometrically connected
scheme separated of finite type over a field.

\section{Cohomological motivic algebras}
\label{CMAsection}

Let $K$ be a commutative ring.

\subsection{}

As in our previous works,
we use $\infty$-categories of mixed motives.
They are obtained from the model (dg, etc) categories of motives
or the $\infty$-categorical version of Voevodsky's construction.
In this paper, we adopt symmetric monoidal
model categories constructed by Cisinski and D\'eglise \cite{CD1}, \cite{CDT}.
Let $X$ be a smooth scheme separated of finite type over 
a perfect field $k$ (or more generally, a noetherian regular scheme).
Let $\Sm_X$ denote the category of smooth schemes separated of finite type
over $X$.
Let $\mathcal{N}^{tr}(X)$ be the Grothendieck abelian category of 
Nisnevich sheaves of $K$-modules with transfers over $X$ (see e.g.
\cite[Example 2.4]{CD1} or \cite{CDT} for this notion).
Let $\Comp(\mathcal{N}^{tr}(X))$ be the symmetric monoidal category
of (possibly unbounded) cochain complexes of $\mathcal{N}^{tr}(X)$.
Then $\Comp(\mathcal{N}^{tr}(X))$
admits a stable symmetric monoidal combinatorial
model category structure,
see \cite[Section 4, Example 4.12]{CD1}.
The construction roughly has two steps: one first defines a certain nice
model structure
whose weak equivalences are quasi-isomorphisms of complexes
of sheaves,
and in the next step one
takes a left Bousfield localization of the model structure
at $\mathbb{A}^1$-homotopy.
Using a generalization of the construction of symmetric spectra,
one can ``stabilize'' the tensor operation with a shifted Tate object
over $X$
and obtains a new category $Sp_{\textup{Tate}}(X)$ from $\Comp(\mathcal{N}^{tr}(X))$
which admits a stable
symmetric monoidal combinatorial model category structure described
in \cite[Proposition 7.13, Example 7.15]{CD1}.
Let $\phi:Y\to X$ be a morphism of smooth schemes.
It gives rise to a Quillen adjunction
\[
\phi^*:Sp_{\textup{Tate}}(X)\rightleftarrows Sp_{\textup{Tate}}(Y):\phi_*
\]
where $\phi^*$ is a symmetric monoidal left Quillen functor.
We suppose further 
that $\phi$ is smooth separated of finite type,
then there is a Quillen
 adjunction
\[
\phi_\sharp:Sp_{\textup{Tate}}(Y)\rightleftarrows Sp_{\textup{Tate}}(X):\phi^*.
\]
In this case, $\phi^*$ is both a left Quillen functor and a right Quillen
functor. Thus, it preserves (trivial) fibrations and (trivial) cofibrations.
Moreover, by using Ken Brown's lemma we see that $\phi^*$ preserves arbitrary
weak equivalences.

We let $\DM^\otimes_{eff}(X)$ be the symmetric monoidal
stable presentable $\infty$-category, which is obtained from
the full subcategory of cofibrant objects $\Comp(\mathcal{N}^{tr}(X))^c$
by inverting weak equivalences. We refer to it as the symmetric monoidal
$\infty$-category of effective mixed motives over $X$.
Similarly, $\DM^\otimes(X)$
is defined to be the symmetric monoidal
stable presentable $\infty$-category
obtained from $Sp_{\textup{Tate}}(X)^c$
by inverting weak equivalences.
We call $\DM^\otimes(X)$ the symmetric monoidal
stable persentable $\infty$-category of mixed motives over $X$.
We refer to $K$ as the coefficient ring of $\DM^\otimes(X)$.
We write $\uni_X$ for a unit object of $\DM^\otimes(X)$.
We write $\uni_X(n)$ for the Tate object for $n\in \ZZ$.
Given an object $M$ of $\DM(X)$, we usually write $M(n)$ for
the tensor product $M\otimes \uni_X(n)$ in $\DM(X)$.
The tensor product $\DM(X)\times \DM(X)\to \DM(X)$
on $\DM^\otimes(X)$
preserves small colimits separately in each variable.
The detail construction can be found in \cite[Section 5.1]{Tan}
(the notation is slightly different,
and $X$ is assumed to be the Zariski spectrum of
a perfect field in \cite{Tan}, but it works for a noetherian regular scheme $X$). The homotopy category of the
full subcategory of $\DM(\Spec k)$ spanned by compact objects
can be identified with the triangulated category of geometric motives constructed by Voevodsky \cite{Voe}.

Let $f:X \to \Spec k$ be the structure morphism.
Since
we have the restriction of the symmetric monoidal left Quillen functor
$f^*:Sp_{\textup{Tate}}(\Spec k)^c\rightarrow Sp_{\textup{Tate}}(X)^c$
between full subcategories of cofibrant objects,
inverting weak equivalences
we have a symmetric monoidal colimit-preserving functor
\[
f^*:\DM^\otimes(k):=\DM^\otimes(\Spec k)\to \DM^\otimes (X).
\]
By abuse of notation, we use the same notation
for the induced functor between $\infty$-categories.
By relative adjoint functor theorem \cite[7.3.2.6, 7.3.2.13]{HA}, there is
the right adjoint functor $f_*:\DM(X)\to \DM(k)$ that is lax
symmetric monoidal.
It induces an adjunction
\[
f^*:\CAlg(\DM^\otimes(k))\rightleftarrows \CAlg(\DM^\otimes(X)):f_*.
\]
In particular, $f_*$ carries a commutative algebra object $M$ to a commutative
algebra object $f_*(M)$ in $\DM^\otimes(k)$.
For any smooth scheme $X$, 
$\CAlg(\DM^\otimes(X))$
is a presentable $\infty$-category (cf. \cite[3.2.3.5]{HA}).
There is another left Quillen functor
$f_\sharp:Sp_{\textup{Tate}}(X)\to Sp_{\textup{Tate}}(\Spec k)$.
The restriction $Sp_{\textup{Tate}}(X)^c\to Sp_{\textup{Tate}}(\Spec k)^c$
to cofibrant objects preserves weak equivalences, and therefore
inverting weak equivalences induces
$f_\sharp:\DM(X)\to \DM(k)$. It determines an adjunction
\[
f_\sharp:\DM(X)\rightleftarrows  \DM(k):f^*.
\]
We put $M(X):=f_\sharp f^*(\uni_k)$ where $\uni_k$ is the
unit of $\DM(k)$.

Let us consider the unit object $\uni_X=f^*(\uni_k)$ in $\DM^\otimes(X)$
which we regard as a commutative algebra object in $\DM^\otimes(X)$.
The image $f_*(\uni_X)=f_*f^*(\uni_k)$ is a commutative algebra object in
$\DM^\otimes(k)$, namely, $f_*(\uni_X)$ in $\CAlg(\DM^\otimes(k))$.

\begin{Definition}
We define $M_X$ in $\CAlg(\DM^\otimes(k))$ to be $f_*(\uni_X)$.
We shall refer to $M_X$ as
the {\it cohomological motivic algebra} of $X$ with
coefficients in $K$.
\end{Definition}

\begin{Remark}
This algebra $M_X$ will play a role of a motivic analogue of
the singular cochain complex $C^*(S,K)$ of a topological space $S$
that is endowed with a structure of an $E_\infty$-algebra.
Our principle is that one may consider $M_X$ to be a motivic homotopy
type of $X$ with coefficients in $K$, that occurs in the title of this paper.
On the other hand, $M(X)$ is
a motivic counterpart of the singular chain complex
$C_*(S,K)$.
To our knowledge, contrary to homotopy theory,
there has been little 
attention being paid to invariants represented by highly structured algebras
in the theory of motives.
\end{Remark}

\subsection{}
\label{cohomologicalmot}
We consider functoriality of motivic cohomological algebras.
Let $f:X\to \Spec k$ and
$g:Y\to \Spec k$ be two smooth scheme separated of finite type
over $k$.
Let $\phi:Y\to X$ be a morphism over $k$.
As above, there is an adjunction $\phi^*:\CAlg(\DM^\otimes(X))\rightleftarrows \CAlg(\DM^\otimes(Y)):\phi_*$.
If we write $M_Y$ for $g_*(\uni_Y)$ we have a morphism
\[
M_X=f_*(\uni_X)\to f_*\phi_*\phi^*(\uni_X)\simeq g_*(\uni_Y)=M_Y
\]
in $\CAlg(\DM^\otimes(k))$
where the first map is induced by the unit map $\uni_X\to \phi_*\phi^*(\uni_X)\simeq \phi_*(\uni_Y)$. Thus, the assignment $X\mapsto M_X$ is contravariantly functorial
with respect to $X$.
We will write $\phi^*:M_X\to M_Y$ for this morphism in $\CAlg(\DM^\otimes(k))$
or in the underlying category $\DM(k)$.
Unfortunately, the notation
$\phi^*$ in $\phi^*:M_X\to M_Y$ overlaps
with $\phi^*:\DM(X)\to \DM(Y)$ or $\phi^*:\CAlg(\DM^\otimes(X)) \to \CAlg(\DM^\otimes(Y))$ though these have different meanings.
We hope that it causes no confusion.
The assignment $X\mapsto M(X)$ is covariantly functorial.
For $\phi:Y\to X$, consider the unit map $u:\uni_X\to f^*f_\sharp(\uni_X)$.
We then have
\[
M(Y)=g_\sharp (\uni_Y)\simeq g_\sharp \phi^*(\uni_X)\stackrel{g_\sharp\phi^*(u)}{\longrightarrow} g_\sharp \phi^* f^*f_\sharp(\uni_X)\simeq g_\sharp  g^* f_\sharp(\uni_X)\to f_\sharp(\uni_X)=M(X)
\]
where the final arrow is induced by the counit $g_\sharp  g^*\to \textup{id}$.
Let $\textup{Sm}_k$ be the nerve of the category of smooth schemes
separated of finite type over $k$.
We will give a functorial construction $X\mapsto M_X$
as a functor $\textup{Sm}_k^{op}\to \CAlg(\DM^\otimes(k))$.
The result is summarized as follows:

\begin{Proposition}
\label{cohmotalg}
Let $M(-):\Sm_k\to \DM(k)$ be the functor which carries $X$ to $M(X)$.
We define $\HOM_{\DM(k)}(-,\uni_k):\DM(k)^{op}\to \DM(k)$
to be the functor which carries $M$ to $M^\vee=\HOM_{\DM(k)}(M,\uni_k)$. By $\HOM(-,-)$ we indicate the internal Hom object.
(We will make a construction of these functors below.) 
Let $M(-)^\vee:\Sm_k^{op}\to \DM(k)$ be the composite of the above two functors, which carries $X$ to $M(X)^\vee$.
Then there is a functor
$\Xi:\textup{Sm}_k^{op}\to \CAlg(\DM^\otimes(k))$
which makes the diagram commutative
\[
\xymatrix{
 & \CAlg(\DM^\otimes(k)) \ar[d] \\
\Sm_k^{op} \ar[r]_{M(-)^\vee} \ar[ru]^{\Xi} & \hhh(\DM(k))
}
\]
where the right vertical arrow is the forgetful functor.
\end{Proposition}

We first construct $\Xi:\textup{Sm}_k^{op}\to \CAlg(\DM^\otimes(k))$.
The busy readers are invited to skip the remainder for the time being
and proceed to Section~\ref{somerem} or~\ref{someex}.
We consider the following general situation.
The functor $\Xi$ will appear in Example~\ref{motivicex}
as an example of the following setup.
Let $I$ be the nerve of a category.
Suppose that $I$ has a final object $\star\in I$.
We are mainly interested in the case $I=\textup{Sm}_k$.
Let us consider a family $\{\mathbb{M}(X)\}_{X\in I}$
of symmetric monoidal model categories
indexed by $I$. More precisely, we assign
a combinatorial symmetric monoidal model category $\mathbb{M}(X)$
to any $X\in I$
(we here assume that a unit is cofibrant)
and assign a symmetric monoidal left Quillen functor 
$\phi^*:\mathbb{M}(X)\to \mathbb{M}(Y)$ to any morphism
$Y\to X$ in $I$.
Moreover, suppose that
for $\phi\circ \psi:Z\to Y\to X$ there is
a structural
natural equivalence
$\psi^*\phi^*\simeq  (\phi\circ \psi)^*$.
Main example is the family
$\{Sp_{\textup{Tate}}(X)\}_{X\in \Sm_k}$.
Consider the pair $(\mathbb{M}^c(X), W_X^c)$
such that $\mathbb{M}^c(X)$ is the full subcategory of
cofibrant objects in the model category $\mathbb{M}(X)$,
and $W_X^c$ is the collection of weak equivalences in $\mathbb{M}(X)^c$.
We think of this pair as the nerve of a category $\mathbb{M}^c(X)$
endowed with the collection of morphisms, determined by $W_X^c$.
We apply to the assignment $X\mapsto (\mathbb{M}(X)^c, W_X^c)$
the construction in \cite[Section 4.1.3.1, 4.1.3.2]{HA}
of inverting weak equivalences in symmetric monoidal
categories in the functorial way.
We then get a functor
\[
d:I^{op}\to \CAlg(\wCat)
\]
which carries $X$ to $\mathbb{M}_\infty^\otimes(X):=\mathbb{M}^c(X)[(W_X^c)^{-1}]$. Here $\mathbb{M}^c(X)[(W_X^c)^{-1}]$ is the symmetric monoidal
$\infty$-category obtained from $\mathbb{M}(X)^c$
by inverting $W_X^c$.
The symmetric monoidal structure on $\wCat$ is given by cartesian products,
and $\CAlg(\wCat)$ is naturally identified with the $\infty$-category of
symmetric monoidal (large) $\infty$-categories whose morphisms are
symmetric monoidal functors, cf. \cite{HA}.
Recall that the $\infty$-category $\CAlg(\wCat)$
can be realized as the full subcategory of $\Fun(\Gamma,\wCat)$
spanned by commutative monoid objects,
where $\Gamma$ is the nerve of the category of pointed finite sets.
The functor $d:I^{op}\to \CAlg(\wCat)\subset \Fun(\Gamma,\wCat)$
induces a functor $I^{op}\times \Gamma\to \wCat$.
Applying the relative nerve functor to $I^{op}\times \Gamma\to \wCat$
(cf. \cite[3.2.5]{HTT}),
we have a coCartesian fibration
\[
D:\mathcal{E}\to I^{op}\times \Gamma
\]
such that each restriction $\mathcal{E}_X:=D^{-1}(\{X\}\times \Gamma)\to \{X\}\times \Gamma$
is a symmetric monoidal $\infty$-category equivalent to
$\mathbb{M}^\otimes_\infty(X)$.
Let $P:\overline{\CAlg}(\mathcal{E})\to I^{op}$ be a map of simplicial sets
defined as follows.
For $q:K\to I^{op}$, the set of $K\to \overline{\CAlg}(\mathcal{E})$
over $q$ is defined to be the set of maps $K\times \Gamma\to \mathcal{E}$
extending $q\times \textup{id}:K\times \Gamma\to I^{op}\times \Gamma$.
Namely, it is $\Fun(\Gamma,\mathcal{E})\times_{\Fun(\Gamma,I^{op}\times \Gamma)}I^{op}\stackrel{\textup{pr}_2}{\to} I^{op}$
where $I^{op}\to \Fun(\Gamma,I^{op}\times \Gamma)$
is induced by the identity of $I^{op}\times \Gamma$.
By the stability property \cite[3.1.2.1 (1), 2.4.2.3. (2)]{HTT}
of coCartesian fibrations,
$\overline{\CAlg}(\mathcal{E})\to I^{op}$ is a coCartesian fibration.
Let $\CAlg(\mathcal{E})$
be the largest subcomplex of $\overline{\CAlg}(\mathcal{E})$
that consists of those vertices $v\in \overline{\CAlg}(\mathcal{E})$ such that
$\{P(v)\}\times \Gamma \to \mathcal{E}$
determines a commutative algebra object of $\mathcal{E}_{P(v)}$.
According to \cite[3.1.2.1 (2)]{HTT}
the induced map
$\CAlg(\mathcal{E})\to I^{op}$ is also a coCartesian fibration.
Note that by the construction,
for each $X$ in $I$ the fiber over $X$ is $\CAlg(\mathcal{E}_X)\simeq \CAlg(\mathbb{M}_\infty^\otimes(X))$, and for each $\phi:Y\to X$
in $I$ the induced map $\CAlg(\mathbb{M}_\infty^\otimes(X))\to \CAlg(\mathbb{M}_\infty^\otimes(Y))$
is equivalent to the pullback functor $\phi^*$.
Each $(\mathbb{M}(X)^c,W_X^c)$ admits
a symmetric monoidal functor $(\mathbb{M}^c(\star),W_{\star}^c)\to (\mathbb{M}(X)^c,W_X^c)$ induced by the morphism $X\to \star$, which preserves weak equivalences.
If $d_\star:I^{op}\to \CAlg(\wCat)$ denotes the constant
functor taking value $\mathbb{M}_\infty^\otimes(\star)$,
it gives rise to a natural transformation $d_\star \to d$.
By using the relative nerve functor as above,
one has a map between coCartesian
fibrations
\[
\xymatrix{
I^{op} \times \mathbb{M}_\infty^\otimes(\star) \ar[rr]^{F^*_\circ} \ar[rd]_{\textup{id}\times e}& & \mathcal{E} \ar[ld]^{D} \\
 & I^{op}\times \Gamma
}
\]
where $e:\mathbb{M}_\infty^\otimes(\star)\to \Gamma$ is a coCartesian fibration
that determines the symmetric monoidal $\infty$-category $\mathbb{M}_\infty^\otimes(\star)$.
The horizontal map preserves coCartesian edges.
Apply the same construction of $\CAlg(\mathcal{E})\to I^{op}$
to $I^{op} \times \mathbb{M}_\infty^\otimes(\star)\to I^{op}\times \Gamma$,
we obtain the constant coCartesian
fibration $I^{op}\times \CAlg(\mathbb{M}_\infty^\otimes(\star))\to I^{op}$
and a map of coCartesian fibrations
\[
F^*:I^{op}\times \CAlg(\mathbb{M}_\infty^\otimes(\star)) \to \CAlg(\mathcal{E})
\]
over $I^{op}$.
For each $f:X\to \star$ in $I$, the fiber $\CAlg(\mathbb{M}_\infty^\otimes(\star)) \to \CAlg(\mathcal{E}_X)\simeq \CAlg(\mathbb{M}_\infty^\otimes(X))$ over $X$ is equivalent to $f^*$.
Thus, each fiber admits the right adjoint functor $f_*:\CAlg(\mathbb{M}_\infty^\otimes(X)) \to \CAlg(\mathbb{M}_\infty^\otimes(\star))$.
In addition, $F^*$ preserves coCartesian edges.
Therefore by the relative adjoint functor theorem \cite[7.3.2.6]{HA}
there is a relative right adjoint $F_*:\CAlg(\mathcal{E})\to I^{op}\times \CAlg(\mathbb{M}_\infty^\otimes(\star))$ over $I^{op}$.
(We refer to \cite[7.3.2]{HA} for the notion of relative adjoint functor.)
For each $f:X\to \star$,
the fiber $\CAlg(\mathbb{M}_\infty^\otimes(X))\to \CAlg(\mathbb{M}_\infty^\otimes(\star))$ is equivalent to
$f_*$.

Now we define a functorial assignment $X\mapsto f_*(\uni_{\mathbb{M}(X)})$
where $\uni_{\mathbb{M}(X)}$ is a unit of $\mathbb{M}(X)$
and $f$ is the natural morphism $X\to \star$.
We let $\iota:I^{op}\to \CAlg(\mathbb{M}_\infty^\otimes(\star))$
be the constant functor whose value is the unit $\uni_\star$
of $\mathbb{M}_\infty^\otimes(\star)$.
It yields a section $\textup{id}\times \iota:I^{op}\to I^{op}\times \CAlg(\mathbb{M}_\infty^\otimes(\star))$.
Composing it with $F^*$,
 we obtain a section $S:I^{op}\to \CAlg(\mathcal{E})$ of $\CAlg(\mathcal{E})\to I^{op}$ which carries $X$ to a unit in $\CAlg(\mathcal{E}_X)\simeq \CAlg(\mathbb{M}_\infty^\otimes(X))$ (every edge in $I^{op}$ maps
to a canonical coCartesian edge).
We define $I^{op}\to \CAlg(\mathbb{M}_\infty^\otimes(\star))$
to be the composite
\[
\Phi:I^{op}\stackrel{S}{\to} \CAlg(\mathcal{E})\stackrel{F_*}{\to} I^{op}\times \CAlg(\mathbb{M}_\infty^\otimes(\star)) \stackrel{\textup{pr}_2}{\to} \CAlg(\mathbb{M}_\infty^\otimes(\star)).
\]

\begin{Remark}
\label{sectionadjoint}
We give a little bit more conceptual explanation of $\Phi$.
Let $\mathcal{C}\to \mathcal{O}$ and $\mathcal{D}\to \mathcal{O}$
be categorical fibrations over an $\infty$-category $\OO$. Let
$\alpha:\mathcal{C}\rightleftarrows \mathcal{D}:\beta$ be functors over $\OO$.
Suppose that $\alpha$ is a left adjoint to $\beta$.
Observe that compositions with $\alpha$ and $\beta$ induce an adjoint pair
between functor categories
$\Fun(\OO,\mathcal{C})\rightleftarrows \Fun(\OO,\mathcal{D})$. To see this,
if $\mathcal{M}\to \Delta^1$ is both a coCartesian fibration and a Cartesian
fibration
which represents the adjoint pair $(\alpha,\beta)$ (cf. \cite[5.5.2.1]{HTT}),
the projection
$\Fun(\OO,\mathcal{M})\times_{\Fun(\OO,\Delta^1)}\Delta^1\to \Delta^1$
is both a coCartesian fibration and a Cartesin fibration
that induces an adjoint pair between functor categories,
where $\Delta^1\to \Fun(\OO,\Delta^1)$ is determined by the projection
$\OO\times \Delta^1\to \Delta^1$.
Suppose further that $\alpha$ is a left adjoint to $\beta$ relative to
$\OO$ (cf. \cite[7.3.2.2]{HA}).
The restriction of the above adjunction induces
\[
\textup{Sect}(\alpha):\textup{Sect}_{\OO}(\mathcal{C}):=\Fun_{\OO}(\OO,\mathcal{C}) \rightleftarrows \Fun_{\OO}(\OO,\mathcal{D})=\textup{Sect}_{\OO}(\mathcal{D}):\textup{Sect}(\beta).
\]
We deduce from \cite[7.3.2.5]{HA} that this pair is an adjunction.
We now apply this to 
\[
F^*:I^{op}\times \CAlg(\mathbb{M}_\infty^\otimes(\star)) \rightleftarrows \CAlg(\mathcal{E}):F_*
\]
over $I^{op}$.
We then have the induced adjunction
\[
\textup{Sect}(F^*):\Fun(I^{op},\CAlg(\mathbb{M}_\infty^\otimes(\star))) \simeq \textup{Sect}_{I^{op}}(I^{op}\times \CAlg(\mathbb{M}_\infty^\otimes(\star))) \rightleftarrows \textup{Sect}_{I^{op}}(\CAlg(\mathcal{E})):\textup{Sect}(F_*).
\]
If $\iota\in \Fun(I^{op},\CAlg(\mathbb{M}_\infty^\otimes(\star)))$
is the constant functor with value $\uni_\star$,
the unit transformation $\textup{id}\to \textup{Sect}(F_*)\circ \textup{Sect}(F^*)$ induces $\iota \to \textup{Sect}(F_*)\circ \textup{Sect}(F^*)(\iota)=\Phi$.
\end{Remark}

\begin{Example}
\label{motivicex}
Let $I=\Sm_k$ and $\star=\Spec k$.
Let $\mathbb{M}(X)=Sp_{\textup{Tate}}(X)$.
We define
\[
\Xi:\Sm_k^{op}\to \CAlg(\DM^\otimes(k))
\]
to be $\Phi$.
Unfolding our construction we see that $\Xi$ carries $X$ to $M_X$, and $\phi:Y\to X$ maps to
$\phi^*:M_X\to M_Y$.
\end{Example}

\begin{Remark}
Let $I=\Sm_k$ and $\star=\Spec k$.
Let $\mathbb{M}(X)=\Comp(\mathcal{N}^{tr}(X))$.
In this case, the above construction also works.
But we will not consider this setting: $f_*(\uni_X)$ is not an appropriate object
we want to consider (for example, Theorem~\ref{real1} does not hold).
\end{Remark}

\begin{Example}
\label{topcomplex}
Let $I$ be the category $\textup{Sch}$ of
separated and quasi-compact schemes.
For any $X$ in $\textup{Sch}$, we let $\Comp(X)$ be
the symmetric monoidal category of (possibly unbounded) cochain
complexes of quasi-coherent sheaves on $X$.
According to \cite[Example 2.3, 3.1, 3.2]{CD1}, there is a symmetric monoidal
model structure on $\Comp(X)$ such that weak equivalences are
quasi-isomorphisms, and for any $Y\to X$ in $\textup{Sch}$
the pullback functor $\Comp(X)\to \Comp(Y)$ is a
left Quillen functor.
Put $\Comp(X)=\mathbb{M}(X)$.
One can apply to this setting our construction and obtain
$\textup{Sch}^{op} \to \CAlg(\mathbb{M}_\infty^\otimes(\Spec \ZZ))$.
\end{Example}

Next we define a functor $\Sm_k\to \DM(k)$ which carries $X$ to $M(X)$.
In some sense, the construction is the dual of that of $\Xi$ and is easier.
We continue to work with the family $\{\mathbb{M}(X)\}$.
Assume that for each $f:X\to \star$ in $I$,
$f^*:\mathbb{M}(\star)\to \mathbb{M}(X)$ is also right Quillen functor
(therefore, it preserves arbitrary weak equivalences).
We denote by $f_\sharp:\mathbb{M}(X)\to \mathbb{M}(\star)$
the left adjoint.
Applying the ``dual version'' of the relative nerve functor or the unstraightning functor
to $X\mapsto \mathbb{M}_\infty(X)$, we obtain a Cartesian fibration $\mathcal{F}\to I$. For each $X\in I$, its fiber is equivalent to $\mathbb{M}_\infty(X)$.
Notice that it is not a coCartesian fibration but a Cartesian fibration.
As in the case of $\mathcal{E}\to I^{op}$, the
natural pullback functors $\mathbb{M}(\star)\to \mathbb{M}(X)$
induce a morphism of Cartesian fibrations
\[
\xymatrix{
\mathcal{F}  \ar[rd] & & I \times \mathbb{M}_\infty(\star) \ar[ld] \ar[ll]_{G^*} \\
 & I.
}
\]
where $I \times \mathbb{M}_\infty(\star)\to I$ is the projection
that is regarded as a Cartesian fibration corresponding to the constant functor
$I\to \wCat$ with value $\mathbb{M}_\infty(\star)$.
Each fiber of the horizontal map over $X\in I$ is equivalent to $f^*$
where $f:X\to \star$ is the natural morphism.
Therefore it admits a left
adjoint functor $f_\sharp:\mathbb{M}_\infty(X)\to \mathbb{M}_\infty(\star)$.
Moreover, $G^*$ preserves Cartesian edges.
Thus, by the relative adjoint functor theorem \cite[7.3.2.6]{HA}
there is a left adjoint $G_\sharp:\mathcal{F}\to I \times \mathbb{M}_\infty(\star)$ relative to $I$. (Its fiber over $X\in I$ is equivalent to $f_\sharp$.)
Let $u:I\to I\times \mathbb{M}_\infty(\star)$ be the functor
determined by the identity $I\to I$ and the constant functor
$I\to \mathbb{M}_\infty(\star)$ taking the value $\uni_\star$.
Then $\Psi:I\to \mathbb{M}_\infty(\star)$ is defined to be the composite
\[
I\stackrel{u}{\to} I\times \mathbb{M}_\infty(\star)\stackrel{G^*}{\to}\mathcal{F}\stackrel{G_\sharp}{\to} I\times \mathbb{M}_\infty(\star)\stackrel{\operatorname{pr}_2}{\to}\mathbb{M}_\infty(\star).
\]

\begin{Example}
\label{motivicex2}
Let $I=\Sm_k$ and $\star=\Spec k$.
Let $\mathbb{M}(X)=Sp_{\textup{Tate}}(X)$.
We define $M(-):\Sm_k \to \DM(k)$
to be $\Psi$.
By our construction, it sends $X$ to an object equivalent to
$M(X)$.
\end{Example}

We define a functor
$\HOM_{\DM(X)}(-,\uni_X):\DM(X)^{op}\to \DM(X)$
as follows. We let
\[
\HOM_{Sp_{\textup{Tate}}(X)}(-,\uni_X'):(Sp_{\textup{Tate}}(X)^{c})^{op} \to Sp_{\textup{Tate}}(X)
\]
be the functor given by $M\mapsto \HOM_{Sp_{\textup{Tate}}(X)}(M,\uni_X')$, where $\HOM_{Sp_{\textup{Tate}}(X)}(-,-)$ denotes the
internal Hom object in $Sp_{\textup{Tate}}(X)$, and
$\uni_X'$ is a fibrant model of the unit $\uni_X$.
By the axiom of symmetric monoidal model category,
the functor $\HOM_{Sp_{\textup{Tate}}(X)}(-,\uni_X')$
preserves weak equivalences. We define 
$\HOM_{\DM(X)}(-,\uni_X):\DM(X)^{op}\to \DM(X)$
to be the functor obtained from
$\HOM_{Sp_{\textup{Tate}}(X)}(-,\uni_X')$
by inverting weak equivalences.

\vspace{2mm}

{\it Proof of Proposition~\ref{cohmotalg}.}
We have constructed the functor $\Xi:\Sm_k^{op}\to \CAlg(\DM^\otimes(k))$
and $M(-):\Sm_k\to \DM(k)$
in Example~\ref{motivicex} and~\ref{motivicex2}.
For simplicity, we write $\Xi$ also for the composite $\Sm_k^{op}\stackrel{\Xi}{\to} \CAlg(\DM^\otimes(k))\to \DM(k)$.
We first observe that for $f:X\to \Spec k$ in $\Sm_k$
there is a canonical equivalence $M(X)^\vee\stackrel{\sim}{\to} M_X=f_*(\uni_X)$. Actually, this equivalence follows from the equivalences of mapping spaces
\begin{eqnarray*}
\Map_{\DM(k)}(M,\HOM_{\DM(k)}(f_\sharp\uni_X,\uni_k)) &\simeq& \Map_{\DM(k)}(M\otimes f_\sharp\uni_X,\uni_k) \\
&\simeq& \Map_{\DM(k)}(f_\sharp(\uni_X),\HOM_{\DM(k)}(M,\uni_k)) \\
&\simeq& \Map_{\DM(X)}(\uni_X,f^*\HOM_{\DM(k)}(M,\uni_k)) \\
&\simeq& \Map_{\DM(X)}(\uni_X,\HOM_{\DM(X)}(f^*(M),f^*(\uni_k))) \\
&\simeq& \Map_{\DM(X)}(f^*(M),\uni_X) \\
&\simeq& \Map_{\DM(k)}(M,f_*(\uni_X))
\end{eqnarray*}
for any $M\in \DM(k)$.
We here used adjunctions $(f_\sharp,f^*)$, $(f^*,f_*)$ and
$f^*\HOM_{\DM(k)}(M,\uni_k)\simeq \HOM_{\DM(X)}(f^*(M),f^*(\uni_k))$.
If we take $M=\HOM_{\DM(k)}(f_\sharp\uni_X,\uni_k)=M(X)^\vee$, then
the identity of $M$ corresponds to $M(X)^\vee\stackrel{\sim}{\to} f_*(\uni_X)$.
The equivalence $M(X)^\vee=f_\sharp(\uni_X)^\vee\to f_*(\uni_X)$
comes from the dual of $\uni_X\to f^*f_\sharp(\uni_X)$:
\[
f^*(f_\sharp(\uni_X)^\vee)\simeq (f^*f_\sharp(\uni_X))^\vee  \to \uni_X
\]
and the composition with
$f_\sharp(\uni_X)^\vee\to f_*f^*(f_\sharp(\uni_X)^\vee)$
where $(-)^\vee$ denotes
the weak dual, that is,
$\HOM_{\DM(-)}(-,\uni_{(-)})$.
By the functoriality of adjoint maps,
it is easy to check that $M(X)^\vee=f_\sharp(\uni_X)^\vee\to f_*(\uni_X)$
is functorial with respect to $X\in \Sm_k$ at the level of
homotopy category, namely,
the functor $\Sm_k^{op}\stackrel{M(-)^{\vee}}{\to}\DM(k)\to \hhh(\DM(k))$
is naturally equivalent to
$\Sm_{k}^{op}\stackrel{\Xi}{\to} \DM(k)\to \hhh(\DM(k))$.
\QED

\begin{Remark}
There should be several approaches to a generalization to singular varieties.
One possible way is to use cubical hyperresolutions of singular varieties \cite{GNPP}
when $k$ is a field of characteristic zero,
so that Hironaka's resolution of singularities is available.
Another one is to adopt a formalism, which works for singular varieties,
such as Beilinson motives \cite{CDT}
when the coefficient ring $K$ is a field of characteristic zero.
But we will not treat singular schemes in this paper.
\end{Remark}

\subsection{}
\label{somerem}
We give some remarks about properties of cohomological motivic algebras.

\begin{Remark}
Since $M_X$ is the weak dual $\Hom_{\DM(k)}(M(X),\uni_k)$ of $M(X)$,
one can observe that $X\mapsto M_X$ satisfies $\mathbb{A}^1$-homotopy
invariance and Nisnevich descent property. Namely, for any
projection $X\times \mathbb{A}^1\to X$ with fiber of the affine line $\mathbb{A}^1=\Spec k[x]$, $M_{X}\to M_{X\times \mathbb{A}^1}$ is an equivalence in $\CAlg(\DM^\otimes(k))$.
For any pullback diagram 
\[
\xymatrix{
V\simeq U\times_{X}Y \ar[r] \ar[d] & Y \ar[d]^{f} \\
U \ar[r]^{j} & X
}
\]
in $\Sm_k$ such that $f$ is \'etale, $j$ is an open immersion and
$(Y\backslash V)_{red} \to  (X\backslash U)_{red}$
is an isomorphism, the induced morphism
$M_X \to M_U\times_{M_V}M_Y$ is an equivalence in $\CAlg(\DM^\otimes(k))$.
\end{Remark}

The following is the Kunneth formula for cohomological motivic algebras.

\begin{Proposition}
\label{kunneth}
There exist a canonical equivalence
$\uni_k\simeq \Xi(\Spec k)$.
Suppose that $X$ and $Y$ are projective and smooth over $\Spec k$.
Then there exists
a canonical equivalence $\Xi(X)\otimes \Xi(Y)=M_X\otimes M_Y \simeq M_{X\times Y}=\Xi(X\times Y)$.
\end{Proposition}

\Proof
The first assertion is obvious.
Next we prove the second assertion.
Consider the Cartesian diagram
\[
\xymatrix{
X\times Y \ar[r]^q \ar[d]_p & Y \ar[d]^g \\
X \ar[r]_f & \Spec k.
}
\]
We will prove that
$p^*\otimes q^*:M_{X}\otimes M_Y\to M_{X\times Y}$
induced by
$p^*:M_X\to M_{X\times Y}$ and $q^*:M_Y\to M_{X\times Y}$
is an equivalence.
For this purpose,
we apply the projection formula for the smooth proper morphism $f$
and the base change theorem for smooth proper morphism $g$
\cite[Theorem 1]{CDT}: we have the sequence of morphisms
induced by unit maps and counit maps of adjunctions
\begin{eqnarray*}
f_*(\uni_X)\otimes g_*(\uni_Y)&\to& f_*f^*(f_*(\uni_X)\otimes g_*(\uni_Y)) \simeq f_*(f^*f_*(\uni_X)\otimes f^*g_*g^*(\uni_k))   \\
&\to& f_*(\uni_X\otimes f^*g_*g^*(\uni_k)) \simeq f_*(f^*g_*g^*(\uni_k)) \\
&\to& f_*(p_*p^*f^*g_*g^*(\uni_k)) \simeq f_*(p_*q^*g^*g_*g^*(\uni_k)) \\
&\to& f_*(p_*q^*g^*(\uni_k)) \simeq f_*p_*(\uni_{X\times Y}) \\
\end{eqnarray*}
whose composite is an equivalence since the projection formula and the base change
theorem imply that the above sequence induces
$f_*(\uni_X)\otimes g_*(\uni_X) \simeq  f_*(\uni_X\otimes f^*g_*(\uni_Y))$ and $f^*g_*(\uni_Y)\simeq p_*q^*(\uni_Y)$. It will suffice to check
that this composite coincides with $p^*\otimes q^*$.
It is straightforward to verify that
\[
f^*(\uni_k)\to f^*g_*g^*(\uni_k)\to p_*p^*f^*g_*g^*(\uni_k)=p_*q^*g^*g_*g^*(\uni_k)\to p_*q^*g^*(\uni_k)=p_*p^*f^*(\uni_k)
\]
is equivalent to $f^*(\uni_k)\to p_*p^*f^*(\uni_k)$
induced by the unit map $\textup{id} \to p_*p^*$.
Then we see that $f_*(\uni_X)\otimes \uni_k\to f_*(\uni_X)\otimes g_*(\uni_X)\stackrel{\sim}{\to}f_*p_*p^*(\uni_{X})$ is equivalent to $p^*:M_X=f_*(\uni_X)\to f_*p_*p^*(\uni_X)=M_{X\times Y}$.
Similarly,
$\uni_k\otimes g_*(\uni_Y) \to f_*(\uni_X)\otimes g_*(\uni_X)\stackrel{\sim}{\to}f_*p_*p^*(\uni_{X})$ is equivalent to $q^*:M_Y\to M_{X\times Y}$.
Thus, $p^*\otimes q^*:M_X\otimes M_Y\to M_{X\times Y}$ is an equivalence.
\QED

\subsection{}
\label{someex}
We will study various objects in $\CAlg(\DM^\otimes(k))$
other than $M_X$:

\begin{Example}
\label{motpath}
Let $X\in \Sm_k$.
Let $x:Y=\Spec k\to X$ and $y:Z=\Spec k\to X$ be two $k$-rational points
on $X$.
Then we have the pushout diagram
\[
\xymatrix{
M_X \ar[r]^{x^*} \ar[d]_{y^*} & M_{\Spec k} \ar[d] \\
M_{\Spec k} \ar[r] & M_{\Spec k}\otimes_{M_X}M_{\Spec k}.
}
\]
in $\CAlg(\DM^\otimes(k))$. Keep in mind that pushouts
in $\CAlg(\DM^\otimes(k))$ do not commute with pushouts in $\DM(k)$
through the forgetful functor. 
By Proposition~\ref{kunneth}, $M_{\Spec k}\simeq \uni_k$.
Thus, 
$M_{\Spec k}\otimes_{M_X}M_{\Spec k}\simeq \uni_k\otimes_{M_X}\uni_k$ in $\CAlg(\DM^\otimes(k))$. We call
$P_{X}(x,y):=\uni_k\otimes_{M_X}\uni_k$ the {\it motivic algebra of path torsors} from $x$
to $y$.
\end{Example}

\begin{Example}
\label{motloop}
Consider $M_X\otimes_{M_{X}\otimes M_X}M_X$.
Note that $\CAlg(\DM^\otimes(k))$ is presentable, and thus
$\CAlg(\DM^\otimes(k))$ is tensored over $\SSS$.
There is a canonical equivalence $S^1\otimes M_X\simeq M_X\otimes_{M_{X}\otimes M_X}M_X$
where $S^1$ is the circle which belongs to $\SSS$.
Thus, by the functoriality of the tensor operation,
$\Map_{\SSS}(S^1,S^1)\simeq S^1$ naturally acts on
$S^1\otimes M_X$ (it is a version of Connes operator; the precise formulation is left to the reader).
We refer to $HHM_X:=M_X\otimes_{M_{X}\otimes M_X}M_X$ as the {\it motivic algebra of free loop space} of $X$.
\end{Example}

\subsection{}
\label{hopfsection}

In Example~\ref{motpath}, if one supposes $x=y$, then
$P_X(x,y)$ has an additional structure.
The augmentation $M_X\to \uni_k\simeq M_{\Spec k}$,
induced by $x:\Spec k\to X$, gives rise to
\[
\uni_k\otimes_{M_X} \uni_k\simeq \uni_k\otimes_{M_X}M_X\otimes_{M_X} \uni_k\to \uni_k\otimes_{M_X}\uni_k\otimes_{M_X} \uni_k\simeq (\uni_k\otimes_{M_X}\uni_k)\otimes (\uni_k\otimes_{M_X}\uni_k) 
\]
and $\uni_k\otimes_{M_X}\uni_k\to \uni_k\otimes_{\uni_{k}}\uni_k\simeq \uni_k$
in $\CAlg(\DM^\otimes(k))$. There is also the flip
$\uni_k\otimes_{M_X}\uni_k\simeq \uni_k\otimes_{M_X}\uni_k$.
Informally, these data define a structure of a cogroup object
on $\uni_k\otimes_{M_X}\uni_k$ in $\CAlg(\DM^\otimes(k))$.
Here $\CAlg(\DM^\otimes(k))$ is endowed with the coCartesian monoidal
structure given by coproducts.
The precise formulation of this structure is as follows.
Let $\Delta_+$ be the category of (possibly nonempty)
linearly ordered finite sets.
Objects are the empty set $[-1]$, $[0]=\{0\}, [1]=\{0,1\},[2]=\{0, 1, 2\},\ldots$. Note that $\Delta_+$ without $[-1]$ is $\Delta$.
The morphism $M_X\to \uni_k$ is described by
$\NNNN(\Delta_+^{\le0})=\NNNN(\{[-1]\to [0]\}) \to \CAlg(\DM^\otimes(k))$.
Since $\CAlg(\DM^\otimes(k))$ has small colimits (in fact, presentable),
the map $\NNNN(\Delta_+^{\le0})\to \CAlg(\DM^\otimes(k))$ admits
a left Kan extension
$\NNNN(\Delta_+)\to \CAlg(\DM^\otimes(k))$.
Namely,
\[
\mathcal{G}(X,x):\NNNN(\Delta_+)^{op}\to \CAlg(\DM^\otimes(k))^{op}
\]
is the Cech nerve associated to
$\NNNN(\{[-1]\to [0]\})^{op}\to\CAlg(\DM^{\otimes}(k))^{op}$ (cf. \cite[6.1.2.11]{HTT}). Its evaluation of $\mathcal{G}(X,x)$ at $[1]$ is equivalent to 
$\uni_k\otimes_{M_X}\uni_k$.
The restriction $\NNNN(\Delta)^{op}\to \CAlg(\DM^\otimes(k))^{op}$
is a group object of $\CAlg(\DM^\otimes(k))^{op}$ (i.e.,
a cogroup object in $\CAlg(\DM^\otimes(k))$).
Namely, it determines
a group structure on $\uni_k\otimes_{M_X}\uni_k$
in $\CAlg(\DM^\otimes(k))^{op}$.
We refer to e.g. \cite[7.2.2.1]{HTT} for the notion of group objects.

Next we define an iterated generalization of $\mathcal{G}(X,x)$.
Consider the restriction $\NNNN(\{[1]\to [0]\})^{op} \subset \NNNN(\Delta_+)^{op} \to \CAlg(\DM^\otimes(k))^{op}$
of the above Cech nerve
$\mathcal{G}^{(1)}(X,x):=\mathcal{G}(X,x):\NNNN(\Delta_+)^{op}\to \CAlg(\DM^\otimes(k))^{op}$.
There is a unique isomorphism
$\NNNN(\Delta_+^{\le0})\simeq \NNNN(\{[1]\to [0]\})$.
Consider the composite $\NNNN(\Delta_+^{\le0})^{op} \simeq \NNNN(\{[1]\to [0]\})^{op} \subset \NNNN(\Delta_+)^{op} \to \CAlg(\DM^\otimes(k))^{op}$.
Once again,
take a rigth Kan extension
$\mathcal{G}^{(2)}(X,x):\NNNN(\Delta_+)^{op}\to \CAlg(\DM^\otimes(k))^{op}$
of this composite.
Repeating this process $n$ times, we obtain
\[
\mathcal{G}^{(n+1)}(X,x):\NNNN(\Delta_+)^{op}\to \CAlg(\DM^\otimes(k))^{op}.
\]
By abuse of notation, we write $\mathcal{G}^{(n+1)}(X,x)$ for the
group object defined as the restriction $\NNNN(\Delta)^{op}\subset \NNNN(\Delta_+)^{op}\to \CAlg(\DM^\otimes(k))^{op}$. (Moreover, one can endow $\mathcal{G}^{(n+1)}(X,x)$ with a structure of an $E_{n+1}$-monoid, but we will not use this enhanced structure.)

\section{Realized motivic rational homotopy type}
\label{realization}

We will consider the realizations of algebraic objects
that appears in Section~\ref{CMAsection} such as $M_X$.
The coefficient field $K$ is
a field of characteristic zero.

\subsection{}
There are several mixed Weil cohomology theories: singular (Betti)
cohomology, \'etale cohomology, analytic de Rham cohomology, 
algebraic de Rham cohomology, rigid cohomology, etc
(see \cite[17.2]{CDT} for mixed Weil cohomology).
To a mixed Weil cohomology theory $E$ with coefficient field $K$,
one can associate a symmetric monoidal colimit-preserving functor
\[
\mathsf{R}_E:\DM^\otimes(k)\longrightarrow \mathsf{D}^\otimes (K)
\]
(see \cite[Section 5]{Tan} for details of the construction in the $\infty$-categorical setting) which is called the {\it realization functor} associated to $E$.
Here $\mathsf{D}^\otimes(K)$
is the derived $\infty$-category of $K$-vector
spaces (see Section~\ref{convention}).
By the relative adjoint functor theorem \cite[7.3.2.6, 7.3.2.13]{HA},
the realization functor $\mathsf{R}_E$ induces an adjunction
\[
\CAlg(\mathsf{R}_E):\CAlg(\DM^\otimes (k))\rightleftarrows \CAlg(\mathsf{D}^\otimes (K))\simeq \CAlg_{K}:\mathsf{M}_E
\]
where $\CAlg(\mathsf{R}_E)$ is the functor induced by $\mathsf{R}_E$,
and $\mathsf{M}_E$ is a right adjoint.
We shall refer to $\CAlg(\mathsf{R}_E):\CAlg(\DM^\otimes (k))\to \CAlg_{K}$
as the {\it multiplicative realization functor}.

In this Section, we consider
the realization functor associated to singular cohomology theory:
\[
\mathsf{R}:\DM^\otimes(k) \longrightarrow \mathsf{D}^\otimes(\QQ).
\]
We here suppose that the base field $k$ is embedded into the complex number
field $\CC$, and the coefficient field $K$ is $\QQ$. 
This functor sends the object $M(X)$ to a complex $\mathsf{R}(M(X))$
that is quasi-isomorphic to the singular chain complex $C_*(X^t,\QQ)$
with rational coefficients. Here $X^t$ stands for the underlying
topological space
of the complex manifold $X\times_{\Spec k}\Spec \CC$.
For ease of notation, when no confusion is likely to arise,
we often write $\mathsf{R}$ for the
multiplicative realization functor $\CAlg(\mathsf{R}):\CAlg(\DM^\otimes(k)) \rightarrow \CAlg_{\QQ}$.

\subsection{}

There are several algebraic models that describe rational homotopy types
of topological spaces.
Quillen \cite{Q} uses {\it differential graded (dg) Lie algebras} whereas
Sullivan \cite{S} adopts {\it commutative differential graded (dg) algebras} as models.
Two approaches are related via Koszul duality between dg Lie
algebras and (augmented) commutative dg algebras.
In this paper, we use cochain
algebras of {\it polynomial differential forms} introduced by Sullivan
as algebraic models of the rational homotopy types of topological spaces.

Let us recall the definition of a cochain algebra of 
polynomial differential forms on a topological space $S$,
see \cite[Section 10]{FHT} for the comprehensive reference.
For a simplicial set $P$, we let $A_{PL}(P)$ be the commutative differential
graded (dg) algebra with rational
coefficients of polynomial
differential forms. This commutative dg algebra is defined as follows
(but we will not apparently need this explicit definition).
For each $n\ge 0$,
we let $\Omega_n$ be the commutative
dg algebra of ``polynomial differential forms on the standard $n$-simplex'',
that is,
\[
\Omega_n:=\QQ[u_0,\ldots, u_n,du_0,\ldots, du_n]/(\Sigma_{i=0}^{n}u_i-1, \Sigma_{i=0}^n du_i)
\]
where $\QQ[u_0,\ldots u_n,du_0,\ldots, du_n]$
is the free commutative graded algebra generated by $u_0,\ldots, u_n$
and  $du_0,\ldots, du_n$ with cohomological degrees $|u_i|=0$, $|du_i|=1$ for each $i$,
and the differential carries $u_i$ and $du_i$ to $du_i$ and $0$, respectively.
For any map $f:\Delta^n\to \Delta^m$, the pullback morphism
$f^*:\Omega_m\to \Omega_n$ of commutative dg algebras is defined
in a natural way (see e.g. \cite[Section 10 (c)]{FHT}).
An element of $A_{PL}(P)$ of (cohomological) degree $r$
is data that consists of a collection $\{ w_{\alpha}\}$
indexed by the set of all morphisms $\alpha:\Delta^n\to P$ from
standard simplices
such that
\begin{itemize}
\item each $\omega_{\alpha}$ is an element of $\Omega_n$ of degree $r$,

\item $f^*(w_{\beta})=w_{\alpha}$ for any $\alpha:\Delta^n\to P$, $\beta:\Delta^m\to P$,
and $f:\Delta^n\to \Delta^m$ such that $\beta \circ f=\alpha$.

\end{itemize} 
The multiplication is given by $\{ w_{\alpha}\}\cdot\{ w'_{\alpha}\}=\{ w_{\alpha}w'_{\alpha}\}$, and the differential is given by $d\{ w_{\alpha}\}=\{ dw_{\alpha}\}$.
If $\phi:P\to P'$ is a map of simplicial sets and $\{\omega_{\alpha}\}_{\alpha:\Delta^n\to P'}$ is an element of $A_{PL}(P')$, then
$\phi^*\{\omega_\alpha\}$ is defined to be $\{\omega_{\phi\circ \beta}\}_{\beta:\Delta^n\to P}$. It gives rise to a map $\phi^*:A_{PL}(P')\to A_{PL}(P)$
of commutative dg algebras.
Let $T$ be a topological space.
If we write $S_*(T)$ for the singuar simplicial complex
whose $n$-th term is the set of singular $n$-simplices,
the commutative dg algebra $A_{PL}(T)$ is defined to be
$A_{PL}(S_*(T))$.

The assignment $P\mapsto A_{PL}(P)$ gives rise to
a functor $A_{PL}:\sSet \to (\CAlg_{\QQ}^{dg})^{op}$ to the category $\CAlg_{\QQ}^{dg}$
of commutative dg algebras over $\QQ$.
There exists a canonical equivalence
between the $\infty$-category $\SSS$ of spaces
and the $\infty$-category obtained
from $\sSet$ by inverting weak homotopy equivalences (cf. \cite[1.3.4.21]{HA}).
As observed below,
the functor $A_{PL}$
sends a weak homotopy equivalence in $\sSet$
to a quasi-isomorphism in $\CAlg_{\QQ}^{dg}$.
Therefore, $A_{PL}:\sSet \to (\CAlg_{\QQ}^{dg})^{op}$ induces
\[
A_{PL,\infty}:\SSS \longrightarrow \NNNN(\CAlg_{\QQ}^{dg})[W^{-1}]^{op}\simeq \CAlg_{\QQ}^{op}.
\]
For a topological space $T$, we shall denote by $A_{PL,\infty}(T)$
the image of $A_{PL}(T)$ in $\CAlg_{\QQ}$.

First we will describe the induced functor $A_{PL,\infty}:\SSS\to \CAlg_{\QQ}^{op}$ in an intrinsic way.

\begin{Proposition}
\label{PLproperty}
The followings hold:
\begin{enumerate}
\renewcommand{\labelenumi}{(\theenumi)}

\item The functor $A_{PL}:\sSet \to \CAlg_{\QQ}^{dg}$
sends a weak homotopy equivalence in $\sSet$
to a quasi-isomorphism in $\CAlg_{\QQ}^{dg}$,

\item $A_{PL,\infty}(\Delta^0)\simeq \QQ$,

\item $A_{PL,\infty}:\SSS\to \CAlg_{\QQ}^{op}$ preserves small colimits.

\end{enumerate}
\end{Proposition}

\begin{Remark}
\label{PLremark}
The functor $A_{PL,\infty}$ is uniquely determined by the properties
(2) and (3) in Proposition~\ref{PLproperty}.
Let $\Fun^{\textup{L}}(\SSS,\CAlg_{\QQ}^{op})$ be the full subcategory
of $\Fun(\SSS,\CAlg_{\QQ}^{op})$ spanned by those functors that
preserve small colimits.
Then by left Kan extension \cite[5.1.5.6]{HTT},
the map $p:\Delta^0\to \SSS$ with value $\Delta^0$ (i.e. the
contractible space) induces an equivalence
\[
\Fun^{\textup{L}}(\SSS,\CAlg_{\QQ}^{op}) \stackrel{\sim}{\to}  \Fun(\Delta^0,\CAlg_{\QQ}^{op})\simeq \CAlg_{\QQ}^{op}.
\]
Therefore, the colimit-preserving functor $A_{PL,\infty}$ is uniquely
determined by the value $\QQ$ of the contractible space.
Namely, if $u:\Delta^0\to \CAlg_{\QQ}^{op}$ denotes the map
determined by the object $\QQ$ of $\CAlg_{\QQ}$, then
$A_{PL,\infty}:\SSS\to \CAlg_{\QQ}^{op}$
is a left Kan extension of $u:\Delta^0\to \CAlg_{\QQ}^{op}$
along $p:\Delta^0 \to \SSS$.
\end{Remark}

\Proof
We first prove (1).
Let $\CAlg_{\QQ}^{dg} \to \Comp(\QQ)$ be the forgetful functor to the
category $\Comp(\QQ)$ of complexes of $\QQ$-vector spaces.
It is enough to show that the
composite $\sSet\to (\CAlg_{\QQ}^{dg})^{op}\to \Comp(\QQ)^{op}$ 
preserves quasi-isomorphisms.
According to \cite[Theorem 10.9]{FHT}, there is  the zig-zag of quasi-isomorphisms in $\Comp(\QQ)$
\[
C^*(P) \to B(P) \leftarrow A_{PL}(P)
\]
where $C^*(P)$ is the cochain complex associated to a simplicial set $P$ with
rational coefficients, and
$B(P)$ is a certain ``intermediate'' cochain complex
associated to $P$.
These quasi-isomorphisms are functorial in the sense that
for any map $P\to P'$ of simplicial sets,
they commute with $A_{PL}(P')\to A_{PL}(P)$,
$B(P')\to B(P)$, and $C^*(P')\to C^*(P)$.
Thus, it will suffice to observe that
$C^*:\sSet\to \Comp(\QQ)^{op}$ given by $P\mapsto C^*(P)$
sends weak homotopy equivalences to quasi-isomorphisms.
Let $C_*:\textup{Set}_{\Delta}\to \Comp(\QQ)$
be the functor which carries $P$ to the (normalized) chain complex
$C_*(P)$ with rational coefficients.
Since the dual of any quasi-isomorphism $C_*(P)\to C_*(P')$
is a quasi-isomorphism $C^*(P')\to C^*(P)$,
we are reduced to proving that
$C_*$
sends weak homotopy equivalences to quasi-isomorphisms.
Indeed, it is a well-known fact, but we here describe one of the proofs.
Let $\textup{Vect}_{\Delta}$
denote the category of simplicial objects in the category of $\QQ$-vector
spaces, that is, simplicial $\QQ$-vector spaces.
Consider the adjunction
$\QQ[-]:\sSet\rightleftarrows \textup{Vect}_{\Delta}:U$
where $U$ is the forgetful functor, and $\QQ[-]$ is its left adjoint, that is,
the free functor.
Let us consider $\sSet$ as the Quillen model category whose weak equivalences
are weak homotopy equivalence, and whose cofibrations are monomorphisms.
As in the case of simplicial abelian groups, $\textup{Vect}_{\Delta}$
admits a model category structure in which
$f$ is a weak equivalences (resp. a fibration) if $U(f)$ is a weak equivalence
(resp. a Kan fibration).
Then the pair $(\QQ[-],U)$ is a Quillen adjunction.
Let $N:\textup{Vect}_{\Delta}\stackrel{\sim}{\to} \Comp^{\le 0}(\QQ)$ be
 the Dold-Kan equivalence
which carries a simplicial vector space to its normalized chain complex,
where $\Comp^{\le 0}(\QQ)$ is the full subcategory
 of $\Comp(\QQ)$ spanned by those object $C$ such that $H^i(C)=0$ for $i>0$.
The composite $\sSet\stackrel{\QQ[-]}{\to} \textup{Vect}_{\Delta}\stackrel{N}{\to} \Comp^{\le 0}(\QQ)$
is naturally equivalent to the functor $\sSet\to \Comp^{\le 0}(\QQ)$
which sends $P$ to $C_*(P)$.
The functor $\QQ[-]$ preserves weak equivalences
since every object in $\sSet$ is cofibrant, and $N$ sends weak equivalences
to quasi-isomorphims. Thus, $C_*$ sends weak homotopy equivalences to quasi-isomorphisms.

The equality $A_{PL}(\Delta^0)=\QQ$ is clear from the definition (see \cite[Example 1 in
page 124]{FHT}). Hence (2) follows.

Next we prove (3). Note that the
forgetful functor $\CAlg_{\QQ}\to \Mod_{\QQ}\simeq \mathsf{D}(\QQ)$
preserves limits (cf. \cite[3.2.2.4]{HA}).
Thus, it will
suffice to prove that
$\SSS\stackrel{A_{PL,\infty}}{\to} \CAlg_{\QQ}^{op}\to \Mod_\QQ^{op}$
preserves small coimits; a small colimit diagram in $\SSS$ maps
to a limit diagram in $\Mod_\QQ$.
According to \cite[4.4.2.7]{HTT}, $\SSS\to \Mod_\QQ^{op}$
preserves small colimits if and only if it preserves pushouts
and small coproducts.
It is enough to
show that
$\sSet\stackrel{A_{PL}}{\to} (\CAlg_{\QQ}^{dg})^{op} \to \Comp(\QQ)^{op}$
sends homotopy pushout diagrams and homotopy coproduct diagrams
to homotopy pullback diagrams and homotopy product diagrams in $\Comp(\QQ)$,
respectively.
Here
the second functor is the forgetful functor, and $\Comp(\QQ)$
is endowed with the projective model structure, see Section~\ref{convention}.
As discussed in the proof of (1), we may replace this composite
by $C^*:\sSet\to \Comp(\QQ)^{op}$.
We will observe that $C_*:\sSet \to \Comp(\QQ)$ preserves
homotopy colimits.
We equip $\Comp^{\le0}(\QQ)$
with the projective model
structures (cf. \cite[2.3]{Hov}, \cite[4.1]{SS}).
A morphism $p$ in $\Comp^{\le0}(\QQ)$
is a weak equivalence (resp. a fibration)
if it is a quasi-isomorphism
(resp. surjective in cohomologically
negative degrees).
Cofibrations are monomorphisms
(keep in mind that $\QQ$ is a field).
The free functor $\QQ[-]:\sSet\to \Vect_{\Delta}$
is a left Quillen functor.
The normalization functor
$N:\textup{Vect}_{\Delta}\to \Comp^{\le0}(\QQ)$
is a left Quillen functor (see \cite[Section 4]{SS}).
In addition,
$\Comp^{\le0}(\QQ)\hookrightarrow \Comp(\QQ)$ is a left Quillen functor.
Therefore, we deduce that
$C_*:\sSet \to \Comp(\QQ)$ preserves
homotopy colimits.
Note that $C^*$ is composite $\sSet\stackrel{C_*}{\to} \Comp(\QQ)\to \Comp(\QQ)^{op}$ where the second functor
is given by the hom complex $\mathcal{H}om_{\QQ}(-,\QQ)$.
Then $\mathcal{H}om_{\QQ}(-,\QQ):\Comp(\QQ)\to \Comp(\QQ)^{op}$
preserves homotopy colimits, so that the induces
functor $\Mod_{\QQ}\to \Mod_{\QQ}^{op}$ preserves colimits.
(Indeed, it is enough to check that
it preserves homotopy pushouts and homotopy coproducts.
In $\Comp(\QQ)$, every object is both cofibrant and fibrant.
By the explicit presentation of homotopy pushouts/coproducts cf. \cite[A.2.4.4]{HTT},
we easily see that
$\mathcal{H}om_{\QQ}(-,\QQ)$ sends a homotopy pushout (resp. coproduct)
diagram
to a homotopy pullback (resp. coproduct) diagram.)
Consequently, $\SSS\to \Mod_{\QQ}^{op}$ induced by
$C^*$ preserves small colimits.
\QED

\subsection{}
Let us consider the composite
\[
T:\Sm_k^{op}\stackrel{\Xi}{\to} \CAlg(\DM^\otimes(k))\stackrel{\mathsf{R}}{\to} \CAlg_{\QQ}
\]
See Proposition~\ref{cohmotalg} for $\Xi$. We put $T_X=T(X)=\mathsf{R}(M_X)$.

\begin{Theorem}
\label{real1}
Let $X$ be a smooth scheme separated of finite type over $k\subset \CC$.
Let $X^t$ be the underlying topological space of the complex manifold $X\times_{\Spec k}\Spec \CC$.
There is a canonical equivalence $\mathsf{R}(M_X)=T_X\stackrel{\sim}{\to} A_{PL,\infty}(X^{t})$
in $\CAlg_{\QQ}$. 
\end{Theorem}

\Proof
We first introduce some categories.
Let $\mathsf{D}^{\otimes}(X^{t})$ be the symmetric monoidal
presentable $\infty$-category of complexes of sheaves on $\QQ$-vector
spaces on $X^{t}$.
We define this $\infty$-category by the machinery of
model categories as in Example~\ref{topcomplex}.
Let $\Sh(X^t)$ be the Grothendieck abelian category of sheaves of $\QQ$-vector spaces on $X^t$ and let $\Comp(\Sh(X^t))$ be the category of cochain complexes
of $\Sh(X^t)$. It is endowed with the
symmetric monoidal structure by tensor product.
Thanks to \cite[Theorem 2.5, Example 2.3, Proposition 3.2]{CD1}, there is a symmetric monoidal
combinatorial model category structure on $\Comp(\Sh(X^t))$
in which weak equivalences consists of quasi-isomorphisms.
We then obtain the symmetric monoidal
presentable $\infty$-category $\mathsf{D}^{\otimes}(X^{t})$ from $\Comp(\Sh(X^t))^c$ by inverting weak equivalences (cf. Section~\ref{convention}).
By replacing $X^t$ by the one-point space $*$,
we also have a symmetric monoidal combinatorial model category
$\Comp(\Sh(*))$ which coincides with $\Comp(\QQ)$ endowed with
the projective
model structure.
By abuse of notation we denote the associated symmetric monoidal presentable $\infty$-category by $\mathsf{D}^\otimes(\QQ)$.
The canonical map to the one-point space $f^t:X^t\to *$ induces 
the symmetric monoidal pullback functor
$\Comp(\Sh(*))\to \Comp(\Sh(X^t))$ that is a left Quillen functor
\cite[Theorem 2.14]{CD1}.
It gives rise to a symmetric monoidal colimit-preserving pullback
functor 
$f^{t*}:\mathsf{D}(\QQ)\to \mathsf{D}(X^t)$.
According to relative adjoint functor theorem \cite[7.3.2.6]{HA},
there is a right adjoint functor $f_*^t:\mathsf{D}(X^t)\to \mathsf{D}(\QQ)$
which is lax symmetric monoidal.
We then use the Beilinson motives studied by Cisinski-D\'eglise \cite{CDT}.
Let $\mathbb{M}_{B}(X)$ be a symmetric monoidal combinatorial
model category of Beilinson motives over $X$ with rational
coefficients (see \cite[14.2]{CDT})
and let $\DM_B^\otimes(X)$ be the symmetric monoidal presentable
$\infty$-category obtained from $\mathbb{M}_{B}(X)$.
Since $X$ is regular,
according to \cite[16.1.1, 16.1.4]{CDT}
there is a symmetric monoidal equivalence
$\DM_B^\otimes(X)\stackrel{\sim}{\to} \DM^\otimes(X)$
induced by a symmetric monoidal left Quillen functor 
$\mathbb{M}_{B}(X)\to Sp_{\textup{Tate}}(X)$
(by \cite[16.1.4]{CDT} the induced functor between their homotopy categories
is an equivalence, from which the equivalence of 
stable $\infty$-categories follows, see e.g. \cite[Lemma 5.8]{Tan}).
The equivalences $\DM_B^\otimes(X)\simeq \DM^\otimes(X)$
and $\DM_B^\otimes(\Spec k)\simeq \DM^\otimes(k)$
commute with pullback functors.
Let $\mathsf{R}_X:\DM^\otimes(X)\simeq
\DM_B^\otimes(X)\to \mathsf{D}^\otimes(X^t)$
be the (relative) 
realization functor that is a symmetric monoidal functor.
It is obtained from symmetric monoidal functors
of model categories (cf. Section~\ref{convention}):
as explained in \cite[17.1.7]{CDT} that uses the construction of Ayoub,
there is a diagram of symmetric monoidal functors
$\mathbb{M}_{B}(X)\stackrel{r}{\to} \mathbb{M}(X)\stackrel{p}{\leftarrow} \Comp(\Sh(X^t))$
of model categories where $\mathbb{M}(X)$ is an intermediate
symmetric monoidal model category, $r$ is a symmetric monoidal left Quillen
functor, and $p$ induces an equivalence of symmetric monoidal $\infty$-categories.
Similarly, we have the realization functor $\mathsf{R}:\DM^\otimes(k)\simeq \DM_B^\otimes(\Spec k)\to \mathsf{D}^\otimes(\QQ)$ of singular cohomology theory.
The functors $\mathsf{R}$ and $\mathsf{R}_X$ commute with the pullback
functors (because of the construction).
Therefore, we have the diagram
\[
\xymatrix{
\DM^\otimes(X) \ar[r]^{\mathsf{R}_X} \ar@<5pt>[d]^{f_*} & \mathsf{D}^\otimes(X^{t}) \ar[d]_{f_*^{t}} \\
\DM^\otimes(k) \ar[r]_{\mathsf{R}} \ar[u]^{f^*} & \mathsf{D}^\otimes(\QQ) \ar@<-5pt>[u]_{f^{t*}}.
}
\]
with a canonical equivalence
$\mathsf{R}_X\circ f^*\simeq f^{t*}\circ \mathsf{R}$
of symmetric monoidal functors.
Let $A$ be a commutative algebra object in $\DM^\otimes(X)$,
that is, an object of $\CAlg(\DM^\otimes(X))$.
Consider the canonical exchange map $e:\mathsf{R}(f_*(A))\to f_*^{t}(\mathsf{R}_X(A))$ in $\mathsf{D}(\QQ)$.
This map is the composition of morphisms
\begin{eqnarray*}
\RRR(f_*(A)) &\to& f^{t}_*f^{t*}(\RRR(f_*(A)) \\
&\simeq& f^{t}_*\RRR_Xf^*(f_*(A)) \\
&\simeq& f^{t}_*\RRR_X(f^*f_*)(A) \\
&\to& f^{t}_*\RRR_X(A)
\end{eqnarray*}
where the first arrow is induced by the unit map $\textup{id}\to f^{t}_*f^{t*}$, the second arrow is induced by $\RRR_Xf^*\simeq f^{t*}\RRR$, and the fourth one is induced by the counit map $f^*f_*\to \textup{id}$.
The unit map $\textup{id}\to f^{t}_*f^{t*}$
and the counit map $f^*f_*\to \textup{id}$ are promoted to
a unit map and a counit map for
adjunctions
$\CAlg(\mathsf{D}^\otimes(\QQ))\rightleftarrows \CAlg(\mathsf{D}^\otimes(X^t))$ and $\CAlg(\DM^\otimes(k))\rightleftarrows \CAlg(\DM^\otimes(X))$, respectively.
In particular, $e:\mathsf{R}(f_*(A))\to f_*^{t}(\mathsf{R}_X(A))$ is
promoted to
an morphism in $\CAlg(\mathsf{D}^\otimes(\QQ))\simeq \CAlg_{\QQ}$.
By \cite[17.2.18, 4.4.25]{CDT}, if $A$ is compact in the underlying $\infty$-category
$\DM(X)$, $e$ is an equivalence. In particular, if $A=\uni_X$,
we have a canonocal equivalence
$\mathsf{R}(f_*(\uni_X))=\mathsf{R}(M_X)\stackrel{\sim}{\to} f_*^{t}(\mathsf{R}_X(\uni_X))$ in $\CAlg_{\QQ}$.
Consequently, to prove our assertion it suffices to prove that
$f_*^{t}(\uni_{X^t})$ is equivalent to $A_{PL,\infty}(X^{t})$
where $\uni_{X^t}$ is the unit of $\mathsf{D}^\otimes(X^{t})$,
i.e., the constant sheaf with value $\QQ$.
For this purpose, recall first that since $X\times_{\Spec k}\Spec \CC$ is a complex smooth scheme separated of finite type, the underlying topological space $X^t$ is
a hausdorff paracompact smooth manifold.
Therefore, according to \cite[Theorem 5.1]{BT}, it admits 
a good cover $\mathcal{U}=\{U_\lambda\}_{\lambda\in I}$, that is,
an open cover $\mathcal{U}=\{U_\lambda\}_{\lambda\in I}$ such that
every non-empty finite intersection
$U_{\lambda_0}\cap\ldots \cap U_{\lambda_r}$ is contractible.
Take the augmented simplicial diagram of the Cech nerve $U_\bullet \to U_{-1}:=X^{t}$ associated to the cover.
The $n$-th term $U_n$ of $U_\bullet$ is the disjoint union of intersections of $n+1$ open sets in $\mathcal{U}$.
We denote by $j_{U_n}:U_n\to X^t=U_{-1}$ the canonical map.
If we think of $U_\bullet \to X^t$ as an augmented simplicial diagram
in $\mathcal{S}$, then by Dugger-Isaksen
\cite[Theorem 1.1]{DI}, it is a colimit diagram.
According to Proposition~\ref{PLproperty},
the functor $A_{PL,\infty}:\SSS\to \CAlg_{\QQ}^{op}$ commutes with small
colimits. Thus,
the canonical morphism
$ A_{PL,\infty}(X^t)\to \varprojlim_{[n]\in \Delta} A_{PL,\infty}(U_n)$
is an equivalence where $\varprojlim_{[n]\in \Delta} A_{PL,\infty}(U_n)$
is a limit of the cosimplicial diagram in $\CAlg_\QQ$.
Thus, it is enough to show that $\varprojlim_{[n]\in \Delta} A_{PL,\infty}(U_n)\simeq f_*^{t}(\uni_{X^t})$.
For $i\ge -1$, we let $\textup{Comp}(\Sh(U_i))$ be the category of 
complexes of sheaves of $\QQ$-vector spaces on $U_i$.
As in the case of $\mathsf{D}(X^t)$,
by the model structure in \cite[2.3, 2.5]{CD1}
we have a symmetric monoidal presentable $\infty$-category
$\mathsf{D}^\otimes(U_n)$ from $\textup{Comp}(\Sh(U_n))$.
For each morphism $U_n\to U_m$,
a symmetric monoidal colimit-preserving functor
$\mathsf{D}^\otimes(U_m)\to \mathsf{D}^\otimes(U_n)$.
It gives rise to a cosimplicial diagram of symmetric monoidal
$\infty$-categories
which we denote simply by $\mathsf{D}^\otimes(U_\bullet)$.
It also has the natural
coaugmentation $\mathsf{D}^\otimes(X^{t})\to \mathsf{D}^\otimes(U_\bullet)$.
Let $\Gamma(U_n,-):\mathsf{D}(U_n)\to \mathsf{D}(\QQ)$ be
the (derived) global section functor,
that is a lax symmetric monoidal right adjoint functor
to the pullback functor $\mathsf{D}^\otimes(\QQ)\to \mathsf{D}^\otimes(U_i)$ of $U_i\to *$.
We denote by $\uni_{U_n}$ the unit of
$\mathsf{D}(U_n)$ that corresponds to the constant sheaf with value $\QQ$.
Note $f^t_*(-)=\Gamma(U_{-1},-)=\Gamma(X^t,-)$, and
$\Gamma(U_{n},\uni_{U_n})$ in $\mathsf{D}(\QQ)$ is a complex computing the sheaf cohomology of $U_i$ with coefficients in $\QQ$.
Remember that $U_n$ is a disjoint union of contractible
spaces for $n\ge0$.
For each connected component $V$ of $U_n$,
$\Gamma(V,\uni_{U_n}|_{V})$ in $\CAlg_{\QQ}$ is
an initial object of $\CAlg_{\QQ}$, i.e., $\QQ$
since the unit map $\QQ\to \Gamma(V,\uni_{U_n}|_{V})$
is an equivalence in $\mathsf{D}(\QQ)$, cf. \cite[Corollary 3.2.1.9]{HA}.
By Proposition~\ref{PLproperty}, 
the image of a contractible space under $A_{PL,\infty}$ is $\QQ$.
Therefore, $\Gamma(U_n,\uni_{U_n})\in \CAlg_{\QQ}$ is
equivalent to $A_{PL,\infty}(U_n)$, i.e., $\Gamma(U_i,\uni_{U_i})\simeq \prod_{\pi_0(U_n)}\QQ \simeq A_{PL,\infty}(U_n)$
for $i\ge0$ ($\pi_0(-)$ is the set of connected components). We may consider $\{\Gamma(U_i,\uni_{U_i})\}_{[i]\in \Delta}$
to be a cosimplicial diagram of ordinary commutative algebras (arising
from connected components of $U_\bullet$).
We then have $\varprojlim A_{PL,\infty}(U_n)\simeq \varprojlim\Gamma(U_n,\uni_{U_n})$.
It will suffice to prove that
the canonical morphism $\Gamma(X^{t},\uni_{X^{t}})\to \varprojlim_{[n]\in \Delta}\Gamma(U_n,\uni_{U_n})$
in $\mathsf{D}(\QQ)$ is an equivalence (we may and will disregard their commutative algebra structures).
To this end, we use the descent for hypercovers on $X^t$.
Let $j_{U_n!}:\mathsf{D}(U_n)\to \mathsf{D}(X^t)$ be the left adjoint to the restriction $j_{U_n}^*:\mathsf{D}(X^t)\to \mathsf{D}(U_n)$.
According to \cite[Example 2.3, Theorem 2.5]{CD1}, we see that the canonical
morphism $\varinjlim_{[n]\in \Delta^{op}}j_{U_n!}(\uni_{U_n})\to \uni_{X^t}$ is an equivalence in $\mathsf{D}(X^t)$.
For any $F$ in $\mathsf{D}(X^t)$, it induces an equivalence
$\Gamma(X^t,F)\stackrel{\sim}{\to} \varprojlim_{[n]\in \Delta}\Gamma(U_n,F|_{U_n})$.
In particular, we have a canonical equivalence
$\Gamma(X^{t},\uni_{X^{t}})\stackrel{\sim}{\to} \varprojlim_{[n]\in \Delta}\Gamma(U_i,\uni_{U_i})$.
\QED

\begin{Remark}
\label{real1R}
Let $\phi:Y\to X$ be a morphism in $\Sm_k$.
Then $\phi^*:M_X\to M_Y$ induces $\mathsf{R}(\phi^*):\mathsf{R}(M_X)=T_X\to \mathsf{R}(M_X)=T_Y$.
On the other hand, the associated
continuous map $\phi^t:Y^t\to X^t$ of topological
spaces induces $\phi^{t*}:A_{PL,\infty}(X^t)\to A_{PL,\infty}(Y^t)$ induced
by $A_{PL}(X^t)\to A_{PL}(Y^t)$.
The morphism $\mathsf{R}(\phi^*):T_X\to T_Y$ in $\CAlg_{\QQ}$
is equivalent to $\phi^{t*}:A_{PL,\infty}(X^t)\to A_{PL,\infty}(Y^t)$
through equivalences $T_X\simeq A_{PL,\infty}(X^t)$ and $T_Y\simeq A_{PL,\infty}(Y^t)$ in Theorem~\ref{real1}.

To observe this, note first that by the compatibility of the realization
functor with pushforward functors, 
$\mathsf{R}(\phi^*)$ can be identified with
$f_*^t(\uni_{X^t})\to g_*^t(\uni_{Y^t})$ induced by $\phi^t:Y^t\to X^t$
where $g^t:Y^t\to *$ is the canonical map to one point space. 
Let us unfold the equaivalence given in the proof of 
Theorem~\ref{real1}.
As in the proof,
choose a good cover $\mathcal{U}=\{U_\lambda\}_{\lambda\in I}$ of $X^t$
and take the augmented Cech nerve $U_\bullet\to X^t=U_{-1}$.
We know from the proof of Theorem~\ref{real1} that
there are canonical equivalences
$\Gamma(U_n,\uni_{U_n})\stackrel{\sim}{\to} \prod_{\alpha\in \pi_0(U_n)}\Gamma(U_{n,\alpha},\uni_{U_{n,\alpha}})\stackrel{\sim}{\leftarrow}\prod_{\alpha\in \pi_0(U_n)}\QQ$ in $\CAlg_\QQ$
where 
each $U_n$ is a disjoint union $\sqcup_{\alpha\in \pi_0(U_n)}U_{n,\alpha}$ of contractible spaces.
Similarly, we have canonical equivalences
$A_{PL,\infty}(U_n)\stackrel{\sim}{\to} \prod_{\alpha\in \pi_0(U_n)}A_{PL,\infty}(U_{n,\alpha})\stackrel{\sim}{\leftarrow}\prod_{\alpha\in \pi_0(U_n)}\QQ$.
Both cosimplicial
objects $\{\Gamma(U_n,\uni_{U_n})\}_{[n]\in \Delta}$
and $\{A_{PL}(U_n)\}_{[n]\in \Delta}$ are equivalent to
the cosimplicial ordinary
commutative $\QQ$-algebra, regarded as a
cosimplicial object in $\CAlg_\QQ$,
that is defined by the assignment
$[n] \mapsto \prod_{\alpha\in \pi_0(U_n)}\QQ=\QQ^{\pi_0(U_n)}$
such that for any $[n]\to [m]$, $\QQ^{\pi_0(U_n)}\to \QQ^{\pi_0(U_m)}$
is induced by the map $\pi_0(U_m)\to \pi_0(U_n)$ (by
the superscript we mean cotensor). It gives
rise to $A_{PL,\infty}(X^t) \simeq \varprojlim A_{PL,\infty}(U_n)\simeq \varprojlim \Gamma(U_n,\uni_{U_n})\simeq f_*^t(\uni_{X^t})$.
Taking account of these steps,
we are reduced to checking
 a functoriality of good covers: it suffices to verify
that
if $\mathcal{U}=\{U_\lambda\}_{\lambda\in I}$
is a good cover of $X^t$,
then there is a good over $\mathcal{V}=\{V_{\mu}\}_{\mu \in J}$ of $Y^t$
such that any $V_\mu \to Y^t\to X^t$ factors through some $U_{\lambda}\to X^t$. Actually,
it follows from the proof of the existence of a good cover.
See \cite[Corollary 5.2]{BT} and
the discussion after the proof of \cite[Theorem 5.1]{BT}.
\end{Remark}

It is useful to have a smooth de Rham model of $T_X$.
We will describe $T_X\otimes_{\QQ}\RR$ in terms of smooth differential
forms.
By $X_{\infty}$ we mean the underlying differential manifold of
$X\times_{\Spec k}\Spec \CC$.
Let $A_{X_\infty}$ be the commutative dg algebra
of $C^\infty$ real differential forms on $X_{\infty}$.
We call $A_{X_\infty}$ the smooth de Rham algebra
on $X_\infty$
We think of $A_{X_\infty}$ as an object in $\CAlg_{\RR}$.

\begin{Corollary}
\label{real2}
Consider the base change $T_X\otimes_{\QQ}\RR$ which belongs to $\CAlg_{\RR}$.
There is an equivalence $T_X\otimes_{\QQ}\RR\simeq A_{X_\infty}$
in $\CAlg_{\RR}$.
\end{Corollary}

\Proof
There is a zig-zag of quasi-isomorphisms between
$A_{X_\infty}$ and $A_{PL}(X^t)\otimes_{\QQ}\RR$
(see \cite[Theorem 11.4]{FHT}). Thus, by Theorem~\ref{real1}
we see that $T_X\otimes_{\QQ}\RR\simeq A_{X_\infty}$.
\QED

By using Theorem~\ref{real1} and Remark~\ref{real1R},
we can easily prove the following:

\begin{Proposition}
Let $\CAlg(\DM^\otimes(k))\to \CAlg_{\QQ}$
be the multiplicative realization functor.
Then the image of the motivic algebra of path torsor $P(X,x,y)$ (cf. Example~\ref{motpath}) in $\CAlg_{\QQ}$ is equivalent to
the pushout 
$\QQ\otimes_{A_{PL,\infty}(X^t)}\QQ$ associated to two augmentations
$A_{PL,\infty}(X^t)\to \QQ$ and $A_{PL,\infty}(X^t)\to \QQ$ respectively
induced by
points $x$ and $y$ in $X^t$.
(We remark that $\QQ\otimes_{A_{PL,\infty}(X)}\QQ$
can be obtained by a bar construction 
of $A_{PL}(X)$ with two augmentations, see \cite{Ol}.)

The image of the motivic algebra of free loop space $HHM_X$ (cf. Example~\ref{motloop}) in $\CAlg_{\QQ}$ is
$A_{PL,\infty}(X^t)\otimes_{A_{PL,\infty}(X^t)\otimes A_{PL,\infty}(X^t)}A_{PL,\infty}(X^t)\simeq S^1\otimes A_{PL,\infty}(X^t)$.
(It might be worth mentioning that if $X^t$ is simply connected, then $A_{PL,\infty}(X^t)\otimes_{A_{PL,\infty}(X^t)\otimes A_{PL,\infty}(X^t)}A_{PL,\infty}(X^t)$
is equivalent to
$A_{PL,\infty}(LX^t)$ where  $LX^t$ is the free loop space of $X^t$ \cite[Example 1 in page 206]{FHT}.)
\end{Proposition}

\subsection{}
\label{preDAG}
Before proceeding the next subsection, we introduce some algebro-geometric
notions. Let $K$ be a field of characteristic zero.
Let $\CAlg^{\textup{dis}}_K$ be the full subcategory of $\CAlg_K$ that is spanned by discrete objects $C$,
i.e., $H^i(C)=0$ for $i\neq 0$.
Put another way, we let $\Mod_K^{\DIS}$ be the
(symmetric monoidal) full subcategory of $\Mod_K\simeq \mathsf{D}(K)$
spanned by discrete objects $M$,
i.e., $H^i(M)=0$ for $i\neq 0$.
This full subcategory is nothing else but (the nerve of) the category
of $K$-vector spaces.
Then $\CAlg_K^{\textup{dis}}=\CAlg(\Mod^{\DIS}_K)$.
 The $\infty$-category
$\CAlg^{\textup{dis}}_K$ is naturally equivalent
to the nerve of category of ordinary commutative $K$-algebras.
Let $\Aff_K$ be the opposite category of $\CAlg_K$.
We write $\Spec R$ for an object in $\Aff_K$
that corresponds to $R\in \CAlg_K$.
We shall refer to it as a derived affine scheme (or affine scheme) over $K$.
The Yoneda embedding identifies $\Aff_K$ with the full subcategory of
$\Fun(\CAlg_K,\SSS)$. This embedding preserves
small limits.
The functor $\Spec R:\CAlg_K \to \SSS$
corepresented by $R$ satisfies the sheaf condition with respect to flat
topology, see e.g. \cite{DAG}.
We often regard $\Spec R$ as a sheaf $\CAlg_K\to \SSS$.
We remark that in the literature of derived geometry (see e.g. \cite{DAG}
for its $\infty$-categorical theory),
$\Spec R$ with $R\in \CAlg_K$ is usually called a nonconnective
(derived) affine scheme.
Let $\Aff^{\textup{dis}}_K$
be the full subcategory of $\Aff_K$ that
corresponds to $\CAlg^{\textup{dis}}_K$.
One can naturally identify $\Aff^{\textup{dis}}_K$
with the category of ordinary affine schemes over $K$
(keep in mind that the full subcategories
$\Aff_K^{\DIS}$ are not closed under some constructions;
for example, in general,
fiber products in $\Aff_K^{\DIS}$
are not compatible
with those in $\Aff_K$).

For an $\infty$-category $\mathcal{C}$ that has  
finite products, we write $\Grp(\mathcal{C})$ for the $\infty$-category
of group objects in $\mathcal{C}$.
We shall
call a group object in $\Aff_K$ a derived affine group scheme over $K$.
There is a canonical Yoneda embedding $\Grp(\Aff_K)\hookrightarrow \Fun(\CAlg_K,\Grp(\SSS))$. Therefore, through this functor we often
think of a derived affine group scheme
as a sheaf $\CAlg_K\to \Grp(\SSS)$.
Put another way, $\Spec R$ in $\Grp(\Aff_K)$ amounts to
a commutative Hopf algebra object $R$ in $\Mod_K^\otimes$.
See \cite[Appendix A]{Tan} for details.

\subsection{}
\label{realhopfsection}
\begin{Definition}
\label{realhopfdef}
In Section~\ref{hopfsection},
for a pointed smooth variety $(X,x)$ and a natural number $n\ge1$,
we have defined the group
object
$\mathcal{G}^{(n)}(X,x):\NNNN(\Delta^{op})\to \CAlg(\DM^\otimes(k))^{op}$.
Since the multiplicative realization functor
$\CAlg(\mathsf{R}_E):\CAlg(\DM^\otimes(k))\to \CAlg_{K}$
preserves 
coproducts,
we see that
the composite
\[
G_E^{(n)}(X,x):\NNNN(\Delta^{op})\stackrel{\mathcal{G}^{(n)}(X,x)^{op}}{\longrightarrow} \CAlg(\DM^\otimes(k))^{op} \stackrel{\CAlg(\mathsf{R}_E)^{op}}{\longrightarrow} \CAlg_{K}^{op}=\Aff_{K}
\]
is a group object in $\Aff_{K}$.
Namely, $G_E^{(n)}(X,x)$ is a derived affine group scheme over $K$.
If no confusion is likely to arise,
we often write $G^{(n)}(X,x)$ for $G_E^{(n)}(X,x)$.
\end{Definition}

\begin{Proposition}
\label{cechtop}
Suppose that $k$ is embedded in $\CC$ and
consider the case of singular realization $\mathsf{R}=\mathsf{R}_E$.
The point $x$ on $X$ determines a point of
the associated topological space $X^t$
which we denote also by $x$.
Let $\Spec \QQ\to \Spec T_X$ be a morphism induced by $x$.
Then the derived affine group scheme $G(X,x)=G^{(1)}(X,x)$
is equivalent to the Cech nerve obtained from $\Spec \QQ\to \Spec T_X$.
The iterated group scheme $G^{(n)}(X,x)$ ($n\ge2$) also has
a similar description.
\end{Proposition}

\Proof
By Remark~\ref{real1R},
the map $T_X=\mathsf{R}(M_X)\to \QQ=\mathsf{R}(M_{\Spec k})$
induced by $M_X\to M_{\Spec k}=\uni_k$
can be viewed as the map $T_X\to \QQ$ induced by $x\in X^t$.
Remember that
the opposite of the multiplicative realization functor
$\CAlg(\DM^\otimes(k))^{op}\to \CAlg^{op}_{\QQ}=\Aff_{\QQ}$
preserves small limits.
Therefore, the derived affine group scheme $G(X,x)$ is the Cech nerve
of $\Spec \QQ\to \Spec T_X$ in $\Aff_{\QQ}$.
The second claim is clear from this argument.
\QED

\section{Motivic Galois action}
\label{MGA}

Let $K$ be a field of characteristic zero.
Let $\mathsf{R}_E:\DM^\otimes(k)\to \mathsf{D}^\otimes(K)$ be a
realization functor associated to a mixed Weil cohomology theory $E$.
In \cite{Tan} (see also \cite{Bar}, \cite{DTD}), we constructed
a derived affine group scheme $\mathsf{MG}_E$
over $K$ out of $\mathsf{R}_E$,
which we refer to as the {\it derived motivic Galois group}
with respect to $E$.
It has many favorable properties such as the consistency of motivic conjectures.
The most important property of $\mathsf{MG}_E$ for us is that
it represents the automorphism group
of the symmetric monoidal functor $\mathsf{R}_E$,
see Definition~\ref{realauto} or \cite{Tan} for the formulation.
Besides,
we have the usual affine group scheme
$MG_E$ associated to $\MG_E$
which we call the {\it motivic Galois group} with respect to $E$.
Note that a symmetric monoidal natural equivalence from $\mathsf{R}_E$
to itself induces a natural equivalence from $\CAlg(\mathsf{R}_E):\CAlg(\DM^\otimes(k))\to \CAlg_K$ to itself.
Actually, there is a canonical morphism from the automorphism
group
of $\mathsf{R}_E$ to the automorphism group of $\CAlg(\mathsf{R}_E)$.
Since $M_X$ belongs to $\CAlg(\DM^\otimes(k))$ for a smooth variety $X$,
the automorphism group of $\CAlg(\mathsf{R}_E)$ acts on
the image of $M_X$, e.g. $A_{PL,\infty}(X^t)$ in $\CAlg_{\QQ}$.
Consequently, it gives
rise to an action of the derived affine group scheme $\mathsf{MG}_E$
on the image of $M_X$.
Based on this natural idea,
we will
construct motivic Galois actions,
i.e., actions of $\MG_E$
by using the machinery of $\infty$-categories
in Section~\ref{actionstep1}.
In Section~\ref{diagramaction}, we focus on the case of
a cosimplicial diagram in $\CAlg(\DM^\otimes(k))$. The motivating cases
come from Section~\ref{hopfsection} and Section~\ref{realhopfsection}:
it yields an action of $\MG_E$ on the derived
affine group schemes $G_E^{(n)}(X,x)$ in Definition~\ref{realhopfdef}.

In Sections~\ref{underlyinggroup} to \ref{finalaction},
we turn to study how to obtain an action of the pro-algebraic group
$MG_E$ on the pro-unipotent completions of homotopy groups
and related invariants arising from various cohomology theories.
Our approach is to deduce the actions of $MG_E$ from the actions of $\MG_E$ on $G_E^{(n)}(X,x)$.
For example, under the situation of Proposition~\ref{cechtop},
one can derive the pro-unipotent completion of the fundamental group of $X^t$ from $G^{(1)}(X,x)$: if $G^{(1)}(X,x)=\Spec A$, the pro-unipotent completion
is given by $\Spec H^0(A)$
($A\simeq \QQ\otimes_{A_{PL}(X^t)}\QQ$ in $\CAlg_{\QQ}$).

\subsection{}
\label{actionstep1}

Our first task is to construct motivic Galois actions
on the images of multipilicative realization functors such as $A_{PL}(X^t)$.

\subsubsection{}

\begin{Definition}
\label{autogroup}
Let $\mathcal{I}$ be an $\infty$-category and $D:\mathcal{I}\to \Cat$ a functor. Suppose that $\mathcal{I}$ has an initial object $\xi$.
Let $C$ be an object of $D(\xi)$.
Let $(-)^\simeq:\Cat \to \SSS$ be the functor which carries an $\infty$-category
$\mathcal{C}$ to its largest Kan subcomplex $\mathcal{C}^\simeq$.
Namely, it is the right adjoint to the inclusion $\SSS\to \Cat$.
Let $\mathcal{F}_D\to \mathcal{I}$ be a left fibration
obtained by applying the unstraightening functor or relative nerve functor \cite{HTT} to $\mathcal{I}\to \Cat\to \SSS$.
By \cite[3.3.3.4]{HTT}, a section $\mathcal{I}\to \mathcal{F}_D$
of $\mathcal{F}_D\to \mathcal{I}$ corresponds to an object
in the limit $\varprojlim_{i\in \mathcal{I}}D(i)^\simeq$ in $\SSS$.
We let $s:\mathcal{I}\to \mathcal{F}_D$ be the section
that corresponds to the image of $C$ under the canonoical functor
$D(\xi)^\simeq\to \varprojlim_{i\in \mathcal{I}}D(i)^\simeq$.
Through the correspondence between left fibrations over $\mathcal{I}$
and functors $\mathcal{I}\to \SSS$ (cf. \cite[3.2, 4.2.4.4]{HTT}), 
$\mathcal{F}_D\to \mathcal{I}$ endowed with the section $s$
amounts to the functor $(-)^\simeq \circ D:\mathcal{I}\to \SSS$
with a natural transformation $\ast \to (-)^\simeq \circ D$
from the constant functor
$\ast:\mathcal{I}\to \SSS$ taking the value $\Delta^0$.
By the adjunction, the natural transtransformation is described as
a functor $D_*:\mathcal{I}\to \SSS_{*}:=\SSS_{\Delta^0/}\subset \Fun(\Delta^1,\SSS)$
such that the composition $\mathcal{I}\to \SSS_{*}\to \SSS$
with the forgetful functor is $(-)^\simeq \circ D$.
We shall refer to $D_*$ as the functor extended by $C$.
Let $\Grp(\SSS)$ denote the category of group objects in $\SSS$ (see e.g. \cite[7.2.2.1]{HTT}, \cite[Definition A.2]{Tan}).
Let $\Omega_*:\SSS_*\to \Grp(\SSS)$
be the functor which carries the based space $S$ to the based loop space
$\Omega_*(S)$.
We define the automorphism group functor of $C$ over $\mathcal{I}$
to be the composite
\[
\Aut_{\mathcal{I}}(C):\mathcal{I}\stackrel{D_*}{\longrightarrow}\SSS_*\stackrel{\Omega_*}{\longrightarrow} \Grp(\SSS).
\]
We usually write $\Aut(C)$ for $\Aut_{\mathcal{I}}(C)$.
\end{Definition}

\begin{Remark}
\label{autogroupR}
For any object $i$ in $\mathcal{I}$, the composition 
$\mathcal{I}\stackrel{\Aut_{\mathcal{I}}(C)}{\to} \Grp(\SSS)\to \SSS$
with the forgetful functor sends $i$ to the $\infty$-groupoid (space)
that is equivalent to the mapping space $\Map_{D(i)}(f(C),f(C))$
where $f:\xi\to i$ is the canonical functor from the initial object.
Indeed, the composite $\mathcal{I} \to \SSS$
sends $i$ to the fiber product $\Delta^0\times_{D(i)^\simeq}\Delta^0$
in $\SSS$,
defined by the map $\Delta^0\to D(i)$ determined by $f(C)$.
The fiber product $\Delta^0\times_{D(i)^\simeq}\Delta^0$
is explicitly given by the fiber product $\{f(C)\}\times_{D(i)^\simeq}\Fun(\Delta^1,D(i)^\simeq)\times_{D(i)^\simeq}\{f(C)\}$ of (genuine) simplicial sets, that is a model of the mapping space (cf. \cite[1.2.2, 4.2.1.8]{HTT}).
\end{Remark}

\begin{Definition}
\label{symfuncauto}
Let $\CAlg_{(-)}:\CAlg_K\to \wCat$ be a functor which carries $A$ to
$\CAlg_A$ where $\CAlg_A$ is the $\infty$-category of commutative
ring spectra over $A$, that is, commutative algebra objects
in $\Mod_A^\otimes$ (a morphism $A\to A'$
maps to $\CAlg_A\to \CAlg_{A'}$ given by the base change $\otimes_{A}A'$,
see Section~\ref{pre2} for the formulation). 
Let $C$ be an object of $\CAlg_{K}$.
We apply Definition~\ref{autogroup} to $\CAlg_{(-)}:\mathcal{I}=\CAlg_K\to \wCat$
and $C$ after replacing $\Cat$ and $\SSS$ by $\wCat$ and $\wSSS$, respectively.
We then define $\Aut_{\CAlg_K}(C):\CAlg_{K}\to \Grp(\wSSS)$ 
to be the automorphism group functor of $C$ over $\CAlg_K$.

Let $L$ be an $\infty$-category.
Let $(-)^L:\wCat\to \wCat$ be the functor which carries
$\mathcal{C}$ to $\Fun(L,\mathcal{C})$.
Namely, it is given by cotensoring with $L$.
Let $h:L\to \CAlg_K$ be a functor (which we will consider to be a diagram
in $\CAlg_K$ indexed by $L$).
Consider the composition 
\[
\mu_L:\CAlg_K\stackrel{\CAlg_{(-)}}{\to} \wCat\stackrel{(-)^L}{\to} \wCat.
\]
Applying Definition~\ref{autogroup} to $\mu_L$
and $h$,
we define $\Aut_{\CAlg_K}(h):\CAlg_K\to \Grp(\wSSS)$ to be the automorphism group functor of $h$
over $\CAlg_{K}$.
Notice that $\Aut_{\CAlg_K}(C)$ is the spacial case of $\Aut_{\CAlg_K}(h)$.
We usually write $\Aut(C)$ and $\Aut(h)$
for $\Aut_{\CAlg_K}(C)$ and $\Aut_{\CAlg_K}(h)$, respectively.
\end{Definition}

\begin{Definition}
Let $\Mod_{(-)}:\CAlg_K\to \wCat$ be a functor which carries $A$ to
$\Mod_A$ (a morphism $A\to A'$
maps to $\Mod_A\to \Mod_{A'}$ given by the base change $\otimes_{A}A'$,
see Section~\ref{pre2} for the formulation). 
Let $P$ be an object of $\mathsf{D}(K)\simeq \Mod_K$.
Applying Definition~\ref{autogroup} to $\Mod_{(-)}:\CAlg_K\to \wCat$
and $P$,
we define $\Aut(P)=\Aut_{\CAlg_K}(P):\CAlg_{K}\to \Grp(\wSSS)$ 
to be the automorphism group functor of $P$ over $\CAlg_K$.
\end{Definition}

Let $\mathsf{R}_E:\DM^\otimes(k)\to \mathsf{D}^\otimes(K)=\Mod_K^\otimes$
be the realization functor associate to a mixed Weil cohomology theory $E$
with coefficients in a field $K$ of characteristic zero.
The coefficient field of $\DM(k)$ 
will be $K$, but one can also
adopt the setting where the coefficient
field of $\DM(k)$ is $\QQ$ (one may choose either one
depending on the purpose).
Let $\mathsf{MG}_E=\Spec B$
be a derived affine group
scheme over $K$
which we call
the derived motivic Galois 
group with respect to $E$ (see \cite{Tan}).
Here the fundamental property of $\mathsf{MG}_E$ for us
is that it represents
the automorphism group functor 
$\Aut(\mathsf{R}_E):\CAlg_K\to \Grp(\wSSS)$
of the realization functor $\mathsf{R}_E$
(see Definition~\ref{realauto} for its definition).
Namely, if one regards $\mathsf{MG}_E$ as a functor
$\CAlg_K\to \Grp(\wSSS)$, then we have an equivalence
$\mathsf{MG}_E\simeq \Aut(\mathsf{R}_E)$.

\begin{Proposition}
\label{goodaction}
Let $C$ be an object of $\CAlg(\DM^\otimes(k))$.
There is a (canonical) action of $\mathsf{MG}_E$ on $\CAlg(\mathsf{R}_E)(C)$.
(Recall that $\CAlg(\mathsf{R}_E):\CAlg(\DM^\otimes(k))\to \CAlg_K$
is the multiplicative realization functor, Section~\ref{realization}.)
Namely, there is a morphism $\mathsf{MG}_E\to \Aut(\CAlg(\mathsf{R}_E)(C))$
in $\Fun(\CAlg_K,\Grp(\wSSS))$.
In particular, we have a (canonical) action of $\mathsf{MG}_E$ on $\CAlg(\mathsf{R}_E)(M_X)$. Moreover, the following properties hold:
\begin{enumerate}
\renewcommand{\labelenumi}{(\theenumi)}

\item The actions are functorial in $\CAlg(\DM^\otimes(k))$:
Namely, if we let $p:L\to \CAlg(\DM^\otimes(k))$ be a functor from an $\infty$-category $L$ and let $h:L\to \CAlg(\DM^\otimes(k))\stackrel{\CAlg(\mathsf{R}_E)}{\to} \CAlg_{K}$
be the composition with the multiplicative realization functor, then there is a morphism $\mathsf{MG}_E\to \Aut(h)$.
For a functor $g:M\to L$ of $\infty$-categories,
the action (morphism) $\mathsf{MG}_E\to \Aut(h\circ g)$ 
is naturally equivalent to $\mathsf{MG}_E\to \Aut(h)\to \Aut(h\circ g)$
where the the first arrow is given by the action on $h$, and the second arrow
is induced  by the composition with $M\to L$.

\item The action is compatible with the formation of colimits:
Let $p:L\to \CAlg(\DM^\otimes(k))$ be a functor from a small $\infty$-category, and $\overline{p}:L^\triangleright \to \CAlg(\DM^\otimes(k))$ a colimit diagram of $p$ (here $(-)^\triangleright$ indicates the right cone
\cite{HTT}). Let $C$ be the colimit
in $\CAlg(\DM^\otimes(k))$, that is, the image of the cone point.
Let
$q:L\to \CAlg_K$ and $\overline{q}:L^\triangleright \to \CAlg_K$
be the composites $\CAlg(\mathsf{R}_E)\circ p$ and $\CAlg(\mathsf{R}_E)\circ \overline{p}$, respectively.
Then the (action) morphism $\mathsf{MG}_E\to \Aut(\CAlg(\mathsf{R}_E)(C))$ factors through
the morphism $\mathsf{MG}_{E}\to \Aut(q)$ in the sense that the
restriction to $L$ induces an equivalence $\Aut(\overline{q})\stackrel{\sim}{\to}\Aut(q)$, and the composite
\[
\mathsf{MG}_{E}\to \Aut(q)\simeq  \Aut(\overline{q})\to \Aut(\CAlg(\mathsf{R}_E)(C))
\]
is naturally equivalent to the ``action''
$\mathsf{MG}_E\to \Aut(\CAlg(\mathsf{R}_E)(C))$. Here the final arrow is
induced by the restriction to the cone point of $L^\triangleright$.

\item There is a (canonical) action of $\mathsf{MG}_E$
on $\mathsf{R}_E(C)$, that is a morphism $\mathsf{MG}_E\to \Aut(\mathsf{R}_E(C))$. We here distinguish the underlying module
$\mathsf{R}_E(C)$ in $\mathsf{D}(K)$ from
$\CAlg(\mathsf{R}_E)(C)$ in $\CAlg_K$.
The action on $\CAlg(\mathsf{R}_E)(C)$ is compatible with
that on $\mathsf{R}_E(C)$ in the sense that there is a canonical morphism
$\Aut(\CAlg(\mathsf{R}_E)(C))\to \Aut(\mathsf{R}_E(C))$ induced by the forgetful functor, and $\mathsf{MG}_E\to \Aut(\mathsf{R}_E(C))$
is equivalent to the composite $\mathsf{MG}_E\to \Aut(\CAlg(\mathsf{R}_E)(C))\to \Aut(\mathsf{R}_E(C))$.
\end{enumerate}
\end{Proposition}

\begin{Corollary}
Suppose that $k$ is embedded in $\CC$. Let
$X^t$ be the underlying topological space of $X\times_{\Spec k}\Spec \CC$.
If $\mathsf{MG}$
denotes the derived motivic Galois group with respect to the singular cohomology
theory, there is a canonical
action of $\mathsf{MG}$ on $A_{PL,\infty}(X^t)\simeq T_X$.
\end{Corollary}

\Proof
Combine Proposition~\ref{goodaction} and Theorem~\ref{real1}.
\QED

\begin{Remark}
Let $A\in \CAlg_K$
and let $g:\Delta^0\to \mathsf{MG}_E(A)$ be an ``$A$-valued point''.
Through the equivalence $\mathsf{MG}_E(A)\simeq \Aut(\mathsf{R}_E)(A)$,
$g$ may be viewed as an automorphism of the composite $\DM^\otimes(k)\stackrel{\mathsf{R}_E}{\to} \Mod^\otimes_K\stackrel{\otimes_KA}{\to} \Mod^\otimes_A$.
It gives rise to an automorphism $u$ of the composite
\[
\CAlg(\DM^\otimes(k)) \stackrel{\CAlg(\mathsf{R}_E)}{\to} \CAlg_K\stackrel{\otimes_KA}{\to} \CAlg_A
\]
(see Section~\ref{prepre} below).
The image $\Delta^0\to\Aut(\CAlg(\mathsf{R}_E)(C))(A)$ of $g$
under the ``action'' $\mathsf{MG}_E(A)\to \Aut(\CAlg(\mathsf{R}_E)(C))(A)$
is a class of an equivalence $\CAlg(\mathsf{R}_E)(C)\otimes_KA\stackrel{\sim}{\to}\CAlg(\mathsf{R}_E)(C)\otimes_KA$ in $\CAlg_A$ obtained from
the automorphism $u$
by evaluating at $C$ (composing with the map
$\Delta^0\to  \CAlg(\DM^\otimes(k))$ determined by $C$).
\end{Remark}

\begin{Remark}
\label{motsubcategory}
One can replace $\DM^\otimes(k)=\mathcal{C}^\otimes$ by
a stable subcategory $\mathcal{E}^\otimes\subset \DM^\otimes(k)$ that is closed under small colimits
and is generated by a small set of dualizable objects.
Again by the main result of \cite{Tan} there is a derived affine group
scheme $\mathsf{MG}_{E,\mathcal{E}^\otimes}$ that represents
$\Aut(\mathsf{R}_E|_{\mathcal{E}^\otimes})$, and for $C\in \CAlg(\mathcal{E}^\otimes)$, 
$\mathsf{MG}_{E,\mathcal{E}^\otimes}$ acts on $\CAlg(\mathsf{R}_E(C))$.
In certain good cases, one can obtain
$\mathsf{MG}_{E,\mathcal{E}^\otimes}$ by means of equivariant
bar constructions, see
\cite{Bar}, \cite{DTD}, \cite{Spi}.
\end{Remark}

\subsubsection{}
\label{prepre}
We start with some $\infty$-categorical preliminary constructions.
To make things elementary, we make some efforts to make extensive use of
the machinery of simplicial categories, i.e., simplicially enriched categories, whereas in 
the earlier version of this manuscript in 2016,
many constructions heavily rely on the theory of left/(co)Cartesian fibrations.

\vspace{2mm}

Let $\operatorname{Cat}_{\infty}^{\textup{sMon},\Delta}$ be a simplicial category defined
as follows. The objects of $\operatorname{Cat}_{\infty}^{\textup{sMon},\Delta}$
are symmetric monoidal small
$\infty$-categories $\mathcal{C}^\otimes\to \Gamma$.
Give two symmetric monoidal $\infty$-categories $\mathcal{C}^\otimes \to \Gamma$ and $\mathcal{D}^\otimes\to \Gamma$, we define 
$\Fun_{\Gamma}^\otimes(\mathcal{C}^\otimes,\mathcal{D}^\otimes)$
to be the full subcategory of $\Fun_{\Gamma}(\mathcal{C}^\otimes,\mathcal{D}^\otimes)$ that consists of symmetric monoidal functors (cf. \cite[2.1.2]{HA}).
We define the mapping simplicial set
$\Map^\otimes(\mathcal{C}^\otimes,\mathcal{D}^\otimes):=\Map_{\operatorname{Cat}_{\infty}^{\textup{sMon},\Delta}}(\mathcal{C}^\otimes,\mathcal{D}^\otimes)$ to be the largest Kan subcomplex of 
$\Fun_{\Gamma}^\otimes(\mathcal{C}^\otimes,\mathcal{D}^\otimes)$.
The composition is defined by the restriction of composition of
 function complexes. The $\infty$-category $\operatorname{Cat}_{\infty}^{\textup{sMon}}$ is defined to be the simplicial nerve of $\operatorname{Cat}_{\infty}^{\textup{sMon},\Delta}$.

We let $\operatorname{Cat}_{\infty}^{\Delta}$ be the simplicial category
defined as follows. Objects are $\infty$-categories, and given two $\infty$-categories $\mathcal{C}$ and $\mathcal{D}$, the simplicial set $\Map(\mathcal{C},\mathcal{D})$ is the largest Kan subcomplex of $\Fun(\mathcal{C},\mathcal{D})$.
By definition, the simplicial nerve of $\operatorname{Cat}_{\infty}^{\Delta}$ is $\Cat$.

Let $\operatorname{Kan}^{\Delta}$ be the simplicial full subcategory of
$\operatorname{Cat}_{\infty}^{\Delta}$ that consists of Kan complexes.
For a symmetric monoidal $\infty$-category $\mathcal{C}^\otimes$,
the assignment $\mathcal{D}^\otimes \mapsto \Map^\otimes(\mathcal{C}^\otimes,\mathcal{D}^\otimes)$ determines a simplicial functor $h^\Delta_{\mathcal{C}^\otimes}:\operatorname{Cat}_{\infty}^{\textup{sMon},\Delta} \to \operatorname{Kan}^{\Delta}$ in the natural way.
Taking the simplicial nerve, we obtain $h_{\mathcal{C}^\otimes}:=\NNNN(h_{\mathcal{C}^\otimes}^\Delta):\operatorname{Cat}_{\infty}^{\textup{sMon}}=\NNNN(\operatorname{Cat}_{\infty}^{\textup{sMon},\Delta})\to \NNNN(\operatorname{Kan}^{\Delta})=\SSS$. We remark that it is equivalent to the functor $\operatorname{Cat}_{\infty}^{\textup{sMon}} \to \SSS$ corepresented by
$\mathcal{C}^\otimes$ defined in \cite[5.1.3]{HTT} (in the dual form).
Similarly, for an $\infty$-category $\mathcal{C}$,
the assigment $\mathcal{D} \mapsto \Map(\mathcal{C},\mathcal{D})$ determines a simplicial functor $h^\Delta_{\mathcal{C}}:\operatorname{Cat}_{\infty}^{\Delta} \to \operatorname{Kan}^{\Delta}$.
Taking the simplicial nerve, we obtain $h_{\mathcal{C}}:=\NNNN(h_{\mathcal{C}}^\Delta):\operatorname{Cat}_{\infty}=\NNNN(\operatorname{Cat}_{\infty}^{\Delta})\to \NNNN(\operatorname{Kan}^{\Delta})=\SSS$.

Next we construct a functor
$\CAlg:\smCat\to \Cat$
from the $\infty$-category of symmetric monoidal (small) $\infty$-categories
to the $\infty$-category of $\infty$-categories, which sends $\mathcal{C}^\otimes$ to $\CAlg(\mathcal{C}^\otimes)$. 
For this purpose we construct a simplicial functor
\[
\CAlg^{\Delta}:\operatorname{Cat}_{\infty}^{\textup{sMon},\Delta}\longrightarrow \operatorname{Cat}_{\infty}^{\Delta}
\]
which carries $\mathcal{C}^\otimes\to \Gamma$ to $\CAlg(\mathcal{C}^\otimes)=\Fun^{\textup{lax}}_{\Gamma}(\Gamma,\mathcal{C}^\otimes)$ where
$\Fun^{\textup{lax}}_{\Gamma}(-,-)$ indicates the full subcategory of
$\Fun_{\Gamma}(-,-)$ that consists of lax symmetric monoidal functors.
To do this, given two symmetric monoidal $\infty$-categories
we will define a map of simplicial sets
$\Map^\otimes(\mathcal{C}^\otimes,\mathcal{D}^\otimes)\to \Map(\CAlg(\mathcal{C}^\otimes),\CAlg(\mathcal{D}^\otimes))$.
Let $K$ be a simplicial set and $f:K\to \Map^\otimes(\mathcal{C}^\otimes,\mathcal{D}^\otimes)$ a map of simplcial sets.
The map amounts to a map of marked simplicial sets $\mathcal{C}^\otimes\times K^{\sharp}\to \mathcal{D}^\otimes$ over $\Gamma$ where $K^\sharp$ denotes
the marked simplicial sets such that all edges are marked.
To the map we associate a map
of simplicial sets $\CAlg(\mathcal{C}^\otimes)\times K\to \CAlg(\mathcal{D}^\otimes)$, equivalently $K\to \Fun(\CAlg(\mathcal{C}^\otimes),\CAlg(\mathcal{D}^\otimes))$ as follows.
Note that for a simplicial set $S$, $S\to \Fun_{\Gamma}(\Gamma,\mathcal{C}^\otimes\times_{\Gamma}(\Gamma\times K))$ corresponds to a pair of maps $S\times \Gamma\to \mathcal{C}^\otimes$ over $\Gamma$ and $S\times \Gamma \to K$.
To $S\to \CAlg(\mathcal{C}^\otimes)\times K$ corresponding to
$\phi:S\times \Gamma\to \mathcal{C}^\otimes$ over $\Gamma$ and $\psi:S\to K$
we associate $S\to \Fun_{\Gamma}(\Gamma,\mathcal{C}^\otimes\times_{\Gamma}(\Gamma\times K))$ corresponding to the pair
$\phi:S\times \Gamma\to \mathcal{C}^\otimes$ over $\Gamma$
and $S\times \Gamma\stackrel{\textup{pr}_1}{\to} S\stackrel{\psi}{\to} K$.
It gives rise to a map
\[
r:\Fun_{\Gamma}^{\textup{lax}}(\Gamma,\mathcal{C}^\otimes)\times K \to \Fun_{\Gamma}(\Gamma,\mathcal{C}^\otimes\times K).
\]
Let $c:\Fun_{\Gamma}(\Gamma,\mathcal{C}^\otimes\times K)\times \Fun_{\Gamma}(\mathcal{C}^\otimes\times K,\mathcal{D}^\otimes)
\to \Fun_{\Gamma}(\Gamma,\mathcal{D}^\otimes)$
be composition. Let $\iota:\Delta^0\to \Fun_{\Gamma}(\mathcal{C}^\otimes\times K,\mathcal{D}^\otimes)$ be the map determined by $\mathcal{C}^\otimes\times K \to\mathcal{D}^\otimes$ over $\Gamma$ that corresponds to $f$.
Consider the following composite
\begin{eqnarray*}
\Fun_{\Gamma}^{\textup{lax}}(\Gamma,\mathcal{C}^\otimes)\times K \simeq (\Fun^{\textup{lax}}_{\Gamma}(\Gamma,\mathcal{C}^\otimes)\times K)\times \Delta^0 &\stackrel{r\times \iota}{\longrightarrow}& \Fun_{\Gamma}(\Gamma,\mathcal{C}^\otimes\times K)\times \Fun_{\Gamma}(\mathcal{C}^\otimes\times K,\mathcal{D}^\otimes) \\
&\stackrel{c}{\longrightarrow}& \Fun_{\Gamma}(\Gamma,\mathcal{D}^\otimes).
\end{eqnarray*}
The image of composition is contained in $\Fun_{\Gamma}^{\textup{lax}}(\Gamma,\mathcal{D}^\otimes)$. Therefore we obtain $\CAlg(\mathcal{C}^\otimes)\times K\to \CAlg(\mathcal{D}^\otimes)$ from $f$. According to
the functoriality with respect to $K$
it yields 
\[
\Map^\otimes(\mathcal{C}^\otimes,\mathcal{D}^\otimes)\to \Fun(\CAlg(\mathcal{C}^\otimes),\CAlg(\mathcal{D}^\otimes)).
\]
Since $\Map^\otimes(\mathcal{C}^\otimes,\mathcal{D}^\otimes)$ is a Kan complex,
its image is contained in $\Map(\CAlg(\mathcal{C}^\otimes),\CAlg(\mathcal{D}^\otimes))$. It is straightforward to see that
$\mathcal{C}^\otimes\mapsto \CAlg(\mathcal{C}^\otimes)$ and 
$\Map^\otimes(\mathcal{C}^\otimes,\mathcal{D}^\otimes)\to \Map(\CAlg(\mathcal{C}^\otimes),\CAlg(\mathcal{D}^\otimes))$ determine a simplicial
functor
$\CAlg^\Delta:\operatorname{Cat}_{\infty}^{\textup{sMon},\Delta}\to \operatorname{Cat}_{\infty}^{\Delta}$.
Taking the simplicial nerves we obtain a functor of $\infty$-categories
\[
\CAlg:\smCat \longrightarrow \Cat.
\]

There is another obvious simplicial functor
$\textup{For}^\Delta:\operatorname{Cat}_{\infty}^{\textup{sMon},\Delta}\longrightarrow \operatorname{Cat}_{\infty}^{\Delta}$ which carries any symmetric monoidal $\infty$-category $\pi:\mathcal{C}^\otimes\to \Gamma$
to the fiber $\pi^{-1}(\langle 1 \rangle)$, i.e., the underlying $\infty$-categories $\mathcal{C}$.
There is the forgetful functor $\CAlg(\mathcal{C}^\otimes)\to \mathcal{C}$
which is defined as
$\Fun_{\Gamma}^{\textup{lax}}(\Gamma,\mathcal{C}^\otimes) \to \Fun_{\Gamma}(\{\langle 1 \rangle\},\mathcal{C}^\otimes)$ induced by composition with $\{\langle 1 \rangle\}\to \Gamma$.
It gives rise to a simplicial natural transformation $\CAlg^\Delta\to \textup{For}^\Delta$.

\subsubsection{}
\label{pre2}

Replacing the universe $\mathbb{U}$ by a larger universe $\mathbb{U}\in \mathbb{V}$,
we define
the $\infty$-category $\wCat$ of $\mathbb{V}$-small $\infty$-categories
, the $\infty$-category $\wsmCat$ of symmetric monoidal $\mathbb{V}$-small $\infty$-categories, and $\widehat{\CAlg}:\wsmCat\to \wCat$ instead of $\Cat$, $\smCat$ and $\CAlg:\smCat\to \Cat$. But for simplicity we write $\CAlg$ for $\widehat{\CAlg}$.

Let $\Theta_K:\CAlg_K\to \wsmCat$ be a functor which carries $A$ to
$\Mod_A^\otimes$ where $\Mod_A^\otimes$ is the symmetric monoidal $\infty$-category of $A$-module spectra (see \cite{HA}, \cite[Appendix A.4]{HA} for the precise construction).
Any morphism $A\to A'$ maps to the symmetric monoidal
functor $\Mod_A^\otimes\to \Mod_{A'}^\otimes$ informally given by the base 
change $\otimes_AA'$.
Let $\NNNN(\textup{For}^\Delta):\wsmCat\to \wCat$ be the forgetful functor.
We define $\Mod_{(-)}:\CAlg_K\to \wCat$ to be the composite
of $\Theta_K$ and the forgetful functor.
We define $\CAlg_{(-)}$ to be the composite
$\CAlg_K\stackrel{\Theta_K}{\longrightarrow} \wsmCat\stackrel{\CAlg}{\longrightarrow} \wCat$.

\begin{Definition}
\label{realauto}
Consider the composite $\rho:\CAlg_K \stackrel{\Theta_K}{\to} \wsmCat\stackrel{h_{\DM(k)^\otimes}}{\longrightarrow} \wSSS$, which carries $A$ to $\Map^\otimes(\DM^\otimes(k),\Mod_A^\otimes)$.
Let $\mathsf{R}_E:\DM^\otimes(k)\to \mathsf{D}^\otimes(K)=\Mod_K^\otimes$
be the realization functor.
It may be viewed
as an object of $\Map^\otimes(\DM^\otimes(k),\Mod_K^\otimes)$.
Applying Definition~\ref{autogroup} to $\rho:\CAlg_K \to \wSSS$
and $\mathsf{R}_E$ we define the automorphism group functor
$\Aut(\mathsf{R}_E):\CAlg_K\to \Grp(\wSSS)$ of $\mathsf{R}_E$
over $\CAlg_K$.
\end{Definition}

\begin{Remark}
The definition of $\Aut(\mathsf{R}_E)$ is apparently different from
that in \cite{Tan} because in {\it loc.cit.} we use the full subcategory
$\DM_\vee^\otimes(k)$ spanned by compact (dualizable) objects
instead of $\DM^\otimes(k)$. But this point is neglective.
Since $\DM^\otimes(k)$ is canonically
equivalent to the symmetric monoidal
$\infty$-category $\Ind(\DM_\vee^\otimes(k))$
of Ind-objects, thus by the (symmetric monodial) Kan extension,
we see that there is a canonical equivalence $\Aut(\mathsf{R}_E)\simeq \Aut(\mathsf{R}_E|_{\DM_\vee^\otimes(k)})$ induced by the restriction 
to $\DM_\vee^\otimes(k)\subset \DM^\otimes(k)$.
\end{Remark}

\subsubsection{}

{\it Construction of the action/Proof of Proposition~\ref{goodaction}}.
Let $L$ be an $\infty$-category.
Consider the following three simplicial functors:
\begin{itemize}
\item Put $\alpha^\Delta=h_{\DM^\otimes(k)}^\Delta:\widehat{\operatorname{Cat}}_{\infty}^{\textup{sMon},\Delta}\to \widehat{\operatorname{Kan}}^{\Delta}$.
It sends a symmetric monoidal $\infty$-category
$\mathcal{D}^\otimes$ to the Kan complex
$\Map^\otimes(\DM^\otimes(k),\mathcal{D}^\otimes)$.

\item Let $\beta_{L}^\Delta:\widehat{\operatorname{Cat}}_{\infty}^{\textup{sMon},\Delta}\to \widehat{\operatorname{Kan}}^{\Delta}$ be a simplicial functor that carries $\mathcal{D}^\otimes$ to $\Map(L,\CAlg(\mathcal{D}^\otimes))$. It is defined as the composite
$\widehat{\operatorname{Cat}}_{\infty}^{\textup{sMon},\Delta}\stackrel{\CAlg^\Delta}{\longrightarrow} \widehat{\operatorname{Cat}}_{\infty}^{\Delta}\stackrel{h_{L}}{\longrightarrow} \widehat{\operatorname{Kan}}^{\Delta}$.

\item Let $\gamma_{L}^\Delta:\widehat{\operatorname{Cat}}_{\infty}^{\textup{sMon},\Delta}\to \widehat{\operatorname{Kan}}^{\Delta}$ be a simplicial functor that carries $\mathcal{D}^\otimes$
to $\Map(L,\mathcal{D})$.
It is defined as the composite
$\widehat{\operatorname{Cat}}_{\infty}^{\textup{sMon},\Delta}\stackrel{\operatorname{For}^\Delta}{\longrightarrow} \widehat{\operatorname{Cat}}_{\infty}^{\Delta} \stackrel{h_{L}}{\longrightarrow} \widehat{\operatorname{Kan}}^{\Delta}$.

\end{itemize}
For each $\mathcal{D}^\otimes$, the simplicial functor $\CAlg^\Delta$ induces a map of simplicial sets
\[
\Map^\otimes(\DM^\otimes(k),\mathcal{D}^\otimes)\to 
\Map(\CAlg(\DM^\otimes(k)),\CAlg(\mathcal{D}^\otimes)).
\]
It is easy to check that these maps determine a simplicial natural
transformation $\alpha^\Delta\to \beta_{\CAlg(\DM^\otimes(k))}^\Delta$.
Similarly, $\operatorname{For}^\Delta$ induces a map of simplicial sets
\[
\Map^\otimes(\DM^\otimes(k),\mathcal{D}^\otimes)\to 
\Map(\DM(k),\mathcal{D}).
\]
It gives rise to 
a simplicial natural
transformation $\alpha^\Delta\to \gamma_{\DM(k)}^\Delta$.
Let $L\to \CAlg(\DM^\otimes(k))$ be a functor.
The composition induces simplicial natural transformation
$\beta_{\CAlg(\DM^\otimes(k))}^\Delta \to \beta_{L}^\Delta$.
Also, $L\to \CAlg(\DM^\otimes(k))\stackrel{\textup{forget}}{\to} \DM(k)$
induces
$\gamma_{\DM(k)}^\Delta \to \gamma_{L}^\Delta$.

Now applying the simplicial nerve functor to $\alpha^\Delta$
we obtain $\alpha=h_{\DM^\otimes(k)}:\wsmCat \to \wSSS$.
Similarly, we obtain $\beta_{L},\gamma_L:\wsmCat \to \wSSS$
from $\beta_{L}^\Delta$ and $\gamma_L^\Delta$.
Consider the simplicial natural transformation
$\alpha^\Delta\to \beta_{\CAlg(\DM^\otimes(K))}^\Delta \to \beta_{L}^\Delta$.
It determines a natural transformation from $\alpha$ to $\beta_L$.
In fact,
we think of $\alpha^\Delta\to \beta_{L}^\Delta$
as a simplicial functor $[1]\times \widehat{\operatorname{Cat}}_{\infty}^{\textup{sMon},\Delta} \to \operatorname{Kan}^\Delta$
such that $[1]=\{0,1\}$ is the linearly ordered set
regarded as a (simplicial) category, and
the restriction to $\{0\}\times \widehat{\operatorname{Cat}}_{\infty}^{\textup{sMon},\Delta} \to \operatorname{Kan}^\Delta$ 
(resp. $\{1\}\times \widehat{\operatorname{Cat}}_{\infty}^{\textup{sMon},\Delta} \to \operatorname{Kan}^\Delta$)
is $\alpha^\Delta$
(resp. $\beta_L^\Delta$).
Since the simplicial nerve functor preserves products,
$\Delta^1\times \wsmCat\simeq \NNNN(\{0\to 1\}\times \widehat{\operatorname{Cat}}_{\infty}^{\textup{sMon},\Delta}) \to \NNNN(\widehat{\operatorname{Kan}}^\Delta)=\wSSS$
defines a natural transformation from $\alpha$ to $\beta_L$,
that is, $\Delta^1\times \wsmCat\to \wSSS$ such that
$\{0\}\times \wsmCat\to \wSSS$ is $\alpha$, and $\{1\}\times \wsmCat\to \wSSS$
is $\beta_L$.
Similarly, $\alpha^\Delta\to \gamma^{\Delta}_{\DM(k)}\to \gamma_{L}^\Delta$
determines a natural transformation from $\alpha$ to $\gamma_L$.

\vspace{1mm}

Next for $p:L\to \CAlg(\DM^\otimes(k))$ and $h=\CAlg(\mathsf{R}_E)\circ p$,
we construct an action of $\MG_E$ on $\Aut(h)$ (cf. Definition~\ref{symfuncauto}).
If $C$ is an object of $\CAlg(\DM^\otimes(k))$,
the automorphism group functor $\Aut(\CAlg(\mathsf{R}_E(C)))$ of $\CAlg(\mathsf{R}_E(C))$ over $\CAlg_k$
is nothing but $\Aut(h)$ where $L=\Delta^0$, and the functor $p:\Delta^0\to \CAlg(\DM^\otimes(k))$ is determined by $C$.
Let $\Delta^1\times \wsmCat\to \wSSS$
be the natural transformation from $\alpha$ to $\beta_L$
defined above. Composing with $\Theta_K$, we have
$\Delta^1\times \CAlg_K \to \wSSS$, that is a natural transformation
from $\rho=\alpha\circ \Theta_K$ to $(-)^\simeq\circ \mu_L=\beta_L\circ \Theta_K$
(we here use the notation in Definition~\ref{symfuncauto}, \ref{realauto}).
Remember that $\mathsf{R}_E$ is an object of $\Map^\otimes(\DM^\otimes(k),\Mod_K^\otimes)$. Thus, as in Definition~\ref{autogroup},
both $\alpha$ and $\beta_L$ are respectively  promoted to functors
$\alpha_*,\beta_{L*}:\wsmCat\to \wSSS_*$ extended by $\mathsf{R}_E$
and $h\in \Map(L,\CAlg_K)$, and
$\Delta^1\times \CAlg_K \to \wSSS$ is promoted to
a natural transformation $\Delta^1\times \wsmCat \to \wSSS_*$
from $\alpha_*$ to $\beta_{L*}$.
Composing $\Omega_*:\wSSS_*\to \Grp(\wSSS)$ and $\Theta_K$,
we obtain
\[
\Delta^1\times \CAlg_K\to \Delta^1\times \wsmCat\to \wSSS_*\to \Grp(\wSSS)
\]
that is a natural tranformation from $\Aut(\mathsf{R}_E)$ to
$\Aut(h)$ (cf. Definition~\ref{symfuncauto}, \ref{realauto}).
Since we have the equivalence $\MG_E\simeq \Aut(\mathsf{R}_E)$,
it defines a morphism $\MG_E\simeq \Aut(\mathsf{R}_E)\to \Aut(h)$
in $\Fun(\CAlg_K,\Grp(\wSSS))$. An action of $\MG_E$ on $h$
is defined to be this morphism.

\vspace{1mm}

We prove the property (1).
For a map $g:M\to L$, there is a simplicial natural
transformation $\beta^\Delta_L\to \beta_M^\Delta$ induced by the composition
with $g$. Therefore, by our construction
the functoriality is obvious.

Next we prove the property (2).
Let $\mathcal{K}$ be the full subcategory
of $\Fun(L^\triangleright,\CAlg_A)$, that consists of
those functors $F:L^\triangleright\to \CAlg_K$
such that the image of the cone point of $L^{\triangleright}$ is a colimit of the restriction $F|_{L}$.
Then by taking account of left Kan extensions \cite[4.3.2.15]{HTT} (keep in mind that $\CAlg_A$
admits small colimits),
the map $\Fun(L^\triangleright,\CAlg_A)\to \Fun(L,\CAlg_A)$ given by
the restriction induces 
an equivalence $\mathcal{K}\stackrel{\sim}{\to} \Fun(L,\CAlg_A)$
of $\infty$-categories.
Note that $\overline{p}:L^\triangleright\to \CAlg(\DM^\otimes(k))$ is a colimit diagram
(of $L\to \CAlg(\DM^\otimes(k))$).
The composite 
$\overline{q}:L^\triangleright\to \CAlg(\DM^\otimes(k))\stackrel{\CAlg(\mathsf{R}_E))}{\to} \CAlg_K$ is also a colimit diagram
because $\CAlg(\mathsf{R}_E)$ is a left adjoint.
Also, the base change $\otimes_KA:\CAlg_K\to \CAlg_A$ is a left
adjoint.
Thus, the composite $L^\triangleright \stackrel{\overline{q}}{\to} \CAlg_K\to \CAlg_A$ belongs to $\mathcal{K}$.
By these observations, we see that $\Aut(\overline{q})\to \Aut(q)$
induced by the restriction is an equivalence in $\Fun(\CAlg_K,\Grp(\wSSS))$.
By the functoriality (1), we have the desired factorization of the action.

Finally, we prove (3).
On can define $\MG_E\to \Aut(\mathsf{R}_E(C))$
by using
$\alpha^\Delta \to \gamma_{\DM(k)}^\Delta \to \gamma_L^\Delta$
and $\mathsf{R}_E$ in the same way as we constructed
$\MG_E\to \Aut(\CAlg(\mathsf{R}_E)(C))$
from
$\alpha^\Delta \to \beta_{\DM^\otimes(k)}^\Delta \to \beta_L^\Delta$
and $\mathsf{R}_E$.
There is a simplicial natural transformation
$\beta_L^\Delta\to \gamma_L^{\Delta}$
which is given by $\Map(L,\CAlg(\mathcal{D}^\otimes))\to \Map(L,\mathcal{D}))$
induced by the composition with the forgetful functor
$\CAlg(\mathcal{D}^\otimes)\to \mathcal{D}$ for each $\mathcal{D}^\otimes$.
By the simplicial nerve fucntor and the construction in Definiton~\ref{autogroup}, it gives rise to $\Aut(\CAlg(\mathsf{R}_E)(C))\to \Aut(\mathsf{R}_E(C))$.
Note that the simplicial natural transformation
$\beta_L^\Delta\to \gamma_L^{\Delta}$ commutes with $\alpha^\Delta\to \beta_L^\Delta$ and $\alpha^\Delta\to \gamma_L^\Delta$.
By this commutativily we see that $\MG_E\to \Aut(\CAlg(\mathsf{R}_E)(C))\to \Aut(\mathsf{R}_E(C))$ is naturally equivalent to
$\MG_E\to \Aut(\mathsf{R}_E(C))$.
\QED

\begin{Remark}
\label{freemot}
Let $M$ be an object of $\DM(k)$.
Let $\Free_{\DM(k)}(M)$ in $\CAlg(\DM^\otimes(k))$
be the free commutative algebra
object generated by $M$ (see Definition~\ref{freealg}).
Let us observe that the action of $\MG_E$ on $\CAlg(\mathsf{R}_E)(\Free_{\DM(k)}(M))$ is essentially determined by
the action of of $\MG_E$ on $\mathsf{R}_E(M)$.
Since the realization functor is a left adjoint,
there is a canonical equivalence
$\Free_{K}(\mathsf{R}_E(M))\simeq \CAlg(\mathsf{R}_E)(\Free_{\DM(k)}(M))$
where $\Free_K:=\Free_{\Mod_K}$ is the free functor $\Mod_K\to \CAlg_K$,
i.e., the left adjoint to the forgetful functor.
Let $S$ be a space that belongs to $\SSS$.
Let $f:S\to \MG_E(K)\simeq \Aut(\mathsf{R}_E)(K)$ be a morphism (in $\SSS$).
Let $g:S\to \Aut(\CAlg(\mathsf{R}_E)(\Free_{\DM(k)}(M)))(K)\simeq \Map_{\CAlg_K}(\CAlg(\mathsf{R}_E)(\Free_{\DM(k)}(M)),\CAlg(\mathsf{R}_E)(\Free_{\DM(k)}(M)))$
be a class of the map induced by the action of $f$.
The forgetful functor induces morphisms
\begin{eqnarray*}
\Map_{\CAlg_K}(\CAlg(\mathsf{R}_E)(\Free_{\DM(k)}(M)),\CAlg(\mathsf{R}_E)(\Free_{\DM(k)}(M)))\ \ \ \ \ \ \ \ \ \ \ \ \ \ \ \ \ \ \ \ \  \\
\ \ \ \ \ \ \ \ \ \ \ \ \ \ \ \ \ \ \ \to \Map_{\Mod_K}(\CAlg(\mathsf{R}_E)(\Free_{\DM(k)}(M))^{\sharp},\CAlg(\mathsf{R}_E)(\Free_{\DM(k)}(M))^{\sharp}) \\
\simeq \Map_{\Mod_K}(\Free_{K}(\mathsf{R}_E(M))^{\sharp},\Free_{K}(\mathsf{R}_E(M))^{\sharp})\ \ \ \ \ \ \ \ \ \ \ \ \ \ \ \ \ \ \ \ \ \ \ \ \ \ 
\end{eqnarray*}
in $\SSS$ where $(-)^\sharp$ here indicates the underlying object.
By the compatibility (3) in Proposition~\ref{goodaction},
the image of $g$
is equivalent to the map
\[
h:S\to \Aut(\Free_{K}(\mathsf{R}_E(M))^{\sharp})(K)\simeq \Map_{\Mod_K}(\Free_{K}(\mathsf{R}_E(M))^{\sharp},\Free_{K}(\mathsf{R}_E(M))^{\sharp})
\]
that is determined by the action of $f$ on $\Free_{K}(\mathsf{R}_E(M))^{\sharp}$.
The composition with the canonical (unit) map $\mathsf{R}_E(M)\to \Free_{K}(\mathsf{R}_E(M))^{\sharp}$ yields the morphisms
\begin{eqnarray*}
\Map_{\Mod_K}(\Free_{K}(\mathsf{R}_E(M))^{\sharp},\Free_{K}(\mathsf{R}_E(M))^{\sharp}) &\to& \Map_{\Mod_K}(\mathsf{R}_E(M),\Free_{K}(\mathsf{R}_E(M))^{\sharp}) \\ &\stackrel{i}{\leftarrow}& \Map_{\Mod_K}(\mathsf{R}_E(M),\mathsf{R}_E(M))
\end{eqnarray*}
in $\SSS$. By the functoriality similar to (1) in Proposition~\ref{goodaction},
we see that the image of $h$ in $\Map_{\Mod_K}(\mathsf{R}_E(M),\Free_{K}(\mathsf{R}_E(M))^{\sharp})$ is equivalent to the image of
$r:S\to \Aut(\mathsf{R}_E(M))(K)\simeq \Map_{\Mod_K}(\mathsf{R}_E(M),\mathsf{R}_E(M))$ that is determined by the action of $f$ on $\mathsf{R}_E(M)$.
Note that by the adjunction, 
the composition gives an equivalence 
\[
\Map_{\CAlg_K}(\CAlg(\mathsf{R}_E)(\Free_{\DM(k)}(M)),\CAlg(\mathsf{R}_E)(\Free_{\DM(k)}(M)))\stackrel{\sim}{\to} \Map_{\Mod_K}(\mathsf{R}_E(M),\Free_{K}(\mathsf{R}_E(M))^{\sharp})
\]
in $\SSS$.
Also, the left arrow $i$ is a fully faithful functor
since $\mathsf{R}_E(M)\to \Free_{K}(\mathsf{R}_E(M))^{\sharp}$ defines a direct summand
of $\Free_{K}(\mathsf{R}_E(M))^{\sharp}$.
The image of $g$ in $\Map_{\Mod_K}(\mathsf{R}_E(M),\Free_{K}(\mathsf{R}_E(M))^{\sharp})$ lies in the essential image of $i$.
The image of $g$ is equivalent to
the image of $r$ under $i$. One can adopt this argument
to not only $K$ but arbitrary $A\in\CAlg_K$.
We remark that any object of $\CAlg(\DM^\otimes(k))$
is constructed from free commutative algebra objects
by forming colimits, see Section~\ref{Sull1}.
\end{Remark}

\subsection{}
\label{diagramaction}
Let $\Fun(\NNNN(\Delta^{op}),\Aff_K)$ be the $\infty$-category
of simplicial diagrams in $\Aff_K$. The $\infty$-category of
group objects in $\Aff_K$, i.e., derived affine group schemes,
is its full subcategory consisting of those simplicial diagram
satisfying the condition of group objects (cf. Section~\ref{preDAG}, see also
\cite[Appendix]{Tan}).
In Section~\ref{diagramaction}
we focus on actions on such objects.
We continue to use the notation in Section~\ref{actionstep1}.

Let $\mathfrak{C}:\NNNN(\Delta) \to \CAlg(\DM^\otimes(k))$ be
a functor which we regard as a cosimplicial diagram of commutative algebra objects in $\DM^\otimes(k)$.
Suppose that $\mathfrak{C}^{op}:\NNNN(\Delta)^{op} \to \CAlg(\DM^\otimes(k))^{op}$ is a group object.
One of our main examples is the opposite of the group object
$\mathcal{G}^{(n+1)}(X,x):\NNNN(\Delta)^{op}\to \CAlg(\DM^\otimes(k))^{op}$ introduced in Section~\ref{hopfsection}.
The multiplicative realization functor $\CAlg(\mathsf{R}_E)$ preserves coproducts and
sends a unit to $K\in \CAlg_K$.
It follows that
the composite 
\[
G_\bullet:\NNNN(\Delta)^{op} \stackrel{\mathfrak{C}^{op}}{\to} \CAlg(\DM^\otimes(k))^{op} \stackrel{\CAlg(\mathsf{R}_E)}{\longrightarrow} \CAlg_K^{op}=\Aff_K
\]
is a group object, that is, a derived affine group scheme over $K$.
We denote it simply by
$G$.

Invoking Proposition~\ref{goodaction} (see also Definition~\ref{symfuncauto})
to the opposite of the group object
$G_\bullet^{op}=h:\NNNN(\Delta)=L\to \CAlg_K$,
we get a morphism
\[
\mathsf{MG}_E\to \Aut(G_\bullet^{op})
\]
in $\Fun(\CAlg_K,\Grp(\wSSS))$, that is, an action of $\mathsf{MG}_E$
on $G_\bullet^{op}$.
Put
$\Aut(G):=\Aut(G_\bullet^{op})$.
Thus we have

\begin{Proposition}
\label{actiongroup}
Let $\mathfrak{C}^{op}:\NNNN(\Delta)^{op}\to \CAlg(\DM^\otimes(k))^{op}$
be a group object.
Let $G$ be the derived affine group scheme over $K$ that is induced by
$\mathfrak{C}^{op}$.
Then there is a (canonical) action of $\mathsf{MG}_E$ on $G$, that is a morphism $\mathsf{MG}_E\to \Aut(G)$
in $\Fun(\CAlg_K,\Grp(\wSSS))$.
\end{Proposition}

\begin{Remark}
We remark that informally $\Aut(G_\bullet^{op})$  is the automorphism group
of the cosimplicial object $G_\bullet^{op}$
in $\CAlg_K$.
Therefore, by our convention $\Aut(G)=\Aut(G_\bullet^{op})$, the morphism
$\mathsf{MG}_E\to \Aut(G)$ should be viewed as the ``right action'' on $G$
that corresponds to the ``left action'' on $G_\bullet^{op}$.
\end{Remark}

\begin{Remark}
The action is functorial with repect to a morphism of derived 
affine group schemes. 
Let
$\mathfrak{C}^{'op}:\NNNN(\Delta)^{op}\to \CAlg(\DM^\otimes(k))^{op}$
be another group object and $G'_{\bullet}:\NNNN(\Delta)^{op}\to \CAlg_K^{op}=\Aff_K$ the derived affine group scheme induced by the composition with
the multiplicative realization functor.
Suppose that there is a morphism (i.e., a natural
transformation) $\mathfrak{C}^{op}\to \mathfrak{C}^{'op}$.
It gives rise to $\theta:\Delta^1\times \NNNN(\Delta)\to \CAlg(\DM^\otimes(k))\to \CAlg_K$, such that $\{0\}\times \NNNN(\Delta) \to \CAlg_K$ is $G_{\bullet}^{'op}$,
and $\{0\}\times \NNNN(\Delta) \to \CAlg_K$ is $G_{\bullet}^{op}$.
By (1) of Proposition~\ref{goodaction}
the actions of $\MG_E$ on $G_{\bullet}^{op}$ and 
$G_{\bullet}^{'op}$ are simultaneously promoted to
an action on $\Aut(\theta)$, i.e., $\MG_E\to  \Aut(\theta)$.

\end{Remark}

\begin{Example}
\label{exgroup}
Let $(X,x:\Spec k\to X)$ be a pointed smooth variety over $k$.
As discussed in Section~\ref{realhopfsection} it gives rise to
a derived affine group scheme $G_E^{(n)}(X,x):\NNNN(\Delta)^{op}\to \Aff_{K}$.
Therefore, $\mathsf{MG}_E$ acts on $G_E^{(n)}(X,x)$.
\end{Example}

\subsection{}
\label{underlyinggroup}
In \cite{Tan} we defined the motivic Galois group
$MG_E$ of $\DM(k)$ (with respect to $E$)
to be a usual affine group scheme over $K$
(i.e., a pro-algebraic group) obtained from $\MG_E$.
Also, we can construct a usual affine group scheme
$\GG_E^{(n)}(X,x)$ from $G_E^{(n)}(X,x)$, Example~\ref{exgroup}.
In general,
if $G$ is a derived affine group scheme over the field of characteristic zero
$K$,
one can obtain a usual affine group scheme (i.e., pro-algebraic group)
$\overline{G}$
over $K$
from $G$, which we will call the underlying affine
group scheme (cf. \cite{Tan}).
We briefly review the procedure.

\vspace{1mm}

Let $\CAlg_K^{dg}$ be the category of
commutative dg algebras $C$ over $K$ (cf. Section~\ref{convention}).
Let $\CAlg_K^{dg,\ge0}$ be the full subcategory
of $\CAlg_K^{dg}$ that consists of
those objects $C$ such that
$H^i(C)=0$ for $i<0$.
It admits a combinatorial model category structure
such that a morphism $f:C\to C'$
is a weak equivalence (resp.
fibration) if the underlying map is a quasi-isomorphism
(resp. a surjective in each degree), see
\cite[Proposition 5.3]{Ol} or \cite[Theorem 6.2.6]{Fre}. Any object
is fibrant.
Any ordinary commutative algebra over $K$ is
a cofibrant object in $\CAlg_K^{dg,\ge0}$
when it is regarded as a commutative
dg algebra placed in degree zero.
The inclusion $\CAlg_K^{dg,\ge0}\hookrightarrow \CAlg_K^{dg}$
is a right Quillen functor. Its left adjoint $\tau:\CAlg_K^{dg}\to \CAlg_K^{dg,\ge0}$ carries
$C$ to the quotient of $C$ by the differential graded ideal
generated by elements $x\in C^i$ for $i<0$.
Namely, we have a Quillen adjunction $\tau:\CAlg_K^{dg}\rightleftarrows \CAlg_K^{dg,\ge0}$. We shall write
$\CAlg_K^{\ge0}$ for the $\infty$-category obtained from
the full subcatgory of cofibrant objects in $\CAlg_K^{dg,\ge0}$
by inverting weak equivalences.
The Quillen adjunction induces an adjunction of $\infty$-categories
\[
\tau:\CAlg_K\rightleftarrows \CAlg_K^{\ge0}
\]
\cite{Maz} where by ease of notation we write $\tau$ also for
the induced left adjoint functor
$\CAlg_K\rightarrow \CAlg_K^{\ge0}$.
We put $G=\Spec C$ with $C\in \CAlg_K$.
The functor $\tau$ preserves colimits, especially coproducts.
We put $\Aff_K^{\ge0}=(\CAlg_K^{\ge0})^{op}$.
We write $\Spec R$ for the object in $\Aff_K^{\ge0}$ corresponding to
$R\in \CAlg_K^{\ge0}$.
Then $\Spec \tau C$ inherits a group structure from $G=\Spec C$.
Namely, $\Spec \tau C$ is a group object
in $\Aff_K^{\ge0}$.
There is a fully faithful left adjoint $\CAlg_K^{\DIS}\to \CAlg_K^{\ge0}$
induced by the natural inclusion from the category of ordinary
commutative $K$-algebras to $\CAlg_K^{dg,\ge0}$.
Its right adjoint
$\CAlg_K^{\ge0} \to \CAlg_K^{\DIS}$
is given by taking the cohomology $C \mapsto H^0(C)$.
The inclusion $\CAlg_K^{\DIS}\to \CAlg_K$ is canonically
equivalent to
the composite $\CAlg_K^{\DIS}\to \CAlg_K^{\ge0}\to \CAlg_K$.
Also, the left adjoint $\tau$ is compatible with
inclusions
$\CAlg_K^{\DIS}\subset \CAlg_K$ and $\CAlg_K^{\DIS}\subset \CAlg^{\ge0}_K$
(use the fact that any object $C$ in $\CAlg_K^{\ge0}$ is the limit of a
cosimplicial diagram of ordinary $K$-algebras).
Consider $G$ to be the functor $\CAlg_K\to \Grp(\SSS)$.
Its restriction $G^\circ:=G|_{\CAlg_K^{\DIS}}:\CAlg_K^{\DIS}\to \Grp(\SSS)$
is naturally equivalent to the functor given by $A\mapsto \Map_{\CAlg_K^{\ge0}}(\tau C,A)$.
Take the cohomology $H^0(\tau C)$.
The structure of a commutative Hopf ring spectrum on $\tau C$
over $K$ (that is, the ``dual'' of the group structure
on $\Spec \tau C$ in $\Aff_K^{\ge0}$, see \cite[Appendix]{Tan}) gives the structure of a commutative Hopf ring on $H^0(\tau C)$ over $K$.
Namely, the comultiplication $\tau C\to \tau C\otimes_K\tau C$,
the counit $\tau C\to K$ and the antipode give rise to the structure of comultiplication
$H^{0}(\tau C)\to H^0(\tau C\otimes_K\tau C)\simeq H^0(\tau C)\otimes_K H^0(\tau C)$ of $H^0(\tau C)$, etc. We denote the associated affine group scheme
by $\overline{G}=\Spec H^0(\tau C)$.
We shall refer to $\overline{G}$ as the {\it underlying affine group
scheme} of $G$ (or the {\it coarse moduli space} for $G$ as in \cite{Tan}).
The assignment $G\mapsto \overline{G}$
is functorial and we actually have a
functor $\Grp(\Aff_K)\to \Grp(\Aff_K^{\DIS})$ which sends $G$ to
the associated affine group scheme $\overline{G}$.
By the adjunction,
the natural morphism $\pi:\Spec \tau C\to \overline{G}=\Spec H^0(\tau C)$ is universal among morphisms to ordinary
affine schemes over $K$ in $\hhh(\Aff_K^{\ge0})$
(note that $\Aff_K^{\ge0}$ contains $\Aff_K^{\DIS}$ as a full subcategory).
Namely, if $\phi:\Spec \tau C \to H$ is a
morphism to an ordinary affine scheme
$H$ in $\hhh(\Aff_K^{\ge0})$, there is a unique morphism $\psi:\overline{G}\to H$
such that $\phi=\psi\circ \pi$. 
In addition, $H$ is an affine group scheme and $\phi:\Spec \tau C \to H$ is a
homomorphism to the affine group scheme over $K$,
then there is a unique homomorphism $\psi:\overline{G}\to H$ in $\hhh(\Grp(\Aff_K^{\ge0})$
such that $\phi=\psi\circ \pi$.

\vspace{2mm}

As mentioned above,
we define $MG_E$ to be the underlying affine group scheme of $\MG_E$.
For the properties of $MG_E$ we refer to \cite{Tan}, \cite{Bar}, \cite{DTD},
\cite{PM}.

We define $\GG^{(n)}(X,x):=\GG_E^{(n)}(X,x)$ to be the underlying affine group scheme of $G_E^{(n)}(X,x)$ (cf. Section~\ref{realhopfsection}).

We consider a geometric interpretation of $\GG^{(n)}(X,x)$.
Suppose that $K=\QQ$ and
the base field $k$ is embeded in $\CC$.
We consider the case when the realization functor
is  associated to singular cohomology theory.

\begin{Proposition}
\label{undcomp}
Let $(X,x:\Spec k \to X)$ be a pointed smooth variety over $k$.
Let $\pi_i(X^t,x)$ be the homotopy group of the underlying topological
space $X^t=X\times_{\Spec k}\Spec \CC$.
For any $n\ge1$, the affine group schemes $\overline{G}^{(n)}(X,x)$
is a unipotent group scheme (i.e., a pro-unipotent algebraic group).
Moreover, $\overline{G}^{(1)}(X,x)$ is the pro-unipotent completion of
 $\pi_1(X^t,x)$ over $K=\QQ$.
Suppose further that the topological space
$X^t$ is nilpotent and of finite type (e.g. simply connected smooth varieties).
Then $\overline{G}^{(n)}(X,x)$ is a pro-unipotent completion of $\pi_n(X^t,x)$
for $n\ge2$.
\end{Proposition}

Before proceeding the proof,
we briefly recall
the notion of affinization (affination in French) studied in \cite{AS}
(in \cite{AS}, cosimplicial algebras are used instead of dg algebras,
see \cite[6.4]{Fre} for the comparison as a Quillen
equivalence between the model category of cosimplicial algebras
and $\CAlg_K^{dg,\ge0}$).
Let $\widehat{\CAlg}^{\ge0}_K$ be the $\mathbb{V}$-version of
$\CAlg^{\ge0}_K$ (cf. Section~\ref{underlyinggroup}).
Write $\widehat{\textup{Aff}}^{\ge0}_K:=(\widehat{\CAlg}^{\ge0}_K)^{op}$.
We write $\Spec R$ for an object of $\widehat{\textup{Aff}}^{\ge0}_K$
corresponding to $R\in \widehat{\CAlg}^{\ge0}_K$.
There is an adjunction
\[
\mathcal{O}:\Fun(\CAlg^{\textup{dis}}_K,\widehat{\SSS}) \rightleftarrows \widehat{\Aff}^{\ge0}_K
\]
where $\mathcal{O}$ is a left Kan extension of
the inclusion
$\Aff^{\textup{dis}}_K\hookrightarrow \widehat{\Aff}^{\ge0}_K$
along the Yoneda embedding $\Aff^{\textup{dis}}_K\to \Fun(\CAlg^{\textup{dis}}_K,\widehat{\SSS})$
(cf. \cite[Section 2.2]{AS}).
The right adjoint sends $R\in \CAlg_K^{\ge0}$ to the functor $h_R:\CAlg^{\textup{dis}}_K\to \widehat{\SSS}$ informally given by
$A\mapsto \Map_{\widehat{\CAlg}^{\ge0}_K}(R,A)$.
The restriction $\Aff^{\ge0}_K=(\CAlg_K^{\ge0})^{op}\to \Fun(\CAlg^{\textup{dis}}_K,\widehat{\SSS})$ of the right adjoint is fully faithful.
Let $F:\CAlg_K^{\DIS}\to \wSSS$ be a functor.
If $\mathcal{O}(F)$ belongs to $\Aff_{K}^{\ge0}$ (not to $\widehat{\Aff}_{K}^{\ge0}$), we refer to $\mathcal{O}(F)$
as the affinization of $F$.
An object $P$ in $\SSS$
can be viewed as the constant functor $\CAlg_{K}^{\DIS}\to \widehat{\SSS}$
with value $P$.
One can consider the affinization of the space $P\in \SSS$.
The composite
$\SSS=\Fun(\Delta^0,\SSS)\to \Fun(\CAlg_K^{\DIS},\widehat{\SSS}) \to \widehat{\Aff}^{\ge0}_K$
preserves small colimits and sends a contractible space to $\Spec K$,
where the first arrow is the functor given by the composition
with $\CAlg_K^{\DIS}\to \Delta^0$.
Consequently, the composite carries the space $S \in\SSS$ to
$\Spec K^{S}$ where $K^S$ is the cotensor with the space $S$.
By Proposition~\ref{PLproperty} and Remark~\ref{PLremark},
we conclude that $\SSS\to \Aff_K^{\ge0}\to \Aff_K$
is equivalent to $A_{PL,\infty}$.
By Theorem~\ref{real1}, $\Spec T_X$ in $\Aff_{K}$ is the affinization of $X^t$.

\vspace{2mm}

\Proof
There are several ways to prove the assertion, and we will give one
of them.
We treat the case $n=1$.
Let $G^{(1)}(X,x)^{\circ}:\CAlg_{K}^{\DIS}\subset\CAlg_K\to \Grp(\wSSS)$
denote the restriction.
It carries $A$ to $\Omega_*\Spec T_X(A)$, where $\Spec T_X(A)$
is the space of $A$-valued points on $\Spec T_X$, and $\Omega_*\Spec T_X(A)$
is its base loop space (the base point comes from $x_A:\Spec A\to \Spec K\to \Spec T_X$). 
We let $G_\circ^{(1)}(X,x):\CAlg_K^{\DIS}\to \Grp(\textup{Set})$
be the sheaf of groups with respect to fpqc topology associated to
the presheaf $A \mapsto \pi_0(\Omega_*\Spec T_X(A))\simeq \pi_1(\Spec T_X(A),x_A)$.
Then according to \cite[2.4.5]{AS}
(or \cite[4.4.8]{DAG}), $G_\circ^{(1)}(X,x)$
is represented by a unipotent affine group scheme
(i.e., a pro-unipotent algebraic group).
(We remark that
there is a canonical equivalence $\Map_{\CAlg_K}(T_X, A)\simeq \Map_{\CAlg_K^{\ge0}}(K^{X^t},A)$ for any $A\in \CAlg_K^{\DIS}$,
see \cite[7.2]{Ol}.)
Note that the natural morphism $G^{(1)}(X,x)^{\circ} \to G_\circ^{(1)}(X,x)$
is universal among morphisms to sheaves of groups on $\CAlg_K^{\DIS}$.
On the other hand, there is the natural map $G^{(1)}(X,x)^{\circ}\to \overline{G}^{(1)}(X,x)$
(recall that if $G^{(1)}(X,x)=\Spec C$, the restriction
$G^{(1)}(X,x)^{\circ}$ is represented by $\Spec \tau C$).
Consequently, by the universal property
there is a natural morphism $G_\circ^{(1)}(X,x)\to\overline{G}^{(1)}(X,x)$ of affine group schemes over $K$.
We wish to show that it is an isomorphism.
Since $K=\QQ$ is characteristic zero and
$G_\circ^{(1)}(X,x)\to \overline{G}^{(1)}(X,x)$
is a morphism as affine group schemes over $K$,
it is enough to prove that
for any algebraically closed field $L$, the induced map 
$G_\circ^{(1)}(X,x)(L)\to \overline{G}^{(1)}(X,x)(L)$
of sets of $L$-valued points is bijective.
In fact, according to \cite[Theorem 5.17]{Tan}
(its proof that works also for $G^{(1)}(X,x)$ instead of $\MG_E$) and \cite[VIII 4.4.8]{DAG},
we see that
$G_\circ^{(1)}(X,x)(L)\to \overline{G}^{(1)}(X,x)(L)$ is bijective.
It follows that $\overline{G}^{(1)}(X,x)$
is a unipotent affine group scheme.
By \cite[2.4.11]{AS} and Theorem~\ref{real1}, the group
scheme
$G_\circ^{(1)}(X,x)\simeq \overline{G}^{(1)}(X,x)$ is naturally isomorphic to a pro-unipotent completion of $\pi_1(X^t,x)$ (that is
endowed with the morphism
form the constant functor with value $\pi_1(X^t,x)$).
The case of $n\ge 2$ is similar.
If $G^{(n)}(X,x)^{\circ}:\CAlg_{K}^{\DIS}\subset\CAlg_K\to \Grp(\wSSS)$
denotes the restriction of $G^{(n)}(X,x)$,
it carries $A$ to the $n$-fold
loop space $\Omega_*^{n}\Spec T_X(A)$.
As in the case of $n=1$ (\cite[2.4.5]{AS}),
we observe that 
the sheaf associated to
the presheaf $A \mapsto \pi_n(\Spec T_X(A),x_A)$ is
isomorphic to $\overline{G}^{(n)}(X,x)$.
Then the final assertion follows from \cite[2.5.3]{AS}.
\QED

 \subsection{}
\label{finalaction}

We will construct an action of the motivic Galois group $MG_E$ on
the affine group scheme $\GG^{(n)}(X,x):=\GG_E^{(n)}(X,x)$.
Unfortunately, if one does not assume motivic conjectures
that imply the existence of a motivic $t$-structure,
it seems difficult to obtain an action of $MG_E$ on $\GG^{(n)}(X,x)$
from that of $\MG_E$ on $G^{(n)}(X,x)$ in a purely categorical way.
To overcome this issue, we use a method of homological algebras,
which yields
a natural action of $MG_E$ on $\GG_E^{(n)}(X,x)$.

For a usual affine group scheme $H$ over $K$,
we let $\Gamma(H)$ be the (ordinary) coordinate ring
on $H$, that is a commutative Hopf ring over $K$.
We let $\Aut(H):\CAlg_K^{\DIS}\to \Grp(\textup{Set})$
be the functor which assigns $A$ to the group of automorphisms
of the commutative Hopf ring
$\Gamma(H)\otimes_KA \stackrel{\sim}{\to} \Gamma(H)\otimes_KA$ over $A$.

\begin{Theorem}
\label{actiononhomotopy}
Let $(X,x)$ be a pointed smooth variety over $k$.
Then there is a (canonical)
morphism $MG_E\to \Aut(\overline{G}^{(n)}(X,x))$ in $\Fun(\CAlg_K^{\DIS},\Grp(\textup{Set}))$, that is,
an action of $MG_E$ on $\GG^{(n)}(X,x)$.
In other words, the action is described as an action on the scheme $\GG^{(n)}(X,x)$
\[
\GG^{(n)}(X,x)\times MG_E \to \GG^{(n)}(X,x)
\]
which is compatible with the group structure. Moreover,
the following properties hold:
\begin{enumerate}
\renewcommand{\labelenumi}{(\theenumi)}

\item The action is functorial:
Let $\phi:(X,x)\to (Y,y)$ be a morphism of smooth varieties
over $k$ that sends $x$ to $y$.
Let $\phi_*:\overline{G}^{(n)}(X,x)\to \overline{G}^{(n)}(Y,y)$
be the induced morphism of group schemes.
Then the action of $MG_E$ commutes with $\phi_*$.

\item
The action has a moduli theoretic interpretation in a coarse sense
(see Remark~\ref{coarsemoduli}).
\end{enumerate}
\end{Theorem}

\begin{Corollary}
\label{actiononhomotopy2}
Suppose that $k$ is embedded in $\CC$ and consider the case of singular
realization.
Let $\pi_i(X^t,x)_{uni}$ be the pro-unipotent completion
of $\pi_i(X^t,x)$ over $\QQ$.
Then 
we have a canonical action
\[
\pi_1(X^t,x)_{uni} \times MG \to \pi_1(X^t,x)_{uni}.
\]
If $X^t$ is nilpotent and of finite type, there
is a canonical action of $MG$ on $\pi_n(X^t,x)_{uni}$
for $n\ge2$.
\end{Corollary}

\Proof
It follows from Theorem~\ref{actiononhomotopy} and Proposition~\ref{undcomp}.
\QED

{\it Construction of an action/Proof of Theorem~\ref{actiononhomotopy}.}
Let $G^{(n)}(X,x):\NNNN(\Delta)^{op}\to \Aff_K$
be the derived affine group scheme over $K$, associated
to $(X,x)$ (see Section~\ref{realhopfsection}).
Let $\Gamma(G^{(n)}(X,x))$ be the image of
$[1]$ under $G^{(n)}(X,x)^{op}:\NNNN(\Delta)\to \CAlg_K$.
(Namely, $\Gamma(G^{(n)}(X,x))$ is the underlying algebra of commutative Hopf algebra
object $G^{(n)}(X,x)^{op}$ in $\CAlg_K$.)
Let $\MG_E=\Spec C$.
The identity
$\MG_E\to \MG_E$ determines a component of the space
($\infty$-groupoid) $\MG_E(C)$.
The action on $G^{(n)}(X,x)$
(cf. Proposition~\ref{actiongroup} and Example~\ref{exgroup})
induces its image in
$\Aut(G^{(n)}(X,x))(C)$.
The equivalence class of the image in $\Aut(G^{(n)}(X,x))(C)$
gives rise to a morphism 
$\Gamma(G^{(n)}(X,x))\otimes_K C\stackrel{\sim}{\to} \Gamma(G^{(n)}(X,x))\otimes_K C$
in $\CAlg_C$ (cf. Remark~\ref{autogroupR}).
Composing with the unit $K\to C$, we have
\[
\theta:\Gamma(G^{(n)}(X,x))=\Gamma(G^{(n)}(X,x))\otimes_K K \to \Gamma(G^{(n)}(X,x))\otimes_K C\stackrel{\sim}{\to} \Gamma(G^{(n)}(X,x))\otimes_K C.
\]
The composite is a coaction of $C$ on $\Gamma(G^{(n)}(X,x))$
at the level of the homotopy category $\hhh(\CAlg_K)$.
Namely, if we think of $C$
as a coalgebra in $\hhh(\CAlg_K)$
determined by the class of comultiplication
$C\to C\otimes_KC$ and the unit $C\to K$,
then 
$\Gamma(G^{(n)}(X,x))\to\Gamma(G^{(n)}(X,x))\otimes_K C$
is an (associative) coaction on $\Gamma(G^{(n)}(X,x))$
in the obvious sense. Also, it commutes with the structure of
coalgebra on $\Gamma(G^{(n)}(X,x))$ at the level of homotopy category.
Let $B:=\tau C$ (see Section~\ref{underlyinggroup} for $\tau$.
Applying $\tau$ to $\theta$ we obtain
\[
\rho:\tau(\Gamma(G^{(n)}(X,x)))\to \tau(\Gamma(G^{(n)}(X,x)))\otimes_K B\stackrel{\sim}{\to} \tau(\Gamma(G^{(n)}(X,x)))\otimes_K B.
\]
Taking the cohomology in the $0$-th term we have
\begin{eqnarray*}
\xi:H^0(\tau(\Gamma(G^{(n)}(X,x))))\to H^0(\tau(\Gamma(G^{(n)}(X,x)))\otimes_K B) &\stackrel{\sim}{\to}& H^0(\tau(\Gamma(G^{(n)}(X,x)))\otimes_K B) \\
&\simeq& H^0(\tau(\Gamma(G^{(n)}(X,x)))\otimes H^0(B).
\end{eqnarray*}
Recall that the commutative Hopf
ring $\Gamma(\overline{G}^{(n)}(X,x))$ of the coordinate ring on $\overline{G}^{(n)}(X,x)$
is $H^0(\tau(\Gamma(G^{(n)}(X,x))))$ equipped with the structure of
commutative Hopf ring that comes from the structures on $\Gamma(G^{(n)}(X,x))$.
Moreover, $MG_E=\Spec H^0(B)$.
The morphism $\xi$ is a coaction of $H^0(B)$ on the commutative $K$-algbera $H^0(\tau(\Gamma(G^{(n)}(X,x))))=\Gamma(\overline{G}^{(n)}(X,x))$ which is
compatible with the structure of coalgebra
on $H^0(\tau(\Gamma(G^{(n)}(X,x))))$.
It gives rise to an action
\[
\GG^{(n)}(X,x)\times MG_E \to \GG^{(n)}(X,x).
\]
The functoriality (1) is obvious from the construction.
\QED

\begin{Remark}
\label{coarsemoduli}
The affine group scheme $MG_E$ is a coarse moduli
space for $\MG_E$.
It has a coarse moduli theoretic interpretation:
for any field $L$ over $K$, $\MG_E^\circ \to MG_E$ induces
an isomorphism $\pi_0(\MG_E(L))\stackrel{\sim}{\to} MG_E(L)$
of sets where $\pi_0(\MG_E(L))$ is the set of connected components,
i.e., the set of equivalence classes of $L$-valued points on $\MG_E$
(cf. \cite[Theorem 1.3]{Tan}).
By $\MG_E\simeq \Aut(\mathsf{R}_E)$,
the set $MG_E(K)$ is naturally identified with
the set of equivalence classes of the automorphism
of $\mathsf{R}_E:\DM^\otimes(k)\to\Mod^\otimes_K$.
Suppose that $q\in MG_E(K)$ corresponds to an automorphism $\sigma$
of $\mathsf{R}_E$.
The automorphism of $\GG^{(n)}(X,x)$ induced by $q$
is the automorphism induced by $\sigma$.
Recall that $\sigma$ induces an automorphism of the multiplicative
realization functor
$\CAlg(\mathsf{R}_E):\CAlg(\DM^\otimes(k))\to\CAlg_K$ (cf. Section~\ref{actionstep1}).
It gives rise to an automorphism
 on $G^{(n)}(X,x)^{op}:\NNNN(\Delta)\to \CAlg_K$
(cf. Section~\ref{diagramaction}).
The induced automorphism $\Gamma(G^{(n)}(X,x))\stackrel{\sim}{\to}\Gamma(G^{(n)}(X,x))$
gives rise to
$a:H^0(\tau\Gamma(G(X,x)))\stackrel{\sim}{\to}H^0(\tau\Gamma(G(X,x)))$.
By our construction, the action of $q$ is equal to $a$.
This interpretation holds also for any field $L$ over $K$.
\end{Remark}

\section{Sullivan models and computational results}
\label{Sullivanmodel}

In rational homotopy theory, an inductive construction of a Sullivan model
is quite powerful.
Let $S$ be a topological space and $A_{PL}(S)$ the
commutative dg algebra of polynomial differential forms.
As in Section~\ref{realization}
we write $A_{PL,\infty}(S)$ for the image of $A_{PL}(S)$ in $\CAlg_{\QQ}$.
Let $\Free_{\QQ}$ denote the free functor $\Mod_{\QQ}\simeq \mathsf{D}(\QQ)\to \CAlg_{\QQ}$ which is defined to be a left adjoint to the forgetful functor
$\CAlg_{\QQ}\to \mathsf{D}(\QQ)$.
Contrary to genuine commutative dg algebras,
in the setting of $\CAlg_{\QQ}$ it is nonsense to say 
what a underlying graded algebra is.
But in the language of $\CAlg_{\QQ}$, the inductive construction describes
$A_{PL,\infty}(S)$ as a colimit of a sequence
\[
\QQ\simeq A_0\to A_1\to \cdots \to A_n\to A_{n+1}\to \cdots
\]
such that for any $n\ge 0$, $A_{n+1}$ fits in the pushout diagram of the form
\[
\xymatrix{
\Free_{\QQ}(V) \ar[r] \ar[d] & A_n \ar[d] \\
\QQ \ar[r] & A_{n+1}
}
\]
in $\CAlg_{\QQ}$ where $V$ is a $\ZZ$-graded vector space over $\QQ$
regarded as an object in $\mathsf{D}(\QQ)$, and the vertical arrow
is $\Free_{\QQ}(V)\to \Free_{\QQ}(0)\simeq \QQ$ induced by $V\to 0$.
Note that $\Free_{\QQ}(V)\to A_n$ is determined by a morphism $V\to A_n$
in $\mathsf{D}(\QQ)$.
Suppose that $V$ is concentrated in a fixed positive degree $n$, i.e.,
$V^i=0$ for $i\ne n$, and the $\QQ$-vector space $V^n$ is finite dimensional.
Then $\Free_{\QQ}(V)$ is
the commutative dg algebra that corresponds to the rational homotopy type
of the Eilenberg-MacLane space $K((V^n)^\vee,n)$.
Informally, the above sequence may be thought of as a presentation
of $A_{PL}(S)$ as a ``successive extension'' of ``simple pieces''
of the form $\Free_{\QQ}(V[1])$.

We will apply this approach to $\CAlg(\DM^\otimes(k))$
and study motivic cohomological algebras.
Free commutative algebra objects in $\DM(k)$ play the role of
free commutative dg algebras.
Actually, from the tannakian viewpoint,
such free objects are quite ``simple'' objects, see Remark~\ref{freemot}.
Put another way,
presentations of successive extensions by free objects
is useful for computations of a motivic counterpart
of rational homotopy groups.
We will inroduce the notion of
cotangent motives in Section~\ref{cotangenthomotopy}.
We then apply
the study of structures of cohomological motivic algebras in this Section
to obtain explicit descriptions of cotangent motives (Theorem~\ref{motivichomotopyexp}).

In this Section, we work with rational coefficients,
but $\QQ$ can be replaced by any field $K$ of characteristic zero.

\subsection{}
\label{Sullivanmodel0}
We will study some ``elementary examples'' such as projective spaces.
We also 
hope that the reader will get the feeling of 
the idea of the
constructions of ``Sullivan models'' of motivic cohomological algebras
in $\CAlg(\DM^\otimes(k))$.

Recall free commutative algebras in a general setting.

\begin{Definition}
\label{freealg}
Let $\mathcal{C}^\otimes$ be a symmetric monoidal
$\infty$-category
that has small colimits and
the tensor product $\otimes:\mathcal{C}\times \mathcal{C}\to \mathcal{C}$
preserves small colimits separately in each variable.
Let $u_{\mathcal{C}}:\CAlg(\mathcal{C}^\otimes)\to \mathcal{C}$
be the forgetful functor.
By \cite[3.1.3]{HA}, there exists a left adjoint
\[
\Free_{\mathcal{C}}:\mathcal{C} \longrightarrow \CAlg(\mathcal{C}^\otimes)
\]
to $u_\mathcal{C}$, which we shall call the free functor of
$\mathcal{C}^\otimes$ (\cite{HA} treats a broader setting).
Given $C\in \mathcal{C}$ we refer to $\Free_{\mathcal{C}}(C)$
as the free commutative algebra (object) generated by $C$.
We often omit the notation $u_{\mathcal{C}}$.

For $A\in\CAlg(\mathcal{C}^\otimes)$,
by the adjunction,
a morphism $f:\Free_{\mathcal{C}}(C)\to A$
corresponds to the composite
$\alpha:C\stackrel{\textup{unit}}{\to} u_{\mathcal{C}}(\Free_{\mathcal{C}}(C))\stackrel{u_{\mathcal{C}}(f)}{\to} u_{\mathcal{C}}(A)$ in $\mathcal{C}$. 
We say that $f:\Free_{\mathcal{C}}(C)\to A$
is classified by $\alpha$.

According to \cite[3.1.3.13]{HA}, the underlying object
$\Free_{\mathcal{C}}(C)$ is equivalent to the coproduct
$\sqcup_{n\ge0} \Sym^n(C)$
in $\mathcal{C}$, where $\Sym_{\mathcal{C}}^n(C)$ is the $n$-fold symmetric product (we usually omit the subscript when the setting is obvious).
If $\mathcal{D}^\otimes$ is a symmetric monoidal
$\infty$-category having the same property and 
$F:\mathcal{C}^\otimes\to \mathcal{D}^\otimes$
is a colimit-preserving functor, then there is a canonical
equivalence 
$\Free_{\mathcal{D}}(F(C))\stackrel{\sim}{\to} F(\Free_{\mathcal{C}}(C))$
for any $C\in \mathcal{C}$.
\end{Definition}

\vspace{2mm}

\subsubsection{}
We start with results that are useful for computations.
Let $\GL_{d}$ be the general linear algebraic group over $\QQ$.
Let $\Vect^\otimes(\GL_d)$ be the symmetric monoidal
abelian category of (possibly
infinite dimensional) representations of $\GL_d$, that is,
$\QQ$-vector spaces
with action of $\GL_d$.
The symmetric monoidal category
$\Comp(\GL_d):=\Comp(\Vect(\GL_d))$ of (possibly unbounded) cochain complexes 
admits a proper combinatorial symmetric monoidal model structure 
such that (i) $f:C\to C'$ is a weak equivalence
if a quasi-isomorphism, (ii) every object is cofibrant,
and (iii) $\{\iota_M:S^{n+1}M \hookrightarrow D^nM\}_{M\in I \atop \\ n\in \ZZ}$ is 
a set of generating cofibrations consisting of natural inclusions,
where $I$ is the set of irreducible
representations of $\GL_d$, and 
$S^nM$ (reps. $D^nM$) in $\Comp(\GL_d)$ defined by $(S^nM)^n=M$ and $(S^nM)^m=0$
for $m\ne n$ (resp. $(D^nM)^{n}=(D^{n}M)^{n+1}=M$, $D^mM=0$
for $m\ne n, n+1$, and $d:(D^{n}M)^{n}\to (D^{n}M)^{n+1}$ is the identity),
see \cite[Section 2.3]{DTD}, \cite[Corollary 3.5]{CD1} for details.
Let $\Rep^\otimes(\GL_d)$ be the symmetric
monoidal $\infty$-category, which
is obtained from $\Comp(\GL_d)$
by inverting quasi-isomorphisms.
Let $\CAlg(\Rep^\otimes(\GL_d))$ be the $\infty$-category
of commutative algebra objects in $\Rep^\otimes(\GL_d)$.

\begin{Lemma}
\label{GLcommalg}
We denote by $\CAlg(\Comp(\GL_d))$ the category of commutative algebra
objects in $\Comp(\GL_d)$. (We may think of an object
as a commutative dg algebra equipped with action of $\GL_d$.)
Then there is a combinatorial model structure on
$\CAlg(\Comp(\GL_d))$ such that
a morphism $f:A\to A'$ in $\CAlg(\Comp(\GL_d))$
is a weak equivalence (resp. a fibration) if
$f$ is a weak equivalence (reps. a fibration) in the underlying category
$\Comp(\GL_d)$.
In addition, if $\CAlg(\Comp(\GL_d))[W^{-1}]$
denotes the $\infty$-category obtained from the full subcategory
of cofibrants in $\CAlg(\Comp(\GL_d))$
by inverting weak equivalences, then the canonical functor
\[
\CAlg(\Comp(\GL_d))[W^{-1}]\to \CAlg(\Rep^\otimes(\GL_d))
\]
is an equivalence of $\infty$-categories.
\end{Lemma}

\Proof
Thanks to \cite[4.5.4.4, 4.5.4.6, 4.5.4.7]{HA}, it is enough
to prove that every cofibration in $\Comp(\GL_d)$ is a power cofibration
in the sense of \cite[4.5.4.2]{HA}.
To this end, we first observe that
a morphism $f:C\to C'$ in $\Comp(G):=\Comp(\Vect(G))$
is a cofibration if and only if $f$ is a monomorphism
when $G$ is either $\GL_d$ or a symmetric group $\Sigma_n$.
Let $M$ be an irreducible representation of $G$.
By the representation theory of $\GL_d$ or $\Sigma_n$,
$\Vect(G)$ is semi-simple and
$\Hom_{\Vect(G)}(M,M)=\QQ$ for any irreducible representation $M$ of
$G$. Let $\xi_M:\Vect(G)\to \Vect$ be the functor to the category of $\QQ$-vector spaces, that is given by $N\mapsto \Hom_{\Vect(G)}(M,N)$.
Taking the product indexed by the set $I(G)$ of irreducible representations
of $G$,
we have $\sqcap_{M\in I(G)} \xi_M:\Vect(G)\to \sqcap_{I(G)}\Vect$.
Note that this functor is an equivalence of categories
and induces an equivalence $\Comp(G)\to \sqcap_{I(G)}\Comp(\QQ)$
in the obvious way.
For an irreducible representation
$P$,  $S^{n+1}P \to D^nP$ corresponds to a morphism $\{f_M\}_{M\in I(G)}$ in 
$\sqcap_{I(G)}\Comp(\QQ)$ such that
$f_P:S^{n+1}\QQ\to D^n\QQ$ and $f_M=0$ if $M\ne P$ through this equivalence.
Therefore, it will suffice to show that
the smallest weakly saturated class containing $\{S^{n+1}\QQ\to D^n\QQ\}_{n\in \ZZ}$
coincides with a collection of monomorphisms in $\Comp(\QQ)$.
In fact, $\{S^{n+1}\QQ\to D^n\QQ\}_{n\in \ZZ}$ is a set of
generating cofibrations
in the projective model structure of $\Comp(\QQ)$.
Since $\QQ$ is a field, a morphism in $\Comp(\QQ)$ is a cofibration
with respect to the projective model structure exactly when it is 
a monomorphism. Thus, we conclude that
a morphism $f:C\to C'$ in $\Comp(G)$
is a cofibration if and only if $f$ is a monomorphism.

Next we prove that a cofibration $f:C\to C'$ of $\Comp(\GL_d)$
is a power cofibration.
We say that $f$ is a power cofibration if
a $\Sigma_n$-equivariant map
$\wedge^n(f):\square^n(f)\to (C')^{\otimes n}$
is a cofibration in $\Comp(\GL_d)^{\Sigma_n}$ for any $n\ge 0$.
Here $\Comp(\GL_d)^{\Sigma_n}$ is the category of objects
in $\Comp(\GL_d)$
endowed with action of the symmetric group $\Sigma_n$,
which is equipped with the projective model structure.
We refer to \cite[4.4.4.1]{HA} for these definitions and notations.
Let $U:\Comp(\GL_d)\to \Comp(\QQ)$ be the forgetful functor, that
is a symmetric monoidal
left adjoint. It follows that
$\wedge^n U(f)\simeq U(\wedge^n(f))$.
Suppose that $f$ is a cofibration.
Then $U(f)$ is a cofibration with respect to the projective model
structure because it is a monomorphism. According to \cite[7.1.4.7]{HA},
$U(f)$ is a power cofibration. Thus
by the above consideration $U(\wedge^n(f))\simeq \wedge^n U(f)$ is a monomorphism.
Then $\wedge^n(f)$ is a monomorphism in $\Comp(\GL_d)$.
Note that there is a canonical
equivalence $(\sqcap_{I}\Comp(\QQ))^{\Sigma_n}
\stackrel{\sim}{\to} \sqcap_{I}(\Comp(\QQ)^{\Sigma_n})=\sqcap_{I}\Comp(\Sigma_n)$.
The image of $\wedge^n(f)$ in $\sqcap_{I}(\Comp(\QQ)^{\Sigma_n})$
is a monomorphism. Again by the above consideration, the image is
a cofibration in $\sqcap_{I}(\Comp(\QQ)^{\Sigma_n})$ endowed with
the projective structure.
Therefore, $\wedge^n(f)$ has the left lifting property
with respect to epimorphic quasi-isomorphisms
in $\Comp(\GL_d)^{\Sigma_n}$, namely, it is a cofibration.
\QED

Let $u:\CAlg(\Comp(\GL_d))\to \Comp(\GL_d)$
be the forgetful functor.
By the definition of the model structure on $\CAlg(\Comp(\GL_d))$
in Lemma~\ref{GLcommalg}, it is a right Quillen functor.
Let $\Free_{\Comp(\GL_d)}:\Comp(\GL_d)\to \CAlg(\Comp(\GL_d))$
be a left Quillen functor to $u$.
It is the free functor
of $\Comp(\GL_d)$. Since every object in $\Comp(\GL_d)$
is cofibrant, thus $\Free_{\Comp(\GL_d)}$ preserves weak equivalences;
that is to say, it is ``derived''.
Let $u_\infty:\CAlg(\Rep(\GL_d))\to \Rep(\GL_d)$ be the forgetful
functor
of $\infty$-categories. We write $\Free_{\Rep(\GL_d)}:\Rep(\GL_d)\to\CAlg(\Rep(\GL_d))$ for the free functor of $\Rep^\otimes(\GL_d)$.
The following Lemma guarantees compatibility between
$\Free_{\Rep(\GL_d)}$ and $\Free_{\Comp(\GL_d)}$.

\begin{Lemma}
\label{freecomp}
Let $C$ be an object in $\Comp(\GL_d)$.
By abuse of notation, we write $C$ (resp. $\Free_{\Comp(\GL_d)}(C)$)
for the images of the cofibrant object $C$ (resp. $\Free_{\Comp(\GL_d)}(C)$) in $\Rep(\GL_d)$
(resp. $\CAlg(\Rep(\GL_d))$).
Then there is a canonical equivalence
$\Free_{\Comp(\GL_d)}(C) \simeq \Free_{\Rep(\GL_d)}(C)$ in $\CAlg(\Rep(\GL_d))$,
which commutes with $C\to u_{\infty}(\Free_{\Comp(\GL_d)}(C))$ and
$C\to u_\infty(\Free_{\Rep(\GL_d)}(C))$.
\end{Lemma}

\Proof
The forgetful functors $u$ and $u_\infty$ commute with
canonical maps $\CAlg(\Comp(\GL_{d}))\to \CAlg(\Rep(\GL_d))$
and $\Comp(\GL_d)\to \Rep(\GL_d)$.
By Lemma~\ref{GLcommalg} we  identify the induced functor $\hhh(u_\infty):\hhh(\CAlg(\Rep(\GL_d)))\to \hhh(\Rep(\GL_d))$ of homotopy categories with the right adjoint
\[
\overline{u}:\hhh(\CAlg(\Comp(\GL_d))[W^{-1}])\to \hhh(\Rep(\GL_d))
\]
of homotopy categories
induced by the right Quillen functor $u$.
Thus, we can identify the left adjoint
$\hhh(\Free_{\Rep(\GL_d)}):\hhh(\Rep(\GL_d))\to \hhh(\CAlg(\Rep(\GL_d)))$
with the left adjoint $\hhh(\Rep(\GL_d))\to \hhh(\CAlg(\Comp(\GL_d))[W^{-1}])$
induced by $\Free_{\Comp(\GL_d)}$.
\QED

\begin{Proposition}
\label{hopushout}
Let $A$ be a cofibrant object in $\CAlg(\Comp(\GL_d))$
and let $\alpha:C\to u(A)$ be a morphism in $\Comp(\GL_d)$.
Let $\phi_\alpha:\Free_{\Comp(\GL_d)}(C)\to A$ be the morphism
classified by $\alpha$. Let $\iota:S^{0}\QQ \hookrightarrow D^{-1}\QQ$
be the cofibration in $\Comp(\GL_d)$, where $\QQ$ here denotes the
unit object in $\Comp(\GL_d)$ (we abuse notation).
Let $\Free_{\Comp(\GL_d)}(C)\to \Free_{\Comp(\GL_d)}(C\otimes (D^{-1}\QQ))$ be the morphism induced by $C\otimes\iota:C\simeq C\otimes (S^0\QQ) \to C\otimes(D^{-1}\QQ)$.
Let $A\langle \alpha \rangle$ be the pushout of the following diagram in
$\CAlg(\Comp(\GL_d))$:
\[
\xymatrix
{
\Free_{\Comp(\GL_d)}(C)\ar[r]^{\phi_{\alpha}} \ar[d] & A \ar[d] \\
\Free_{\Comp(\GL_d)}(C\otimes (D^{-1}\QQ)) \ar[r] & A\langle \alpha \rangle.
}
\]
Then this diagram is a homotopy pushout.
See Remark~\ref{explicitdg} for the explicit presentation of $A\langle \alpha \rangle$.
\end{Proposition}

\begin{Remark}
\label{explicitdg}
The commutative algebra object $A\langle \alpha \rangle$
is regarded as a commutative dg algebra endowed with an action of $\GL_d$.
The explicit presentation of $A\langle \alpha \rangle$ is described as follows
(see the proof of Proposition~\ref{hopushout}).
For simplicity, we suppose that differential of $C$ is zero and
we view it as a graded vector space with an action of $\GL_d$.
This assumption is not essential in practice because $\Vect(\GL_d)$ is
semi-simple. Let $\overline{A}$ be the underlying
graded algebra
of $A$ obtained by forgetting the differential.
The underlying graded algebra of $A\langle \alpha \rangle $ is 
given by the tensor product
$\overline{A}\otimes \Free_{\Comp(\GL_d)}(C[1])$
of commutative graded algebras with the action of $\GL_d$.
If one forgets the action of $\GL_d$ on $\Free_{\Comp(\GL_d)}(C[1])$,
then it is the free commutative graded algebra generated by
the underlying graded algebra of $C[1]$.
The differential
on $\overline{A}\otimes \Free_{\Comp(\GL_d)}(C[1])$
is given by the differential on $A$ and $\partial|_{C}=\alpha$.
When $\GL_d$ is the trivial, i.e., $d=0$ or one forgets the action of
$\GL_d$,
then the construction of $A\langle \alpha \rangle$ is
classical, see \cite[2.2.2]{Hin}.
\end{Remark}

\begin{Example}
\label{Sullex}
Let $\mathbb{G}_m=\GL_1$ and let $\chi_{i}$ in $\Comp(\mathbb{G}_m)$
be one dimensional representation of $\mathbb{G}_m$ of weight $i$
placed in degree zero.
Let $A=\Free_{\Comp(\mathbb{G}_m)}(\chi_1[-2])$ be the 
free commutative algebra generated by $\chi_1[-2]$.
The underlying cochain complex is $\oplus_{i\ge0}\chi_i[-2i]$
with zero differential.
Let $\alpha:\chi_{n+1}[-2n-2]\to \oplus_{i\ge0}\chi_i[-2i]=A$
be the canonical inclusion.
Let us consider $A\langle \alpha \rangle$.
Note that $\Free_{\Comp(\mathbb{G}_m)}(\chi_{n+1}[-2n-1])$ is the trivial square zero extension $\chi_0\oplus \chi_{n+1}[-2n-1]$
by $\chi_{n+1}[-2n-1]$ (since the generator is in the odd degree).
The underlying graded algebra is $(\oplus_{i\ge0}\chi_i[-2i])\otimes (\chi_0\oplus \chi_{n+1}[-2n-1])$. The non-zero part of differential
is given by ``identities'' $\chi_i[-2i]\otimes \chi_{n+1}[-2n-1]\to \chi_{i+n+1}[-2i-2n-2]\otimes \chi_0$ for $i\ge0$.
\end{Example}

The standard consequence of Proposition~\ref{hopushout} is

\begin{Corollary}
The image of the square diagramin Proposition~\ref{hopushout}
in $\CAlg(\Rep(\GL_d))$ is a pushout diagram.
We remark that the image of $\Free_{\Comp(\GL_d)}(C)$ 
and $\Free_{\Comp(\GL_d)}(C\otimes (D^{-1}K))$ in $\CAlg(\Rep(\GL_d))$ are
equivalent to $\Free_{\Rep(\GL_d)}(C)$ and the unit algebra, respectively
(Lemma~\ref{freecomp}).
\end{Corollary}

{\it Proof of Proposition~\ref{hopushout}.}
Let $B$ be a pushout of $C\otimes D^{-1}\QQ \leftarrow C\to u(A)$ in $\Comp(\GL_d)$, that is the standard mapping cone $(u(A)\oplus C[1], d)$
of $\alpha:C\to u(A)$. Since $u(A)$ is cofibrant
and $C\otimes S^0K\to C\otimes D^{-1}K$
is a cofibration, $B$ is a homotopy pushout, see e.g. \cite[A.2.4.4]{HTT}.
Then we have the commutative diagram
\[
\xymatrix
{
\Free_{\Comp(\GL_d)}(C) \ar[r] \ar[d] & \Free_{\Comp(\GL_d)}(u(A)) \ar[r] \ar[d]& A \ar[d] \\
\Free_{\Comp(\GL_d)}(C\otimes D^{-1}K) \ar[r] &  \Free_{\Comp(\GL_d)}(B)\ar[r]  & A\langle \alpha \rangle
}
\]
that consists of pushout squares. The upper right horizontal map is the
counit map.
Since $A$ is cofibrant and 
the left vertical arrow is an cofibration, 
again by \cite[A.2.4.4]{HTT} both left and right
(and the outer) squares are homotopy pushouts,
as claimed. The explicit structure of $A\langle \alpha \rangle$ in Remark~\ref{explicitdg}
can easily be seen from the right pushout.
\QED

\subsubsection{}
\label{cps}
We will consider the $n$-dimensional projective space $\PP^n$
over a perfect field $k$.

We denote by $\Free_{\DM(k)}:\DM(k)\to \CAlg(\DM^\otimes(k))$
the free functor of $\DM^\otimes(k)$.
For ease of notation, we put $\Free:=\Free_{\DM(k)}$.

By projective bundle theorem, there is a decomposition 
\[
M_{\PP^n}\simeq M(\PP^n)^\vee\simeq \uni_k\oplus \uni_k(-1)[-2]\oplus \ldots \oplus \uni_k(-n)[-2n]=\oplus_{i=0}^n\uni_k(-i)[-2i]
\]
in $\DM(k)$, see e.g. \cite[Lec.15]{MVW}.
Consider the inclusion $\iota:\uni_k(-1)[-2]\hookrightarrow M_{\PP^n}\simeq \oplus_{i=0}^n\uni_k(-i)[-2i]$ that is a morphism in $\DM(k)$.
It gives rise to a morphism
\[
f:\Free(\uni_k(-1)[-2])\to M_{\PP^n}
\]
in $\CAlg(\DM^\otimes(k))$,
 that is classified by $\iota$.
We note that $\Free(\uni_k(-1)[-2])\simeq \oplus_{i\ge0} \uni(-i)[-2i]$ in $\DM(k)$.
Observe that for $j>n$, the composite
\[
\uni_k(-j)[-2j]\hookrightarrow  \oplus_{i\ge0} \uni(-i)[-2i]\simeq \Free(\uni_k(-1)[-2])\to M_{\PP^n}
\]
is null homotopic. Indeed, $\uni_k(-j)[-2j]\to \uni_k(-i)[-2i]$ is null homotopic for $0\le i \le n$
since $\uni_k(j)[2j]\otimes(\uni_k(-j)[-2j]\to \uni_k(-i)[-2i])$
corresponds to
an element of motivic cohomology $H_{M}^{2j-2i}(\Spec k,j-i)\simeq \CH^{j-i}(\Spec k)=0$. Here $\CH^p(-)$ denotes the $p$-th
Chow group, and the comparison isomorphism between motivic cohomology 
and (higher) Chow groups is due to Voevodsky.
Next
we let
\[
g:\Free(\uni_k(-n-1)[-2n-2])\to \Free(\uni_k(-1)[-2])
\]
be a morphism that is classified by the inclusion $\uni_k(-n-1)[-2n-2]\hookrightarrow \Free(\uni_k(-1)[-2])$. Consider the
morphism $h:\Free(\uni_k(-n-1)[-2n-2])\to \Free(0)\simeq \uni_k$ induced by $\uni_k(-n-1)[-2n-2]\to 0$.
Take a pushout 
\[
S_{\PP^n}:=\Free_{\DM(k)}(\uni_k(-1)[-2])\otimes_{\Free_{\DM(k)}(\uni_k(-n-1)[-2n-2])}\uni_k
\]
along $h$ in $\CAlg(\DM^\otimes(k))$.
Note that $f\circ g$ factors through $h:\Free(\uni_k(-n-1)[-2n-2])\to \Free(0)\simeq \uni_k$ because
$\uni_k(-n-1)[-2n-2]\to M_{\PP^n}$ is null homotopic.
Consequently, by the universal property of the pushout
we obtain the induced morphism
\[
S_{\PP^n}\to M_{\PP^n}.
\]

\begin{Proposition}
\label{cpsprop}
The morphism $S_{\PP^n}\to M_{\PP^n}$
is an equivalence in $\CAlg(\DM^\otimes(k))$.
\end{Proposition}

\Proof
We first claim that
$\oplus_{i\ge0} \uni_k(-i)[-2i]\simeq \Free(\uni_k(-1)[-2])\to M_{\PP^n}\simeq \oplus_{i=0}^n \uni_k(-i)[-2i]$ induces an equivalence
$\Free(\uni_k(-1)[-2])\supset \uni_i(-i)[-2i]\stackrel{\sim}{\to}  \uni_i(-i)[-2i]\subset M_{\PP^n}$ for $0\le i\le n$.
As discussed before this Proposition, 
$\uni_i(-i)[-2i]\subset \Free(\uni_k(-1)[-2])\to M_{\PP^n}$ is null homotopic if $i>n$ because $\Hom_{\hhh(\DM(k))}(\uni_k(a)[2a],\uni_k(b)[2b])$
is $\QQ$ (resp. $0$) if $a=b$
(resp. $a\ne b$).
Consider the dual $M(\PP^n)\simeq \oplus_{i=0}^n \uni_k(i)[2i]$
of the isomorphism $M_{\PP^n}\simeq \oplus_{i=0}^n \uni_k(-i)[-2i]$.
Recall that the Chow ring $\CH^*(\PP^n)$ is isomorphic to
$\ZZ[H]/(H^{n+1})$ where $H\in \CH^1(\PP^n)$ is a class of
a hyperplane.
The projection
$M(\PP^n)\to \uni_k(i)[2i]$ corresponds to a generator of Chow group
$\QQ=\CH^i(\PP^n)\otimes_{\ZZ}\QQ\simeq H_{M}^{2i}(X,i)\simeq \Hom_{\hhh(\DM(k))}(M(\PP^n),\uni_k(i)[2i])$. Using scalar multiplication (if necessary),
we may and will assume that $M(\PP^n)\to \uni_k(i)[2i]$
corresponds to $H^i$.
Now we prove our claim by induction on $i$.
By the construction, the case of $i=1$ is clear.
We suppose that the case $i(< n-1)$ is true. We will show the case $i+1$.
By Lemma~\ref{homotopyfree}, 
$\Free(\uni_k(-1)[-2])$ in the homotopy category $\hhh(\DM(k))$
is also regarded as the free commutative algebra object lying in $\CAlg(\hhh(\DM^\otimes(k))$
generated by $\uni_k(-1)[-2]$ in $\hhh(\DM^\otimes(k))$.
Thus, the multiplication
map $\Free(\uni_k(-1)[-2])\otimes \Free(\uni_k(-1)[-2])\to \Free(\uni_k(-1)[-2])$ induces an isomorphism from the component $\uni_k(-a)[-2a]\otimes \uni_k(-b)[-2b]$ in the domain to $\uni_k(-a-b)[-2a-2b]$
in the target. Therefore, by the induction hypothesis
and the compatibility of multiplication maps, 
if the multiplication $M_{\PP^n}\otimes M_{\PP^n}\to M_{\PP^n}$
induces an isomorphism of the composite
\[
\xi:\uni_k(-1)[-2]\otimes \uni_k(-i)[-2i]\hookrightarrow M_{\PP^n}\otimes M_{\PP^n}\to M_{\PP^n}\to \uni_k(-i-1)[-2i-2],
\]
 then $\Free(\uni_k(-1)[-2])\to M_{\PP^n}$ induces an isomorphism
from the component $\uni_k(-i-1)(-2i-2)$ in the domain to $\uni_k(-i-1)[-2i-2]\subset M_{\PP^n}$ (namely, the case $i+1$ holds). 
Note that the dual $M(\PP^n)\to \uni_k(i)[2i]$
of $\uni_k(-i)[-2i]\to M_{\PP^n}$
corresponds to the element $H^i\in \CH^i(\PP^n)$ (for any $i$).
Observe that the dual $M(\PP^n)\to \uni_k(i+1)[2i+2]$ of the composite
$l:\uni_k(-1)[-2]\otimes \uni_k(-i)[-2i]\hookrightarrow M_{\PP^n}\otimes M_{\PP^n}\to M_{\PP^n}$ corresponds to the intersection product
$H^{i+1}=H\cdot H^i\in \CH(\PP^n)$.
To see this, recall that the product of motivic cohomology
\begin{eqnarray*}
\Hom_{\hhh(\DM(k))}(M(\PP^n),\uni_k(1)[2])\otimes \Hom_{\hhh(\DM(k))}(M(\PP^n),\uni_k(i)[2i])\ \ \ \ \ \ \ \  \\
\ \ \ \ \ \ \ \ \to \Hom_{\hhh(\DM(k))}(M(\PP^n),\uni_k(i+1)[2i+2])
\end{eqnarray*}
is induced by the composition with $M(\PP^n)\to M(\PP^n)\otimes M(\PP^n)$ defined by the diagonal map.
By Lemma~\ref{multcomult} below, the multiplication
$M_{\PP^n}\otimes M_{\PP^n}\to M_{\PP^n}$ is the dual of
$M(\PP^n)\to M(\PP^n)\otimes M(\PP^n)$.
In addition, the product structure on motivic cohomology
is compatible with that of
(higher) Chow groups via the comparison isomorphism \cite{KY}.
Therefore, we conclude that
the dual of $l$ corresponds to $H^{i+1}\in \CH^{i+1}(\PP^n)$.
It follows that $\xi$ is an isomorphism.

Next, by \cite[Theorem 3.1]{DTD}
there is a colimit-preserving symmetric monoidal
functor $F:\Rep^\otimes(\mathbb{G}_m)\to \DM^\otimes(k)$
which sends one dimensional representation $\chi_1$
of weight one placed in degree zero to $\uni_k(1)$.
Here $\mathbb{G}_m:=\GL_1$ and we denote by $\chi_p$ one dimensional
representation of weight $p$.
Let $\Free_{\Comp(\mathbb{G}_m)}(\chi_{-1}[-2])$
and $\Free_{\Comp(\mathbb{G}_m)}(\chi_{-n-1}[-2n-2])$
be the free commutative algebra in $\Comp(\mathbb{G}_m)$
generated by $\chi_{-1}[-2]$ and $\chi_{-n-1}[-2n-2]$, respectively.
Let $\Free_{\Comp(\mathbb{G}_m)}(\chi_{-n-1}[-2n-2])\to \Free_{\Comp(\mathbb{G}_m)}(\chi_{-1}[-2])$ be the morphism classified by the inclusion
$\alpha:\chi_{-n-1}[-2n-2]\hookrightarrow \Free_{\Comp(\mathbb{G}_m)}(\chi_{-1}[-2])$.
Take a homotopy pushout $ \Free_{\Comp(\mathbb{G}_m)}(\chi_{-1}[-2])\langle \alpha \rangle$, see Proposition~\ref{hopushout}.
By Remark~\ref{explicitdg}, the easy computation shows that
$ \Free_{\Comp(\mathbb{G}_m)}(\chi_{-1}[-2])\langle \alpha \rangle \simeq \oplus_{i=0}^n\chi_{-i}[-2i]$ in $\hhh(\Rep(\mathbb{G}_m))$
and the natural map
\[
\Free_{\Comp(\mathbb{G}_m)}(\chi_{-1}[-2])\simeq \oplus_{i\ge0}\chi_{-i}[-2i]\to \Free_{\Comp(\mathbb{G}_m)}(\chi_{-1}[-2])\langle \alpha \rangle\simeq \oplus_{i=0}^n\chi_{-i}[-2i]
\]
is the projection
(cf. Example~\ref{Sullex}).
By abuse of notation,
we will write $\chi_i$, $\Free_{\Comp(\mathbb{G}_m)}(\chi_{-1}[-2])$ and likes
also for their images in $\Rep(\mathbb{G}_m)$ or $\CAlg(\Rep(\mathbb{G}_m))$.
Note that $F$ sends the $\chi_i$ to $\uni_k(i)$ in $\DM(k)$.
The left adjoint functor
$\CAlg(F):\CAlg(\Rep(\mathbb{G}_m))\to \CAlg(\DM^\otimes(k))$
sends $\Free_{\Comp(\mathbb{G}_m)}(\chi_{-n-1}[-2n-2])\to \Free_{\Comp(\mathbb{G}_m)}(\chi_{-1}[-2])$
to $g$.
Then since $\CAlg(F)$ preserves a pushout,
$\Free_{\Comp(\mathbb{G}_m)}(\chi_{-1}[-2])\to \Free_{\Comp(\mathbb{G}_m)}(\chi_{-1}[-2])\langle \alpha \rangle$ maps to the canonical morphism
$\Free(\uni_k(-1)[-2]) \to S_{\PP^n}$.
We see that the composite
\[
\oplus_{i=0}^n\uni_k(-i)[-2i]\hookrightarrow \oplus_{i\ge0}\uni_k(-i)[-2i]\simeq \Free(\uni_k(-1)[-2]) \to S_{\PP^n}\simeq \oplus_{i=0}^n\uni_k(-i)[-2i]
\]
is an equivalence.
Taking account of the first claim of this proof, we see that
the underlying morphism $S_{\PP^n}\to M_{\PP^n}$ in $\DM(k)$ is 
an equivalence. Thus, $S_{\PP^n}\to M_{\PP^n}$ in $\DM(k)$ is an equivalence
in $\CAlg(\DM^\otimes(k))$.
\QED

\begin{Remark}
Suppose that the base field $k$ is embedded in $\CC$.
Let $\mathsf{R}:\CAlg(\DM^\otimes(k))\to \CAlg_{\QQ}$
be
the multiplicative realization functor considered in Section~\ref{realization}.
The multiplicative realization functor commutes with
free functors and the formulation of colimits.
Then the above construction of $S_{\PP^n}$ and the equivalence $S_{\PP^n}\simeq M_{\PP^n}$ is compatible with the classical construction of a Sullivan
model of $A_{PL}(\CC \PP^n)$ where $\CC \PP^n$ is the complex projective space.
The morphism $\mathsf{R}(\Free(\uni_k(-1)[-2]))\simeq \Free_{\QQ}(\QQ[-2])\to \mathsf{R}(M_{\PP^n})\simeq A_{PL,\infty}(\CC \PP^n)$ induced by $f$
is determined by a morphism $\QQ[-2]\to A_{PL,\infty}(\CC \PP^n)$
defined by a generator of $H^2(\CC\PP^n,\QQ)=\QQ$. This is
the first step of the construction of a Sullivan model.
The subsequent steps are also compatible. See e.g. \cite{He}.
Also, we remark that $\pi_{i}(\CC\PP^n)\otimes_{\ZZ}\QQ=\QQ$ if $i=2, 2n+1$,
and $\pi_{i}(\CC\PP^n)\otimes_{\ZZ}\QQ=0$ if otherwise.
See also Theorem~\ref{motivichomotopyexp} and Remark~\ref{motcpsR}.
\end{Remark}

\begin{Remark}
The object $M_{\PP^n}$ lies in the full subcategory of mixed Tate motives
in $\DM(k)$. But the above argument works for arbitrary perfect base
fields and does not need a (conjectural) motivic $t$-structure.
\end{Remark}

\subsubsection{}
\label{minus1}
Let $\mathbb{A}^n$ denote the $n$-dimensional affine space over a perfect
field $k$.
Let $X=\mathbb{A}^n-\{p\}$ be the open subscheme
of $\mathbb{A}^n$
that is obtained by removing a $k$-rational point
$p$.
Let $j:X\to \mathbb{A}^n$ be the open
immersion.
By the dual of the Gysin triangle \cite[14.5]{MVW},
we have a distinguished triangle
\[
\uni_k(-n)[-2n]\to M_{\mathbb{A}^n}\stackrel{j^*}{\to} M_{X}
\]
in the triangulated category $\hhh(\DM(k))$. Note that
$M_{\mathbb{A}^n}\simeq \uni_k$ and $\uni_k(-n)[-2n]\to M_{\mathbb{A}^n}$
is null homotopic (see the case in \ref{cps}).
Hence we have an equivalence $M_{X}\simeq \uni_k\oplus \uni_k(-n)[-2n+1]$
in $\DM(k)$.
We let $\Free(\uni_k(-n)[-2n+1])\to M_{X}$ be a morphism in
$\CAlg(\DM^\otimes(k))$, that is classified by the inclusion $\uni_k(-n)[-2n+1]\hookrightarrow \uni_k\oplus \uni_k(-n)[-2n+1]\simeq M_{X}$.

\begin{Proposition}
\label{sphere}
The morphism 
$\Free_{\DM(k)}(\uni_k(-n)[-2n+1])\to M_{X}$ is an equivalence.
\end{Proposition}

\Proof
We continue to use the notation in the proof of Proposition~\ref{cpsprop}
and the
colimit-preserving symmetric monoidal functor
$F:\Rep^\otimes(\mathbb{G}_m)\to \DM^\otimes(k)$.
Let $\Free_{\Comp(\mathbb{G}_m)}(\chi_{-n}[-2n+1])$
be the free algebra that belongs to $\CAlg(\Comp(\mathbb{G}_m))$
(keep in mind that it can be viewed as
a commutative dg algebra endowed with an action of $\mathbb{G}_m$).
Since the generator is concentrated in
the odd degree $2n-1$, by the Koszul sign rule
there is an isomorphism $\Free_{\Comp(\mathbb{G}_m)}(\chi_{-n}[-2n+1])\simeq \chi_0\oplus \chi_{-n}[-2n+1]$ as objects in $\Comp(\mathbb{G}_m)$.
The functor $F$ carries $\Free_{\Comp(\mathbb{G}_m)}(\chi_{-n}[-2n+1])$
to $\Free(\uni_k(-n)[-2n+1])$ in $\CAlg(\DM^\otimes(k))$.
Thus, the underlying object of $\Free(\uni_k(-n)[-2n+1])$
is equivalent to $\uni_k\oplus \uni_k(-n)[-2n+1]$.
Moreover, the canonical inclusion (unit map)
$\uni_k(-n)[-2n+1]\to \Free(\uni_k(-n)[-2n+1])$ is compatible with $\uni_k(-n)[-2n+1]\hookrightarrow \uni_k\oplus \uni_k(-n)[-2n+1]$.
Using these facts we deduce that $\Free(\uni_k(-n)[-2n+1])\simeq \uni_k\oplus \uni_k(-n)[-2n+1] \to M_{X}\simeq \uni_k\oplus \uni_k(-n)[-2n+1]$
is an equivalence, as desired.
\QED

\begin{Remark}
Suppose that the base field $k$ is embedded in $\CC$.
Then the complex manifold $X\times_{\Spec k}\Spec \CC$
is homotopy equivalent to the $(2n-1)$-dimensional sphere $S^{2n-1}$.
Proposition~\ref{sphere} is a motivic generalization
of the fact that the free commutative dg algebra generated
by one dimensional vector space placed in (cohomological) degree $2n-1$ is
a Sullivan model of $A_{PL}(S^{2n-1})$ (cf. \cite[Example 1 in page 142]{FHT}).
\end{Remark}

\subsubsection{}
\label{minus2}
Proposition~\ref{cpsprop} and~\ref{sphere}
gives explicit ``models'' $S_{\PP^n}$, $\Free_{\DM(k)}(\uni_k(-n)[-2n+1])$
of motivic cohomological algebras.
The constructions of models have only finitely many steps.
As in the classical rational homotopy theory,
an inductive construction often consists of infinite steps.
The following is such an example.

Let $Y=\mathbb{A}^n-\{p\}-\{q\}$ be the open subscheme
of $\mathbb{A}^n$ that is obtained by removing two $k$-rational points
$p,q$. Let $s:Y\to \Spec k$ denote the structure morphism.

\begin{Proposition}
\label{minus2prop}
Let $A_0=\uni_k$ be the unit algebra in $\CAlg(\DM^\otimes(k))$
and let $A_0=\uni_k \to M_{Y}$ be a unique morphism from the initial
object $\uni_k$ in $\CAlg(\DM^\otimes(k))$.
Then there is a refinement of $A_0\to M_{Y}$
\[
\uni_k=A_0\to A_1\to A_2\to \cdots \to A_i\to A_{i+1} \to \cdots\to M_Y
\]
that satisfies the following properties:
\begin{enumerate}
\renewcommand{\labelenumi}{(\theenumi)}

\item
The canonical morphism
$\varinjlim_{i\ge0}A_i\to M_Y$ is an equivalence.
Here $\varinjlim_{i}A_i$ be a colimit of the sequence in
$\CAlg(\DM^\otimes(k))$.

\item Let $V_i$ be the kernel (homotopy fiber) of
$A_i\to M_Y$ in $\DM(k)$ for any $i\ge0$. Then
for each $i\ge0$, $A_i\to A_{i+1}$ is of the form
$A_i\to A_i\otimes_{\Free(V_i)}\uni_k$ given by the pushout of 
$A_{i}\leftarrow \Free(V_i) \rightarrow \uni_k$
where $\Free(V_i)\to A_i$ is classified by $V_i\to A_i$.
\end{enumerate}
Moreover, for $n\ge2$, one can explicitly compute each $A_i$
 in the sense explained below.
\end{Proposition}

The first statement is a consequence of a
more general fact, see Lemma~\ref{Sullivanmodel1} below.
We explain the second statement, that is,
the procedure of an explicit computation.
We will compute the lower degrees $A_1$, $A_2$, $A_3$.
We can apply the same procedure and arguments also to higher
degrees and we leave it to the interested reader.

We continue to use the notation in Section~\ref{cps}, \ref{minus1}.
As in the case of $X=\mathbb{A}^n-\{p\}$,
applying the dual of Gysin triangle to the open immersion
$Y\hookrightarrow \mathbb{A}^n$, we see that
there is an equivalence
$M_Y \simeq \uni_k\oplus \mathsf{1}(-n)[-2n+1]^{\oplus 2}$
in $\DM(k)$.
The morphism $s^*:\uni_k\to \uni_k\oplus
\mathsf{1}(-n)[-2n+1]^{\oplus 2}\simeq M_Y$
induces an equivalence $\uni_k\stackrel{\sim}{\to} \uni_k\hookrightarrow
\uni_k\oplus
\uni_k(-n)[-2n+1]^{\oplus 2}$. Thus,
$V_0\simeq \uni_k(-n)[-2n]^{\oplus 2}$.
We then find that
\[
A_1=\Free(0)\otimes_{\Free(\uni_k(-n)[-2n]^{\oplus 2})}\Free(0) \simeq \Free(0\sqcup_{\uni_k(-n)[-2n]^{\oplus 2}}0)\simeq \Free(\uni_k(-n)[-2n+1]^{\oplus 2}).
\]
The induced morphism $f:A_1=\Free(0)\otimes_{\Free(\uni_k(-n)[-2n]^{\oplus 2})}\Free(0)\simeq \Free(\mathsf{1}_k(-n)[-2n+1]^{\oplus 2})\to M_Y$ is
classified by the inclusion $\iota:\uni_k(-n)[-2n+1]^{\oplus 2}\hookrightarrow
\uni_k\oplus \uni_k(-n)[-2n+1]^{\oplus 2}\simeq M_Y$.
Let $F:\Rep^\otimes(\mathbb{G}_m)\to \DM^\otimes(k)$
be the colimit-preserving symmetric monoidal functor
which carries $\chi_1$ to $\uni_k(1)$ (cf. the proof of Proposition~\ref{cpsprop}).
Consider $\Free_{\Comp(\mathbb{G}_m)}(\chi_{-n}[-2n+1]^{\oplus 2})$.
The underlying object in $\Comp(\mathbb{G}_m)$
is isomorphic to $\uni_k\oplus \chi_{-n}[-2n+1]^{\oplus 2} \oplus \Sym^2(\chi_{-n}[-2n+1]^{\oplus 2})\simeq \uni_k\oplus \chi_{-n}[-2n+1]^{\oplus 2} \oplus \chi_{-2n}[-4n+2]$. The image of $\Free_{\Comp(\mathbb{G}_m)}(\chi_{-n}[-2n+1]^{\oplus 2})$ under $\CAlg(\Rep^\otimes(\mathbb{G}_m))\to \CAlg(\DM^\otimes(k))$ is equivalent to $A_1$.
The composite
$\uni_k\oplus \uni_k(-n)[-2n+1]^{\oplus 2}\hookrightarrow
\uni_k\oplus \uni_k(-n)[-2n+1]^{\oplus 2}\oplus \uni_k(-2n)[-4n+2]\simeq \Free(\mathsf{1}_k(-n)[-2n+1]^{\oplus 2})\to M_Y$
is an equivalence.
Note that
a morphism $\uni_k(-2n)[-4n+2]\to \uni_k\oplus \uni_k(-n)[-2n+1]^{\oplus 2}$
is null homotopic because it corresponds to
an element in $\Hom_{\hhh(\DM(k))}(\uni_k,\uni_k(n)[2n-1])^{\oplus 2})\oplus \Hom_{\hhh(\DM(k))}(\uni_k,\uni_k(2n)[4n-2])\simeq (\CH^n(\Spec k,1)^{\oplus 2}\oplus \CH^{2n}(\Spec k,2))\otimes_{\ZZ}\QQ=0$ (we use the condition $n\ge2$). Here $CH^i(-,j)$ is the Bloch's higher Chow group. Hence $V_1=\uni_k(-2n)[-4n+2]$
and $V_1\to A_1\simeq \uni_k\oplus \uni_k(-n)[-2n+1]^{\oplus 2}\oplus \uni_k(-2n)[-4n+2]$ may be viewed as the canonical inclusion.
We see that
\[
A_2=\Free(\uni_k(-n)[-2n+1]^{\oplus 2})\oplus_{\Free(\uni_k(-2n)[-4n+2])}\uni_k.\]
Consider $\Free_{\Comp(\mathbb{G}_m)}(\chi_{-2n}[-4n+2])\to \Free_{\Comp(\mathbb{G}_m)}(\chi_{-n}[-2n+1]^{\oplus 2})$ classified by
the inclusion $\alpha:\chi_{-2n}[-4n+2]\hookrightarrow \Free_{\Comp(\mathbb{G}_m)}(\chi_{-n}[-2n+1]^{\oplus 2})$.
Let $\Free_{\Comp(\mathbb{G}_m)}(\chi_{-n}[-2n+1]^{\oplus 2})\langle \alpha \rangle$ be the homotopy pushout, see Proposition~\ref{hopushout}.
Note that the image of $\Free_{\Comp(\mathbb{G}_m)}(\chi_{-n}[-2n+1]^{\oplus 2})\langle \alpha \rangle$ in $\CAlg(\DM^\otimes(k))$ (under $F$)
is equivalent to $A_2$.
By the computation using Remark~\ref{explicitdg},
we see that
\[
\Free_{\Comp(\mathbb{G}_m)}(\chi_{-n}[-2n+1]^{\oplus 2})\langle \alpha \rangle\simeq \chi_0\oplus \chi_{-n}[-2n+1]^{\oplus 2}\oplus \chi_{-3n}[-6n+2]^{\oplus 2}\oplus \chi_{-4n}[-8n+3]
\]
in $\Rep(\mathbb{G}_m)$. Hence $A_2\simeq \uni_k\oplus \uni_k(-n)[-2n+1]^{\oplus 2}\oplus \uni_k(-3n)[-6n+2]^{\oplus 2}\oplus \uni_k(-4n)[-8n+3]$.
By the argument similar to the case of $V_1$, we see that $V_2=\uni_k(-3n)[-6n+2]^{\oplus 2}\oplus \uni_k(-4n)[-8n+3]$ and $V_2\to A_2$ may be viewed as the canonical inclusion. We thus find
\[
A_3=A_2\otimes_{\Free(\uni_k(-3n)[-6n+2]^{\oplus 2}\oplus \uni_k(-4n)[-8n+3])}\uni_k.
\]

\subsubsection{}
\label{Sull1}

 Let $\mathcal{C}^\otimes$ be a stable presentable
 $\infty$-category endowed with a symmetric monoidal
 structure whose tensor
 operation
$\mathcal{C}\times \mathcal{C} \to \mathcal{C}$
preserves small colimits separately in each variable.
Let $\Free_{\CCC}:\CCC\to \CAlg(\CCC^\otimes)$ be the free functor of $\mathcal{C}^\otimes$.
Let $A$ and $B$ be commutative algebra objects in $\CAlg(\CCC^\otimes)$
and $f:A\to B$ be a morphism in $\CAlg(\CCC^\otimes)$.
Let $V$ be the kernel of $f$ in the stable $\infty$-category $\CCC$,
i.e.,
the pullback $A\times_{B}\{0\}$.
Let
$\sigma:\Free_{\CCC}(V)\to A$ in $\CAlg(\CCC^\otimes)$
be the morphism classified by $V\to A$.
Let $\epsilon:\Free_{\CCC}(V)\to \uni_{\CCC}=\Free_{\CCC}(0)$
be the morphism induced by $V\to 0$ where $\uni_{\CCC}$
is the unit algebra in $\CAlg(\CCC^\otimes)$.
We have a pushout diagram
\[
\xymatrix{
\Free_{\CCC}(V) \ar[r]^\sigma \ar[d]_\epsilon & A \ar[d] \\
\uni_{\CCC} \ar[r] & A(f)
}
\]
in $\CAlg(\CCC^\otimes)$.
Note that the composite $\Free_{\CCC}(V)\to A \stackrel{f}{\to} B$ factors
through $\Free_{\CCC}(V)\to \uni_{\CCC}$.
We have a factorization 
\[
A\to A(f)\stackrel{f'}{\to} B
\]
 of $f$. Applying this procedure to $f':A(f)\to B$ 
we obtain a refined factorization $A\to A(f)\to A(f,f'):=A(f)(f')\to B$.
Repeating it in the inductive way we have a sequence in $\CAlg(\CCC^\otimes)_{/B}$
described as
\[
A=A_0\to A_1\to A_2\to \cdots \to A_n\to A_{n+1} \to \cdots
\]
where $A_1=A(f)$, $A_2=A(f,f')\ldots$.
We denote by $f_n:A_n\to B$ the structural morphism.
We shall refer to this sequence as the inductive sequence associated to
$A\to B$.

\begin{Lemma}
\label{Sullivanmodel1}
Let $\varinjlim_{n}A_n$ be a colimit of the sequence in $\CAlg(\CCC^\otimes)$.
Then the canonical morphism
$\varinjlim_{n}A_n \to B$ is an equivalence in $\CAlg(\CCC^\otimes)$.
\end{Lemma}

\Proof
According to \cite[3.2.3.1]{HA},
the forgetful functor $\CAlg(\CCC^\otimes)\to \CCC$
preserves filtered colimits.
Hence it is enough to prove that
a colimit $\varinjlim_{n}A_n$ in $\CCC$ (by abuse of notation
we continue to use the same symbol)
is naturally equivalent to $B$ in $\CCC$.
If $V_n$ denotes the kernel of $f_n:A_n\to B$ in $\CCC$, then
$V_n\to \Free_{\CCC}(V_n)\to A_n\to A_{n+1}$ is null-homotopic.
Thus, $A_n\to A_{n+1}$ factors as composition
$A_n\to \Coker(V_n\to A_n)\to A_{n+1}$ in $\CCC$ where
$\Coker(-)$ stands for cokernel (cofiber/cone) in $\CCC$.
The sequence $A\to A_{1}\to A_2\to\cdots$ in $\CCC$
is refined as
\[
A\to A_1\to \Coker(V_1\to A_1)\to A_2\to \Coker(V_2\to A_2)\to A_3\to\cdots.
\]
By cofinality, the colimit of this sequence is naturally equivalent
to $\varinjlim_{n}A_n$.
Notice that $\Coker(V_n\to A_n)\simeq B$ in $\CCC$
for any $n\ge1$.
Hence we deduce that 
$\varinjlim_{n}A_n \to B$ is an equivalence in $\CCC$.
\QED

\begin{Remark}
\label{Sull1R}
Let $\mathcal{D}^\otimes$ be another stable presentable
 $\infty$-category endowed with a symmetric monoidal
 structure whose tensor
 operation
$\mathcal{D}\times \mathcal{D} \to \mathcal{D}$
preserves small colimits separately in each variable.
Let $F:\mathcal{C}^\otimes\to \mathcal{D}^{\otimes}$ be 
a symmetric monoidal functor that preserves small colimits.
Our main example of interest is 
the realization functor $\mathsf{R}:\DM^\otimes(k) \to \mathsf{D}^\otimes(\QQ)$.
Let $A=A_0\to A_1 \to \cdots\to B$
be the inductive sequence
associated to $f:A\to B$.
Note that $\mathcal{C}\to \mathcal{D}$ is an exact functor of stable $\infty$-categories, and $\CAlg(F):\CAlg(\mathcal{C}^\otimes)\to \CAlg(\mathcal{D}^\otimes)$
preserves small colimits.
Then the sequence $F(A_0)\to F(A_1)\to \cdots \to F(B)$
is canonically equivalent to the inductive sequence
associated to $F(A)\to F(B)$ as a diagram in $\CAlg(\mathcal{D}^{\otimes})_{/F(B)}$.
\end{Remark}


\subsection{}
\label{semiab}
Let $G$ be a semi-abelian variety over $k$.
There is a (canonical) equivalence
\[
M(G)\stackrel{\sim}\to \oplus_{n\ge0}M_n(G)
\]
in $\DM(k)$ such that $M_n(G)=\Sym^n(M_1(G))$.
This is a result of Ancona-EnrightWard-Huber \cite{AEH},
which builds upon the works
of Shermenev, Deniger-Murre and K\"unnemann
on a decomposition of the motives of an abelian variety, see \cite{Kun}
and references therein.
If $G$ is an extension of a $g$-dimesional abelian
variety by a torus of rank $r$,
then $\Sym^n(M_1(G))\simeq 0$ for $n> 2g+r$.
The direct summand $M_1(G)$ is represented,
as an object in $\Comp(\mathcal{N}^{tr}(X))$, by 
the \'etale sheaf of $\QQ$-vector spaces
given by $S\mapsto \Hom_{\Sm_k}(S,G)\otimes_{\ZZ}\QQ$
which is promoted to a sheaf with transfers (see
e.g. \cite[Section 2.1]{AEH}).

By using their work, we prove the following:

\begin{Theorem}
\label{abelianvar}
Let $M_1(G)^\vee$ be the dual of $M_1(G)$ in $\DM(k)$ ($M_1(G)$
is a dualizable object).
Let $\Free_{\DM(k)}(M_1(G)^\vee)$ be a free commutative algebra object
in $\DM(k)$ generated by $M_1(G)^\vee$.
Then there is an equivalence
\[
\Free_{\DM(k)}(M_1(G)^\vee) \stackrel{\sim}{\longrightarrow} M_G
\]
in $\CAlg(\DM^\otimes(k))$.
\end{Theorem}

\begin{Remark}
Let $\mathbb{G}$ be a connected compact Lie group.
A theorem of Hopf says that 
there are elements $x_1,\ldots,x_n$ of odd degrees in
$H^*(\mathbb{G},\QQ)$ such that
$H^*(\mathbb{G},\QQ)$ is a free commutative graded algebra
generated by $x_1,\ldots,x_n$.
One can deduce from this theorem
that a Sullivan model of $A_{PL}(\mathbb{G})$ is
given by a free commutative
graded algebra generated by some graded vector space, see \cite[Example 3 in
page 143]{FHT}.
Theorem~\ref{abelianvar} may be thought of as a generalization of
this homotopical statement to $\CAlg(\DM^\otimes(k))$ for semi-abelian varieties.
\end{Remark}

\begin{Lemma}
\label{homotopyfree}
Let $\CCC^\otimes$ be a symmetric monoidal presentable
$\infty$-category
whose tensor
 operation
$\mathcal{C}\times \mathcal{C} \to \mathcal{C}$
preserves small colimits separately in each variable.
Suppose that $\CCC^\otimes$ is $K$-linear, namely,
it is endowed with a colimit-preserving
symmetric monoidal functor $\Mod_K^\otimes\to \CCC^\otimes$
($K$ is a field of characteristic zero).
Let $\hhh(\CCC)^\otimes$ be the homotopy category of $\CCC$
endowed with a symmetric monoidal structure induced by
that of $\CCC^\otimes$.
The canonical functor $\pi:\CCC\to \hhh(\CCC)$
can be promoted to a symmetric monoidal functor.
Let $\pi':\CAlg(\CCC^\otimes)\to \CAlg(\hhh(\CCC)^\otimes)$
be the ``projection'' induced by the symmetric monoidal functor $\pi$.
In this Lemma we use the temporaty notation
$\Free:=\Free_{\CCC}:\CCC\to \CAlg(\CCC^\otimes)$
be a free functor of $\CCC$.
Let
$\Free^h:=\Free_{\hhh(\CCC)}:\hhh(\CCC)\to \CAlg(\hhh(\CCC)^\otimes)$
be a free functor of $\hhh(\CCC)$.
Let $\theta:\CAlg(\CCC^\otimes)\to \CCC$ and $\theta^h:\CAlg(\hhh(\CCC)^\otimes)\to \hhh(\CCC)$ be forgetful functors.
Let $C$ be an object in $\CCC$.
The unit map $C\to \theta(\Free(C))$ induces
$\pi(C)\to \pi(\theta(\Free(C)))=\theta^h(\pi'(\Free(C)))$.
By the adjunction
$(\Free^h,\theta^h)$, it gives rise to $\sigma:\Free^h(\pi(C))\to \pi'(\Free(C))$.
Then the canonical morphism $\sigma$ is an equivalence.
\end{Lemma}

\Proof
Let $A=\theta^h(\pi'(\Free(C)))$.
The $n$-fold multiplication $A^{\otimes n}\to A$
induces $\Sym_{\hhh(\CCC)}^n (A)\to A$ where $\Sym_{\hhh(\CCC)}^n(-)$ is the $n$-fold symmetric
product in the $K$-linear idempotent complete category $\hhh(\CCC)$.
The map $\pi(C)\to A$ induces $\Sym_{\hhh(\CCC)}^n(\pi(C))\to \Sym_{\hhh(\CCC)}^n(A)\to A$.
Taking its coproduct we have
\[
\tau:\oplus_{n\ge0} \Sym_{\hhh(\CCC)}^n(\pi(C))\to A.
\]
Taking account of the canonical
equivalence $\oplus_{n\ge0} \Sym_{\hhh(\CCC)}^n(\pi(C))\simeq \Free^h(\pi(C))$,
it will suffice to show that $\tau$ is an isomorphism in $\hhh(\CCC)$.
By \cite[3.1.3.13]{HA}, there is an equivalence $\oplus_{n\ge 0}\Sym_{\CCC}^n(C)\to \theta(\Free(C))$ where each $\Sym_{\mathcal{C}}^n(C)\to \theta(\Free(C))$
is induced by the composition of 
$C^{\otimes n}\to \Free(C)^{\otimes n}$
and the $n$-fold multiplication $\Free(C)^{\otimes n} \to \Free(C)$.
Here $\Sym_{\CCC}^n(C)$ is the symmetric product in $\CCC$.
Therefore,
it is enough to prove that
the natural morphism $\Sym_{\hhh(\CCC)}^n(\pi(C))\to \pi(\Sym_{\CCC}^n(C))$
is an isomorphism.
Note that for any $D$ in $\CCC$,
the set $\Hom_{\hhh(\CCC)}(\Sym_{\hhh(\CCC)}^n(\pi(C)),\pi(D))$
is the invariant part $\Hom_{\hhh(\CCC)}(\pi(C)^{\otimes n},\pi(D))^{\Sigma_n}$
of $\Hom_{\hhh(\CCC)}(\pi(C)^{\otimes n},\pi(D))$
with the permutation action of the symmetric group $\Sigma_n$.
On the other hand, the hom complex
$\Hom_{\CCC}(\Sym_{\CCC}^n(C),D)$ in $\Mod_K$ is a limit $\Hom_{\CCC}(C^{\otimes n},D)^{\Sigma_n}$ of
$\Hom_{\CCC}(C^{\otimes n},D)$ with permutation action of $\Sigma_n$.
(By definition, the hom complex $\Hom_{\CCC}(C,D)$
is given by the image of $D$ under the right adjoint $\Hom_{\CCC}(C,-)$
to the colimit preserving functor $(-)\otimes C:\Mod_K\to\CCC$.)
Since $K$ is a field of characteristic zero (the semi-simplicity of
representations of finite groups), we have
$H^0(\Hom_{\CCC}(C^{\otimes n},D)^{\Sigma_n})=H^0(\Hom_{\CCC}(C^{\otimes n},D))^{\Sigma_n}=\Hom_{\hhh(\CCC)}(\Sym_{\hhh(\CCC)}^n(\pi(C)),\pi(D))$.
Thus, we see that $\Sym_{\hhh(\CCC)}^n(\pi(C))\to \pi(\Sym_{\CCC}^n(C))$
is an isomorphism.
\QED

\begin{Remark}
\label{homotopyfreeR}
By the proof, 
if we define the canonical
functor $\pi:\DM(k)\to \hhh(\DM(k))$,
then we have
a canonical isomorphism $\Sym_{\hhh(\DM(k))}^n(\pi(C))\simeq \pi(\Sym_{\DM(k)}^n(C))$.
Namely, $\pi$ commutes with
the formulation of symmetric products.
By this canonical isomorphism, we often abuse notation
by writing $\Sym^n(C)$
for both $\Sym_{\hhh(\DM(k))}^n(\pi(C))$ and $\Sym_{\DM(k)}^n(C)$.
\end{Remark}

{\it Proof of Theorem~\ref{abelianvar}.}
Let $\alpha_G:M(G)\to M_1(G)$ be the morphism described in \cite[2.1.4]{AEH}
(in {\it loc. cit.},
$\alpha_G$ is a morphism $\DM_{eff}(k)$, but we here regard
it as a morphism in $\DM(k)$).
We remark also that in \cite[2.1.4]{AEH} \'etale motives
are empolyed, but
$\DM^\otimes(k)$ agrees with the \'etale version
since $K$ is a field of characteristic zero, cf. \cite{MVW}, \cite[1.6.1]{AEH}.
Let $\alpha^\vee_G:M_1(G)^\vee\to M(G)^\vee$ be a dual of $\alpha_G$.
Since $M_G=M(G)^\vee$ in $\DM(k)$,
$\alpha_G^\vee$ induces
a morphism $\Free(M_1(G)^\vee)\to M_G$ in
$\CAlg(\DM^\otimes(k))$.
We will prove that it is an equivalence.
To see this, it is enough to show that
$\pi'(\Free(M_1(G)^\vee))\to \pi'(M_G)$ is an isomorphism
where $\pi':\CAlg(\DM^\otimes(k))\to \CAlg(\hhh(\DM(k))^\otimes)$
is the canonical functor (we continue to use the notation
in Lemma~\ref{homotopyfree}).
Lemma~\ref{homotopyfree} guarantees that 
$\Free^h(\pi(M_1(G)^\vee))\stackrel{\sim}{\to} \pi'(\Free(M_1(G)^\vee))$.
The composite $\Free^h(\pi(M_1(G)^\vee))\to \pi'(M_G)$ is
induced by $\pi(M_1(G)^\vee)\to \pi(\theta(M_G))=\theta^h(\pi'(M_G))$.
The proof is reduced to showing that this morphism
$\Free^h(\pi(M_1(G)^\vee))\to \pi'(M_G)$ is an isomorphism in $\hhh(\DM(k))$.
The each factor $\phi_n:\Sym^n(\pi(M_1(G)^\vee)\to \pi'(M_G)$ of $\oplus_{n\ge0}\Sym^n(\pi(M_1(G)^\vee) \stackrel{\sim}{\to} \Free^h(\pi(M_1(G)^\vee))\to \pi'(M_G)$
is induced by $\pi(M_1(G)^\vee)^{\otimes n} \stackrel{(\alpha_G^\vee)^{\otimes n}}{\to} \theta^h(\pi'(M_G))^{\otimes n}\to \theta^h(\pi'(M_G))$
where the second morphism is the $n$-fold multiplication.
In the following Lemmata, we will observe that $\phi_n$ is a dual of
the projection $\psi_n:M(G)\to \Sym^n(M_1(G))$ of
the equivalence
$M(G)\stackrel{\sim}{\to} \oplus_{0\le n\le 2g+r}\Sym^n(M_1(G))$
proved in \cite[Theorem 7.1.1]{AEH}. It will finish the proof.
\QED

\begin{Lemma}
\label{algcoalg}
$\phi_n:\Sym^n(\pi(M_1(G)^\vee))\to \theta^h(\pi'(M_G))$
is a dual of $\psi_n:M(G)\to \Sym^n(M_1(G))$.
\end{Lemma}

\Proof
We first recall $\psi_n$.
We work with the homotopy category
$\hhh(\DM(k))$.
By abuse of notation,
we put $M(G)=\pi(M(G))$, $M_1(G)=\pi(M_1(G))$,
$\Sym^n(M_1(G)^\vee)=\Sym^n(\pi(M_1(G))^\vee)$,
$\Sym^n(M_1(G))=\Sym^n(\pi(M_1(G)))$,  $M_G=\theta^h(\pi'(M_G))$, etc.
These idetifications are harmless (cf. Lemma~\ref{homotopyfree} and
Remark~\ref{homotopyfreeR}).
The morphism $\psi_n:M(G)\to  \Sym^n(M_1(G))$
is the composite
$M(G)\to M(G)^{\otimes n}\stackrel{\alpha_G^{\otimes n}}\to M_1(G)^{\otimes n}\to \Sym^n(M_1(G))$ where the first morphism is the $n$-fold comultiplication
and the the third morphism is the canonical projection.
By ease of notation,
we let $f_\sharp:\hhh(\DM(G))\rightleftarrows \hhh(\DM(k)):f^*$
be the adjoint pair induced by $f_\sharp:\DM(G)\rightleftarrows \DM(k):f^*$
where $f:G\to \Spec k$ is the structure morphism.
The colax monoidal functor $f_\sharp$
induces the coalgebra structure on $M(G)=f_\sharp(\uni_X)$ in $\hhh(\DM(k))$:
the comultiplication is given by the composition
\[
f_\sharp(\uni_G)=f_\sharp(\uni_G \otimes \uni_G)\to f_\sharp(f^*f_\sharp(\uni_G)\otimes f^*f_\sharp(\uni_G))\simeq f_\sharp f^*(f_\sharp(\uni_G)\otimes f_\sharp(\uni_G)) \to f_\sharp(\uni_G)\otimes f_\sharp(\uni_G)
\]
where the left arrow is induced by the counit of $(f_\sharp,f^*)$,
and the right arrow
is induced by the unit. The counit $M(G)\to \uni_k$ is
$f_\sharp f^*(\uni_k)\to \uni_k$.
If one regards $\Sym^n (M_1(G))$ as a direct summand of $M_1(G)^{\otimes n}$,
then $M(G)\to M(G)^{\otimes n}\stackrel{\alpha_G^{\otimes n}}\to M_1(G)^{\otimes n}$ factors as
$M(G)\stackrel{\psi_n}{\to}  \Sym^n(M_1(G))\to M(G)^{\otimes n}$.
On the other hand,
$\phi_n:\Sym^n(M_1(G)^{\vee})\to M_G$ is induced by
$\Sigma_n$-equivariant morphism $(M_1(G)^{\vee})^{\otimes n}\stackrel{(\alpha_G^\vee)^{\otimes n}}{\to} M_G^{\otimes n } \to M_G$ where $M_G$ has the trivial
action, and the second arrow is the $n$-fold multiplication.
To prove our assertion of this Lemma, it is enough to show 
the following general fact:

\begin{Lemma}
\label{multcomult}
Let $X$ be a smooth scheme separated of finite type over $k$.
The multiplication morphism $M_X\otimes M_X\to M_X$ is a
dual of the comultiplication morphism $M(X)\to M(X)\otimes M(X)$
given by the diagonal
in $\hhh(\DM(k))$ through the isomorphism $M_X\simeq M(X)^\vee$.
(We remark that $M(X)$ is also dualizable in $\DM^\otimes(k)$
since we work with coefficients
of characteristic zero.)
\end{Lemma}

\Proof
We here write $\uni:=\uni_X$ and the structure morphism $f:X\to \Spec k$.
Remember that the multiplication $M_X\otimes M_X\to M_X$
is given by the composition
\[
f_*(\uni)\otimes f_*(\uni)\to f_*f^*(f_*(\uni)\otimes f_*(\uni))\simeq f_*(f^*f_*(\uni)\otimes f^*f_*(\uni))
\to f_*(\uni\otimes \uni)\simeq f_*(\uni)
\]
such that the left arrow is induced by the counit of the adjunction
$(f^*,f_*)$, and the right arrow is induced  by its counit.
The canonical isomorphism $\eta:M(X)^\vee \stackrel{\sim}{\to}M_X$
is defined as follows (see the proof of Proposition~\ref{cohmotalg}). 
For $M \in \DM(X)$, consider the unit $M\to f^*f_\sharp(M)$.
Taking the dual and $f_*$, we have
$f_*f^*(f_\sharp (M))^\vee \simeq f_*(f^*f_\sharp (M))^\vee\to f_*(M^\vee)$.
Composing with the unit
$(f_\sharp (M))^\vee \to f_*f^*(f_\sharp (M))^\vee$
we obtain $\eta_M:(f_\sharp (M))^\vee\to f_*(M^\vee)$
which determines an isomorphism $\eta=\eta_{\uni_X}:M(X)^\vee \stackrel{\sim}{\to}M_X$.
We will check that
the dual of $f_\sharp(\uni)\to f_\sharp(\uni)\otimes f_\sharp(\uni)$ is $f_*(\uni)\otimes f_*(\uni)\to f_*(\uni)$
through $\eta:f_\sharp (\uni)^\vee\simeq f_*(\uni)$.
By using the counit of $(f^*,f_*)$ and its counit-unit equations,
we see that
the dual $f^*f_\sharp(M)^\vee\to M^\vee$ of $M\to f^*f_\sharp(M)$
is $f^*(f_\sharp(M)^\vee)\stackrel{f^*(\eta_M)}{\to} f^*f_*(M^\vee)\to M^\vee$
where the final arrow is the counit of $(f^*,f_*)$.
When $M=\uni$, we deduce that the unit $s:\uni\to f^*f_\sharp(\uni)$ is the dual of the counit $t:f^*f_*(\uni)\to \uni$
through the isomorphism $\eta:f_\sharp(\uni)^\vee\simeq f_*(\uni)$. It follows that its tensor product
$s\otimes s:\uni\otimes \uni \to f^*f_\sharp(\uni)\otimes f^*f_\sharp(\uni)$
is the dual of $t\otimes t:f^*f_*(\uni)\otimes f^*f_*(\uni)\to \uni\otimes \uni$
through the isomorphism through the isomorphism $\eta:f_\sharp(\uni)^\vee\simeq f_*(\uni)$.
Thus the triangle in the following diagram commutes.
\[
\xymatrix{
(f_\sharp(\uni)\otimes f_\sharp(\uni))^\vee \ar[r]^(0.45){a} \ar[d]_{\simeq} & (f_\sharp f^*(f_\sharp(\uni)\otimes f_\sharp(\uni)))^\vee \ar[r]^(0.7)b \ar[d]_{\eta_{f^*(f_\sharp(\uni)\otimes f_\sharp(\uni))}} & f_\sharp(\uni)^\vee \ar[dd]^{\eta} \\
f_\sharp(\uni)^\vee\otimes f_\sharp(\uni)^\vee \ar[r]^(0.45){c} \ar[d]_{\eta\otimes \eta} & f_*f^*(f_\sharp(\uni)^\vee\otimes f_\sharp(\uni)^\vee) \ar[rd]^{f_*(s^\vee \otimes s^\vee)} \ar[d]_{f_*f^*(\eta\otimes \eta)} &  \\
f_*(\uni)\otimes f_*(\uni) \ar[r]^(0.45){d} & f_*f^*(f_*(\uni)\otimes f_*(\uni))  \ar[r]_(0.7){f_*(t\otimes t)} & f_*(\uni)
}
\]
Here $a$ is induced by
the dual of $f_\sharp f^*\to \textup{id}$, and $b$ is induced by the
dual of $s\otimes s:\uni \otimes \uni \to f^*f_\sharp (\uni)\otimes f^*f_\sharp(\uni)=f^*(f_\sharp(\uni)\otimes f_\sharp(\uni))$.
Note that the composite of the upper horizontal arrows is the dual of comultiplication
$M(X)\to M(X)\otimes M(X)$.
Both $c$ and $d$ is induced by the unit $\textup{id}\to f_*f^*$.
The composite of lower horizontal arrows is the multiplication
$M_X\otimes M_X\to M_X$.
The commutativity of other squares follow from
the contravariant functoriality of $\eta_M$ with respect to $M$,
the functoriality/naturality of $\textup{id}\to f_*f^*$,
and the counit-unit equations for the adjunction $(f_\sharp,f^*)$.
Thus, we have a commutativity of the outer square, which completes the proof.
\QED

\begin{Remark}
The unit map $\uni_k\to M_X$ is nothing but a dual
of the morphism $M(X)\to M(\Spec k)=\uni_k$ induced by $f$.
\end{Remark}


\subsection{}

\subsubsection{}
\label{curveexp}
We consider a once-punctured smooth proper curve,
that is, $C=\overline{C}-\{p\}$ obtained by removing
a $k$-rational point $p$ from
a connected smooth proper curve $\overline{C}$
over the perfect field $k$.
Let $j:C\to \overline{C}$
be the open immersion.
The genus of $\overline{C}$ is $g\ge1$.
If $k=\CC$, the fundamental
group of the underlying topological space
is the free group generated by $2g$ elements.

Let $M(\overline{C})\simeq M_0(\overline{C})\oplus M_1(\overline{C})\oplus M_{2}(\overline{C})$
be a (Chow-K\"unneth) decomposition of $M(\overline{C})$ 
such that $M_0(\overline{C})\simeq \uni_k$
and $M_2(\overline{C})\simeq \uni_k(1)[2]$ (see the first paragraph
of the proof of Lemma~\ref{lemmasqzero} for the precise
formulation).
(In this case, $M(C)$ is equivalent to
$M_0(\overline{C})\oplus M_1(\overline{C})$ as an object in $\DM(k)$.)
We put $M^i_{\overline{C}}:=M_i(\overline{C})^\vee$
so that $M_{\overline{C}}\simeq \oplus_{i=0}^2M^i_{\overline{C}}$.
Let
\[
A_0=\uni_k\to A_1\to A_2\to \cdots \to A_n\to A_{n+1}\to \cdots
\]
be the inductive sequence in $\CAlg(\DM^\otimes(k))_{/M_C}$
associated to the unique morphism
$\uni_k \to M_C$ in $\CAlg(\DM^\otimes(k))$ (cf. Section~\ref{Sull1}).
By Lemma~\ref{Sullivanmodel1}, the colimit $\varinjlim_{n\ge0}A_n\simeq M_C$.

\begin{Theorem}
\label{curve3}
The first three terms $A_1$, $A_2$ and $A_3$ are computed as follows:

\begin{enumerate}
\renewcommand{\labelenumi}{(\theenumi)}

\item $A_1$ is $\Free(M^1_{\overline{C}})$, and $f_1:A_1\to M_C$ is classified by $M^1_{\overline{C}} \hookrightarrow M^0_{\overline{C}}\oplus M^1_{\overline{C}}\oplus M^2_{\overline{C}}\simeq M_{\overline{C}}\stackrel{j^*}{\to} M_C$.

\item $f_1:A_1\to M_C$ is the composite $M_{\Alb_{\overline{C}}}\stackrel{u^*}{\to} M_{\overline{C}}\stackrel{j^*}{\to} M_C$ up
to an equivalence $\Free(M^1_{\overline{C}})\simeq M_{\Alb_{\overline{C}}}$, where $u:\overline{C}\to \Alb_{\overline{C}}$ is the Albanese (Abel-Jacobi) morphism into the Albanese
variety, which carries $p$ to the origin.

\item Let $W_1$ be $\oplus_{i=2}^{2g}\Sym^i(M_{\overline{C}}^1)$.
Let $\Free(W_1)\to A_1$ be the morphism in $\CAlg(\DM^\otimes(k))$
that is classified by the inclusion
$W_1\to A_1\simeq \oplus_{i\ge0}\Sym^i(M_{\overline{C}}^1)$
in $\DM(k)$. Then $A_2$ is the pushout $A_1\otimes_{\Free(W_1)}\uni_k$,

\item Let $W_2$ be the object in $\DM(k)$
which will be defined just before the proof (Section~\ref{curve3proof}).
Then $A_3$ has the form of the pushout $A_2\otimes_{\Free(W_2)}\uni_k$.

\end{enumerate}

\end{Theorem}

\begin{Remark}
The symmetric product $\Sym^N(M^1_{\overline{C}})$ is zero for $N>2g$
(see the proof of Lemma~\ref{lemmasqzero}).
Thus $A_1=\oplus_{i\ge0}\Sym^i(M^1_{\overline{C}})=\oplus_{i=0}^{2g}\Sym^i(M^1_{\overline{C}})$.
\end{Remark}

\begin{Remark}
As mentioned in Introduction,
the sequence $\{A_n\}_{n\ge0}$ can be
viewed as a step-by-step description of the ``non-abelian structure'' of $C$.
To give a feeling for this, let us make the following observation.
Suppose that $c$ is a $k$-rational point on $C$, and let
$M_C\to \uni_k$ be the induced augmentation. The sequence
$\{A_n\}_{n\ge0}$ is promoted to a sequence
to $\CAlg(\DM^\otimes(k))_{/\uni_k}$ in the obvious way.
By applying the construction in Section~\ref{hopfsection} to $A_n\to \uni_k$,
we obtain the sequence 
of cogroup objects in $\CAlg(\DM^\otimes(k))$, which we denote
by $\{\uni_k\otimes_{A_n}\uni_k\}_{n\ge0}$ (we abuse notation
since $\uni_k\otimes_{A_n}\uni_k$
is the underlying object).
Now consider the ``topological aspect'' of this sequence.
For this purpose,
suppose further that $k\subset \CC$. Let $R_n\to \QQ$ be the image of $A_n\to \uni_k$
under the singular realization functor. Let $R\to \QQ$ be the image of $M_C\to \uni_k$.
Then
$\{R_n\}_{n\ge0}$ is the inductive sequence associated to $\QQ\to R$ (cf.
Remark~\ref{Sull1R}), and
the image of $\{\uni_k\otimes_{A_n}\uni_k\}_{n\ge0}$ under the realization
is $\{\QQ\otimes_{R_n}\QQ\}_{n\ge0}$ (we abuse notation again).
Taking the $0$-th cohomology, we have the sequence of the pro-unipotent
algebraic groups
\[
\cdots \to \Spec H^0(\QQ\otimes_{R_n}\QQ)\to \cdots \to \Spec H^0(\QQ\otimes_{R_1}\QQ)\to \Spec H^0(\QQ\otimes_{R_0}\QQ)\simeq \Spec \QQ.
\]
Define $G_n:=\Spec H^0(\QQ\otimes_{R_n}\QQ)$.
In this case, by Theorem~\ref{curve3} (i), $G_1$ is a commutative
unipotent group of rank $2g$ (in fact, $\Free(M_{\overline{C}}^1[1])$ maps to
$\QQ\otimes_{R_1}\QQ\simeq \Free_{\QQ}(\QQ^{\oplus 2g})\simeq H^0(\QQ\otimes_{R_1}\QQ)$).
Recall that $G:=\Spec H^0(\QQ\otimes_{R}\QQ)$ is the pro-unipotent completion
$\pi_1(C^t,c)$ of $C^t$ (cf. Section~\ref{MGA}).
By a standard argument in rational homotopy theory,
each morphism $G_{n+1}\to G_n$ is a surjective morphism
with a commutative kernel, and the canonical morphism
$\pi_1(C^t,c)_{uni}\simeq G\to \varprojlim_{n\ge0}G_n$ is an isomorphism of pro-unipotent algebraic groups.

\end{Remark}

\begin{Example}
Let $\overline{C}$ be an elliptic curve
and let $C=\overline{C}-\{0\}$ be the open curve obtained by removing
the origin $0$.
Then by Theorem~\ref{curve3}, one can easily see that
$A_1=\Free(M^1_{\overline{C}})$, $A_2=\Free(M^1_{\overline{C}})\otimes_{\Free(\uni_k(-1)[-2])}\uni_k$, and $A_2$ is equivalent to $\uni_k\oplus M^1_{\overline{C}}\oplus M^1_{\overline{C}}(-1)[-1]\oplus \uni_k(-2)[-3]$
as an object in $\DM(k)$.
We have $W_2=M^1_{\overline{C}}(-1)[-1]\oplus \uni_k(-2)[-3]$,
and the third term $A_3$ is of the form
$A_2\otimes_{\Free(W_2)}\uni_k$.
\end{Example}

\subsubsection{}
\begin{Lemma}
\label{lemmasqzero}
The multiplication map $M_C\otimes M_C\to M_C$ in the homotopy category
$\hhh(\DM(k))$ is
\[
m:(\uni_k\oplus M^1_{\overline{C}})^{\otimes 2} \simeq \uni_k\otimes \uni_k\oplus \uni_k\otimes M^1_{\overline{C}}\oplus M^1_{\overline{C}}\otimes \uni_k\oplus (M^1_{\overline{C}})^{\otimes 2} \to \uni_k\oplus M^1_{\overline{C}}
\]
defined as a coproduct
of $m|_{(M^1_{\overline{C}})^{\otimes2}}=0$
and ``identities''
$\uni_k\otimes \uni_k\to \uni_k$, $\uni_k\otimes M^1_{\overline{C}}\to M^1_{\overline{C}}$,
$M^1_{\overline{C}}\otimes \uni_k\to M^1_{\overline{C}}$.
Namely, $M_C$ is the trivial square zero extension of $\uni_k$
by $M_{\overline{C}}^1$
in $\hhh(\DM(k))$.
\end{Lemma}

\begin{Remark}
The unit $\uni_k\to M_C$ may be identified with the
morphism $\uni_k=M_{\Spec k}\to M_{C}$
determined by the structure morphism $C\to \Spec k$.
By the construction of the decomposition (see below),
it is the inclusion $\uni_k\hookrightarrow \uni_k\oplus M^1_{\overline{C}}$. Thus the non-trivial part of the Lemma is $m|_{(M^1_{\overline{C}})^{\otimes2}}=0$.
\end{Remark}

In the proof of the above Lemma,
we 
discuss decompositions of motives
and use the category $\Chow_k$ of Chow motives
with rational coefficients, cf. \cite[Section 1]{Sc}, \cite[Section 2.2]{MNP}.
We choose the contravariant Chow motives since we will refer to \cite{Sc}
and \cite{MNP} in which the authors adopt the contravariant formulation.
But $\DM(k)$ is a covariant theory in the sense that there is
the canonical covariant functor $\Sm_k\to \DM(k)$ given by $X\mapsto M(X)$
while $\Chow_k$ has a contravariant functor $\textup{SmPr}_k\to \Chow_k$
given by $X\mapsto \Ch(X)$. Here
$\textup{SmPr}_k$ is the category of connected 
smooth projective varieties over $k$,
and 
following \cite{MNP} we denote by $\Ch(X)$ the Chow motive of
$X$ (that is $h(X)$ in \cite{Sc}). The relation between $\DM(k)$
and $\Chow_k$ is quite well-known, but the difference between
the covariant and the contravariant formulations is likely
to cause unnecessary confusion.
We thus give some remarks.
There is a fully faithful $\QQ$-linear functor
$\Chow_k^{op}\to \hhh(\DM(k))$ that is symmetric monoidal,
\cite[20.1, 20.2]{MVW}, \cite[9.3.6]{MNP}.
It carries $\Ch(X)$ to $M(X)$.
The Lefschetz motive $\mathbb{L}$ maps to $\uni_k(1)[2]$.
As the level of hom sets,
\begin{eqnarray*}
\Hom_{\Chow_k}(\Ch(Y),\Ch(X))=\CH^{d}(Y\times X) &\stackrel{\textup{transpose}}{\to}& \CH^{d}(X\times Y) \\
&\simeq& \Hom_{\hhh(\DM(k))}(M(X\times Y),\uni_k(d)[2d]) \\
&\simeq& \Hom_{\hhh(\DM(k))}(M(X),M(Y)^\vee\otimes \uni_k(d)[2d]) \\
&\simeq& \Hom_{\hhh(\DM(k))}(M(X),M(Y))
\end{eqnarray*}
where $d$ and $e$ are the dimensions of $Y$ and $X$, respectively.
For $f:X\to Y$, we write $f^*:\Ch(Y)\to \Ch(X)$ for the class of the transposed graph
${}^t\Gamma_f$ in $\CH^{d}(Y\times X)$. We also use $f_*:\Ch(X)\to \Ch(Y)\otimes \mathbb{L}^{\otimes e-d}$ that
corresponds to the class of $\Gamma_f$ in $\CH^{e+(d-e)}(X\times Y)$.
The functor $\Chow_k^{op}\to \hhh(\DM(k))$ carries
$f^*:\Ch(Y)\to \Ch(X)$ to $M(X)\to M(Y)$ induced by the graph of $f$,
that is the dual of $f^*:M_Y\to M_X$ (the final $f^*$ is defined in Section~\ref{cohomologicalmot}).

{\it Proof of Lemma~\ref{lemmasqzero}.}
For ease of notation, $X=\overline{C}$.
We first recall the decomposition $\Ch^0(X)\oplus \Ch^1(X)\oplus \Ch^2(X)\simeq  \Ch(X)$. We define the retract $\Ch(\Spec k)\stackrel{s^*}{\to} \Ch(X)\stackrel{p^*}{\to} \Ch(\Spec k)$ given by $p:\Spec k=\{p\}\to X$ and the structure morphism $s:X\to \Spec k$.
There is also the retract defined by
$\mathbb{L} \stackrel{p_*\otimes \mathbb{L}}{\to} \Ch(X)\stackrel{s_*}{\to} \mathbb{L}$. The components $\Ch^0(X)=\uni_k$ 
and $\Ch^2(X)$ arise from the first retract and the second retract, respectively. Here we abuse notation by writing $\uni_k$ for the unit object in $\Chow_k$
because it corresponds to the unit object in $\DM(k)$.
The component $\Ch^1(X)$ can be described as the Picard motives
in the sense of J.P. Murre \cite[Section 6.2]{MNP}, \cite[Section 4]{Sc}.
Let $\Alb_X$ be the Albanese variety of $X$ and let $\Pic_X$ be
the Picard variety of $X$.
Note that
\[
\CH^1(X\times X)\supset \Hom^*((X,p),(\Pic_X,0))\simeq \Hom_{AV}(\Alb_X,\Pic_X)
\]
where $\Hom^*$ indicates the set of morphisms that preserve base points,
and $\Hom_{AV}$ indicates the set of morphisms of abelian varieties.
Here we implicitly use the Albanese morphism $(X,p)\to (\Alb_X,0)$.
The set $\Hom^*((X,p),(\Pic_X,0))$ corresponds to 
the subgroup of $\CH^1(X\times X)$,
that consists of those classes of divisors
$D\in \CH^1(X\times X)$ such that $(\textup{id}_X\times p)^*(D)=0$
and
$(p\times \textup{id}_X)^*(D)=0$ in $\CH^1(X)$.
We will call such divisors $p$-normalized divisors
and denote by $\CH^1_{(p)}(X\times X)$ the subgroup of $p$-normalized divisors.
Consider the isomorphism $\theta:\Alb_X\stackrel{\sim}{\to} \Pic_X$
defined by the theta divisor.
By \cite[Lemma 6.2.6]{MNP}
the element $\pi_1$ in $\Hom_{\Chow_k}(\Ch(X),\Ch(X))=\CH^1(X\times X)\otimes_{\ZZ}\QQ$
corresponding to $\theta$ is an idempotent morphism of $\Ch(X)$.
We define $\Ch^{1}(X)$ to be the object corresponding to $\pi_1$,
namely, $\Ch^{1}(X)\to \Ch(X)$ is $\Ker(\textup{id}-\pi_1)\to \Ch(X)$.
Let $M_0(X)\oplus M_1(X)\oplus M_2(X)\simeq M(X)$
be the decomposition that arises from the decomposition
of $\Ch(X)\simeq \Ch^1(X)\oplus \Ch^1(X)\oplus \Ch^2(X)$.
Put $M^i_X=M_i(X)^\vee$ and let
$M^0_X\oplus M^1_X\oplus M^2_X\simeq M_X$ be the decomposition
obtained by taking the dual. We remark that $M^0_X\simeq \uni_k$
and $M^2_X\simeq \uni_k(-1)[-2]$.

Next we construct a Picard motive $\Ch^1(\Alb_X)$ by using the
Albanese (Abel-Jacobi) map $u:X\to \Alb_X$ which carries $p$ to the origin.
Consider the isomorphisms
$\Pic_{\Alb_X}\stackrel{\sim}{\to} \Pic_X\stackrel{\theta}{\leftarrow}\Alb_X\stackrel{\sim}{\to} \Alb_{\Alb_X}$. 
The third morphism induced by the functoriality
is an isomorphism because of the universal property of Albanese vatieties,
and the first morphism is its dual.
Let $\sigma:\Alb_{\Alb_X}\to \Pic_{\Alb_X}$
be the inverse of the composite.
If we denote by $\CH^1_{(0)}(\Alb_X\times \Alb_X)$ the subgroup
of $0$-normalized divisors ($0$ is the origin),
we have the canonical isomorphisms $\Hom_{AV}(\Alb_{\Alb_X},\Pic_{\Alb_X})\simeq \Hom_{AV}(\Alb_X,\Pic_{\Alb_X})\simeq \CH^1_{(0)}(\Alb_X\times \Alb_X)$.
Let $Z$ be the divisor that corresponds to $\sigma$
and let $\phi:\Ch(\Alb_X)\otimes \mathbb{L}^{\otimes 1-g}\to \Ch(\Alb_X)$
be the morphism defined by $Z\in \CH^1(\Alb_X\times\Alb_X)$.
Let $\omega_1:\Ch(\Alb_X)\to \Ch(\Alb_X)$ be the composite
\[
\Ch(\Alb_X)\stackrel{u^*}{\to} \Ch(X)\stackrel{u_*}{\to} \Ch(\Alb_X)\otimes \mathbb{L}^{\otimes 1-g}\stackrel{\phi}{\to} \Ch(\Alb_X).
\]
We can apply the proof of \cite[Lemma 6.2.6]{MNP}
and see that $\omega_1$ is an idempotent map.
We define $\Ch^1(\Alb_X)$ to be $\Ker(\textup{id}-\omega_1)$.

We now claim that the composite
\[
\Ch^1(\Alb_X)\to \Ch(\Alb_X)\stackrel{u^*}{\to} \Ch(X)\to \Ch^1(X)
\]
is an isomorphism where the first arrow is the canonical monomorphism
and the final arrow is the ``projection''.
As observed in \cite[proof of Lemma 6.2.6]{MNP},
the equality $(E):\phi\circ u_*\circ u^* \circ \phi \circ u_*=\phi\circ u_*$
holds (indeed, $\omega_1\circ \omega_1=\omega_1$ is a direct consequence
of $(E)$).
Thus, $u^*\circ \phi\circ u_*\circ u^* \circ \phi \circ u_*=u^*\circ \phi\circ u_*$. Namely, $u^*\circ \phi\circ u_*:\Ch(X) \to \Ch(X)$ is an idempotent
morphism. The morphism $\pi_1$ coincides with $u^*\circ \phi\circ u_*$.
Actually, again by the observation in \cite[proof of Lemma 6.2.6]{MNP},
$u^*\circ \phi\circ u_*$ corresponds to the composite
$\Alb_X\stackrel{\sim}{\to} \Alb_{\Alb_X}\stackrel{\sigma}{\to} \Pic_{\Alb_X}\stackrel{\sim}{\to} \Pic_X$, that is $\theta$, through
$\CH^1(X\times X)\supset \Hom^*((X,p),(\Pic_X,0))\simeq \Hom_{AV}(\Alb_X,\Pic_X)$.
Let $F:=(u^*\circ \phi\circ u_*)\circ u^*\circ (\phi \circ u_* \circ u^*)$
and $G:=(\phi \circ u_* \circ u^*)\circ \phi\circ u_*\circ (u^*\circ \phi\circ u_*)$. To prove our claim,
it will suffice to show $F\circ G=u^*\circ \phi \circ u_*$
and $G\circ F=\phi\circ u_*\circ u^*$. These equalities follow from 
$(E)$. We also see that $\phi\circ u_*$
induces $\Ch^1(X)\to\Ch(X)\stackrel{\phi\circ u_*}{\to} \Ch(\Alb_X)\to \Ch^{1}(\Alb_X)$ is an isomorphism.

Let $\Free_{\Chow_k}(\Ch^1(\Alb_X))=\oplus_{n\ge0}\Sym^n(\Ch^1(\Alb_X))$
be the free commutative algebra object in $\Chow_k$ and
let $h:\Free_{\Chow_k}(\Ch^1(\Alb_X))\to \Ch(\Alb_X)$ be the morphism
of commutative algebra objects that is classified by $\Ch^1(\Alb_X)\to \Ch(\Alb_X)$. Here $\Ch(\Alb_X)$ admits the commutative algebra structure defined by
$\Ch(\Spec k)\to \Ch(\Alb_X)$ induced by the structure morphism
and $\Ch(X)\otimes \Ch(X)\to \Ch(X)$ induced by the diagonal.
We will show that $h$ is an isomorphism.
Let $R_l:\Chow_k\to \textup{GrVect}$ be the (symmetric monoidal)
$l$-adic realization functor to the category of $\ZZ$-graded $\QQ_l$-vector
space (the symmetric monoidal structure on $\textup{GrVect}$
adopts the Koszul rule). For a projective smooth variety $U$,
it carries $\Ch(U)$ to the $\ZZ$-graded $\QQ_l$-vector space $\textup{H}^*_{\textup{\'et}}(U\times_{k}\overline{k},\QQ_l)$ of \'etale cohomology 
($\overline{k}$ is a separable closure).
Then $\textup{H}^*_{\textup{\'et}}(\Alb_{X}\times_{k}\overline{k},\QQ_l)$
is the free commutative graded algebra generated by
$\textup{H}^1_{\textup{\'et}}(\Alb_{X}\times_{k}\overline{k},\QQ_l)$ placed in
degree one.
By \cite[Theorem 6.2.1]{MNP}, $R_l(\Ch^1(\Alb_X))$ is $\textup{H}^1_{\textup{\'et}}(\Alb_X\times_{k}\overline{k},\QQ_l)$ placed in
degree one, and $R_l(\Ch^1(\Alb_X)\to \Ch(\Alb_X))$ is
$\textup{H}^1_{\textup{\'et}}(\Alb_X\times_{k}\overline{k},\QQ_l)\hookrightarrow
\textup{H}^*_{\textup{\'et}}(\Alb_X\times_{k}\overline{k},\QQ_l)$.
We then conclude that $R_l(h)$ is an isomorphism.
Since $\Sym^{N}(\Ch^1(\Alb_X))=\Sym^N(\Ch^1(X))=0$ for $N>2g$ (see e.g. \cite{Kun}),
$\Free_{\Chow_k}(\Ch^1(\Alb_X))=\oplus_{n=0}^{2g}\Sym^n(\Ch^1(\Alb_X))$.
Both $\Ch(\Alb_X)$ and $\Free_{\Chow_k}(\Ch^1(\Alb_X))$
are Kimura finite (see e.g. \cite{MNP} for this notion).
Thanks to Andr\'e-Kahn \cite[Proposition 1.4.4.(b), Theorem 9.2.2]{AK} (explained also in \cite[Theorem 1.3.1]{AEH}),
we deduce from the isomorphism $R_l(h)$ that $h$ is an isomorphism.
(We remark that $h$ is not necessarily compatible with the equivalence
in Theorem~\ref{abelianvar}.)

Next consider the composition
\[
\psi:\Free_{\Chow_k}(\Ch^1(\Alb_1))=\oplus_{i=0}^{2g}\Sym^{i}(\Ch^1(\Alb_X))
\stackrel{h}{\simeq}\Ch(\Alb_{X})\stackrel{u^*}{\to} \Ch(X)\stackrel{\pi_1}{\to} \Ch^1(X).
\]
We will
show that for $i\ne1$, $\Sym^{i}(\Ch^1(\Alb_X))\hookrightarrow \Ch(\Alb_{X})\stackrel{\psi}{\to} \Ch^1(X)$ is zero. Let $\omega_i:\Ch(\Alb_X)\to \Sym^{i}(\Ch^1(\Alb_X))\to \Ch(\Alb_X)$ denote the idempotent map arising from the direct summand
$\Sym^{i}(\Ch^1(\Alb_X))$. Note that $\omega_1:\Ch(\Alb_X)\to \Ch(\Alb_X)$
equals to $\Ch(\Alb_X)\stackrel{u^*}{\to}\Ch(X)\stackrel{\pi_1}{\to}\Ch^1(X)\hookrightarrow \Ch(X)\stackrel{\phi\circ u_*}{\to} \Ch(\Alb_X)$.
Indeed, the equality $(E)$ implies that $\phi\circ  u_*\circ  \pi_1\circ  u^*=\phi \circ u_* \circ u^*=\omega_1$.
Suppose that $\Sym^{i}(\Ch(\Alb_X))\to \Ch^1(X)$ induced by $\psi$ is not
zero. It follows that $w_1\circ w_i$ is not zero
because $\phi\circ u_*$ induces the isomorphism $\Ch^1(X)\simeq \Ch^1(\Alb_X)\subset\Ch(\Alb_X)$.
For $i\ne 1$, this contradicts the orthogonality $w_1\circ w_i=0$. 
Hence $\Sym^{i}(\Ch(\Alb_X))\to \Ch^1(X)$ is zero for $i\ne 1$.
Remember that $\Ch(X)$ has a commutative algebra structure (in $\Chow_k$)
that is defined by the structure morphism and the diagonal in the
same way as $\Ch(\Alb_X)$. In addition, $u^*$ is 
a homomorphism of commutative algebras.
The homomorphism $u^*$ induces an isomorphism $\Ch^1(\Alb_X)\simeq \Ch^1(X)$.
Taking account of the compatibility of multiplications,
we see that the multiplication $\Ch(X)\otimes \Ch(X)\to \Ch(X)\simeq \Ch^1(X)\oplus \Ch^1(X)\oplus \Ch^2(X)$
sends $\Ch^1(X)\otimes \Ch^1(X)$ to the direct summand $\Ch^2(X)\subset \Ch(X)$.
Namely, the composition $\Ch^1(X)\otimes \Ch^1(X)\hookrightarrow \Ch(X)\otimes \Ch(X)\to \Ch(X)\to \Ch^1(X)$ is zero (by the compatibility of the unit maps
the projection to $\Ch^0(X)$ is also zero).
Now move to $\hhh(\DM(k))$.
The commutative algebra $\Ch(X)$ corresponds to the cocommutative coalgebra
$M(X)$ whose coalgebra structure is determined by the structure morphism
and the diagonal. Take the dual of $M(X)$, that is, $M_X$ in $\hhh(\DM(k))$.
According to Lemma~\ref{multcomult},
the algebra structure of $M_X$ (in $\hhh(\DM(k))$)
given by the coalgebra structure of $M(X)$ coincides with that of
$M_X$ induced by the cohomological motivic algebra.
We have proved that the multiplication $M_X\otimes M_X\to M_X$
sends the component $M_X^1\otimes M_X^1$ to the direct summand $M^2_X$.
If $j:C\to X$ denote the open immersion, then we have
the dual of the Gysin distinguished triangle \cite[14.5]{MVW}
\[
\uni_k(-1)[-2] \stackrel{\eta}{\to} M_X\stackrel{j^*}{\to} M_C\to 
\]
in $\hhh(\DM(k))$. By the exact sequence
$\CH^0(\Spec k)\stackrel{p_*}{\to} \CH^1(X)\stackrel{j^*}{\to} \CH^1(C)\to 0$,
the composite 
$\uni_k(-1)[-2]\stackrel{e}{\to}M_X\stackrel{j^*}{\to} M_C$ is zero
where $e$ is the dual of the morphism
$M(X)\to \uni_k(1)[2]$ corresponding to the class of $p$ in $\CH^1(X)$.
Since $e$ is non-trivial and $\End_{\hhh(\DM(k))}(\uni_k(-1)[-2])=\QQ$, thus we may suppose that $e=\eta$.
Consequently, $\eta$ is the canonical inclusion
$\uni_k(-1)[-2]\simeq M^2_X\hookrightarrow M^0_X\oplus M_X^1\oplus M_X^2$.
It follows that
$j^*$ is identified with the projection $M_X\simeq M^0_X\oplus M_X^1\oplus M_X^2\to M_X^0\oplus M_X^1\simeq M_C$ (with respect this decomposition).
Therefore, the multiplication $M_C\otimes M_C\to M_C$ sends
the component $M_X^1\otimes M_X^1\subset M_C\otimes M_C$
to zero.
\QED

\begin{Remark}
One can ask
whether or not $M_C$ is a trivial square zero extension of $\uni_k$
by some motive $M$ in $\DM(k)$
(not only at the level of  $\hhh(\DM(k))$).
It would be an interesting problem.
We refer to \cite[7.3.4]{HA} for the notion of trivial square zero extensions
in $\infty$-categorical setting.
If $M_C$ is the trivial square zero extension
at the level of $\DM(k)$,
it should be regarded as {\it formality} of $M_C$.
Suppose that we are given a connected affine smooth
curve $C$ over $\CC$.
Write $C^t$ for the underlying topological space of $C$.
Then $A_{PL}(C^t)$ is equivalent to
the square zero extension $\QQ\oplus H^1(C^t,\QQ)[-1]$
of $\QQ=H^0(C^t,\QQ)$ by $H^1(C^t,\QQ)[-1]$ in $\CAlg_{\QQ}$.
(Namely, $A_{PL}(C^t)$ is formal.)
The problem about formality of $M_C$ makes sense for arbitrary
(geometrically connected) affine smooth curves.
\end{Remark}

\subsubsection{}
\label{curve3proof}
Before the proof of Theorem~\ref{curve3}, we will define $W_2$.
Let $K$ be the standard representation of $\GL_{2g}$, that is,
the $2g$-dimensional vector space $V$ endowed with the canonical
action of $\Aut(V)=\GL_{2g}$.
We usually consider $K$ to be the complex
concentrated in degree zero, that belongs to either $\Comp(\GL_2)$
or $\Rep(\GL_2)$.
Let $\Free_{\Comp(\GL_2)}(K[-1])$ is the free commutative algebra object
in $\Comp(\GL_{2g})$, that is isomorphic to
$\oplus_{i=0}^{2g}\Sym^i(K[-1])$ as an object of $\Comp(\GL_{2g})$.
Put $U_1=\oplus_{i=2}^{2g}\Sym^i(K[-1])$ and consider the inclusion
$\alpha:U_1\hookrightarrow \oplus_{i=0}^{2g}\Sym^i(K[-1])$.
Let $\phi_{\alpha}:\Free_{\Comp(\GL_{2g})}(U_1)\to \Free_{\Comp(\GL_{2g})}(K[-1])$ be the morphism classified by $\alpha$
and let $\Free_{\Comp(\GL_{2g})}(K[-1])\langle \alpha\rangle$
be the homotopy pushout (cf. Proposition~\ref{hopushout}).
Consider 
$\Free_{\Comp(\GL_{2g})}(K[-1])\langle \alpha\rangle$
as an object in $\Comp(\GL_{2g})$. Then by the explicit presentation
in Remark~\ref{explicitdg}, we find that
its $0$-th cohomogy is the unit, and the first cohomology is $K$.
Thus, $\Free_{\Comp(\GL_{2g})}(K[-1])\langle \alpha\rangle$
is equivalent to $\uni\oplus K[-1]\oplus U_2$ in $\Rep(\GL_{2g})$
where $\uni$ is a unit object in $\Rep(\GL_{2g})$, and $U_2$ is concentrated
in the degrees larger than one.
(We remark that in practice one can compute $U_2$ explicitly
by means of the representation theory of $\GL_{2g}$.)
Note that the wedge product $\wedge^{N}(M^1_{\overline{C}}[1])=\Sym^N(M^1_{\overline{C}})[N]$
is zero exactly when $N>2g$
because $M^1_{\overline{C}}$ is equivalent to
the dual of the direct summand $M_1(\Alb_{\overline{C}})$
arising from $\Ch^1(\Alb_{\overline{C}})$
(see the proof of Lemma~\ref{lemmasqzero}),
and $\Sym^i(M_1(\Alb_{\overline{C}})^\vee)=0$
for $N>2g$, see e.g. \cite{Kun} for this vanishing.
By \cite[Theorem 3.1, Proposition 6.1]{DTD}, 
there is a
colimit-preserving symmetric monoidal functor $F:\Rep^\otimes(\GL_{2g})\to \DM^\otimes(k)$ which carries $K$ to $M^1_{\overline{C}}[1]$.
Indeed, $F$ is a unique $\QQ$-linear symmetric monoidal functor having
this property (see \cite{DTD} for the detail of the formulation).
We define $W_2$ to be $F(U_2)$.

{\it Proof of Theorem~\ref{curve3}.}
We first prove (1). For simplicity, we put $X=\overline{C}$.
Taking account of the construction of the
decomposition $\Ch(C)\simeq \Ch^0(X)\oplus\Ch^1(X)$ in the proof of Lemma~\ref{lemmasqzero}, $\uni_k=M_{\Spec k}\to M_C$ induced by the
structure morphism $C\to \Spec k$ is
identified with the inclusion $\uni_k=M^0_{X}\hookrightarrow M^0_{X}\oplus M^1_{X}$.
Note that the unit algebra $\uni_k$ is an initial object. Thus,
the kernel of $A_0=\uni_k\to M_{C}$ in $\DM(k)$ is $M^1_{X}[-1]$,
that is, $V_0=M^1_{X}[-1]$.
Therefore $A_1=\uni_k\otimes_{\Free(M^1_{X}[-1])}\uni_k\simeq \Free(0\sqcup_{M^1_{X}[-1]}0)\simeq \Free(M^1_{X})$ ($\sqcup$ indicates the pushout).
The composite $M_X\to \Free(M^1_{X})\to M_C$
is equivalent to the inclusion
$M^1_X\simeq 0\sqcup_{M^1_{X}[-1]}0\to \uni_k\sqcup_{M^1_{X}[-1]}0\simeq \uni_k\oplus M^1_X$ where the second arrow is induced by $0\to \uni_k$.
Next we prove (2).
Let $M(\Alb_X)\to M_1(\Alb_X)$ be the morphism
arising from $\Ch^1(\Alb_X)\to \Ch(\Alb_X)$ (see the second paragraph of
the proof of Lemma~\ref{lemmasqzero}).
If one takes its dual
$M^1_{\Alb_{X}}=M_1(\Alb_X)^\vee\to M_{\Alb_X}$, then
by Lemma~\ref{homotopyfree},
$\Free(M^1_{\Alb_{X}})\to M_{\Alb_X}$ classified by it is an equivalence
in $\CAlg(\DM^\otimes(k))$.
By the isomorphism $\Ch^1(\Alb_X)\simeq \Ch^1(X)$ in the third
paragraph of the proof of Lemma~\ref{lemmasqzero},
the composite $\Free(M^1_{\Alb_{X}})\to M_{\Alb_X}\stackrel{u^*}{\to} M_{X}\stackrel{j^*}{\to} M_C$ induces an equivalence
$M^1_{\Alb_X}\to \Free(M^1_{\Alb_{X}})\to M_C\to M^1_X$.
Also, $M^1_{\Alb_X}\to M_C\to M^0_X$ is null homotopic.
Consider $\Free(M^1_{X})\simeq \Free(M^1_{\Alb_X})$ induced by
$\Ch^1(X)\simeq \Ch^1(\Alb_X)$. Then
$\Free(M^1_{X})\simeq \Free(M^1_{\Alb_X})\simeq M_{\Alb_X}\stackrel{j^*u^*}{\to} M_C$ is equivalent to $A_1\to M_C$.

Next we prove (3).
Let $V_1$ is the kernel of $A_1=\Free(M^1_{X})\to M_C$.
Then $M^1_{X} \to\Free(M^1_{X})\to M_C\simeq \uni_k\oplus M^1_{X}$
may be viewed as the inclusion.
In addition, Lemma~\ref{lemmasqzero} shows that
$M_C\otimes M_C\to M_C$ kills $M^1_X\otimes M^1_X$.
Thus, taking account of the commutative algebra
structure of $\Free(M^1_{X})$ in $\hhh(\DM(k))$
we deduce that $\Free(M^1_{X})=\oplus_{i=0}^{2g}\Sym^i(M^1_X)\to M_C\simeq \uni_k\oplus M^1_{X}$ can be identified with the projection.
Hence $V_1\to \Free(M^1_{X})$ is $\oplus_{i=2}^{2g}\Sym^i(M^1_X)\hookrightarrow \oplus_{i=0}^{2g}\Sym^i(M^1_X)$.
Let $\Free(V_1)\to \Free(M^1_{X})$ is the morphism classified by $V_1\to\Free(M^1_{X})$. Thus $A_2=\Free(M^1_{X})\otimes_{\Free(V_1)}\uni_k$.

Next we prove that (4).
Note that we already defined an ``explicit'' model of $A_2$ before this proof.
Namely, $A_2$ is equivalent to the image of
$\Free_{\Comp(\GL_{2g})}(K[-1])\langle \alpha\rangle$
under $\CAlg(F):\CAlg(\Rep^\otimes(\GL_{2g}))\to \CAlg(\DM^\otimes(k))$.
Thus $A_2\simeq \uni_k\oplus M^1_{X}\oplus W_2$.
Moreover, using the sequence $A_1\to A_2\to M_C$
we find that the composite $r:\uni_k\oplus M^1_{X}\hookrightarrow \uni_k\oplus M^1_{X}\oplus W_2\simeq A_2\to M_C\simeq \uni_k\oplus M^1_{X}$ is an equivalence.
Put $h:W_2\hookrightarrow\uni_k\oplus M^1_{X}\oplus W_2\simeq A_2\to M_C\simeq \uni_k\oplus M^1_{X}$.
Then $H=(-r^{-1}\circ h) \oplus \textup{id}_{W_2}:W_2\to (\uni_k\oplus M^1_{X})\oplus W_2\simeq A_2$
is the kernel of $A_2\to M_C$ (we expect that $h$ is zero).
Let $\Free(W_2)\to A_2$ be the morphism classified by $H$.
Then $A_3=A_2\otimes_{\Free(W_2)}\uni_k$.
\QED

\section{Cotangent complexes and homotopy groups}
\label{cotangenthomotopy}

We introduce a cotangent motive of a pointed
(smooth) scheme $(X,x)$.
Under a suitable condition,
the dual of rationalized homotopy groups will appear as the realization of the cotangent
motive. The notion of cotagent motives is inspired by
Sullivan's description of homotopy groups
in terms of the space of indecomposable elements of a minimal Sullivan model.
We may think of cotangent motive as motives
of (dual of) rationalized homotopy groups.
In this Section, the coefficient ring of $\DM(k)$ is $\QQ$.

\subsection{}
Let $(X,x:\Spec k\to X)$ be a pointed smooth scheme
over $k$. It gives rise to an augmented
object $x^*:M_X\to \uni_k=M_{\Spec k}$.
We will define an object of $\DM(k)$ by means
of cotangent complexes for $\CAlg(\DM^\otimes(k))$.
For this purpose, we use the theory of cotangent complexes
for presentable $\infty$-categories, developed in \cite[Section 7.3]{HA}.
This theory is a vast generalization of cotangent complexes
(topological Andr\'e-Quillen homology) for $E_\infty$-algebras.
Let us briefly recall some definitions about cotangent complexes for
the reader's convenience.
Let $\CCC$ be a presentable $\infty$-category
and let $A$ be an object in $\CCC$.
Let $\SP(\CCC_{/A})$ be the stabilization (stable envelope) of $\CCC_{/A}$ (cf. \cite[1,4]{HA}).
Let $(\CCC_{/A})_{*}$ denote the $\infty$-category of pointed objects
of $\CCC_{/A}$: one may take $(\CCC_{/A})_{*}=(\CCC_{/A})_{A/}$.
Then $\SP(\CCC_{/A})$ is defined to be the limit of the sequence of $\infty$-categories
\[
\cdots \stackrel{\Omega_*}{\to} (\CCC_{/A})_{*} \stackrel{\Omega_*}{\to} (\CCC_{/A})_{*} \stackrel{\Omega_*}{\to} (\CCC_{/A})_{*},
\]
where $\Omega_*$ is informally given by $S\mapsto *\times_{S}*$
($*$ is a final object). The stable $\infty$-category $\SP(\CCC_{/A})$ is also
presentable.
Another presentation of $\SP(\CCC_{/A})$
is the $\infty$-category of spectrum objects of $\CCC$, see \cite[1.4.2]{HA}.
There is a canonical functor $\Omega^{\infty}:\SP(\CCC_{/A})\to (\CCC_{/A})_*\to \CCC_{/A}$ where the first arrow
is the projection to $(\CCC_{/A})_*$ placed in the right end in the above sequence, and the second arrow is the forgetful functor.
Let $\Sigma_{+}^{\infty}:\CCC_{/A} \to \SP(\CCC_{/A})$ be a left adjoint to
$\Omega^\infty$, whose existence is ensured by adjoint functor theorem
since $\Omega^\infty$ preserves small limits and is accessible.
An absolute cotangent complex $L_{A}$ of $A$ is defined to be $\Sigma^{\infty}_+(A\stackrel{\textup{id}}{\to} A)$.
If $A$ is an initial object, then $L_A$ is a zero object.
We now take $\CCC$ to be $\CAlg(\DM^\otimes(k))$. 
By \cite[7.3.4.13]{HA}, there is a canonical equivalence
$\SP(\CAlg(\DM^\otimes(k))_{/A})\simeq \Mod_{A}(\DM(k))$
of $\infty$-categories. Here $\Mod_{A}(\DM(k))$ denotes
the $\infty$-category of $A$-module objects in $\DM(k)$.
We refer to \cite[3.3.3, 4.5]{HA} for the notion of module
objects over a commutative algebra object.
We have the adjunction
\[
\Sigma_+^{\infty}:\CAlg(\DM^\otimes(k))_{/A}\rightleftarrows \SP(\CAlg(\DM^\otimes(k))_{/A})\simeq \Mod_{A}(\DM(k)):\Omega^{\infty}.
\]
We regard $L_A=\Sigma_+^{\infty}(A\stackrel{\textup{id}}{\to} A)$ as an object of $\Mod_{A}(\DM(k))$.
Let $\phi:A\to B$ be a morphism in $\CAlg(\DM^\otimes(k))$.
Let $(-)\otimes_AB:\Mod_{A}(\DM(k))\to \Mod_B(\DM(k))$ denote
a left adjoint to the forgetful functor $\Mod_B(\DM(k))\to\Mod_{A}(\DM(k))$
induced by $A\to B$. Then as in the classical theory of cotangent complexes,
there is a canonical morphism $L_A\otimes_{A}B \to L_B$;
indeed, $L_A\otimes_{A}B\simeq \Sigma_+^{\infty}(A \stackrel{\phi}{\to} B)$
when $A \to B$ is thought of as an object of $\CAlg(\DM^\otimes(k))_{/B}$
(see \cite[7.3.2.14,
7.3.3, 7.3.4.18]{HA} and Remark~\ref{basefree}).
We define the relative cotangent complex $L_{B/A}$ of $A\to B$ to be
a cokernel (cofiber) of $L_A\otimes_{A} B \to L_B$ in $\Mod_B(\DM(k))$.

\begin{Definition}
\label{cotangentmotivedef}
Let $(X,x)$ be a pointed smooth scheme separated of finite type over $k$.
Let $x^*:M_X\to \uni_k=M_{\Spec k}$ be the morphism induced by $x$.
We define $LM_{(X,x)}$  to be
$L_{M_X}\otimes_{M_X}\uni_k$ in $\DM(k)$.
Here $L_{M_X}$ belongs to $\Mod_{M_X}(\DM(k))$, and $(-)\otimes_{M_X}\uni_k:\Mod_{M_X}(\DM(k))\to \Mod_{\uni_k}(\DM(k))\simeq \DM(k)$
is induced by $x^*$.
We shall refer to $LM_{(X,x)}$
as the {\it cotangent motive} of $X$ at $x$.
For $i\in \ZZ$
and $j\in\ZZ$,
we define $\prod_{i,j}(X,x):=\Hom_{\hhh(\DM(k))}(LM_{(X,x)},\uni_k(-j)[-i])$.
\end{Definition}

\begin{Remark}
There is a canonical equivalence $LM_{(X,x)}=L_{M_X}\otimes_{M_X}\uni_k\simeq L_{\uni_k/M_X}[-1]$.
Indeed, there is the distinguished triangle (cofiber sequence) arising from
$\uni_k\to M_X\to \uni_k$:
\[
L_{M_X}\otimes_{M_X}\uni_k\to L_{\uni_k/\uni_k}\to L_{\uni_k/M_X}\to
\]
in the homotopy category of $\DM(k)$, see \cite[7.3.3.5]{HA}.
In addition, $L_{\uni_k/\uni_k}\simeq 0$.
It follows that
$L_{M_X}\otimes_{M_X}\uni_k\simeq  L_{\uni_k/M_X}[-1]$.
\end{Remark}

\begin{Remark}
The definition of the cotangent motives makes sense also when
we work with an arbitrary coeffiecient ring $K$ of $\DM(k)$.
\end{Remark}

The main result of this Section is the following:

\begin{Theorem}
\label{motivichomotopy}
Let $X$ be a smooth variety over $k$ and let $x$ be a $k$-rational point.
Suppose that $k$ is embedded in $\CC$
and the underlying topological space
$X^t$ of $X\times_{\Spec k}\Spec \CC$ is simply connected.
Then the (singular) realization functor $\mathsf{R}:\DM(k)\to \mathsf{D}(\QQ)$ carries
$LM_{(X,x)}$ to $\oplus_{2\le i }(\pi_i(X^t,x)\otimes_{\ZZ}\QQ)^{\vee}[-i]$
in $\mathsf{D}(\QQ)\simeq \Mod_{\QQ}$.
Namely,
there is an isomorphism
\[
H^i(\mathsf{R}(LM_{(X,x)}))\simeq (\pi_i(X^t,x)\otimes_{\ZZ}\QQ)^{\vee}
\]
for $i\ge 2$, and $H^i(\mathsf{R}(LM_{(X,x)}))=0$ for $i<2$.
Here $(\pi_i(X^t,x)\otimes_{\ZZ}\QQ)^{\vee}$
is the dual $\QQ$-vector space of $\pi_i(X^t,x)\otimes_{\ZZ}\QQ$.
\end{Theorem}

\begin{Remark}
Under the hypothesis of Theorem~\ref{motivichomotopy},
the cohomology $H^i(X^t,\QQ)$ is finite dimensional for any $i\ge 0$.
Indeed, the simply connectedness is not necessary for this finiteness.
In general, if $S$ is a simply connected topological space
whose cohomoogy $H^i(S,\QQ)$ is finite dimensional for any $i\ge0$, then
$\pi_{i}(S,s)\otimes_{\ZZ}\QQ$ is finite dimensional for any $i\ge2$.
\end{Remark}

\subsection{}
 The proof proceeds in several steps.

\begin{Lemma}
\label{realmorita}
Let $\mathsf{R}_E:\DM^\otimes(k)\to \mathsf{D}^\otimes(K)$ be the symmetric monoidal realization functor associated to mixed Weil theory $E$ with coefficients in a field $K$ of characteristic zero.
Let $G$ be a right adjoint to $\mathsf{R}_E$, that is lax symmeytric monoidal.
Let $G(K)$ be the commutative algebra object (i.e., an object of $\CAlg(\DM^\otimes(k))$) where $K$ is the unit algebra in $\mathsf{D}(K)$.
Consider the composition of symmetric monoidal functors
\[
 \Mod_{G(K)}^{\otimes}(\DM(k))\to \Mod_{\mathsf{R}_E(G(K))}^\otimes(\mathsf{D}(K))\to \Mod_{K}^\otimes(\mathsf{D}(K)) \simeq \mathsf{D}^\otimes(K)
\]
where the first arrow is induced by $\mathsf{R}_E$, and 
the second arrow is given by the base change
$(-)\otimes_{\mathsf{R}_E(G(K))}K$ induced by the counit map
$\mathsf{R}_E(G(K))\to K$.
Then the composite is an equivalence, and $\mathsf{R}_E$
is
equivalent to the base change functor
$(-)\otimes_{\uni_k}G(K):\DM^\otimes(k)\to \Mod_{G(K)}^{\otimes}(\DM(k))\simeq \mathsf{D}^\otimes(K)$.
\end{Lemma}

\Proof
If we verify two conditions
\begin{itemize}
\item there is a set $\{M_{\lambda}\}_{\lambda\in \Lambda}$ of compact and
dualizable objects of $\DM(k)$ such that
the whole category $\DM(k)$ is the smallest
stable subcategory which contains $\{M_{\lambda}\}_{\lambda\in \Lambda}$
and is closed under small coproducts
(that is to say, $\{M_{\lambda}\}_{\lambda\in \Lambda}$ is a generator
of $\DM(k)$),

\item each $\mathsf{R}_E(M_{\lambda})$ is compact, and there is some
$\mu\in I$ such that $\mathsf{R}_E(M_{\mu})\simeq K$,
\end{itemize}
then our assertion follows from \cite[Proposition 2.1]{DTD}.
For $X\in \Sm_k$ and $n\in \ZZ$, $M(X)(n)$ is compact in $\DM(k)$,
and the set $\{M(X)(n)\}_{X\in \Sm_k,n\in \ZZ}$ is a generator of $\DM(k)$.
In addition, $M(X)$ is dualizable
because it holds if $X$ is projective,
and 
we work with rational coefficients, so that we can use 
the standard argument based on de Jong's
alteration (or one can directly apply a very general result in
\cite[4.4.3, 4.4.17]{CDT}).
Since $\mathsf{R}_E$ is symmetric monoidal and $M(X)(n)$
is dualizable, $\mathsf{R}_E(M(X)(n))$ is also dualizable.
In $\mathsf{D}(K)$, an object is compact if and only if it is dualizable.
Finally, $\mathsf{R}_E(M(\Spec k))=\mathsf{R}_E(\uni_k)\simeq K$
since $\mathsf{R}_E$ is symmetric monoidal.
Hence the above two conditions are satisfied.
\QED

According to Lemma~\ref{realmorita}, under the setting of Theorem~\ref{motivichomotopy},
we write $P:=G(\QQ)$ and identify the (singular) realization functor
$\mathsf{R}$ with $(-)\otimes_{\uni_k}P:\DM^\otimes(k)\to \Mod_P^\otimes(\DM(k))\simeq \mathsf{D}^\otimes(\QQ)$.
The multiplicative realization functor
$\CAlg(\mathsf{R}):\CAlg(\DM^\otimes(k))\to \CAlg_{\QQ}$
can naturally be identified with
\[
\CAlg(\DM^\otimes(k))\longrightarrow \CAlg(\Mod_P^\otimes(\DM(k))\simeq \CAlg(\DM^\otimes(k))_{P/}
\]
which sends $A$ to $P\simeq \uni_k\otimes P \to A\otimes P$.
For the right equivalence, see \cite[3.4.1.7]{HA}.

We focus on cotangent complexes of commutative dg algebras,
that is, objects of $\CAlg_{\QQ}$.
Let $C$ be an object of $\CAlg_{\QQ}$.
If we take $\CCC$ to be $\CAlg_{\QQ}$ in the above formalism of cotangent
complexes, we have the adjunction
\[
\Sigma_+^{\infty}:(\CAlg_{\QQ})_{/C}\rightleftarrows \SP((\CAlg_{\QQ})_{/C})\simeq \Mod_{C}(\mathsf{D}(\QQ)):\Omega^{\infty}.
\]
Here we abuse notation by using $\Sigma_+^{\infty}, \Omega^{\infty}$ again.
We define the absolute
cotangent complex
$L_C$ of $C$ to be $\Sigma_+^{\infty}(C\stackrel{\textup{id}}{\to} C)$.
Given a morphism $C\to D$ we define $L_{D/C}$ to be a cokernel of
$D\otimes_CL_C\to L_D$ in $\Mod_{D}(\mathsf{D}(\QQ))$.

\begin{Remark}
\label{basefree}
Let $\mathcal{C}$ be either $\CAlg(\DM^\otimes(k))$ or
$\CAlg_{\QQ}=\CAlg(\mathsf{D}^\otimes(\QQ))$.
More generally, $\mathcal{C}$ could be a presentable
$\infty$-category $\CAlg(\DDD^\otimes)$
such that $\DDD^\otimes$ is a symmetric monoidal stable
presentable $\infty$-category
whose tensor product $\DDD\times \DDD\to \DDD$ preserves small colimits
separately in each variable.
Let $A$ and $B$ be objects of $\CCC$. Let $B\to A$ be a morphism.
Consider the adjunction
\[
\Sigma_+^{\infty}:\CCC_{/A} \rightleftarrows \SP(\CCC_{/A})\simeq \Mod_A(\DDD):\Omega^{\infty}.
\]
If we regard $B\to A$ as an object of $\CCC_{/A}$, then
$\Sigma_+^{\infty}$ sends $B\to A$ to $L_{B}\otimes_BA$,
where $(-)\otimes_BA:\Mod_B(\DDD)\to \Mod_A(\DDD)$ denotes the base change
functor. It is a direct consequence of a functorial construction
of cotangent complexes by using the notion
of a tangent bundle in \cite[7.3.2.14]{HA}
and a presentation of the tangent bundle by a presentable fibration
of module
categories \cite[7.3.4.18]{HA}.

Suppose that $A$ is an initial object (that is, a unit algebra).
The above adjunction is extended to 
\[
\DDD \rightleftarrows \CCC_{/A} \rightleftarrows \SP(\CCC_{/A})\simeq \DDD
\]
where the left arrow
$j:\CCC_{/A}=\CAlg(\DDD^\otimes)_{/A} \to \DDD$
is the functor which carries $\epsilon:B\to A$ to the kernel (fiber) $\Ker(\epsilon)$
of $B\to A$ in $\DDD$. The left adjoint $\DDD\to \CCC_{/A}$ to $j$
sends $M\in \DDD$ to
$\Free_{\DDD}(M)\to \Free_{\DDD}(0)\simeq A$ determined by $M\to 0$,
where $\Free_{\DDD}:\DDD\to \CCC$ is the free functor,
see Definition~\ref{freealg}.
By the construction of $\SP(\CCC_{/A})\simeq \DDD$ (cf. \cite[7.3.4.13]{HA}),
the composite $\DDD\simeq \SP(\CCC_{/A}) \stackrel{\Omega^\infty}{\to} \CCC_{/A}\stackrel{j}{\to} \DDD$ is naturally equivalent to the identity functor.
Thus $\Sigma^{\infty}_+$ carries $\Free_{\DDD}(M)\to  A$
to $M$. Namely, $L_{\Free_{\DDD}(M)}\otimes_{\Free_{\DDD}(M)}A\simeq M$.
\end{Remark}

\begin{Remark}
If one considers $x^*:M_X\to \uni_k$ to be an object
of $\CAlg(\DM(k))_{/\uni_k}$, then
its image under $\Sigma^\infty_+:\CAlg(\DM(k))_{/\uni_k}\to \DM(k)$
is $LM_{(X,x)}$
(cf. Remark~\ref{basefree}). The right adjoint $\Omega^\infty:\DM(k)\to \CAlg(\DM(k))_{/\uni_k}$
sends $LM_{(X,x)}$ to a square zero extension of $\uni_k$ by $LM_{(X,x)}$,
which is informally given by $\uni_k\oplus LM_{(X,x)}\stackrel{\textup{pr}_1}{\to} \uni_k$, see \cite[7.3.4]{HA} for square zero extensions.
By the adjunction, we have the unit map
$u:M_X\to \uni_k\oplus LM_{(X,x)}$ in $\CAlg(\DM(k))_{/\uni_k}$.
Let $\overline{M}_X$ be the kernel (fiber) of $M_X\to \uni_k$ in $\DM(k)$. 
It gives rise to a morphism in $\DM(k)$
\[
h:\overline{M}_X\to LM_{(X,x)}
\]
induced by $u$.
This morphism is a motivic version of dual Hurewicz map.
\end{Remark}

\begin{Lemma}
\label{nicebasechange}
Let $\epsilon:A\to \uni_k$ be a morphism in $\CAlg(\DM^\otimes(k))$, that is,
an augmented commutative algebra object in $\DM(k)$.
Let $B:=\mathsf{R}(A)\to \mathsf{R}(\uni_k)= \QQ$ be a
the image of
$\epsilon$ in $\CAlg_{\QQ}$ under the multiplicative realization functor.
Let $L_{B}$ be the (absolute) cotangent complex of $B$
and let $L_B\otimes_{B}\QQ$ be the base change that lies in $\mathsf{D}(\QQ)$.
Then there is a canonical equivalence 
\[
\mathsf{R}(L_{A}\otimes_{A}\uni_k)\simeq L_B\otimes_{B}\QQ
\]
in $\mathsf{D}(\QQ)$.
\end{Lemma}

\Proof
As explained above, Lemma~\ref{realmorita}
allows us to identify the multiplicative realization functor with
$\CAlg(\DM^\otimes(k))\to \CAlg(\DM^\otimes(k))_{P/}\simeq \CAlg_{\QQ}$.
Then we have the pushout diagram
\[
\xymatrix{
\uni_k \ar[r] \ar[d] & P \ar[d] \\
A \ar[r] & A\otimes P
}
\]
in $\CAlg(\DM^\otimes(k))$, and
$B$ corresponds to the right vertical arrow $P \to A\otimes P$
which we regard as an object of $\CAlg(\DM^\otimes(k))_{P/}$.
By \cite[7.3.3.8, 7.3.3.15]{HA},
the absolute cotangent complex of $P \to A\otimes P$ regarded as an object
of $\CAlg(\DM^\otimes(k))_{P/}$ is equivalent to the relative cotangent complex 
$L_{A\otimes P/P}$ of the morphism 
$P\to A\otimes P$ in $\CAlg(\DM(k))$.
It follows that $L_{A\otimes P/P}\simeq L_B$ under the canonical
equivalence $\Mod_{A\otimes P}(\DM(k))\simeq \Mod_{B}(\mathsf{D}(\QQ))$.
The final equivalence is induced by
\begin{eqnarray*}
\Mod_{A\otimes P}(\DM(k)) &\simeq& \SP(\CAlg(\DM^\otimes(k))_{/A\otimes P}) \\
&\simeq& \SP((\CAlg(\DM^\otimes(k))_{P/})_{/A\otimes P}) \\
&\simeq&  \SP((\CAlg_{\QQ})_{/B}) \\
&\simeq& \Mod_{B}(\mathsf{D}(\QQ))
\end{eqnarray*}
where the first and final equivalences follow from \cite[7.3.4.13]{HA},
and the second one follows from \cite[.3.3.9]{HA}.
Since $A\otimes P\to \uni_k\otimes P\simeq P$ corresponds to
$B\to \QQ$, we see that $L_{B}\otimes_{B}\QQ$ corresponds to
$L_{A\otimes P/P}\otimes_{A\otimes P}P$ in $\Mod_{P}(\DM(k))\simeq \mathsf{D}(\QQ)$.
By the base change formula for cotangent complexes \cite[7.3.3.7]{HA},
$L_{A\otimes P/P}\simeq L_{A}\otimes_{A}(A\otimes P)$.
Therefore, we obtain
\[
L_{A\otimes P/P}\otimes_{A\otimes P}P\simeq L_{A}\otimes_{A}(A\otimes P)\otimes_{A\otimes P}P\simeq (L_A\otimes_A\uni_k)\otimes P.
\]
Note that 
$\mathsf{R}(L_{A}\otimes_{A}\uni_k)\simeq (L_A\otimes_A\uni_k)\otimes P$
in $\Mod_P(\DM(k))$. Hence our assertion follows.
\QED

The following is a theorem of Sullivan \cite[Section 8]{S},
reformulated
in terms of cotangent complexes.

\begin{Lemma}
\label{cotangSul}
Let $(S,s)$ be a simply connected topological space $S$
with a point $s$.
Assume that
 the cohomology $H^i(S,\QQ)$ is a
 finite dimensional $\QQ$-vector space for any $i\ge0$.
Let $A_{PL,\infty}(S)$ be the image of $A_{PL}(S)$
in $\CAlg_{\QQ}$ (see Section~\ref{realization}).
Let $A_{PL,\infty}(S)\to \QQ$ be the augmentation induced by $s$.
Then $L_{A_{PL,\infty}(S)}\otimes_{A_{PL,\infty}(S)}\QQ
\simeq \oplus_{2\le i }(\pi_i(S,s)\otimes_{\ZZ}\QQ)^{\vee}[-i]$ in $\mathsf{D}(\QQ)$.
\end{Lemma}

\Proof
For ease of notation, we may assume that
$S$ is a rational space, so that $\pi_i(S,s)$ is a $\QQ$-vector space
for each $i\ge 2$.
Consider a Postnikov tower
\[
S=S_{\infty} \to \cdots \to S_n\to S_{n-1}\to \cdots \to S_2\to S_1.
\]
We first show our assertion in the case of $S_{n}$.
The case of $n=1$ is trivial because $S_1$ is contractible
and $L_{A_{PL,\infty}(S_1)}\simeq 0$. We suppose that
our assertion holds for $S_{n-1}$.
Consider the diagram
\[
\xymatrix{
K(\pi_{n}(S,s),n) \ar[r] \ar[d] & S_n  \ar[d] \\
\ast \ar[r] & S_{n-1}
}
\]
where $\ast$ is a contractible space, $K(\pi_{n}(S,s),n)$
is an Eilenberg-MacLane space, and we here think of the diagram with
a pullback square in $\SSS$.
By a computation for the Eilenberg-MacLane space
\cite[Section 15 Example 3, Section 12, Example 2]{FHT},
$A_{PL,\infty}(K(\pi_{n}(S,s),n))\simeq \Free_{\QQ}(\pi_{n}(S,s)^\vee[-n])$
where $\pi_{n}(S,s)^\vee[-n]$ is the dual $\QQ$-vector space placed in cohomological
degree $n$, that we consider to be an object of $\mathsf{D}(\QQ)$,
and
$\Free_{\QQ}:\mathsf{D}(\QQ)\to \CAlg_{\QQ}$ is the free functor, see Definition~\ref{freealg}.
By \cite[Theorem 15.3]{FHT},
$A_{PL,\infty}(K(\pi_{n}(S,s),n))$ is a pushout of $A_{PL,\infty}(S_{n})\leftarrow A_{PL,\infty}(S_{n-1})\to \QQ\simeq A_{PL,\infty}(\ast)$ in $\CAlg_{\QQ}$
(the result found in \cite{FHT} shows that it is a homotopy pushout in
$\CAlg_{\QQ}^{dg}$).
When $S_{n}$ and $S_{n-1}$ are equipped with (compatible) base points,
$A_{PL,\infty}(K(\pi_{n}(S,s),n))$ is promoted to a pushout in $(\CAlg_{\QQ})_{/\QQ}$. Note that $\Sigma_+^\infty:(\CAlg_{\QQ})_{/\QQ}\to \SP((\CAlg_{\QQ})_{/\QQ})\simeq \mathsf{D}(\QQ)$ preserves small colimits.
Taking account of Remark~\ref{basefree} we have a pushout diagram
\[
\xymatrix{
L_{A_{PL,\infty}(S_{n-1})}\otimes_{A_{PL,\infty}(S_{n-1})}\QQ \ar[r] \ar[d] & L_{A_{PL,\infty}(S_{n})}\otimes_{A_{PL,\infty}(S_{n})}\QQ  \ar[d] \\
0 \ar[r] & \pi_{n}(S,s)^\vee[-n]
}
\]
in $\mathsf{D}(\QQ)$.
By the assumption, $L_{A_{PL,\infty}(S_{n-1})}\otimes_{A_{PL,\infty}(S_{n-1})}\QQ\simeq \oplus_{2\le i \le n-1}\pi_i(S,s)^{\vee}[-i]$.
Then $L_{A_{PL,\infty}(S_{n})}\otimes_{A_{PL,\infty}(S_{n})}\QQ$ is a cokernel (cofiber) of 
$\pi_{n}(S,s)^\vee[-n-1]\to \oplus_{2\le i \le n-1}\pi_i(S,s)^{\vee}[-i]$. Thus the case of $S_{n}$ follows.
Next we show the case of $S$.
For simplicity, $A:=A_{PL,\infty}(S)$ and $A_{n}:=A_{PL,\infty}(S_{n})$.
As the above proof reveals,
$L_{A_{n-1}}\otimes_{A_{n-1}}\QQ \to L_{A_{n}}\otimes_{A_{n}}\QQ$
can be identified with the inclusion
$\oplus_{2\le i \le n-1}\pi_i(S,s)^{\vee}[-i]\to \oplus_{2\le i \le n}\pi_i(S,s)^{\vee}[-i]$.
It will suffice to prove that the canonical morphism
$\varinjlim_{n} L_{A_{n}}\otimes_{A_{n}}\QQ\to L_{A}\otimes_{A}\QQ$
is an equivalence in $\mathsf{D}(\QQ)$.
Since $\Sigma_+^{\infty}$ preserves colimits,
it is enough to show that the canoncial morphism
$\varinjlim_{n}A_n\to A$ is an equivalence $\CAlg_{\QQ}$.
For this we are reduced to proving
the canonical map $\varinjlim_{n}H^{i}(S_n,\QQ)=\varinjlim_{n}H^{i}(A_n)\to H^{i}(S,\QQ)=H^i(A)$
is bijective for $i\ge0$.
By applying Serre spectral sequence
to the fiber sequence $F_{m,n}=*\times_{S_n}S_m \to S_{m}\to S_{n}$
for $n\le m \le \infty$, we see that
$H^{n}(S_n,\QQ)\simeq H^n(S_{n+1},\QQ)\simeq  \ldots \simeq H^n(S,\QQ)$,
so that $\varinjlim_{n}H^{i}(S_n,\QQ)\simeq H^{i}(S,\QQ)$.
\QED

\vspace{2mm}

{\it Proof of Theorem~\ref{motivichomotopy}.}
By Theorem~\ref{real1} and Remark~\ref{real1R},
the image of $M_X\to \uni_k$
can be identified with $A_{PL,\infty}(X^t)\to \QQ$
induced by the point $x$ on $X^t$.
Write $B:=A_{PL,\infty}(X^t)$.
Taking account of Lemma~\ref{nicebasechange},
we see that $\mathsf{R}(LM_{(X,x)})\simeq L_{B}\otimes_{B}\QQ$.
Now our assertion follows from Lemma~\ref{cotangSul}.
\QED

We would like to relate cotangent motives with
partial data of fundamental groups.

\begin{Theorem}
\label{motivichomotopy2}
Let $(X,x)$ be a pointed smooth variety over $k$.
Suppose that
$k$ is embedded in $\CC$.
Let $\pi_i(X^t,x)_{uni}$ be the pro-unipotent completion
of the fundamental group $\pi_i(X^t,x)$ of $X^t$ over $\QQ$.
Then
the $\QQ$-vector space
$H^1(\mathsf{R}(LM_{(X,x)}))$ gets identified with
the cotangent space of the unipotent affine scheme $\pi_1(X^t,x)_{uni}$
at the origin.
\end{Theorem}

\Proof
As in the proof of Theorem~\ref{motivichomotopy},
the image of $M_X\to \uni_k$
can be identified with $A_{PL,\infty}(X^t)\to \QQ$
induced by the point $x$ on $X^t$.
Write $B:=A_{PL,\infty}(X^t)$.
The image of $\uni_k\otimes_{M_X}\uni_k$ under the multiplicative
realization functor can naturally be identified with $\QQ\otimes_{B}\QQ$.
According to Hochschild-Kostant-Rosenberg (HKR)
theorem for $B \in \CAlg_{\QQ}$,
we have $\QQ\otimes_{B}\QQ\simeq \QQ\otimes_{B}B\otimes_{B\otimes B}B\simeq \Free_{\QQ}((L_B\otimes_B\QQ)[1])$
(see e.g. \cite[Prop. 4.4]{BN} for HKR theorem: strictly speaking,
the connectivity on $B$ is assumed in loc. cit., but
its proof shows that the nonconnective affine case holds).
It follows that $H^0(\QQ\otimes_{B}\QQ)\simeq H^0(\Free_{\QQ}((L_B\otimes_B\QQ)[1]))$ (keep in mind that
the dual of the base point $H^0(\QQ\otimes_{B}\QQ)\to H^0(\QQ\otimes_{\QQ}\QQ)\simeq  \QQ$ is identified with $H^0(\Free_{\QQ}((L_B\otimes_B\QQ)[1]))\to H^0(\Free_{\QQ}(0))\simeq \QQ$ induced by $L_B\otimes_B\QQ\to L_{\QQ}\simeq 0$).
Remember that $H^0(\QQ\otimes_{B}\QQ)$ is isomorphic to
the coordinate ring of the pro-unipotent completion
$\pi_1(X^t,x)_{uni}$ of $\pi_1(X^t,x)$
over $\QQ$, cf. Proposition~\ref{undcomp}.
By Lemma~\ref{nicebasechange},
$\mathsf{R}(LM_{(X,x)})\simeq L_B\otimes_{B}\QQ$.
Let us observe that
$H^0(\Free_{\QQ}(\mathsf{R}(LM_{(X,x)})[1]))\simeq 
H^0(\Free_{\QQ}((L_B\otimes_B\QQ)[1]))$
is naturally isomorphic to the free
ordinary commutative $\QQ$-algebra $\Free_{ord}(H^1(L_B\otimes_B\QQ))$
generated by the $\QQ$-vector space
$H^1(L_B\otimes_B\QQ)\simeq H^0(\mathsf{R}(LM_{(X,x)})[1])$.
Taking account of Lemma~\ref{homotopyfree},
we are reduced to showing that $H^i(L_B\otimes_B\QQ)=0$ for $i<1$.
Thus, it will suffice to prove the following Lemma.
\QED

\begin{Lemma}
Let $B$ be a commutative dg algebra over $\QQ$, which
we regard as an object of $\CAlg_{\QQ}$.
Suppose that we are given an augmentation $B\to \QQ$.
Assume that $H^0(B)=\QQ$, and $H^i(B)=0$ for $i<0$.
Then $H^i(L_B\otimes_B\QQ)=0$ for $i<1$.
\end{Lemma}

\Proof
(This fact or equivalent
versions is well-known, but we prove it for the completeness.)
Let $B_0=\QQ\to B_1\to \cdots \to B_n \to\cdots \to  B$
be the inductive sequence associated to
the canonical morphism $\QQ\to B$ in $\CAlg_{\QQ}$,
see Section~\ref{Sull1}.
By Lemma~\ref{Sullivanmodel1},
$\varinjlim_{n}B_n\simeq B$. It follows that $\varinjlim_{n}L_{B_n}\otimes_{B_n}\QQ\simeq L_{B}\otimes_B\QQ$.
Therefore, it is enough
to show that $H^i(L_{B_n}\otimes_{B_n}\QQ)=0$ for $i<1$.
We will prove, by induction on $n\ge 0$,
that (i) $H^0(B_n)=\QQ$, $H^i(B_n)=0$ for $i<0$, 
(ii) $H^1(B_n)\to H^1(B)$ is injective, and
(iii) $H^i(L_{B_n}\otimes_{B_n}\QQ)=0$ for $i<1$.
For $n=0$, this is obvious.
Assume therefore that all (i), (ii), (iii) hold for $n$.
Let $M$ be the kernel (fiber) of $B_n\to B$ in $\mathsf{D}(\QQ)$.
Then $H^i(M)=0$ for $i<2$ by the inductive assumptions (i) and (ii).
By definition, $B_{n+1}=B_{n}\otimes_{\Free_{\QQ}(M)}\QQ$.
By the explicit presentation of the homotopy pushout
$B_{n}\otimes_{\Free_{\QQ}(M)}\QQ$ (Propsition~\ref{hopushout}
and Remark~\ref{explicitdg}),
(i) holds for $B_{n+1}$.
In addition, again by the explicit homotopy pushout,
we have an exact sequence $0\to H^1(B_n)\to H^1(B_{n+1})\to H^2(M)$.
Comparing it with the exact sequence
$0\to H^1(B_n)\to H^1(B)\to H^2(M)$ induced by the
cofiber sequence $M\to B_n\to B$,
we see that $H^1(B_{n+1})\to H^1(B)$ is injective.
Note that $L_{B_{n+1}}\otimes_{B_{n+1}}\QQ$ is a cokernel (cofiber) of
$L_{\Free_{\QQ}(M)}\otimes_{\Free_{\QQ}(M)}\QQ\to L_{B_{n}}\otimes_{B_{n}}\QQ$. By Remark~\ref{basefree},
$L_{\Free_{\QQ}(M)}\otimes_{\Free_{\QQ}(M)}\QQ\simeq M$. 
Taking account of the inductive assumption (iii) for $B_n$,
we conclude that (iii) holds for $B_{n+1}$.
\QED

\subsection{}
We use the explicit computations of
cohomological motivic algebras in Section~\ref{Sullivanmodel}
to obtain explicit presentations of cotangent motives.

\begin{Theorem}
\label{motivichomotopyexp}
We have the following explicit presentations:

\begin{enumerate}
\renewcommand{\labelenumi}{(\theenumi)}

\item
Let $\mathbb{P}^n$ be the projective space over a perfect field $k$
and let $x$ be a $k$-rational point, see Section~\ref{cps}.
Then $LM_{(\mathbb{P}^n,x)}\simeq \uni_k(-1)[-2]\oplus \uni_k(-n-1)[-2n-1]$.

\item Let $X=\mathbb{A}^n-\{p\}$ and let $x$ be a $k$-rational
point, see Section~\ref{minus1}.
Then $LM_{(X,x)}\simeq \uni_{k}(-n)[-2n+1]$.

\item Let $Y=\mathbb{A}^n-\{p\}-\{q\}$ ($n\ge2$)
and let $y$ be a $k$-rational
point, see Section~\ref{minus2}.
Then $LM_{(Y,y)}\simeq \uni_{k}(-n)[-2n+1]^{\oplus 2}\oplus \uni_k(-2n)[-4n+3]\oplus \uni_k(-3n)[-6n+3]^{\oplus 2}\oplus \ldots$.

\item Let $G$ be a semi-abelian variety
and let $o$ be the origin, see Section~\ref{semiab}.
Then $LM_{(G,o)}\simeq M_1(G)^{\vee}$.

\end{enumerate}
\end{Theorem}

\Proof
We show (1). We use the notation in Section~\ref{cps}.
By Propsition~\ref{cpsprop}, $M_{\mathbb{P}^n}\simeq \Free(\uni_{k}(-1)[-2])\otimes_{\Free(\uni_k(-n-1)[-2n-2])}\uni_k$. Let $x^*:M_{\mathbb{P}^n}\to \uni_k$
be the morphism induced by the $k$-rational point $x$.
Note that
$\Free(\uni_{k}(-1)[-2])\to M_{\mathbb{P}^n}\stackrel{x^*}{\to}\uni_k$
is equivalent to
$\Free(\uni_{k}(-1)[-2])\to \Free(0)\simeq \uni_k$
determined by $\uni_{k}(-1)[-2]\to 0$. Indeed, the
morphism $\Free(\uni_{k}(-1)[-2])\to \uni_k$ in $\CAlg(\DM^\otimes(k))$
is classified by the composite $\uni_{k}(-1)[-2]\hookrightarrow \Free(\uni_{k}(-1)[-2])\to \uni_k$ in $\DM(k)$, which is null-homotopic because
$\textup{CH}^1(\Spec k)=0$. Similarly, $\Free(\uni_k(-n-1)[-2n-2])\to M_{\mathbb{P}^n}\to \uni_k$ is equivalent to
$\Free(\uni_k(-n-1)[-2n-2])\to \Free(0)$ determined by $\uni_k(-n-1)[-2n-2]\to 0$.
Since $\Sigma_+^{\infty}:\CAlg(\DM^\otimes(k))_{/\uni_k}\to \SP(\CAlg(\DM^\otimes(k))_{/\uni_k})\simeq \DM(k)$ preserves small colimits, thus
we have the pushout diagram
\[
\xymatrix{
   L_{\Free(\uni_{k}(-n-1)[-2n-2])}\otimes_{\Free(\uni_{k}(-n-1)[-2n-2])}\uni_k  \ar[d] \ar[r] & L_{\Free(\uni_{k}(-1)[-2])}\otimes_{\Free(\uni_{k}(-1)[-2])}\uni_k  \ar[d] \\ 
L_{\uni_k} \ar[r] & LM_{(\mathbb{P}^n,x)}
}
\]
in $\DM(k)$,
cf. Remark~\ref{basefree}. 
Moreover, again by Remark~\ref{basefree},
the upper left term (resp. the upper right term)
is equivalent to $\uni_{k}(-n-1)[-2n-2]$ (resp. $\uni_{k}(-1)[-2]$).
The/any morphism $\uni_{k}(-n-1)[-2n-2]\to \uni_{k}(-1)[-2]$
is null-homotpic.
Combining this consideration with $L_{\uni_k}\simeq0$,
we conclude that $LM_{(\mathbb{P}^n,x)}\simeq \uni_k(-1)[-2]\oplus \uni_k(-n-1)[-2n-1]$. The cases (2) and (4) are easier than (1)
(cf. Proposition~\ref{sphere} and Proposition~\ref{abelianvar}), and the case (3) is similar
to (1) (cf. Proposition~\ref{minus2prop}).
\QED

\begin{Remark}
\label{motcpsR}
Let us consider a meaning of the presentation of the case of projective
spaces. In light of Theorem~\ref{motivichomotopy}, if $k\subset \CC$,
we have
\[
\mathsf{R}(\uni_k(-1))\simeq (\pi_2(\CC\mathbb{P}^n,x)\otimes_{\ZZ}\QQ)^\vee, \ \ \ \ \mathsf{R}(\uni_k(-n-1))\simeq (\pi_{2n+1}(\CC\mathbb{P}^n,x)\otimes_{\ZZ}\QQ)^\vee.
\]
Thus, it is natural to think that $\uni_k(1)$ is a motive for $\pi_2(\CC\mathbb{P}^n,x)\otimes_{\ZZ}\QQ$,
and $\uni_k(n+1)$ is a motive for $\pi_{2n+1}(\CC\mathbb{P}^n,x)\otimes_{\ZZ}\QQ$.
\end{Remark}

\begin{Remark}
According to Theorem~\ref{motivichomotopyexp} (4),
the cotangent motives may also be viewed as
a generalization of (the dual of)
$1$-motives of semi-abelian varieties.
\end{Remark}


\section{Motivic homotopy exact sequence for algebraic curves}
\label{mhes}

Let $X$ be a geometrically connected scheme
of finite type over a perfect field $k$
and let $X_{\bar{k}}$ be the base change to a separable closure $\bar{k}$.
Let $G_k$ denote the absolute Galois group $\Gal(\bar{k}/k)=\pi_1^{\textup{\'et}}(\Spec k,\Spec \bar{k})$.
We write $\pi_1^{\textup{\'et}}(-,a)$ for the \'etale fundamental group of ``$(-)$'' with a base point $a$. Let $\bar{x}:\Spec \bar{k} \to X_{\bar{k}}$ be a geometric
point
and let $x:\Spec \bar{k}\to X$ be the composite.
There is an exact sequence of 
profinite groups
\[
1\to \pi_1^{\textup{\'et}}(X_{\bar{k}},\bar{x})\to \pi_1^{\textup{\'et}}(X,x) \to G_k \to 1
\]
induced by $X_{\bar{k}}=X\times_{\Spec k}\Spec \bar{k} \to X \to \Spec k$.
It is usually called the homotopy exact sequence
because it can be thought of as a fairly precise analogue of
the long exact sequence that comes from
a homotopy fiber sequence of topological spaces.
The higher homotopy groups of \'etale homotopy type of $\Spec k$
in the sense of Artin-Mazur are trivial, and
the above exact sequence may be
understood as a part of a long exact sequence.
In this Section, 
combining the results
of this paper with
the tannakian theory developed in \cite{DTD}
we formulate and prove a motivic counterpart of a homotopy exact sequence 
when $X$ is a smooth curve (Proposition~\ref{MHES}).
The coefficient field of $\DM(k)$ 
and its full subcategories will be $\QQ$,
whereas $K$ will be a coeffcient field of Weil cohomology theory.

Let $C$ be a smooth curve, that is,
a connected one dimensional smooth
scheme
separated of finite type over a perfect field $k$.
Let $j:C\hookrightarrow \overline{C}$ be a smooth compactification
of $C$.
Namely, $\overline{C}$ is a smooth proper curve over $k$, and
$j$ is an open immersion with a dense image.
Let $Z$ denote the complement $\overline{C}- C$, that is
a finite set of
closed points $Z= p_0 \sqcup p_1 \sqcup \ldots \sqcup p_m$.
For simplicity, we assume that $\overline{C}$ admits a $k$-rational point.

We begin by the definition of a symmetric monoidal
full subcategory of $\DM^\otimes(k)$ that is
``smaller'' and more tractable than $\CAlg(\DM^\otimes(k))$.

\begin{Lemma}
\label{subcategory}
Let $A$ be an abelian variety over $k$ and let $l$ be a finite Galois
extension of $k$.
Let $\DM^\otimes(A,l/k)$
be the smallest symmetric monoidal stable
full subcategory of $\DM^\otimes(k)$
which is closed under colimits and contains $M(A)$,
the dual $M(A)^\vee$, $M(\Spec l)$
and Tate objects $\uni_k(n)$ for any $n\in \ZZ$.
(We remark that the symmetric monoidal structure on $\DM^\otimes(A,l/k)$ inherits from that of $\DM^\otimes(k)$, and
$\DM(A,l/k)$ is presentable.)

Let $C$ be a smooth curve over $k$.
Let $k'$ be a Galois field extension of $k$ such that
for any $0\le i \le m$, the residue field $k_i\supset k$
of $p_i$ can be embedded into $k'$.
Let $J_{\overline{C}}$ be the Jacobian variety of $\overline{C}$.
Then $M_C$ lies in $\CAlg(\DM^\otimes(J_{\overline{C}},k'/k))$.
\end{Lemma}

\Proof
Since the underlying object $M_C \in \DM(k)$ is a dual of $M(C)$,
it suffices to prove that $M(C)^\vee$
belongs to $\DM(J_{\overline{C}},k'/k)$.
We note a decomposition $M(J_{\overline{C}})\simeq \oplus_{i=0}^{2g}M_i(J_{\overline{C}})$
for the Jacobian variety $J_{\overline{C}}$ such that
$M_i(J_{\overline{C}})\simeq \Sym^i(M_1(J_{\overline{C}}))$ (see Section~\ref{semiab}).
Here $g$ is the genus of $\overline{C}$.
Also, there is an isomorphism
$M(\overline{C})\simeq \uni_k\oplus M_1(J(\overline{C}))\oplus \uni_k(1)[2]$ in $\DM(k)$. Thus both $M(\overline{C})$ 
and $M(\overline{C})^\vee\simeq M(\overline{C})\otimes \uni_k(-1)[-2]$ lie in $\DM(J_{\overline{C}},k'/k)$.
By Gysin triangle (see \cite[14.5]{MVW}),
there is a distinguished triangle
\[
M(C) \to  M(\overline{C}) \to  M(Z)(1)[2]\to
\]
in the triangulated categories
$\hhh(\DM(k))$.
Therefore, we are reduced to showing that
$M(Z)^\vee\simeq \oplus_{0\le i \le m} M(\Spec k_i)^\vee\simeq \oplus_{0\le i \le m} M(\Spec k_i)$ lies in
$\DM(J_{\overline{C}}, k'/k)$.
Using the functoriality with respect to finite correpondences,
we deduce that each $M(\Spec k_i)$ is a direct summand of $M(\Spec k')$
(since we work with rational coefficients).
\QED

The symmetric monoidal stable presentable $\infty$-category
$\DM^\otimes(A,l/k)$ is a nice property: it is an {\it algebraic fine
tannakian $\infty$-category}.
This notion has been introduced and studied in our work \cite{DTD}.

\begin{Proposition}
We follow the notation in Lemma~\ref{subcategory}.
Let $M_1(A)$ be the direct summand of $M(A)$ in the decomposition in
Section~\ref{semiab}.
Then $M=M_1(A)[-1] \oplus \uni_k(1)\oplus M(\Spec l)$
is a wedge-finite object. Namely, there is an natural number $n$
such that the wedge product $\wedge^{n+1} M$ is zero, and $\wedge^n M$
is an invertible object, see \cite[Section 1]{DTD}.
Consequently, the symmetric monoidal $\infty$-category $\DM^\otimes(A,l/k)$ is an algebraic fine tannakian $\infty$-category,
see \cite[Definition 4.4, Theorem 4.1]{DTD}.
\end{Proposition}

\Proof
By \cite[Proposition 6.1]{DTD} and the fact that $\Hom_{\hhh(\DM(k))}(\uni_k,\uni_k)\simeq \QQ$, it is enough to prove that the wedge product
$\wedge^N M$ is
zero for $N>>0$. To this end, we are reduced to proving that
$\wedge^N (M_1(A)[-1])=0$, $\wedge^N \uni_k(1)=0$, and $\wedge^N M(\Spec l)=0$
for $N>>0$.
By the well-known Kimura finiteness
(see \cite{Kun}, \cite[Thereom7.1.1]{AEH}),
$\wedge^{2e+1} (M_1(A)[-1])\simeq (\Sym^{2e+1}M_1(A))[-2e-1]\simeq 0$
where $e$ is the dimension of $A$.
Also, $\wedge^2 \uni_k(1)=0$
and $\wedge^{d+1}M(\Spec l)=0$. Here $d=[l:k]$.
The final claim follows from the definition of $\DM^\otimes(A,l/k)$
and the definition of algebraic fine tannakian $\infty$-category.
\QED

We define a derived stack from $\DM^\otimes(A,l/k)$ and $M=M_1(A)[-1] \oplus \uni_k(1)\oplus M(\Spec l)$.
By a derived stack over a field $K$, we mean a sheaf $\CAlg_{\QQ}\to \wSSS$
which satisfies a certain geometric condition.
The $\infty$-category $\textup{AlgSt}_{K}$ of derived stacks is defined to be
the full subcategory of $\Fun(\CAlg_{K},\wSSS)$ that consists of derived stacks.
A typical example is a derived affine scheme $\Spec R:\CAlg_{K}\to \wSSS$, that is corepresented by $R\in \CAlg_{K}$. Thus
there is a natural fully faithful embedding $\Aff_{K}\subset \textup{AlgSt}_{K}$.
Another main example for us is a quotient stack $[\Spec R/G]$
that arises from an action of an algebraic affine group scheme 
$G$ on $\Spec R$.
We refer to \cite[Section 2.1]{DTD}
for conventions and terminology concerning derived stacks.

\vspace{2mm}

Applying \cite[Theorem 4.1]{DTD} to $\DM^\otimes(A,l/k)$ with the wedge-finite object
$M$ we obtain

\begin{Corollary}
\label{MGS}
Let $n$ be the natural number such that $\wedge^{n+1}M\simeq 0$ and $\wedge^nM$
is an invertible object.  (Actually, one can see that
$n=2e+d+1$ if $e$ is the dimension of $A$, and $d=[l:k]$.)
There exist a derived stack $\mathcal{X}_{A,l}$ over $\QQ$ such that
$\mathcal{X}_{A,l}$ has a presentation as a quotient stack of the form
$[\Spec V_{A,l}/\GL_n]$ where
$V_{A,l}$ is in $\CAlg_{\QQ}$,
and a symmetric monoidal $\QQ$-linear equivalence
\[
\phi:\QC^\otimes(\mathcal{X}_{A,l})\simeq \DM^\otimes(A,l/k).
\]
Here $\GL_n$ is the general linear group over $\QQ$ that acts on
$V_{A,l}$,
and $\QC^\otimes(\mathcal{X}_{A,l})$ is the symmetric monoidal
$\QQ$-linear presentable $\infty$-category of quasi-coherent complexes
on $\mathcal{X}_{A,l}$.
We shall call $\mathcal{X}_{A,l}$
the {\it motivic Galois stack} associated to $\DM^\otimes(A,l/k)$ and $M$.
For the definition of $\QC(-)$,
we refer to either \cite[Section 2.3]{DTD} or Remark~\ref{MGSR1}.
\end{Corollary}

\begin{Corollary}
\label{trivialcor}
We continue to use the notation in Lemma~\ref{subcategory}.
Then $M_C$ can be naturally regarded as a commutative
object in $\CAlg(\QC^\otimes(\mathcal{X}_{J(\overline{C}),k'}))$.
\end{Corollary}

\Proof
Combine Lemma~\ref{subcategory} and Corollary~\ref{MGS}.
\QED

\begin{Remark}
\label{MGSR1}
For a quotient stack $[\Spec V/G]$ such that $G$ is an algebraic affine group
scheme, $\QC^\otimes([\Spec V/G])$ can be described in the following way.
The action of $G$ on $\Spec V$ can be defined by a simplicial diagram of
derived affine schemes which is informally given by
$[i]\mapsto \Spec V\times G^{\times i}$.
If we put $\Spec R^i=\Spec V\times G^{\times i}$,
then $\QC^\otimes([\Spec V/G])$
is defined to be $\varprojlim_{[i]}\Mod_{R^i}^\otimes$.
The limit of the cosimplicial
diagram $\{\Mod_{R^n}^\otimes\}_{[n]\in \Delta}$
is taken in the $\infty$-category of symmetric monoidal $\infty$-categories.
\end{Remark}

\begin{Remark}
\label{MGSR2}
The stack $\mathcal{X}_{A,l}\simeq [\Spec V_{A,l}/\GL_n]$
is defined as follows (see \cite{DTD} for details):
Let $\Rep^\otimes(\GL_{n})$ be the symmetric monoidal stable $\infty$-category
of representations of $\GL_n$
(cf. Section~\ref{Sullivanmodel}).
There is a canonical equivalence $\QC^\otimes ([\Spec \QQ/\GL_n])\simeq\Rep^\otimes(\GL_{n})$.
Since $[\Spec V_{A,l}/\GL_n]$ is affine over $B\GL_n:=[\Spec \QQ/\GL_n]$,
$\Spec V_{A,l}$ with action of $\GL_n$
can be identified with an object in $\CAlg(\Rep^\otimes(\GL_{n}))$.
By \cite[Theorem 3.1]{DTD}
we have a symmetric monoidal
colimit-preserving functor
$p:\Rep^\otimes(\GL_{n})\to \DM^\otimes(A,l/k)$ which carries
the standard representation of $\GL_n$ placed in degree zero to
$M$. 
By the relative adjoint functor theorem, this functor admits
a lax symmetric monoidal right adjoint
$q:\DM^\otimes(A,l/k)\to \Rep^\otimes(\GL_{n})$.
Thus, $q$ carries a unit object $\uni_{\DM^\otimes(A,l/k)}$
to a commutative algebra object $U_{A,l}:=q(\uni_{\DM^\otimes(A,l/k)})\in \CAlg(\Rep^\otimes(\GL_{n}))$.
This object $U_{A,l}$
amounts to $V_{A,l}$ endowed with action of $\GL_n$,
i.e., data of $[\Spec V_{A,l}/\GL_n]$.
The commutative algebra $V_{A,l}$ in $\CAlg_{\QQ}$ is
the image of $U_{A,l}$ in $\CAlg_{\QQ}$.
We remark that there is a canonical equivalence $\QC^\otimes([\Spec V_{A,l}/\GL_n])\simeq \Mod_{U_{A,l}}^\otimes(\Rep^\otimes(\GL_n))$
where $\Mod_{U_{A,l}}^\otimes(\Rep^\otimes(\GL_n))$
is the symmetric monoidal $\infty$-category of $U_{A,l}$-module objects
in $\Rep^\otimes(\GL_n)$. This equivalence
makes the diagram
\[
\xymatrix{
\QC^\otimes (B\GL_n) \ar[r] \ar[d]_{\simeq} & \QC^\otimes([\Spec V_{A,l}/\GL_n])\ar[d]^{\simeq} \ar[rd]^{\phi} &  \\
\Rep^\otimes(\GL_n) \ar[r]_{\otimes U_{A,l}} & \Mod_{U_{A,l}}^\otimes(\Rep^\otimes(\GL_n)) \ar[r] & \DM^\otimes(A,l/k).
}
\]
commute up to homotopy, where the top horizontal arrow is the
pullback functor of the projection $[\Spec V_{A,l}/\GL_n]\to B\GL_n$.
The equivalence
$\Mod_{U_{A,l}}^\otimes(\Rep^\otimes(\GL_n)) \to \DM^\otimes(A,l/k)$
is defined to be the composite
\[
\Mod_{U_{A,l}}^\otimes(\Rep^\otimes(\GL_n))\to \Mod_{p(U_{A,l})}^\otimes(\DM^\otimes(A,l/k))\to \Mod_{\uni_{\DM^\otimes(A,l/k)}}^\otimes(\DM^\otimes(A,l/k))\simeq \DM^\otimes(A,l/k)
\]
where the first functor is induced  by $p$, and the second functor is induced by the base change along the counit map $p(U_{A,l})=pq(\uni_{\DM^\otimes(A,l/k)})\to \uni_{\DM^\otimes(A,l/k)}$.
The composite of lower horizontal arrows is equivalent to $p$.
\end{Remark}

\begin{Remark}
There is the following uniqueness.
Let $(\mathcal{Y},N)$ be a pair that consists of
a derived stack $\mathcal{Y}$ over $\QQ$, and $N$ is a vector bundle
on $\mathcal{Y}$.
Here by a vector bundle we mean an object $N$
in $\QC(\mathcal{Y})$
such that for any $f:\Spec R\to \mathcal{Y}$,
the restriction $f^*(N)$ is equivalent to
a direct summand of
some finite coproduct $R^{\oplus m}$.
The stack $\mathcal{X}_{A,l}\simeq [\Spec V_{A,l}/\GL_n]$ has a vector bundle $N_{A,l}$
that is defined to be the pullback of the tautological vector bundle
on $BGL_n=[\Spec \QQ/\GL_n]$.
So we have such a pair $(\mathcal{X}_{A,l},N_{A,l})$.
By the diagram in Remark~\ref{MGSR2}, the equivalence
$\phi:\QC^\otimes(\mathcal{X}_{A,l/k})\simeq \DM^\otimes(A,l/k)$
sends $N_{A,l}$ to $M$.
Assume that there
is a symmetric monoidal $\QQ$-linear equivalence
$\QC^\otimes(\mathcal{Y})\simeq \DM^\otimes(A,l/k)$
which sends $N$ to $M$.
Then there is an equivalence
$\mathcal{Y}\simeq \mathcal{X}_{A,l}$ such that
the induced equivalence
$\QC^\otimes(\mathcal{Y})\simeq \QC^{\otimes}(\mathcal{X}_{A,l})$ sends $N$ to $N_{A,l}$.
This uniqueness will not be necessary in this paper,
so that we will not present the proof.
But one can prove it by using arguments in \cite{DTD}.
\end{Remark}

We say that a morphism $\mathcal{X}\to \mathcal{Y}$ of derived stacks over $K$
is affine
if, for any $\Spec R \to \mathcal{Y}$ from a derived affine scheme,
the fiber product $\mathcal{X}\times_{\mathcal{Y}}\Spec R$ belongs to
$\Aff_{K}$.
Let $\Aff_{\mathcal{Y}}$ be the full subcategory of the overcategory $(\textup{AlgSt}_{K})_{/\mathcal{Y}}$ that consists of affine morphisms
$\mathcal{X}\to \mathcal{Y}$.
There is a canonical equivalence $\Aff_{\mathcal{Y}}\simeq \CAlg(\QC^\otimes(\mathcal{Y}))^{op}$
(cf. \cite[Section 2.3]{DTD}, this is a direct generalization of the analogous
fact in the usual scheme theory).

\begin{Definition}
\label{MRHTMGS}
By Corollary~\ref{trivialcor}, let us consider
$M_C$ as an object
in $\CAlg(\QC^\otimes(\mathcal{X}_{J(\overline{C}),k'}))$.
Let $\mathcal{M}_{C}\to \mathcal{X}_{J(\overline{C}),k'}$
be a derived stack affine over $\mathcal{X}_{J(\overline{C}),k'}$
that corresponds to $M_C$ through
the equivalence $\Aff_{\mathcal{X}_{J(\overline{C}),k'}}\simeq \CAlg(\QC^\otimes(\mathcal{X}_{J(\overline{C}),k'}))^{op}$.
\end{Definition}

Let $\mathsf{R}_E:\DM^\otimes(k)\to \mathsf{D}^\otimes(\KK)\simeq \Mod_{\KK}^\otimes$
be the realization functor associated to a mixed Weil Theory $E$
with coefficients in a field $\KK$ of characteristic zero.
By abuse of notation we write $\mathsf{R}_E$ also for the restriction
$\DM^\otimes(A,l/k) \to \mathsf{D}^\otimes(\KK)$.
Suppose that $\mathsf{R}_E(M)$ is concentrated in degree zero $\mathsf{D}(\KK)$
(all known mixed Weil theories satisfy this condition).
As discussed in \cite[Section 4.1]{PM} or \cite[Remark 6.13]{DTD}, it gives rise to a morphism
\[
\rho_E:\Spec \KK\to \mathcal{X}_{A,l}.
\] 
We refer to this morphism as the base point of $\mathsf{R}_E$.

We briefly recall the construction of $\rho_E$.
Let $p:\QC^\otimes(B\GL_n)\simeq  \textup{Rep}^\otimes(\GL_n)\to \DM^\otimes(A,l/k)$ be the sequence contained in the diagram in
Remark~\ref{MGSR2}.
Note that this functor carries the standard representation of $\GL_n$
placed in degree zero to $M$, and the realization functor carries
$M$ to the $n$-dimensional vector space placed in degree zero in $\mathsf{D}(\KK)$. Therefore, by the universal property of $\textup{Rep}^\otimes(\GL_n)$
\cite[Theorem 3.1]{DTD}, the composite 
$\QC^\otimes(B\GL_n)\simeq  \textup{Rep}^\otimes(\GL_n)\to \DM^\otimes(A,l/k)\to \mathsf{D}^\otimes(\KK)$
is equivalent to the pullback functor 
$\QC^\otimes(B\GL_n)\to \mathsf{D}^\otimes(\KK)\simeq \QC^\otimes(\Spec \KK)$
along $\Spec \KK\to \Spec \QQ\to B\GL_n$.
Let $u:\mathsf{D}^\otimes(\KK)\to \Rep^\otimes(\GL_n) \simeq \QC^\otimes(B\GL_n)$
be the lax symmetric monoidal right adjoint to
$\QC^\otimes(B\GL_n)\to \QC^\otimes(\Spec \KK)\simeq \mathsf{D}^\otimes(\KK)$,
whose existence is ensured by the relative adjoint functor theorem.
Then this right adjoint induces $\CAlg(\mathsf{D}^\otimes(\KK))\simeq \CAlg_{\KK}\to \CAlg(\Rep^\otimes(\GL_n))$ which carries the unit algebra $\KK$
to $u(\KK)\simeq \Gamma(\GL_n)\otimes_{\QQ}\KK \in \CAlg(\Rep^\otimes(\GL_n))\simeq \CAlg(\QC^\otimes(B\GL_n))$. Here, write $\Gamma(\GL_n)$ for
the (ordinary) coordinate ring of the general linear group $\GL_n$
which is endowed with the natural action of $\GL_n$. The symbol $\KK$ in $\Gamma(\GL_n)\otimes_{\QQ}\KK$ is understood as the $\QQ$-algebra $\KK$ with the trivial
action of $\GL_n$.
Note that there is a natural morphism $U_{A,l}\to u(\KK)\simeq \Gamma(\GL_n)\otimes_{\QQ}\KK$ in $\CAlg(\QC^\otimes(B\GL_n))$. In fact, if $v:\CAlg_{K}\to \CAlg(\DM^\otimes(A,l/k))$
denotes the right adjoint to
the restricted multiplicative realization functor
$\CAlg(\DM^\otimes(A,l/k))\to\CAlg_{\KK}$, then
there is a unit map $\uni_{\DM(A,l/k)}\to v(\KK)$
that induces $U_{A,l}=q(\uni_{\DM(A,l/k)})\to qv(\KK)=u(\KK)$, as claimed (for the functor $q$, see Remark~\ref{MGSR2}).
By using the equivalence $\Aff_{B\GL_n}\simeq \CAlg(\QC^\otimes(B\GL_n))^{op}$,
we obtain
$\rho_E:\Spec K\simeq [\Spec \Gamma(\GL_n)\otimes_{\QQ}K/\GL_n]\to \mathcal{X}_{A,l}=[\Spec V_{A,l}/\GL_n]$.

\begin{Remark}
\label{geomreal}
By this construction and Remark~\ref{MGSR2}, we see that the diagram
\[
\xymatrix{
\QC^\otimes (\mathcal{X}_{A,l}) \ar[r]^{\rho^*_E} \ar[d]_{\simeq}^\phi & \QC^{\otimes}(\Spec \KK) \ar[d]^\simeq \\
\DM^\otimes(A,l/k) \ar[r]^{\mathsf{R}_E} & \mathsf{D}^\otimes(\KK)
}
\]
commutes up to homotopy,
where $\rho^*_E$ is the pullback functor (cf. \cite[Section 2.3]{DTD}),
the right vertical arrow is a canonical equivalence.
\end{Remark}

One can associate to
the base point $\rho_E:\Spec \KK\to \mathcal{X}_{A,l}$
a derived affine group scheme over $\KK$.
Namely, we take the Cech nerve $G:\NNNN(\Delta_+)^{op}\to \textup{AlgSt}_\KK$
of $\rho_E\times\textup{id}:\Spec \KK\to \mathcal{X}_{A,l}\times_{\Spec \QQ}\Spec \KK$, which is defined to be the right Kan extension
$\NNNN(\Delta_+)^{op}\to \textup{AlgSt}_{\KK}$ of $\NNNN(\Delta_+^{\le0})^{op}=\NNNN(\{[-1]\to [0]\})^{op}\to \textup{AlgSt}_{\KK}$ determined by $\rho_E\times\textup{id}$.
The evaluation $G([1])$ is equivalent to 
$\Spec \KK\times_{\mathcal{X}_{A,l}\times \Spec \KK}\Spec \KK$ which is affine because
the diagonal $[\Spec V_{A,l}/\GL_n]\to [\Spec V_{A,l}/\GL_n]\times [\Spec V_{A,l}/\GL_n]$ is affine.
Thus the restriction of $G$ defines a group object $\NNNN(\Delta)^{op}\to \Aff_K$,
whose underlying derived affine scheme is $\Spec \KK\times_{\mathcal{X}_{A,l}\times \Spec \KK}\Spec \KK$. We write $\Omega_{\rho_E}\mathcal{X}_{A,l}$ for this
derived affine group scheme over $\KK$.
The derived group scheme $\Omega_{\rho_E}\mathcal{X}_{A,l}$
is related to the derived motivic Galois
group:

\begin{Proposition}
\label{loopgp}
Let $\mathsf{MG}_{E,\DM^\otimes(A,l/k)}$ be the derived motivic Galois
group which represents the automorphism group functor $\Aut(\mathsf{R}_E|_{\DM^\otimes(A,l/k)}):\CAlg_{\KK}\to \Grp(\wSSS)$, cf. Remark~\ref{motsubcategory}.
Then $\Omega_{\rho_E}\mathcal{X}_{A,l}$ is naturally equivalent to
$\mathsf{MG}_{E,\DM^\otimes(A,l/k)}$.
\end{Proposition}

\Proof
By Remark~\ref{geomreal}, we have
$\Aut(\mathsf{R}_E|_{\DM^\otimes(A,l/k)})\simeq \Aut(\rho_E^*)$
where 
$\rho^*_E:\QC^\otimes (\mathcal{X}_{A,l}) \to \QC^{\otimes}(\Spec \KK)$.
It will suffice to show that $\Omega_{\rho_E}\mathcal{X}_{A,l}\simeq \Aut(\rho_E^*)$.
This equivalence follows from \cite[Proposition 4.6]{Bar}.
\QED

\begin{Remark}
By the representability of automorphism groups,
the restriction to $\DM^\otimes(A,l/k)$
induces $\MG_E\to \mathsf{MG}_{E,\DM^\otimes(A,l/k)}\simeq \Omega_{\rho_E}\mathcal{X}_{A,l}$.
The action of $\MG_E$ on $\mathsf{R}_E(M_C)$ described in Proposition~\ref{goodaction} factors
through
$\MG_E\to \Omega_{\rho_E}\mathcal{X}_{A,l}$.
\end{Remark}

\vspace{2mm}

Now we are ready to prove 
the following motivic generalization of homotopy exact sequence.

\begin{Proposition}
\label{MHES}
Let $\mathcal{M}_C\to \mathcal{X}_{J(\overline{C}),k'/k}$ be
the affine morphism defined in Definition~\ref{MRHTMGS}.
Let us consider the pullback diagram of derived stacks
\[
\xymatrix{
\mathcal{F}_E \ar[r] \ar[d] &  \mathcal{M}_C \ar[d] \\
\Spec \KK \ar[r]^{\rho_E} & \mathcal{X}_{J(\overline{C}),k'}
}
\]
in $\textup{AlgSt}_{\QQ}$.
(One may think of this diagram as a Cartesian diagram in $\Fun(\CAlg_{\QQ},\wSSS)$.)
Then 
the fiber $\mathcal{F}_E$ is naturally equivalent to $\Spec \mathsf{R}_E(M_C)$,
where $\mathsf{R}_E(M_C)$ in $\CAlg_{\KK}$ is the image of $M_C$
under the multiplicative realization functor $\mathsf{R}_E:\CAlg(\DM^\otimes(k))\to \CAlg_{\KK}$.
In particular, when $E$ is the singular cohomology theory,
by Theorem~\ref{real1} we have
a Cartesian diagram
\[
\xymatrix{
\Spec A_{PL,\infty}(C^t) \ar[r] \ar[d] &  \mathcal{M}_C \ar[d] \\
\Spec \QQ \ar[r]^{\rho_E} & \mathcal{X}_{J(\overline{C}),k'}.
}
\]
\end{Proposition}

\begin{Remark}
The morphism 
$\mathcal{M}_C \to \mathcal{X}_{J(\overline{C}),k'}$
should be thought of as a motivic counterpart of the delooping of
$\pi_1^{\textup{\'et}}(X,x) \to G_k$.
By Proposition~\ref{loopgp}
we can obtain the derived motivic Galois group $\MG_{E,\DM^\otimes(J(\overline{C}),k'/k)}\simeq \Omega_{\rho_E}\mathcal{X}_{J(\overline{C}),k'}$ from the base stack $\Spec \KK\to \mathcal{X}_{J(\overline{C}),k'}$
by using the construction of the base loop space.
The fiber $\mathcal{F}_E$ shoud be understood as a role of
the delooping of $\pi_1^{\textup{\'et}}(X_{\bar{k}},\bar{x})$.
Consider the situation that $k$ is a subfield of $\CC$.
Then
$\pi_1^{\textup{\'et}}(X_{\bar{k}},\bar{x})$ is isomorphic to
the profinite completion of the (topological)
fundamental group $\pi_1(X^t,\bar{x})$
of the underlying topological space $X^t$
of $X\times_{\Spec k}\Spec \CC$.
On the other hand, if we fix a $k$-rational point $c$,
the unipotent group scheme
$\overline{G}^{(1)}(C,c)\simeq \Spec H^0(\QQ\otimes_{A_{PL}(C^t)}\QQ)$ is
the pro-unipotent completion of the topological fundamental group
$\pi_1(C^t,c)$.
\end{Remark}

\Proof
We have already done almost things.
By Remark~\ref{geomreal}, one can identify the multiplicative realization functor
$\CAlg(\DM^\otimes(J(\overline{C}),k'/k)) \to \CAlg_{\KK}$
with 
$\CAlg(\QC^\otimes(\mathcal{X}_{J(\overline{C}),k'}))\to \CAlg_{\KK}\simeq \CAlg(\QC^\otimes(\Spec \KK))$ induced by the pullback functor $\rho_E^*$.
Then we use the observation that the canonical equivalences
\[
\CAlg(\QC^\otimes(\mathcal{X}_{J(\overline{C}),k'}))^{op}\simeq \Aff_{\mathcal{X}_{J(\overline{C}),k'}}\ \  \textup{and}\ \  \CAlg(\QC^\otimes(\Spec \KK))^{op}\simeq \Aff_{\Spec \KK}
\]
are compatible
with pullback functors. Namely, through these canonical equivalences,
the opposite functor $\CAlg(\QC^\otimes(\mathcal{X}_{J(\overline{C}),k'}))^{op}\to \CAlg(\QC^\otimes(\Spec \KK))^{op}$ can be identified with
$\Aff_{\mathcal{X}_{J(\overline{C}),k'}}\to \Aff_{\Spec \KK}=\Aff_{\KK}$
given by $\{\mathcal{Z}\to \mathcal{X}_{J(\overline{C}),k'}\} \mapsto \{\textup{pr}_2:\mathcal{Z}\times_{\mathcal{X}_{J(\overline{C}),k'}}\Spec \KK\to \Spec \KK\}$.
Therefore, we see that
$\mathcal{F}_E$ is equivalent to $\Spec \mathsf{R}_E(M_C)$
via these identifications.
\QED

\appendix

\renewcommand{\theTheorem}{A.\arabic{Theorem}}

\section{Comparison results}

We will compare the motivic algebra of path torsors with an approach by 
Deligne-Goncharov \cite{DG}.

\subsection{}
Suppose that $k$ is a number field. We work with rational coefficients.
We begin by reviewing the category of mixed Tate motives
over $k$.
Let $\DTM:=\DTM(k)$ be the smallest stable subcategory of $\DM(k)$
that is closed under small colimits
and consists of $\uni_k(n)$ for any $n\in \ZZ$.
The stable subcategory $\DTM$ inherits a symmetric monoidal structure
from $\DM(k)$.
We refer to it as the symmetric monoidal stable $\infty$-category of mixed
Tate motives and denote it by $\DTM^\otimes$.
The stable $\infty$-category $\DTM$ is compactly generated.
Let $\DTM_\vee$ denote the stable subcategory
spanned by compact objects.
In particular, $\Ind(\DTM_\vee)\simeq \DTM$ where
$\Ind(-)$ indicates the Ind-category.
 The full subcategory $\DTM_\vee$ coincides with
the stable subcategory consisting of dualizable objects.
Let $(\mathsf{D}(\QQ)_{\ge0},\mathsf{D}(\QQ)_{\le0})$ be
the standard $t$-structure on $\mathsf{D}(\QQ)$
such that $C$ belongs to $\mathsf{D}(\QQ)_{\ge0}$ (resp. $\mathsf{D}(\QQ)_{\le0}$) if and only if
$H^{-i}(C)=H_i(C)=0$ for $i<0$ (resp. $i>0$).
For our conventions on (motivic) $t$-structures, we refer to \cite{HA} and \cite[Section 7]{Bar}.
Under the setting where $k$ is a number field,
there is a nondegenerate bounded $t$-structure on $\DTM_\vee$ given by
\[
\DTM_{\vee,\ge0}:=\mathsf{R}_T^{-1}(\mathsf{D}(\QQ)_{\ge0})\cap \DTM_\vee,\ \  \DTM_{\vee,\le0}:=\mathsf{R}_T^{-1}(\mathsf{D}(\QQ)_{\le0})\cap \DTM_\vee
\]
where $\mathsf{R}_T:\DTM^\otimes\to \mathsf{D}^\otimes(\QQ)$ is the singular realization functor. We call it the motivic $t$-structure on $\DTM_\vee$.
The realization functor  $\DTM_\vee\to \mathsf{D}(\QQ)$ is $t$-exact
and conservative. The both
categories $\DTM_{\vee,\ge0}$ and $\DTM_{\vee,\le0}$
are closed under tensor products.
Let $\TM^\otimes$ be the heart $\DTM_{\vee,\ge0}\cap \DTM_{\vee,\le0}$
which is a symmetric monoidal (furthermore tannakian) abelian category.
We refer to $\TM^\otimes$ as the abelian category of mixed Tate motives.

\subsection{}
The construction in Deligne-Goncharov \cite{DG} employs the idea in
Wojtkowiak \cite{W} that uses cosimplicial schemes.
Let $X$ be a smooth variety over $k$.
Let $x:\Spec k\to X$ and $y:\Spec k\to X$ be two $k$-rational points.
To $(X,x,y)$ we associate a cosimplicial smooth scheme, i.e.,
a functor $P^{\Delta}(X,x,y):\Delta\to \Sm_k:[n]\mapsto X^n$
whose cofaces are defined by
\begin{eqnarray*}
d^0(x_1,\ldots,x_n)=(x_1,\ldots,x_n,x),\ \ d^{n+1}(x_1,\ldots,x_n)=(y,x_1,\ldots,x_n), \\
d^i(x_1,\ldots,x_n)=(x_1,\ldots,x_{n-i+1},x_{n-i+1},\ldots,x_n),\ \ (0<i<n),
\end{eqnarray*}
$d^0,d^1:X^0=\Spec k\rightrightarrows X^1=X$ is given by $x$ and $y$.
The codegeneracy are given by projections.
Recall the functor $\Xi:\Sm_k^{op}\to \CAlg(\DM^\otimes(k))$
from Section~\ref{cohomologicalmot}.
By abuse of notation we write $\Xi$ for the composite
$\Sm_k^{op}\stackrel{\Xi}{\to} \CAlg(\DM^\otimes(k))\to \DM(k)$
where the second functor is the forgetful functor.
Consider the simplicial object in $\DM(k)$ given by the composition
\[
\mathcal{M}_{\Delta}(X,x,y):\NNNN(\Delta)^{op}\stackrel{P^{\Delta}(X,x,y)^{op}}{\longrightarrow} \Sm_k^{op}\stackrel{\Xi}{\to} \DM(k).
\]
Let $\Delta_{s}$ be the subcategory of $\Delta$ whose
objects coincide with that of $\Delta$, and whose morphisms
are injective maps.
The inclusion $\NNNN(\Delta_s)^{op}\hookrightarrow \NNNN(\Delta)^{op}$
is cofinal \cite[6.5.3.7]{HTT}.
It follows that a colimit of $\mathcal{M}_{\Delta}(X,x,y)$ is naturally equivalent to that of the restriction $\mathcal{M}_{\Delta}(X,x,y)|_{\NNNN(\Delta_s)^{op}}:\NNNN(\Delta_s)^{op}\to \DM(k)$.
Let $\Delta_{s,\le n}$ be the full subcategory of $\Delta_{s}$ spanned by
$\{[0],\ldots,[n]\}$ and let $\mathcal{M}_{\Delta_{s,\le n}}(X,x,y):\NNNN(\Delta_{s,\le n})^{op}\to \DM(k)$
denote the restriction of $\mathcal{M}_{\Delta}(X,x,y)$.
Let $\mathcal{M}(X,x,y)$  denote a colimit of $\mathcal{M}_{\Delta}(X,x,y)|_{\NNNN(\Delta_s)^{op}}$ (or equivalently $\mathcal{M}_{\Delta}(X,x,y)$).
Let $\mathcal{M}_n(X,x,y)$
denote a colimit of $\mathcal{M}_{\Delta_{s,\le n}}(X,x,y)$ in $\DM(k)$.
The colimits $\mathcal{M}_n(X,x,y)$ naturally constitute a sequence $\mathcal{M}_0(X,x,y)\to \mathcal{M}_1(X,x,y)\to \cdots$, and there is a canonical equivalence $\varinjlim_{n}\mathcal{M}_n(X,x,y)\simeq \mathcal{M}(X,x,y)$
(cf. \cite[4.2.3]{HTT}).
Now suppose that $M(X)$ belongs to $\DTM_\vee$. Then $M_{X^r}\simeq (M(X)^{\otimes r})^{\vee}$ lies in $\DTM_\vee$. Consequently, the finite colimit
$\mathcal{M}_n(X,x,y)$ belongs to $\DTM_\vee$.
Take the $0$-th cohomology
$H^0(\mathcal{M}_n(X,x,y))$ with respect to motivic $t$-structure.
We let
\[
M_{DG}(X,x,y):=\varinjlim_n H^0(\mathcal{M}_n(X,x,y))
\]
be the filtered colimit
in $\Ind(\TM)$. We refer to it as the Deligne-Goncharov motive
associated to $(X,x,y)$.
According to \cite[7.4]{Bar}, $\DTM\simeq \Ind(\DTM_\vee)$ has
a $t$-structure defined by $(\Ind(\DTM_{\vee,\ge0}),\Ind(\DTM_{\vee,\le0}))$. 
Passing to the $0$-th cohomology (with respect to
$t$-structure) commutes with filtered colimits so
that $M_{DG}(X,x,y)=\varinjlim_n H^0(\mathcal{M}_n(X,x,y))\simeq H^0(\mathcal{M}(X,x,y))$.
Therefore
$M_{DG}(X,x,y)$
is nothing else but the $0$-th cohomology
of a colimit of the simplicial diagram $\mathcal{M}_{\Delta}(X,x,y)$.

\begin{Remark}
Taking advantage of a functorial assignment $X\mapsto M_X$ (see Proposition~\ref{cohmotalg}), we here give the
cohomological construction of $M_{DG}(X,x,y)$
while the homological one is described
in \cite[3.12]{DG}.
Thus, procedures are dual to one another.
In {\it loc. cit.},
one considers the diagram $\NNNN(\Delta_{s,\le n}) \to \DM(k):[r]\mapsto M(X^{r})$ induced by the restricted 
diagram $P^{\Delta_{s,\le n}}(X,x,y):\NNNN(\Delta_{s,\le n}) \to \Sm_k:[r]\mapsto X^{r}$ instead of $\mathcal{M}_{\Delta_{s,\le n}}(X,x,y)$ (see \cite[3.12]{DG}). Then take a finite limit of the diagram in $\DM(k)$
 by means of
Moore complexes.
The pleasant feature of cohomological construction is that
it is not necessary to take the family of the restricted diagrams
(though we take trouble to take them):
one can directly define it to be the $0$-th cohomology of a colimit
of the simplicial
diagram $\mathcal{M}_{\Delta}(X,x,y)$.
\end{Remark}

\begin{Remark}
One can consider a larger subcategory that consists of
Artin-Tate motives. This category contains
not only Tate motives but also motives of the form
$M(\Spec k')$ such that $k'$ is a finite separable extension field of $k$.
We can treat this category 
by using a main result of 
\cite{FI} and \cite[Section 8]{Bar}. But we will not
pursue a generalizaton to this direction.
\end{Remark}

\subsection{}
We will think of $\TM^\otimes$ as a neutral tannakian category over $\QQ$,
which is endowed with the (symmetric monoidal) singular
realization functor to the category of vector spaces over $\QQ$
\[
R_T:\TM^\otimes\to \Vect^\otimes_\QQ.
\]
The Tannaka dual $MTG$ with respect to this functor
is a pro-algebraic group over $\QQ$ which represents the automorphism group
of this symmetric monoidal functor $R_T$.
For any $M\in \TM$
$MTG\simeq \Aut(R_T)$ naturally acts on $R_T(M)$. It gives rise to
a $\QQ$-linear symmetric monoidal equivalence $\TM^\otimes\simeq \Rep^\otimes(MTG)_\vee$
where $\Rep^\otimes(MTG)_\vee$ is the symmetric monoidal abelian category
of finite dimensional representations of $MTG$.
Recall from \cite{Bar} the relation of tannakization
and $MTG$.

\begin{Proposition}[cf. Theorem 7.16 in \cite{Bar}]
Let $\mathsf{MTG}$ be the derived affine group scheme which represents the automorphism group of $\mathsf{R}_T:\DTM^\otimes\to \mathsf{D}^\otimes(\QQ)$, that is, the tannakziation of $\mathsf{R}_T:\DTM_\vee^\otimes\to \mathsf{D}^\otimes(\QQ)$ in the sense of \cite{Tan}.
Then there is a natural isomorphism between
$MTG$ and
the underlying group scheme of $\mathsf{MTG}$.
\end{Proposition}

\begin{Remark}
There are approaches to $\mathsf{MTG}$ by means of
bar constructions, see Spitzweck's
derived tannakian presentation of $\DTM^\otimes$ \cite{Spi}, 
(see also \cite{Bar}, \cite{DTD}).
If we suppose futhermore that $k$ is a number field, then
by Borel's computation of rational motivic cohomology groups
of number fields, it is not difficult to prove that
$\mathsf{MTG}\simeq MTG$.
\end{Remark}

Let $X$ be a smooth variety and assume that $M(X)$ belongs to $\DTM_\vee$
(thus $M_X$ also lies in $\DTM_\vee$).
Let $x,y:\Spec k\rightrightarrows X$ be two $k$-rational points
on $X$.
Recall the motivic algebra of path torsors
\[
P_X(x,y)=\uni_k\otimes_{M_X}\uni_k
\]
in $\CAlg(\DTM^\otimes)\subset \CAlg(\DM^\otimes(k))$ from Example~\ref{motpath}.
Take the cohomology $H^0(\uni_k\otimes_{M_X}\uni_k)$
with respect to the $t$-structure $(\Ind(\DTM_{\vee,\ge0}),\Ind(\DTM_{\vee,\le0}))$.

\begin{Proposition}
\label{algstr}
The cohomology
$H^0(P_X(x,y))$
inherits the structure of commutative algebra object
in $\Ind(\TM)$ from $P_X(x,y)$.
(The construction is described in the proof below.)
\end{Proposition}

\Proof
Note first that $M_X$ is the dual of $M(X)$ in $\DTM$,
and $M(X)$ belongs to $\DTM_{\vee,\ge0}$.
Since $\mathsf{R}_T(M_X)$ is the dual
of $\mathsf{R}_T(M(X))\in \mathsf{D}(\QQ)_{\ge0}$, 
thus $M_X$ lies in $\DTM_{\vee,\le0}$.
Remember that
$\mathsf{R}_T:\CAlg(\DTM^\otimes)\to \CAlg(\mathsf{D}^\otimes(\QQ))$
is a left adjoint (in particular, it preserves colimits).
It follows that
$\mathsf{R}_T(\uni_k\otimes_{M_X}\uni_k)\simeq \QQ\otimes_{T_X}\QQ$.
The pushout $\QQ\otimes_{T_X}\QQ$ lies in $\mathsf{D}(\QQ)_{\le0}$
(for example, compute it by the standard bar construction).

\vspace{1mm}

Now we recall the left completion of $\DTM$ with respect to
$(\Ind(\DTM_{\vee,\ge0}),\Ind(\DTM_{\vee,\le0}))$.
In a nutshell, the left completion of $\DTM$ is 
a symmetric monoidal $t$-exact colimit-preserving functor $\DTM^\otimes\to \oDTM^\otimes$
to the ``left completed'' stable presentable symmetric monoidal $\infty$-category $\oDTM^\otimes$ 
(we refer the reader to \cite[Section 7.2]{Bar} and references therein
for the notions of left completeness and left completion).
The $\infty$-category $\oDTM$ can be described the limit of the diagram indexed by $\ZZ$
\[
\cdots \to \DTM_{\le n+1} \stackrel{\tau_{\le n}}{\to} \DTM_{\le n} \stackrel{\tau_{\le n-1}}{\to} \DTM_{\le n-1} \stackrel{\tau_{\le n-2}}{\to} \cdots
\]
of $\infty$-categories, where $\tau_{\le n}$ are the truncation
functors (we use the homological indexing following \cite{HA}).
According to \cite[3.3.3]{HTT}
the $\infty$-category $\oDTM$ can be identified with
the full subcategory of $\Fun(\NNNN(\ZZ),\DTM)$
spanned by functors $\phi:\NNNN(\ZZ)\to \DTM$ such that
\begin{itemize}
\item for any $n\in \ZZ$, $\phi([n])$ belongs to $\DTM_{\le -n}$,

\item for any $m\le n\in \ZZ$, the associated map $\phi([m])\to \phi([n])$
gives an equivalence $\tau_{\le -n}\phi([m])\to \phi([n])$.
\end{itemize}
Let $\oDTM_{\ge0}$ (resp. $\oDTM_{\le0}$) be the full subcategory of $\oDTM$
spanned by $\phi:\NNNN(\ZZ)\to \DTM$ such that $\phi([n])$ belongs to
$\DTM_{\ge0}$ (resp. $\DTM_{\le0}$) for each $n\in \ZZ$.
The functor $\DTM\to \oDTM$ induces an equivalence $\DTM_{\le0}\to \oDTM_{\le0}$.
The pair $(\oDTM_{\ge0}, \oDTM_{\le0})$ is an accessible, left complete and right complete $t$-structure of $\oDTM$. The functor
$\DTM \to \oDTM$ carries $M$ to $\{\tau_{\le r}M\}_{r\in \ZZ}$.
Since the $t$-structure on $\mathsf{D}(\QQ)$ is left complete,
thus the realization functor $\DTM^\otimes \to \mathsf{D}^\otimes(\QQ)$
factors as $\DTM^\otimes \to \oDTM^\otimes \stackrel{\overline{\mathsf{R}}_T}{\to} \mathsf{D}^\otimes(\QQ)$
such that $\overline{\mathsf{R}}_T:\oDTM^\otimes \to \mathsf{D}^\otimes(\QQ)$ is conservative by
\cite[Corollary 7.3]{Bar}.

\vspace{1mm}

Return to the proof. Since $\DTM^\otimes\to \oDTM^\otimes$
is $t$-exact, we may and will work with $\oDTM$
instead of $\DTM$. By abuse of notation, we write
$\uni_k\otimes_{M_X}\uni_k$ for the image in $\oDTM$.
It follows from the conservativity of $\overline{\mathsf{R}}_T$
that $\uni_k\otimes_{M_X}\uni_k$ belongs to
$\oDTM_{\le0}$. Consider the adjunction $\oDTM_{\ge0}\rightleftarrows \oDTM:\tau_{\ge0}$ where the left adjoint is the symmetric monoidal fully faithful functor.
Thus the right adjoint $\tau_{\ge0}:\oDTM\to \oDTM_{\le 0}$ is lax symmetric monoidal.
For any $M\in \CAlg(\oDTM)$, $\tau_{\ge0}(M)$ is a commutative algebra
object in $\oDTM_{\ge0}^\otimes$.
Consequently, $H^0(\uni_k\otimes_{M_X}\uni_k)=\tau_{\ge0}(\uni_k\otimes_{M_X}\uni_k)$ inherits a commutative algebra structure
\begin{eqnarray*}
H^0(\uni_k\otimes_{M_X}\uni_k)\otimes H^0(\uni_k\otimes_{M_X}\uni_k)\to H^0(\uni_k\otimes_{M_X}\uni_k),\ \ 
H^0(1_k)\to H^0(\uni_k\otimes_{M_X}\uni_k)
\end{eqnarray*}
in $\Ind(\TM)$.
\QED

We put $M(X,x,y):=H^0(\uni_k\otimes_{M_X}\uni_k)$.
By Proposition~\ref{algstr}
we regard it as a commutative algebra in $\Ind(\TM)\simeq \Rep(MTG)$.

\begin{Remark}
We can think of $M(X,x,y)$ also as a commutative $\QQ$-algebra
$H^0(\QQ\otimes_{T_X}\QQ)$ with the canonical
action of $MTG\simeq \Aut(R_T)$.
This action of $MTG$ on $H^0(\QQ\otimes_{T_X}\QQ)$ can be identified with the action in Section~\ref{MGA}, Theorem~\ref{actiononhomotopy}.
As discussed in Section~\ref{actionstep1}, Section~\ref{finalaction},
$\mathsf{MTG}\simeq \Aut(\mathsf{R}_T)$ acts on $\QQ\otimes_{T_X}\QQ\simeq \mathsf{R}_T(\uni_k\otimes_{M_X}\uni_k)$.
It gives rise to an action of the underlying group scheme $MTG$
on $H^0(\QQ\otimes_{T_X}\QQ)$
(but we treated only the case $x=y$).
By \cite[Theorem 7.16]{Bar} and its proof, there is a canonical equivalence
$\Aut(\mathsf{R}_T)\simeq \Aut(\overline{\mathsf{R}}_T)$
as functors $\CAlg_{\QQ}^{\DIS}\to \Grp(\SSS)$ (note that the domain is not $\CAlg_{\QQ}$ but $\CAlg_{\QQ}^{\DIS}$).
In addition, by \cite[Proposition 7.13, 7.12]{Bar}
$(\oDTM^\otimes,\oDTM_{\ge0},\oDTM_{\le0})$ is a locally dimensional $\infty$-category in the sense of Lurie \cite[VIII, Section 5]{DAG}. Therefore, the heart is the tannakian category
$\Rep^\otimes(MTG)$ of (not necessarily finite dimensional)
representations of $MTG$, and the natural morphism
$\mathsf{MTG} \to MTG$ in $\Fun(\CAlg_{\QQ}^{\DIS},\Grp(\SSS))$
can naturally be
identified with $\Aut(\mathsf{R}_T)\to \Aut(R_T)$ induced by the restriction
of natural equivalences to the heart.
Let $L$ be the function field of $MTG$.
Taking account of Theorem~\ref{actiononhomotopy} (2),
the action of the group of $L$-valued point $MTG(L)$ on $H^0(\QQ\otimes_{T_X}\QQ)\otimes_{\QQ}L$
in Theorem~\ref{actiononhomotopy}
coincides with the canonical action of $MTG(L)\simeq \Aut(R_T)(L)$ on
$R_T(H^0(\uni_k\otimes_{M_X}\uni_k))\otimes_{\QQ}L\simeq H^0(\QQ\otimes_{T_X}\QQ)\otimes_{\QQ}L$. Since $MTG$ is integral, the coordinate ring on $MTG$
is a subring of $L$. We then deduce that the action of the group scheme $MTG$ on $H^0(\QQ\otimes_{T_X}\QQ)$
in Theorem~\ref{actiononhomotopy}
coincides with the natural action of $MTG\simeq \Aut(R_T)$.
\end{Remark}

\subsection{}
\label{compthmsection}
\begin{Theorem}
\label{compthm}
There is an isomorphism
\[
M_{DG}(X,x,y)\simeq M(X,x,y)
\]
in $\Ind(\TM)$.
\end{Theorem}

\begin{Lemma}
\label{coproductKan}
Let $\operatorname{Fin}$ be the category of (possibly empty) finite sets.
Let $\mathcal{C}$ be an $\infty$-category which has finite coproducts.
Then $\Fun^+(\operatorname{Fin},\mathcal{C})$ be the full subcategory
of $\Fun(\operatorname{Fin},\mathcal{C})$ spanned by those functors
that preserve finite coproducts.
Let $\Delta^0\to \operatorname{Fin}$ be the map determined by the set having
one element.
Then the composition induces an equivalence $\Fun^+(\operatorname{Fin},\mathcal{C})\to \Fun(\Delta^0,\mathcal{C})=\mathcal{C}$
of $\infty$-categories.
\end{Lemma}

\Proof
We here denote by $*$ the set having one element.
Since $\mathcal{C}$ has finite
coproducts,
any functor $\Delta^0\to \mathcal{C}$
admits a left Kan extension
along the inclusion $\Delta^0 =\{*\}  \to \operatorname{Fin}$.
Moreover, $F:\Fin\to \mathcal{C}$ is a left Kan extension
of $F|_{\{*\}}$ if and only if $F$ preserves finite coproducts.
Thus, by \cite[4.3.2.15]{HTT} 
$\Fun^+(\operatorname{Fin},\mathcal{C})\to \Fun(\Delta^0,\mathcal{C})=\mathcal{C}$ is an equivalence.
\QED

\begin{Example}
\label{coproductKanE}
Let $X\in \Sm_k$.
Let $\langle X \rangle$ be the subcategory of $\Sm_k$
defined as follows:
Objects are finite products of $X$, that is,
$\{\Spec k,X,X^2,\ldots,X^n,\ldots\}$.
A morphism $f:X^n\to X^m$ in $\Sm_k$ is a morphism in $\langle X \rangle$
if and only if $f$ is of the form $X^{n}\to X^m,\ (x_1,\ldots x_n)\mapsto (x_{i_1},\ldots x_{i_m})$ for some $\{i_1,\ldots,i_m\}\subset \{1,\ldots, n\}$.
Then there is an equivalence $\langle X \rangle^{op}\simeq \operatorname{Fin}$
which carries $X^n$ to the set having $n$ elements.
\end{Example}

{\it Proof of Theorem~\ref{compthm}.}
We first prove that
there is a natural isomorphism
$M_{DG}(X,x,y)\simeq M(X,x,y)$ in $\Ind(\TM)$.
Note the equivalence $\uni_k\otimes_{M_X}\uni_k \simeq \uni_k\otimes_{M_{X}\otimes M_X}M_X$
in $\CAlg(\DTM_\vee)$
where the right hand side is determined by $x^*\otimes y^*:M_X\otimes M_X\to \uni_k\otimes \uni_k\simeq \uni_k$ and 
$M_X\otimes M_X\simeq M_{X\times X}\to M_X$ induced by the diagonal $X\to X\times X$.
Here the two projections $X\leftarrow X\times X\to X$
determines a canonical equivalence
$M_X\otimes M_X\to M_{X\times X}$ in $\CAlg(\DTM_\vee^\otimes)$
(one way to see this is to observe that
the conservative realization $\CAlg(\DTM_\vee^\otimes)\to \CAlg_{\QQ}$
sends $M_X\otimes M_X\to M_{X\times X}$ to $T_X\otimes T_X\to T_{X\times X}$
that is an equivalence by K\"unneth formula).
Next we define a certain ``resolution'' of $M_X$ over $M_{X}\otimes M_X$.
For this purpose, let us consider the following cosimplicial
scheme
\[
R^{\Delta}(X):\Delta\to \Sm_k,\ \ \ [n]\mapsto X'\times X^n \times X''
\]
over $X'\times X''=X\times X$. Here, to avoid confusion
we put $X'=X$ and $X''=X$, and $X'\times X''$ is regarded as the constant
cosimplicial scheme. Cofaces are given by
\[
d^i(x_0,x_1,\ldots,x_{n+1})=(x_0,\ldots,x_{n-i+1},x_{n-i+1},\ldots x_{n+1}),\ \ \ 0\le i\le n+1,
\]
and codegeneracies are defined by projections.
If $X\to X'\times X''$ is the diagonal morphism, then
$R^{\Delta}(X)$ has a coaugmentation $X\to R^\Delta(X)$ over $X'\times X''$.
Observe that there is a the fiber product
of cosimplicial schemes
\[
\xymatrix{
P^{\Delta}(X,x,y) \ar[r] \ar[d] & R^{\Delta}(X) \ar[d] \\
\Spec k=(y,x) \ar[r] & X'\times X''
}
\]
where the right vertical map is the projection, and $\Spec k$ is considered to
be the constant cosimplicial scheme.
For each cosimplicial scheme, composing it with $\Xi:\Sm_k^{op}\to \CAlg(\DM^\otimes(k))$
we obtain simplicial objects $\mathcal{M}_\Delta(X,x,y)$, $\mathcal{M}_\Delta(X)$,
$M_{X'}\otimes M_{X''}$, $\uni_k$ in $\CAlg(\DM^\otimes(k))$ respectively from $P^{\Delta}(X,x,y)$,
$R^\Delta(X)$, $X'\times X''$ and $\Spec k$.
Each term of these simplicial objects lies in $\CAlg(\DTM_\vee)$ since
$M_{X^n}\simeq M_X^{\otimes n}$.
Consider the pushout $\uni_k\otimes_{M_{X'}\otimes M_{X''}}\mathcal{M}_\Delta(X)$ of simplicial objects (which consists of termwise pushouts).
There is a natural morphism of simplicial objects
\[
\uni_k\otimes_{M_{X'}\otimes M_{X''}}\mathcal{M}_\Delta(X)\to \mathcal{M}_\Delta(X,x,y).
\]
This morphism is an equivalence. To see this, it will suffice to prove
that the morphism in each term is an equivalence.
The morphism in the $n$-th term is equivalent to
\[
\uni_k\otimes_{M_{X'}\otimes M_{X''}}M_{X'}\otimes M_{X^n}\otimes M_{X''}\to M_{\{y\}\times X^n\times \{x\}}
\]
which is an equivalence.
Let $\mathcal{M}(X)$ be a colimit of $\mathcal{M}_\Delta(X)$
in $\CAlg(\DTM)$.
The coaugmentation $X\to R^\Delta(X)$ over $X'\times X''$
gives rise to $\mathcal{M}(X)\to M_X$ over $M_{X'}\otimes M_{X''}$.
Since $M_{DG}(X,x,y)=H^0(\mathcal{M}(X,x,y))$,
we will show that the induced map
\[
H^0(\uni_k\otimes_{M_{X'}\otimes M_{X''}}\mathcal{M}(X))\to H^0(\uni_k\otimes_{M_{X}\otimes M_X}M_X)
\]
is an isomorphism
in $\Ind(\TM)$.
To this end, 
recall the left completion $\DTM^\otimes \to \oDTM^\otimes$
from the second paragraph of the proof of
Proposition~\ref{algstr}. It is symmetric monoidal, $t$-exact and colimit-preserving.
We may and will replace $\DTM^\otimes$ by $\oDTM^\otimes$.
We show that $\uni_k\otimes_{M_{X'}\otimes M_{X''}}M_X$ is the colimit of $\uni_k\otimes_{M_{X'}\otimes M_{X''}}\mathcal{M}_\Delta(X)$ in $\oDTM$.
The image of $\mathcal{M}_\Delta(X)$ under the realization functor
is
the simplicial diagram in $\CAlg_{\QQ}$
given by 
the composite
\[
s:\Delta^{op} \stackrel{R^{\Delta}(X)^{op}}{\to}\Sm_k^{op} \stackrel{\Xi}{\to} \CAlg(\DM^\otimes(k))\stackrel{\mathsf{R}}{\to}\CAlg_{\QQ},\ \  [n]\mapsto T_{X'\times X^n\times X''}.
\] 
Since $\CAlg(\oDTM) \stackrel{\overline{\mathsf{R}}_T}{\to} \CAlg_{\QQ}$
is conservative and colimit-preserving, we are reduced to proving that
$s:[n]\mapsto T_{X'\times X^n\times X''}$
in $\CAlg_{\QQ}$
has a colimit $T_X$.
We let $F_X:\langle X \rangle^{op}\to \CAlg_{\QQ}$ be the functor
given by $X^m \mapsto T_{X^m}$. The natural projections
induce $T_X^{\otimes m}=T_X\otimes \ldots \otimes T_{X}\stackrel{\sim}{\to}T_{X^m}$, and $T_{\Spec k}\simeq \QQ$.
By Lemma~\ref{coproductKan} and Example~\ref{coproductKanE},
there is a canonical equivalence
$\Fun^+(\langle X \rangle^{op},\CAlg_{\QQ})\simeq \CAlg_{\QQ}$
which carries $F$ to $F(X)$.
Since $F_X$ belongs to $\Fun^+(\langle X \rangle^{op},\CAlg_{\QQ})$,
the functor $F_X$ that preserves
finite coproducts is ``uniquely determined'' by $F_X(X)=T_X$.
Let $A$ be a cofibrant commutative dg algebra over $\QQ$ that represents
$T_X$.
Let $\CAlg_{\QQ}^{dg}\to \CAlg_{\QQ}$ be the canonical functor (see Section~\ref{convention}). Let $f_A:\langle X \rangle^{op}\to \CAlg_{\QQ}^{dg}$ be
the functor given by $X^m \mapsto A^{\otimes m}$,
which corresponds to $A$ through the canonical equivalence
$\Fun^+(\langle X \rangle^{op},\CAlg_{\QQ}^{dg})\simeq \CAlg_{\QQ}^{dg}$.
The composite $F_A:\langle X \rangle^{op}\to \CAlg_{\QQ}$ is
the functor that preserves finite coproducts.
Thus $F_A\in \Fun^+(\langle X \rangle^{op},\CAlg_{\QQ})$.
It follows from $A\simeq T_X$ in $\CAlg_{\QQ}$
that $F_A\simeq F_X$.
Note that $R^\Delta(X)^{op}:\Delta^{op}\to \Sm_k^{op}$ uniquely
factors through the subcategory $\langle X \rangle^{op}\to \Sm_k^{op}$.
The composite
$s:\Delta^{op}\to \langle X \rangle^{op}\stackrel{F_X}{\to} \CAlg_{\QQ}$
is equivalent to $s':\Delta^{op}\to \langle X \rangle^{op}\stackrel{F_A}{\to} \CAlg_{\QQ}$.
We may replace $s$ by $s'$.
By unfolding the definition,
the simplicial commutative dg algebra
$s':\Delta^{op}\to \CAlg_{\QQ}^{dg},\ [n]\mapsto A\otimes A^{\otimes n}\otimes A$ (over $A\otimes A$) is 
the simplicial bar resolution of $A$ over $A\otimes A$:
$[n] \mapsto A\otimes A^{\otimes n} \otimes A$ (see \cite[4.3, 4.4, 4.6]{Ol} or \cite[3.7]{W} for what this means).
The (homotopy) colimit of the simplicial bar resolution
$[n] \mapsto A\otimes A^{\otimes n} \otimes A$ (equivalently the totalization)
is naturally equivalent to $A$.
(We remark that a colimit of a simplicial diagram of commutative algebra
objects is a colimit of simplicial diagram of underlying objects.)
Consequently, $\uni_k\otimes_{M_{X'}\otimes M_{X''}}\mathcal{M}(X)\simeq \uni_k\otimes_{M_{X'}\otimes M_{X''}}M_X$ in $\oDTM$.
Hence we obtain a canonical isomorphism $M_{DG}(X,x,y)\simeq M(X,x,y)$ in $\Ind(\TM)$.
\QED


{99}

\end{document}